\documentclass[11pt]{amsart}

\pdfoutput=1

\usepackage{amsmath,amssymb,amsthm}
\usepackage{mathtools,url}
\usepackage{stmaryrd,xcolor,tikz-cd,textcomp}
\usepackage{enumitem}
\usepackage[colorlinks=true,linkcolor=blue,urlcolor=black,bookmarksopen=true]{hyperref}
\usepackage[open,openlevel=2]{bookmark}
\numberwithin{equation}{section}

\usepackage[top=1.25in, bottom=1in, left=1.25in, right=1.25in]{geometry}
\newcommand{\into}{\hookrightarrow}
\newcommand{\onto}{\twoheadrightarrow}

\newcommand\GHilb{\operatorname{G-Hilb}}

\DeclareMathOperator{\Hom}{Hom}
\DeclareMathOperator{\Ext}{Ext}
\DeclareMathOperator{\End}{End}
\DeclareMathOperator{\Tor}{Tor}

\DeclareMathOperator{\Spec}{Spec}
\DeclareMathOperator{\Proj}{Proj}
\DeclareMathOperator{\img}{im}

\DeclareMathOperator{\Sym}{Sym}
\DeclareMathOperator{\Hilb}{Hilb}
\DeclareMathOperator{\Aut}{Aut}
\DeclareMathOperator{\ch}{ch}
\DeclareMathOperator{\orbch}{\wtilde{\ch}}
\DeclareMathOperator{\tr}{tr}
\DeclareMathOperator{\td}{td}
\DeclareMathOperator{\orbtd}{\wtilde{\td}}
\DeclareMathOperator{\ord}{ord}
\DeclareMathOperator{\Supp}{Supp}
\DeclareMathOperator{\Ind}{Ind}

\DeclareMathOperator{\Quot}{Quot}
\DeclareMathAlphabet{\mathpzc}{OT1}{pzc}{m}{it}

\newcommand{\vphi}{\varphi}
\newcommand{\wtilde}{\widetilde}

\newcommand{\subeq}{\subseteq}

\newcommand{\inprod}[1]{\langle#1\rangle}
\newcommand{\tb}{\textbf}
\newcommand{\orbv}{\wtilde{v}}
\newcommand{\orbw}{\wtilde{w}}

\newcommand{\orbs}{\wtilde{s}}
\newcommand{\orbt}{\wtilde{t}}

\newcommand{\cj}{\overline}
\newcommand{\orbP}{\wtilde{P}}
\newcommand{\orbp}{\wtilde{p}}

\newcommand{\id}{\mathrm{id}}

\newcommand{\D}{\mathrm{D^b}}
\newcommand{\rk}{\mathrm{rk}}
\newcommand{\Pic}{\mathrm{Pic}}

\newcommand{\Qcoh}{\mathrm{Qcoh}}
\newcommand{\Coh}{\mathrm{Coh}}
\newcommand{\Vect}{\mathrm{Vect}}

\newcommand{\Stab}{\mathrm{Stab}}

\newcommand{\eig}{\mathrm{eig}}
\newcommand{\GL}{\mathrm{GL}}
\newcommand{\PGL}{\mathrm{PGL}}
\newcommand{\SL}{\mathrm{SL}}
\newcommand{\reg}{\mathrm{reg}}
\newcommand{\pt}{\mathrm{pt}}
\newcommand{\fix}{\mathrm{fix}}
\newcommand{\mov}{\mathrm{mov}}
\newcommand{\Tot}{\mathrm{Tot}}
\newcommand{\Ad}{\mathrm{Ad}}

\newcommand{\Sch}{\mathrm{Sch}}
\newcommand{\Obj}{\mathrm{Obj}}
\newcommand{\Rep}{\mathrm{Rep}}
\newcommand{\Sh}{\mathrm{Sh}}
\newcommand{\Real}{\mathrm{Re}}

\newcommand{\Amp}{\mathrm{Amp}}
\newcommand{\ob}{\mathrm{ob}}
\newcommand{\fppf}{\mathrm{fppf}}

\newcommand{\gr}{\mathrm{gr}}

\newtheorem{thm}{Theorem}[section]
\newtheorem{cor}[thm]{Corollary}
\newtheorem{prop}[thm]{Proposition}
\newtheorem{lem}[thm]{Lemma}

\theoremstyle{definition}
\newtheorem{defn}[thm]{Definition}
\newtheorem{defn-prop}[thm]{Definition/Proposition}
\newtheorem{ex}[thm]{Example}

\newtheorem{rmk}[thm]{Remark}
\newtheorem{notation}[thm]{Notation}
\newtheorem{qn}[thm]{Question}
\newtheorem{constrn}[thm]{Construction}
\newtheorem{situation}[thm]{Situation}

\newcommand{\CC}{\mathbb{C}}

\newcommand{\PP}{\mathbb{P}}
\newcommand{\QQ}{\mathbb{Q}}
\newcommand{\RR}{\mathbb{R}}

\newcommand{\ZZ}{\mathbb{Z}}

\newcommand{\cA}{\mathcal{A}}

\newcommand{\cD}{\mathcal{D}}
\newcommand{\cE}{\mathcal{E}}
\newcommand{\cF}{\mathcal{F}}

\newcommand{\cH}{\mathcal{H}}

\newcommand{\cK}{\mathcal{K}}
\newcommand{\cL}{\mathcal{L}}
\newcommand{\cM}{\mathcal{M}}
\newcommand{\cN}{\mathcal{N}}
\newcommand{\cO}{\mathcal{O}}
\newcommand{\cP}{\mathcal{P}}

\newcommand{\cT}{\mathcal{T}}

\newcommand{\cV}{\mathcal{V}}
\newcommand{\cW}{\mathcal{W}}
\newcommand{\cX}{\mathcal{X}}
\newcommand{\cY}{\mathcal{Y}}
\newcommand{\cZ}{\mathcal{Z}}

\newcommand{\fg}{\mathfrak{g}}

\newcommand{\tu}[1]{\textup{#1}}

\begin{document}
\title{Orbifold Hirzebruch-Riemann-Roch for Quotient Deligne-Mumford Stacks and Equivariant Moduli Theory on $ K3 $ Surfaces}
\author{Yuhang Chen}

\begin{abstract}
	We study the orbifold Hirzebruch-Riemann-Roch (HRR) theorem for quotient Deligne-Mumford stacks, explore its relation with the representation theory of finite groups, and derive a new orbifold HRR formula via an orbifold Mukai pairing. As a first application, we use this formula to compute the dimensions of $ G $-equivariant moduli spaces of stable sheaves on a $ K3 $ surface $ X $ under the action of a finite subgroup $ G $ of its symplectic automorphism group. We then apply the orbifold HRR formula to reproduce the number of fixed points on $ X $ when $ G $ is cyclic without using the Lefschetz fixed point formula. We prove that under some mild conditions, equivariant moduli spaces of stable sheaves on $ X $ are irreducible symplectic manifolds deformation equivalent to Hilbert schemes of points on $ X $ via a connection between Gieseker and Bridgeland moduli spaces, as well as the derived McKay correspondence.
\end{abstract}

\maketitle
\tableofcontents

\section{Introduction}
Let $ X $ be a complex projective $ K3 $ surface. Its canonical bundle $ \omega_X $ is trivial and admits a nowhere-vanishing global section $ \theta, $ which is a holomorphic symplectic form such that at each point $ p $ of $ X $, there is a non-degenerate skew-symmetric pairing
\begin{equation*}
\theta_p: T_pX \times T_pX \to \CC,
\end{equation*}
where $ T_pX $ is the tangent space of $ X $ at $ p. $ The symplectic form $ \theta $ is unique up to a scalar. An automorphism $ \sigma $ of $ X $ is called symplectic if $ \sigma^*\theta = \theta. $

It's known that (see for example, \cite[Section 15.1]{huybrechts2016lectures}) if a symplectic automorphism $ \sigma $ has finite order $ n, $ then $ n \leq 8, $ and in this case, the number of points fixed by $ \sigma $ is finite and only depends on $ n $ as shown in Table \ref{table_number_fixed_points}. These numbers were computed by Mukai via the Lefschetz fixed point formula in \cite{mukai1988finite}.
\begin{table}
	\centering
	\caption{The number of fixed points $ f_n $ of a symplectic automorphism on a complex projective $ K3 $ surface with finite order $ n $}\label{table_number_fixed_points}
	\begin{tabular}{c|l*{7}{c}}
		$ n $ & \tu{2} & \tu{3} & \tu{4} & \tu{5} & \tu{6} & \tu{7} & \tu{8} \\ \hline
		$ f_n $ & \tu{8} & \tu{6} & \tu{4} & \tu{4} & \tu{2} & \tu{3} & \tu{2}
	\end{tabular}
\end{table}
A symplectic automorphism $ \sigma $ of finite order determines a cyclic group which acts on $ X $ faithfully. In general, we can consider a finite group $ G $ acting on $ X $ faithfully and symplectically, which means $ G $ can be identified with a symplectic automorphism subgroup of $ X $, and hence the possible groups are limited by the intrinsic geometry of the $ K3 $ surface $ X $. 
All such possible symplectic actions were classified by Xiao in \cite{xiao1996galois}. For example, there are in total $ 80 $ possible groups (including the trivial one) with a maximal order of $ 960 $ that can occur as a subgroup of the symplectic automorphism group of a $ K3 $ surface. Therefore, symplectic actions on $ K3 $ surfaces are well understood by now.

On the other hand, moduli spaces of sheaves on $ K3 $ surfaces have been studied extensively by Beauville \cite{beauville1983varietes} , Huybrechts \cite{huybrechts1997birational}, Mukai \cite{mukai1984symplectic}, O'Grady \cite{o1997weight}, and Yoshioka \cite{yoshioka2001moduli} among others. By imposing a suitable stability condition known as Gieseker stability via an ample line bundle $ H $ on $ X $ called a polarization, and choosing an element $ x $ in the numerical Grothendieck ring $ N(X) $, one can construct a projective scheme $ M(X,x) $ from a Quot scheme via the geometric invariant theory (GIT) such that $ M(X,x) $ is a good moduli space of semistable sheaves on $ X $ with numerical class $ x. $ Indeed, such a GIT construction for moduli spaces of sheaves exists for any projective scheme. 

Unfortunately, $ M(X,x) $ is not a coarse moduli space since two semistable sheaves may share a point in $ M(X,x) $ when they are so-called $ S $-equivalent\footnote{Here $ S $ stands for Seshadri who first introduced the notion of $ S $-equivalence for vector bundles over a curve in \cite{seshadri1967space} where he called it strong equivalence.}. The GIT approach produces a quasi-projective scheme $ M^s(X,x) $ which is a coarse moduli space of stable sheaves with numerical class $ x. $ With a suitable choice of $ x $ and $ H, $ we have 
$$ M(X,x) = M^s(X,x), $$ 
i.e., there is no semistable sheaf with numerical class $ x $, and it follows that $ M(X,x) $ is deformation equivalent to a Hilbert scheme of points $ \Hilb^n(Y) $ on a $ K3 $ surface $ Y $, which is an irreducible symplectic manifold\footnote{A smooth complex projective variety $ X $ is an irreducible symplectic manifold if $ H^0(X, \Omega^2) = \CC \omega $ for a holomorphic symplectic form $ \omega $ and $ H^1(X, \cO_X) = 0. $}. Here $ Y $ can be taken to be the original $ K3 $ surface as $ \Hilb^n(X) $ is deformation equivalent to $ \Hilb^n(Y) $ for any $ K3 $ surface $ Y $. Therefore, the deformation type of a moduli space $ M(X,x) $ is determined by its dimension. For an element $ x $ in $ N(X) $, there is an associated Mukai vector $ v(x) = \ch(x) \sqrt{\td_X} $ in the numerical Chow ring $ R(X) $ of $ X $ which is equivalent to the Chern character $ \ch(x) $ and hence to the element $ x $ as well. The intersection product on $ R(X) $ induces an integral Mukai pairing\footnote{Note that the Mukai pairing here differs by a sign from the original one defined by Mukai in \cite{mukai1984moduli}.}
$$ \inprod{v,w} = v_0w_2 - v_1w_1 + v_2w_0, $$
where $ v = v_0 + v_1 + v_2 $ with $ v_i  \in R^i(X) $ and $ w = w_0 + w_1 + w_2 $ with $ w_i \in R^i(X) $ for $ i = 0, 1, 2 $.
The Hirzebruch-Riemann-Roch (HRR) theorem gives a formula
\begin{equation*}\label{HRR}
	\chi(x, y) = \inprod{v(x), v(y)}
\end{equation*}
for all $ x $ and $ y $ in $ N(X). $ It then follows that
\begin{equation*}\label{eq_dim}
	\dim M(X,x) = 2 - \inprod{v(x)^2}.
\end{equation*}
See Appendix \ref{review_mukai} for details on the Mukai pairing and the HRR formula for a proper smooth scheme.


Now, it's natural to ask the following question.
\begin{qn}\label{qn_equivariant_moduli_space}
	What can we say about equivariant moduli spaces of sheaves on a $ K3 $ surface with symplectic automorphisms?
\end{qn}

Let $G$ be a finite subgroup of the symplectic automorphism group of $X.$ The category of $ G $-equivariant sheaves on $ X $ is equivalent to the category of sheaves on the quotient stack $ \cX = [X/G]. $ Therefore, we can consider moduli spaces of sheaves on $ \cX $ instead. By a GIT construction of Nironi in \cite{nironi2009moduli}, there exists a moduli space $ M^{(s)}(\cX,x) $ of semistable (resp. stable) sheaves on $ \cX $ with numerical class $ x $ in $ N(\cX) $. Roughly speaking, one chooses a polarization $ (\cH, \cV) $
on $ \cX $ to impose stability conditions, where $ \cH $ is a line bundle on $ \cX $ which descends to its coarse moduli space, and $ \cV $ is a generating sheaf on $ \cX $ which is a locally free sheaf on $ \cX $ such that it contains all irreducible representations of the automorphism group at each point on $ \cX. $ In our case, the geometric quotient $ X/G $ is a coarse moduli space of $ \cX $, so we can pullback an ample line bundle on $ X/G $ to a line bundle $ \cH $ on $ \cX $. There is also a canonical generating sheaf 
$$ \cV_\reg = \cO_\cX \otimes \rho_\reg $$ 
on $ \cX $, where $ \rho_\reg $ is the regular representation of $ G. $ Therefore, fixing an element $ x $ of $ N(\cX), $ we have a moduli space $ M^{(s)}(\cX,x) $ of semistable (resp. stable) sheaves on $ \cX $ similarly to the case of schemes.

The canonical morphism
$$ p: X \to \cX $$
is smooth, and hence induces an exact pullback functor
\begin{equation*}
p^*: \Coh(\cX) \to \Coh(X), \quad \cE = (E, \phi) \mapsto E.
\end{equation*}
The line bundle $ \cH $ on $ \cX $ pulls back to a $ G $-invariant line bundle $ H = p^*\cH $ on $ X $. 
Let $ N(\cX) $ denote the numerical Grothendieck ring of $ \cX $ and $ X $. 
For a coherent sheaf $ \cE $ (resp. $ E $) on $ \cX $ (resp. $ X $), let $ \gamma(\cE) $ (resp. $ \gamma(E) $) denote its image in $ N(\cX) $ (resp. $ N(X) $).
Then there is a numerical pullback 
$$ p^N: N(\cX) \to N(X), \quad \gamma(\cE) = \gamma(E,\phi) \mapsto \gamma(E). $$ 
The action of $ G $ on $ X $ induces an action on $ N(X) $ and an action on $ R(X) $. Let $ N(X)^G $ and $ R(X)^G $ denote their $ G $-invariant subspaces respectively. 
An element $ x $ in $ N(\cX) $ determines an orbifold Mukai vector $ \orbv(x) $ in the complex numerical Chow ring $ R(I\cX)_\CC $ of the inertia stack $ I\cX. $ An element $ x $ in $ N(\cX) $ determines an orbifold Mukai vector 
$$ \orbv(x) = \orbch(x) \sqrt{\td_{I\cX}}$$ 
in the complex numerical Chow ring $ R(I\cX)_\CC $ of the inertia stack $ I\cX, $ where $ \orbch(x) $ is the orbifold Chern character of $ x $ and $ \td_{I\cX} $ is the Todd class of the inertia stack $ I\cX $.
There is an orbifold Mukai pairing
$$ \inprod{\cdot\ {,}\ \cdot}_{I\cX}: R(I\cX)_\CC \times R(I\cX)_\CC \to \CC $$ 
which is conjugate linear in the first component and linear in the second. The orbifold HRR theorem gives a formula
\begin{equation}\label{eq_orbifold_HRR}
	\chi(x, y) = \inprod{\orbv(x), \orbv(y)}_{I\cX}
\end{equation}
for all $ x $ and $ y $ in $ N(\cX). $ See Section \ref{sec_orbifold_HRR} for details on the orbifold Mukai pairing and the HRR formula for a connected proper smooth quotient Deligne-Mumford (DM) stack. An explicit form of the orbifold Mukai pairing for $ \cX = [X/G] $ is worked out in Section \ref{sec_HRR_K3/G}. By the orbifold HRR formula \eqref{eq_orbifold_HRR}, the moduli space $ M^s(\cX, x) $ is either empty or a smooth quasi-projective scheme with
\begin{equation*}
	\dim M^s(\cX, x) = 2 - \inprod{\orbv(x)^2}_{I\cX}.
\end{equation*}

Now we answer Question \ref{qn_equivariant_moduli_space} in the following theorem.
\begin{thm}\label{thm_main}
	Let $ x $ be an element in $ N(\cX) $ with $ y = p^N x $ in $ N(X)^G $. Suppose $ y $ is primitive with $ \rk(y) > 0 $ and $ H $ is $ y $-generic.
	Then $ M(\cX,x) = M^s(\cX,x), $ which is non-empty if and only if $ \inprod{\orbv(x)^2}_{I\cX} \leq 2. $  If $ M(\cX,x) $ is non-empty, then it is an irreducible symplectic manifold of dimension $ n = 2 - \inprod{\orbv(x)^2}_{I\cX} $ deformation equivalent to $ \Hilb^{n/2}(X) $.
\end{thm}

The proof of Theorem \ref{thm_main} makes use of three notions: 
\begin{enumerate}[font=\normalfont,leftmargin=3em]
	\item Bridgeland stability conditions
	\item Induced stability conditions
	\item Derived McKay correspondence
\end{enumerate}
See Section \ref{sec_proof_main_thm} for a complete proof. We outline the proof here. Let $ M_1 = M(\cX,x). $ Roughly speaking, the argument can be divided into three steps.

\vspace{2pt}

\tb{Step 1.} Our starting point is that the $ G $-invariant ample line bundle $ H $ on the $ K3 $ surface $ X $ determines a real number $ t_0 $ such that for all $ t \geq t_0, $ there exist $ G $-invariant stability conditions $ \sigma_t $ in the distinguished component $ \Stab^\dagger(X) $ in the space of stability conditions on $ \D(X), $ and $ H $-stable sheaves in $ \Coh(X) $ are precisely $ \sigma_t $-stable objects in $ \D(X). $ These $ G $-invariant stability conditions $ \sigma_t $ induce stability conditions $ \wtilde{\sigma}_t $ on $ \D(\cX) $ by a technique developed in \cite{macri2009inducing}. We then show that $ \cH $-stable sheaves in $ \Coh(\cX) $ are precisely $ \wtilde{\sigma}_t $-stable objects in $ \D(\cX). $ This makes it possible to identify the Gieseker moduli space $ M_1 $ with a Bridgeland moduli space $ M_2 $ of stable objects in $ \D(\cX). $

\vspace{2pt}

\tb{Step 2.} The $ G $-Hilbert scheme $ M = \GHilb X $, which is also a $ K3 $ surface, gives the minimal resolution of the surface $ X/G. $ The derived McKay correspondence gives a Fourier-Mukai transform 
$$ \Phi: \D(\cX) \xrightarrow{\sim} \D(M) $$
such that under $ \Phi, $ the stability conditions $ \wtilde{\sigma}_t $ on $ \D(\cX) $ are mapped to stability conditions $ \Phi^S(\wtilde{\sigma}_t) $ on $ \D(M) $, and the Bridgeland moduli space $ M_2 $ on $ \cX $ is then identified with another Bridgeland moduli space $ M_3 $ on $ M. $

\vspace{2pt}

\tb{Step 3.} By \cite[Proposition 6.1]{beckmann2022equivariant}, the stability conditions $ \Phi^S(\wtilde{\sigma}_t) $ are in the distinguished component $ \Stab^\dagger(M). $ Under $ \Phi, $ the primitive element $ x $ in $ N(\cX) $ is mapped to another primitive element $ \Phi^N(x) $ in $ N(M) $. In the same time we can choose $ t $ such that $ \Phi^S(\wtilde{\sigma}_t) $ is generic with respect to the element $ \Phi^N(x) $. The main result of \cite{bottini2024stable} then says that the Bridgeland moduli space $ M_3 $ is deformation equivalent to a Hilbert scheme of points on a $ K3 $ surface, and hence so is our moduli space $ M_1. $

\subsection{Outline}
The paper is structured as follows. In Section \ref{sec_sheaves_on_quotient_stacks}, we review some basic notions about quotient stacks, provide essential definitions, in particular the definition of a point of a quotient stack, study (quasi-)coherent sheaves on a quotient stack $ [X/G], $ and relate them to $ G $-equivariant (quasi-)coherent sheaves on $ X. $

In Section \ref{sec_HRR_formula_stacks}, we review some basic facts about the $ K $-theory of quotient stacks and derive an orbifold HRR formula for connected proper smooth quotient DM stacks. Let $ \cX $ be a connected separated quotient DM stack. We will study the inertia stack $ I\cX $ of $ \cX $ in details. In particular, we will decompose $ I\cX $ into connected components. Such a decomposition is essential in the constructions of the orbifold Chern character $ \orbch(E) $ and the orbifold Todd class $ \orbtd(E) $ of a coherent sheaf $ E $ on $ \cX $ when $ \cX $ is smooth. We will derive explicit formulas for $ \orbch(E) $ and $ \orbtd(E) $. We will define an orbifold Mukai pairing and derive an orbifold HRR formula for a pair of coherent sheaves on $ \cX $ when it is proper and smooth. We will provide plenty of examples along the way. Some of them contain new results. For instance, the orbifold HRR formula for the classifying stack $ BG $ when $ G \cong \ZZ/n\ZZ $ recovers Parseval's theorem for the discrete Fourier transform as shown in Examples \ref{ex_discrete_FT} and \ref{ex_descrete_Parseval}.

In Section \ref{sec_moduli_sheaves_stacks}, we review Gieseker stability conditions for coherent sheaves on polarized projective schemes and stacks, as well as the GIT constructions of the corresponding moduli spaces. We then focus on the case of a projective quotient stack $ [X/G] $ where $X$ is a smooth projective scheme and $G$ is a finite group, and define $ G $-equivariant moduli spaces of stable sheaves on $ X. $ Finally we compare the stabilities of $ G $-equivariant coherent sheaves on $ X $, i.e., coherent sheaves on $ [X/G] $, and their underlying sheaves on $ X $.

In Section \ref{sec_moduli_sheaves_K3/G}, we derive an explicit orbifold HRR formula for a quotient stack $ [K3/G] $ and use it to compute the dimensions of $ G $-equivariant moduli spaces of stable sheaves on a $ K3 $ surface $ X $. As a joyful digression, we will apply the orbifold HRR formula to reproduce the number of fixed points on $ X $ when $ G \cong \ZZ/n\ZZ $ without using the Lefschetz fixed point formula. Next we review the derived McKay correspondence, Bridgeland stability conditions, and induced stability conditions on a quotient stack $ [X/G] $ where $ X $ is a smooth projective variety and $ G $ is a finite group acting on $ X $. We then prove Theorem \ref{thm_main}. In Section \ref{sec_Hilb}, we give the definition of a $ G $-equivariant Hilbert scheme of points on a projective $ G $-scheme for a finite group $G$. We will specialize Theorem \ref{thm_main} in the rank one case to determine the deformation types of $ G $-equivariant Hilbert schemes of points on a $ K3 $ surface. We will compute a few $ G $-equivariant Hilbert schemes of points when $ G \cong \ZZ/2\ZZ $ is generated by a Nikulin involution.

In Appendix \ref{proof_equiv_cat_sheaves}, we prove the equivalence between the category of sheaves on a quotient stack $ \cX = [X/G] $ and that of $ G $-equivariant sheaves on the scheme $ X $. In Appendix \ref{review_mukai}, we review the Mukai pairing and the HRR formula for proper smooth schemes.

\subsection{Conventions}
Let $ S $ be a base scheme. Let $ (\Sch/S) $ be the category of schemes over $ S $. Denote by $ (\Sch/S)_{\fppf} $ the big fppf site of  $ S $. An algebraic space over $ S $ is a sheaf $ X $ of sets on $ (\Sch /S)_\fppf $ such that the diagonal morphism $ X \to X \times_S X $ is representable by schemes over $ S $, and there is a scheme $ U $ over $ S $ and a surjective \'etale morphism $ U \to X $ over $ S. $ A stack over $ S $ is a stack of groupoids over $ (\Sch /S)_\fppf $. An algebraic (resp. DM) stack is a stack $ \cX $ such that the diagonal morphism $ \cX \to \cX \times_S \cX $ is representable by algebraic spaces over $ S $, and there is a scheme $ U $ over $ S $ and a surjective smooth (resp. \'etale) morphism $ U \to \cX $ over $ S. $ Note that we have a chain of strict inclusions of 2-categories over $ S $:
\begin{equation*}
\text{Schemes} \ \subset \ \text{Algebraic Spaces} \ \subset \ \text{DM Stacks} \ \subset \ \text{Algebraic Stacks} \ \subset \ \text{Stacks}
\end{equation*}

If $ S = \Spec k $ for a base field $k,$ then we will use the words ``over $k$" in place of ``over $S$". We will often omit the words ``over $ S $" (resp. ``over $k$") when the base scheme $ S $ (resp. the base field $ k $) is understood. The notation $ \pt = \Spec k $ for a field $ k $ which will be clear in the context, and is usually $ \CC $. The fiber product $ X \times_\CC Y $ of two schemes $ X $ and $ Y $ over $ \CC $ will be denoted by $ X \times Y. $

Let $X$ be a scheme over a field $ k $. A point of $ X $ means a closed point of degree one, which is the same as a $ k $-valued point, i.e., a morphism $ \pt \to X $. The notation $x \in X$ means a point $ x $ of $ X. $

Let $\cX$ be a locally noetherian algebraic stack over $ S $. A sheaf on $ X $ means a coherent sheaf of $ \cO_\cX $-modules. A vector bundle on $ \cX $ means a locally free sheaf on $ \cX $ of a finite constant rank. The categories of quasi-coherent sheaves, coherent sheaves, and vector bundles on $ \cX $ are denoted by $ \Qcoh(\cX), \Coh(\cX), $ and $ \Vect(\cX), $ respectively.

A ring means a commutative ring with unity. For a ring $ R, $ let $ R_\CC = R \otimes_\ZZ \CC $ denote its extension of scalars from $ \ZZ $ to $ \CC, $ and write $ a r$ for a pure tensor $ r \otimes a $ in $ R_\CC; $ let $ R^\times $ denote the multiplicative group of $ R, $ i.e., the group of invertible elements in $ R $.

Let $ G $ be a linear algebraic group over $ \CC $. The identity of $ G $ is denoted by $ 1. $ A scheme over $ \CC $ under a left action of $ G $ is called a $ G $-scheme. A representation of $ G $ means a rational representation of $ G $. The character group of $ G $ is denoted by $ \widehat{G}, $ i.e., $ \widehat{G} = \Hom(G, \CC^*). $ For a positive integer $ n, $ let $ \mu_n $ denote the cyclic group of $ n $th roots of unity in $ \CC. $

\subsection*{Acknowledgements}
I would like to thank my advisor Hsian-Hua Tseng for providing a free research environment and help when necessary. I want to thank Yunfeng Jiang for suggesting the problem of $ [K3/G] $ and constant discussions throughout the research process. I thank Promit Kundu for many discussions on stacks and moduli spaces of sheaves on $ K3 $ surfaces, as well as pointing out a few mistakes in an early version of my paper. I also thank David Anderson for discussions on intersection theory and some useful comments. I am grateful to James Cogdell for discussions on representation theory and many valuable comments. I thank Yonghong Huang for discussions on stacks. I thank Xiping Zhang for discussions on algebraic geometry. I thank Yilong Zhang for discussions on Hilbert polynomials. I also thank Hao Sun for discussions on Bridgeland stability conditions. I am grateful to Jim Bryan for multiple email conversations on equivariant Hilbert schemes of points on $ K3 $ surfaces. I wish to thank Johan de Jong, Dan Edidin, Daniel Huybrechts, Martijn Kool, Georg Oberdieck, Martin Olsson, and Alessandra Sarti for useful email correspondence. I am thankful to Li Li for sharing a draft on Hilbert schemes of points on tame DM stacks with me. I thank Nicolas Addington for pointing out a mistake in an early version of my paper. I also thank Benjamin Call, Deniz Genlik, Jean-Pierre Serre, and Luke Wiljanen for pointing out various typos in earlier versions of this paper.

\section{Sheaves on quotient stacks}\label{sec_sheaves_on_quotient_stacks}

\subsection{Quotient stacks}
In this section we review some basic notions related to quotient stacks and set up a few notations.

We start from the definition of a quotient stack. Let $ S $ be a scheme. 
\begin{defn}\label{defn_quotient_stack}
	Let $ X $ be a separated scheme of finite type over $ S .$ Let $ G $ be an affine smooth group scheme of finite type over $ S. $ Let $ G $ act on $ X $ from the left. Define a category $ [X/G] $ over $ (\Sch/S) $ as follows:
	\begin{enumerate}[font=\normalfont,leftmargin=*]
		\item An object in $ [X/G] $ over a scheme $ Y $ in $ (\Sch/S) $ is a principal $ G $-bundle $ E \to Y $ with a $ G $-equivariant morphism $ E \to X $.
		\item An arrow in $ [X/G] $ over a morphism $ Y \to Z $ in $ (\Sch/S) $ is a morphism of principal $ G $-bundles over $ Y \to Z $ compatible with the $ G $-equivariant morphisms to $ X, $ i.e., a commutative diagram:
		\begin{equation*}
		\begin{tikzcd}[column sep=3em,row sep=3em]
		E \arrow[bend left=40]{rr} \arrow{r} \arrow{d} & F \arrow{d} \arrow{r} & X \\
		Y \arrow{r} & Z &
		\end{tikzcd}
		\end{equation*}
	\end{enumerate}
	\vspace{8pt}
	The category $ [X/G] $ is called the \tb{quotient stack} of $ X $ by $ G $ over $ S $.
\end{defn}
\begin{rmk}
	Let $ [X/G] $ be a quotient stack over $ S $. By definition, $ [X/G] $ is a category fibered in groupoids over $ (\Sch/S). $ Furthermore, $ [X/G] $ is an algebraic stack over $ S $. This is proved in \cite[Example 8.1.12]{olsson2016algebraic} where $ G $ is a smooth group scheme, not necessarily affine. If the stabilizer of every geometric point of $ X $ is finite and reduced, then \cite[Example 7.17]{vistoli1989intersection} shows that $ [X/G] $ is a DM stack over $ S $. In general, the canonical morphism $ X \to [X/G] $ is representable by schemes and is a smooth surjection: a morphism $ Y \to [X/G] $ from a scheme $ Y $ corresponds to an object $ (E \to Y, E \to X) $ in $ [X/G] $ over $ Y, $ which gives a $ 2 $-pullback diagram
	\vspace{8pt}
	\begin{equation*}
	\begin{tikzcd}[column sep=3em,row sep=3em]
	E \arrow{r} \arrow{d} & X \arrow{d} \\
	Y \arrow{r} & \left[X/G\right]
	\end{tikzcd}\vspace{8pt}
	\end{equation*}
	in the 2-category of algebraic stacks over $ S $. Therefore, $ [X/G] $ appears as if $ G $ acts on $ X $ freely such that the canonical morphism $ X \to [X/G] $ is a principal $ G $-bundle.
\end{rmk}
Let $ S = \Spec k $ for a field $ k. $ An important class of quotient stacks over $ k $ are projective stacks including weighted projective stacks and quotients of projective schemes by finite groups. We have a chain of strict inclusions of 2-categories over $ k: $
\begin{equation*}
\text{Projective Schemes} \ \subset \ \text{Projective Stacks} \ \subset \ \text{Proper Tame Quotient DM Stacks}
\end{equation*}
We will define projective stacks over $ k $ within the class of tame stacks over $ k $ in Section \ref{sec_moduli_sheaves_stacks2}. If $ k = \CC, $ then projective stacks over $ k $ have a rather simple description since tameness is automatic in characteristic zero. We will give the definition of a projective stack over $ \CC $ within the class of separated quotient DM stacks over $ \CC $ later in this section.

\vspace{6pt}

From now on, the base scheme $ S = \Spec \CC $ unless otherwise specified. 
\begin{rmk}
	An affine smooth group scheme of finite type is the same as a linear algebraic group $ G \subseteq \GL(n, \CC). $ Such a group $ G $ is also a complex Lie group, so we can talk about its Lie algebra $ \fg, $ i.e., the tangent space of $ G $ at its identity. For example, the multiplicative group
	\begin{equation*}
	\mathbb{G}_m = \Spec \CC[x, x^{-1}]
	\end{equation*}
	over $ \CC $ is the one-dimensional linear algebraic group $ \GL(1,\CC) $ with a Lie algebra $ \CC. $
	The $ \CC $-valued points of $ \mathbb{G}_m $ are the torus $ \CC^*, $ so we will also use the notation $ \CC^* $ for $ \mathbb{G}_m. $ We are often interested in its subgroup $ \mu_n \subset \CC^* $ defined by $ x^n - 1, $ i.e.,
	\begin{equation*}
	\mu_n = \Spec \CC[x]/(x^n - 1).
	\end{equation*}
	The $ \CC $-valued points of $ \mu_n $ are the $ n $-th roots of unity in $ \CC, $ so we can endow the cyclic group $ \ZZ/n\ZZ $ a group scheme structure. In general, every finite group $ G $ can be viewed as a linear algebraic group since we can embed $ G $ into a symmetric group $ S_n, $ which can be represented by permutation matrices in $ \GL(n,\CC). $
\end{rmk}
\begin{notation}
	Let $ [X/G] $ be a quotient stack. There is a diagonal action of $ G $ on $ X \times X $ given by
	\begin{equation*}
	g(x,y) = (gx, gy)
	\end{equation*}
	for $ g \in G $ and $ x, y \in X. $ 
	Since $ G $ acts on itself by conjugation, there is also a diagonal action of $ G $ on $ G \times X $ given by
	\begin{equation*}
	h(g, x) = (hgh^{-1}, hx)
	\end{equation*}
	for $ g, h \in G $ and $ x \in X. $ 
	The action of $ G $ on $ X $ gives a morphism
	\begin{equation*}
	\sigma: G \times X \to X, \quad (g,x) \mapsto gx,
	\end{equation*}
	which restricts to an \tb{orbit} morphism
	\begin{equation*}
	\sigma_x: G \to X, \quad g \mapsto gx,
	\end{equation*}
	for every $ x \in X. $ Let $ p_2: G \times X \to X $ denote the projection to $ X $. Then $ \alpha = (p_2, \sigma) $ gives an \tb{action} morphism, i.e., 
	\begin{equation*}
	\alpha: G \times X \to X \times X, \quad (g,x) \mapsto (x, gx),
	\end{equation*}
	which is $ G $-equivariant with respect to the diagonal $ G $-actions on $ G \times X $ and $ X \times X. $
\end{notation}
\begin{rmk}[The difference between objects and points in a quotient stack]
	Recall that a point of a scheme $ X $ is a morphism $ \pt \to X, $ which corresponds to its image $ x \in X. $ A natural question is:
	\vspace{5pt}
	\begin{center}
		\textit{What is a point of a quotient stack $ \cX = [X/G] ? $}
	\end{center}
	\vspace{5pt}
	This is a subtle question. The desired answer is that a point of the stack $ \cX $ corresponds to a $ G $-orbit in $ X, $ since this would then reduce to the notion of a point of the scheme $ \cX $ when $ G $ acts on $ X $ freely such that $ \cX $ is (i.e., can be represented by) a scheme. With this in mind, now we work out the details. In the language of stacks, $ \pt $ is the category $ (\Sch/\CC) $, and a morphism $ f: \pt \to \cX $ is a functor between two categories. For every $ x \in X, $ an object $ \eta_x $ of the groupoid $ \cX(\CC) $ is a diagram
	\vspace{5pt}
	\begin{equation*}
	\eta_x = \left(
	\begin{tikzcd}[column sep=2.8em,row sep=2.8em]
	G \arrow{r}{\sigma_x} \arrow{d} & X \\
	\pt
	\end{tikzcd}\right),\vspace{8pt}
	\end{equation*}
	where $ G $ is a $ G $-torsor under the left multiplication of $ G $, and $ \sigma_x $ is the orbit morphism associated with the point $ x \in X $. Therefore, we have a set bijection
	\begin{align*}
	X(\CC) & \quad \xrightarrow{\sim} \quad \Obj(\cX(\CC)) \\
	x & \quad \mapsto \quad \eta_x.
	\end{align*}
	The 2-Yoneda lemma says the functor
	\begin{align*}
	\Phi: \Hom(\pt, \cX) & \to \cX(\CC) \\
	f & \mapsto f(\id_{\pt})
	\end{align*}
	is an equivalence of groupoids. The inverse $ \Phi^{-1} $ sends an object $ \eta_x $ of $ \cX(\CC) $ to a morphism $ \eta_x^*: \pt \to \cX, $ which maps a scheme $ Y $ over $ \CC, $ or its structure morphism $ Y \to \pt, $ to be precise, to an object 
	\vspace{5pt}
	\begin{equation*}
	\eta_x^*(Y \to \pt) = \left(
	\begin{tikzcd}[column sep=2.8em,row sep=2.8em]
	G \times Y \arrow{r}{\sigma_x \circ p_1} \arrow{d} & X \\
	Y
	\end{tikzcd}\right) \vspace{8pt}	
	\end{equation*}
	of $ \cX(Y) $ obtained by pulling back $ \eta_x $ along $ Y \to \pt. $
	Since $ \Phi $ is an equivalence of categories rather than an isomorphism of categories, if we define a point of $ \cX $ to be a morphism $ f: \pt \to \cX, $ then there may be more points of $ \cX $ than the objects of $ \cX(\CC), $ which are as many as the points of $ X. $ Observe that there is a set bijection of isomorphism classes of objects
	\begin{equation*}
	\Obj(\Hom(\pt, \cX))/\sim \ \to \Obj(\cX(\CC))/\sim \ , \quad [f] \mapsto \left[f(\id_{\pt})\right].
	\end{equation*}
	Two morphisms $ f_1 $ and $ f_2 $ in $ \Hom(\pt, \cX) $ are isomorphic when there is a natural isomorphism $ \alpha $ between them:
	\vspace{5pt}
	\begin{equation*}
	\begin{tikzcd}[column sep=huge, every label/.append style={font=\normalsize}]
	\pt \arrow[r, bend left=50, "f_1"{name=U, above, outer sep=2pt}]
	\arrow[r, bend right=50, "f_2"{name=D, below, outer sep=2pt}]
	& \cX
	\arrow[shorten <=10pt,shorten >=10pt, Rightarrow, from=U, to=D, "\alpha"{right, outer sep=6pt},"\cong"{left, outer sep=6pt}]
	\end{tikzcd}\vspace{8pt}
	\end{equation*} 
	Two objects $ \eta_x $ and $ \eta_y $ in $ \cX(\CC) $ are isomorphic if and only if $ x $ and $ y $ are in the same orbit, i.e., $ y = gx $ for some $ g \in G, $ which gives a commutative diagram 
	\vspace{5pt}
	\begin{equation*}
	\begin{tikzcd}[column sep=2.8em,row sep=2.8em]
	G \arrow{dr}{R_{g^{-1}}}[swap]{\cong} \arrow[bend left=40]{drr}{\sigma_x} \arrow[bend right=40]{ddr} & [-15 pt] & \\ [-15 pt] 
	& G \arrow{r}{\sigma_{gx}} \arrow{d} & X \\
	& \pt
	\end{tikzcd},\vspace{8pt}
	\end{equation*}
	where $ R_{g^{-1}}: G \to G $ is the right multiplication by $ g^{-1}. $ Therefore, we have a set bijection
	\begin{align*}
	\{G\text{-orbits in }X\} & \quad \xrightarrow{\sim} \quad \Obj(\cX(\CC))/\sim \\
	Gx & \quad \mapsto \quad [\eta_x].
	\end{align*}
	Now it's evident that how a point of $ \cX = [X/G] $ shall be defined such that it corresponds to a $ G $-orbit in $ X. $
	\begin{defn}
		Let $ \cX $ be a quotient stack. A \tb{point} of $ \cX $ is the isomorphism class $ [f] $ of a morphism $ f: \pt \to \cX $ as an object in the groupoid $ \Hom(\pt, \cX). $ A point $ [f] $ of $ \cX $ is called an \tb{orbifold point} if $ \Aut(f) $ is non-trivial; otherwise it's called a \tb{generic point}.
	\end{defn}
\end{rmk}
\begin{rmk}
	Let $ \cX = [X/G] $ be a quotient stack. A point of $ \cX $ corresponds to a $ G $-orbit $ Gx $ in $ X $ which gives a 2-commutative diagram
	\vspace{5pt}
	\begin{equation*}
	\begin{tikzcd}[column sep=2em,row sep=3em]
	Gx \arrow{r} \arrow{d} & X \arrow{d} \\
	\left[Gx/G\right] \arrow{r} & \left[X/G\right].
	\end{tikzcd}\vspace{6pt}
	\end{equation*}
	In other words, a point $ [f] $ of $ \cX $ corresponds to a zero-dimensional substack $ [Gx/G] $ of $ \cX $ which is equivalent to the classifying stack $ BG_x $ where $ G_x $ is the stabilizer of $ x. $ If $ G $ is trivial (resp. acts on $ X $ freely), then $ [X/G] $ is represented by a scheme (resp. an algebraic space) $ Y, $ and hence every zero-dimensional substack $ [Gx/G] $ is equivalent to $ \pt = (\Sch/\CC), $ i.e., the corresponding point $ [f] $ of $ \cX $ is a generic point. 
\end{rmk}
There are two scheme-theoretic notions of orbit spaces defined as follows. Let $ G $ be a linear algebraic group.
\begin{defn}[{\cite[Definition 0.5]{mumford1994geometric}}]
	Let $ X $ be a $ G $-scheme. A morphism $ q: X \to Y $ to a scheme $ Y $ is called a \tb{categorical quotient} of $ X $ (by $ G $) if:
	\begin{enumerate}[font=\normalfont,leftmargin=*]
		\item $ q $ is $ G $-equivariant (here $ Y $ is viewed as a $ G $-scheme with a trivial action).
		\item $ q $ is initial for such morphisms, i.e., if there is another $ G $-equivariant morphism $ q': X \to Y' $ to a scheme $ Y' $ with a trivial $ G $-action, then there is a unique morphism $ f: Y \to Y' $ such that $ q' = f \circ q, $ i.e., the following triangle commutes:
		\vspace{5pt}
		\begin{equation*}
		\begin{tikzcd}[column sep=1.5em,row sep=3em]
		& X \arrow{ld}[swap]{q} \arrow{rd}{q'} & \\
		Y \arrow{rr}{f} & & Y'
		\end{tikzcd}\vspace{8pt}
		\end{equation*}
	\end{enumerate}
\end{defn}
\begin{defn}[{\cite[Definition 0.6]{mumford1994geometric}}]
	Let $ X $ be a $ G $-scheme. A morphism $ q: X \to Y $ to a scheme $ Y $ is called a \tb{geometric quotient} of $ X $ if:
	\begin{enumerate}[font=\normalfont,leftmargin=*]
		\item $ q $ is $ G $-equivariant.
		\item $ q $ is surjective, and every fiber of $ q $ is a $ G $-orbit in $ X $.
		\item A subset $ U \subset Y $ is open if and only if $ q^{-1}(U) \subset X $ is open.
		\item The natural morphism $ \cO_Y \to (q_*\cO_X)^G $ is an isomorphism.
	\end{enumerate}
\end{defn}
Geometric quotients are a stronger notion than categorical quotients.
\begin{prop}[{\cite[Proposition 0.1]{mumford1994geometric}}]
	Let $ X $ be a $ G $-scheme. A geometric quotient of $ X $ is a categorical quotient of $ X, $ and hence is unique up to isomorphism.
\end{prop}
This establishes the uniqueness of a geometric quotient, but it does not always exist.
\begin{ex}
	Let $ \CC^* $ act on $ \CC $ by $ t \cdot x = tx $ for all $ t \in \CC^* $ and $ x \in \CC. $
	Then there is no geometric quotient of $ \CC $ by $ \CC^* $ since the orbit $ \CC^* \cdot 1 $ is not closed in $ \CC. $ However, the geometric quotient of $ \CC - \{0\} $ by $ \CC^* $ exists and is a point which corresponds to the orbit $ \CC^* \cdot 1. $
\end{ex}
If a $ G $-scheme $ X $ admits a geometric quotient, then it will be called ``the" geometric quotient of $ X $ by $ G $ and denoted by $ X/G. $ There are two similar notions to geometric quotients of schemes in the theory of quotient stacks.
\begin{defn}[{\cite[Definition 4.25]{edidin2013equivariant}}]
	Let $ \cX $ be a quotient stack. A morphism $ \pi: \cX \to Y $ to a scheme $ Y $ is called a \tb{coarse moduli scheme} of $ \cX $ if:
	\begin{enumerate}[font=\normalfont,leftmargin=*]
		\item $ \pi $ is initial for morphisms from $ \cX $ to schemes. 
		\item The $ \CC $-valued points of $ \cX $ are in bijection with the $ \CC $-valued points of $ Y $.
	\end{enumerate}
\end{defn}
By the definition of a coarse moduli scheme, it is unique up to isomorphism if it exists. It is straightforward to check the following:
\begin{prop}\label{prop_geom_quot_is_coarse_moduli_scheme}
	Let $ \cX = [X/G] $ be a quotient stack. Then a geometric quotient of $ X $ by $ G $ is a coarse moduli scheme of $ \cX. $
\end{prop}
\begin{defn}
	Let $ \cX $ be a quotient stack. A morphism $ \pi: \cX \to Y $ to an algebraic space $ Y $ is called a \tb{coarse moduli space} of $ \cX $ if:
	\begin{enumerate}[font=\normalfont,leftmargin=*]
		\item $ \pi $ is initial for morphisms from $ \cX $ to algebraic spaces. 
		\item The $ \CC $-valued points of $ \cX $ are in bijection with the $ \CC $-valued points of $ Y $.
	\end{enumerate}
\end{defn}
It's proven in \cite{keel1997quotients} that if $ \cX $ is an algebraic stack locally of finite presentation over a locally noetherian\footnote{The noetherian assumption is removed in \cite{conrad2005keel}.} scheme $ S $ and has a finite inertia, i.e., the canonical morphism $ I\cX \to \cX $ is finite, then $ \cX $ admits a coarse moduli space $ \cX \to Y, $ which is a proper morphism. Therefore, a separated DM stack locally of finite type over a field always has a coarse moduli space.
\begin{defn}
	An algebraic stack $ \cX $ over a scheme $ S $ is \tb{separated} if the diagonal $ \Delta: \cX \to \cX \times_S \cX $ is proper. 
\end{defn}

\begin{rmk}
	The separatedness of an algebraic stack is a generalization of that of a scheme: for a scheme $ X $ over $ S, $ the condition that the diagonal $ \Delta: X \to X \times_S X $ is proper is the same as that the diagonal is a closed immersion.
\end{rmk}

The separatedness of a quotient stack can be checked via the group action.
\begin{prop}[{\cite[Proposition 4.17]{edidin2013equivariant}}]\label{prop_separated}
	A quotient stack $ [X/G] $ is separated if and only if $ G $ acts on $ X $ properly, i.e., the action morphism $\alpha: G \times X \to X \times X $ is proper.
\end{prop}
The following result is immediate.
\begin{cor}
	Let $ [X/G] $ be a separated quotient DM stack. Then $ [X/G] $ is a separated DM stack of finite type.
\end{cor}
\begin{proof}
	By Proposition \ref{prop_separated}, the action morphism $\alpha: G \times X \to X \times X $ is proper. Let $ x $ be a point of $ X. $ The stabilizer $ G_x $ is the pullback of $ \alpha $ along the embedding $ \pt \to X \times X $ with an image point $ (x,x) $ in $ X \times X, $ so $ G_x $ is finite. The stabilizer $ G_x $ is also reduced since every group scheme over a field of characteristic zero is reduced. Therefore, $ [X/G] $ is a separated DM stack. It is of finite type over $ \CC $ because $ X $ is assumed so in Definition \ref{defn_quotient_stack}.
\end{proof}
\begin{ex}
	Let $ [X/G] $ be a quotient stack where $ G $ is a finite group. Then $ [X/G] $ is separated by \cite[Lemma 2.1]{geraschenko2015torus} and the assumption that $ X $ is separated in Definition \ref{defn_quotient_stack}. In particular, the classifying stack $ BG = [\pt /G] $ is separated.
\end{ex}
Consider a quotient stack $ \cX $ over $ \CC $. If a scheme $ Y $ is the coarse moduli space of $ \cX, $ then it is also a coarse moduli scheme of $ \cX. $ But the other direction is not always true, i.e., it can happen that a scheme $ Y $ is the coarse moduli scheme of $ \cX $ but not a coarse moduli space of $ \cX $ (see Example 4.30 in \cite{edidin2013equivariant}). Such phenomenon does not occur when $ \cX $ is separated. We combine Proposition 4.29 and Corollary 4.33 in \cite{edidin2013equivariant} in the following proposition.
\begin{prop}\label{prop_coarse_moduli_scheme_is_geom_quot}
	If a scheme $ Y $ is the coarse moduli scheme of a separated quotient DM stack $ \cX = [X/G], $ then $ Y $ is both the coarse moduli space of $ \cX $ and the geometric quotient of $ X $ by $ G. $
\end{prop}
Propositions \ref{prop_geom_quot_is_coarse_moduli_scheme} and \ref{prop_coarse_moduli_scheme_is_geom_quot} imply the following:
\begin{cor}\label{cor_geom_quot_is_coarse_moduli_space}
	Let $ \cX = [X/G] $ be a separated quotient DM stack. The geometric quotient $ X/G $ of $ X $ by $ G $, if it exists, is both the coarse moduli scheme and the coarse moduli space of $ \cX $.
\end{cor}

\begin{ex}
	The classifying stack $ BG = [\pt/G] $ for a linear algebraic group $ G $ has a coarse moduli space $ \pt, $ which is also the geometric quotient of $ \pt $ by $ G $.
\end{ex}
\begin{ex}
	Consider the separated quotient DM stack $ [\CC^2/\mu_2] $ where the action of $ \mu_2 $ on $ \CC^2 $ is generated by the involution
	\begin{equation*}
	(x,y) \mapsto (-x,-y).
	\end{equation*}
	Then the invariant subring $ \CC[x,y]^{\mu_2} = \CC[x^2,xy,y^2] $ of $ \CC[x,y] $ gives a geometric quotient
	\begin{equation*}
	\CC^2 \to \CC^2/\mu2
	\end{equation*}
	of $ \CC^2 $ by $ \mu_2 $ where $ \CC^2/\mu2 = \Spec \CC[x^2,xy,y^2] $ is a cone in $ \CC^3, $ and also gives a coarse moduli scheme
	\begin{equation*}
	[\CC^2/\mu_2] \to \CC^2/\mu2
	\end{equation*}
	of $ [\CC^2/\mu_2]. $
\end{ex}
In general, it's a difficult question whether a $ G $-scheme $ X $ has a geometric quotient or whether a quotient stack $ [X/G] $ has a coarse moduli scheme. Fortunately, there is a simple criterion when $ G $ is a finite group.
\begin{prop}[{\cite[Chapter III, Theorem 1]{mumford2008abelian}}]\label{prop_geom_quot}
	Let $ G $ be a finite group acting on a scheme $ X $ such that every orbit is contained in an affine open subset of $ X. $ Then there is a geometric quotient $ \pi: X \to X/G $ of $ X $ by $ G $ where the morphism $ \pi $ is finite surjective, and the scheme $ X/G $ is unique up to isomorphism.
\end{prop}
We will need the notion of properness of quotient stacks, e.g., in the definition of the orbifold Euler characteristic.
\begin{defn}[{\cite[Definition 4.21]{edidin2013equivariant}}]
	A quotient stack $ \cX $ is \tb{proper} if it is separated and there is a finite surjective morphism $ Y \to \cX $ where $ Y $ is a proper scheme.
\end{defn}
\begin{rmk}
	A finite surjective morphism $ Y \to [X/G] $ means a finite surjective $ G $-equivariant morphism $ E \to X $ and a $ G $-torsor $ E \to Y $ in the following 2-pullback diagram:
	\vspace{2pt}
	\begin{equation*}
	\begin{tikzcd}[column sep=3em,row sep=3em]
	E \arrow{r} \arrow{d} & X \arrow{d} \\
	Y \arrow{r} & \left[X/G\right]
	\end{tikzcd}\vspace{8pt}
	\end{equation*}
	Such a scheme $ Y $ is also called a finite cover of the stack $ [X/G] $ and always exists for DM stacks by \cite[Theorem 2.7]{edidin2001brauer}. 
	The properness of a quotient stack is a generalization of that of a scheme: when $ G $ is trivial, $ [X/G] $ is proper if and only if $ X $ is separated and there exists a finite surjective morphism $ Y \to X $ from a proper scheme $ Y $; these two conditions are equivalent to $ X $ being proper since the image of a proper scheme in a separated scheme is proper.
\end{rmk}
\begin{ex}
	let $ G $ be a finite group acting on a proper scheme $ X $. Then the quotient stack $ [X/G] $ is a proper DM stack since the canonical morphism $ X \to [X/G] $ is a finite surjective morphism.
\end{ex}
\begin{ex}
	The classifying stack $ BG = [\pt/G] $ for a finite group $ G $ is a proper DM stack.
\end{ex}

We will also need a notion of projectivity of quotient stacks, which is essential for the construction of moduli spaces of sheaves on quotient stacks. 
\begin{defn}\label{defn_projective_stack_over_C}
	A separated quotient DM stack $ \cX $ is called a \tb{projective} (resp. \tb{quasi-projective}) stack if its coarse moduli space is a projective (resp. quasi-projective) scheme.
\end{defn}
As expected, projectivity is stronger than properness for quotient stacks.
\begin{lem}
	A projective stack is proper.
\end{lem}
\begin{proof}
	Let $ \cX $ be a projective stack. Then it is a separated DM stack. Choose a finite cover $ Y \to \cX $. The composition $ Y \to \cX \to M $ is a proper morphism since it's a composition of two proper morphisms. Since $ M $ is a projective scheme and hence proper over $ \CC $, the scheme $ Y $ must be proper over $ \CC $.
\end{proof}
\begin{ex}
	Let $ U = \CC^{n+1} - \{0\}. $ Let $ X = [U/\CC^*] $ where $ \CC^* $ acts on $ U $ by
	\begin{equation*}
	t(x_0, \dots, x_n) = (t^{a_0}x_0, \dots, t^{a_n}x_n)
	\end{equation*}
	for some positive integers $ a_0, \dots, a_n. $ This is a proper action. Therefore, $ X $ is a separated DM stack, known as a \tb{weighted projective stack} and often denoted by $ \PP(a_0, \dots, a_n). $ 
	The morphism
	\begin{equation*}
	U \to U, \quad (x_0, \dots, x_n) \mapsto (x_0^{a_0}, \dots, x_n^{a_n})
	\end{equation*}
	is $ \CC^* $-equivariant and induces a finite surjective morphism
	\begin{equation*}
	\PP^n \to \PP(a_0, \dots, a_n).
	\end{equation*}
	Therefore, $ \PP(a_0, \dots, a_n) $ is a proper DM stack. We also have the corresponding weighted projective space 
	$$ P(a_0, \dots, a_n) = \Proj \CC[x_0, \dots, x_n] $$ 
	where $ x_i $ has degree $ a_i $ for $ i = 0, \dots, n. $ The natural morphism
	\begin{equation*}
	\PP(a_0, \dots, a_n) \to P(a_0, \dots, a_n)
	\end{equation*}
	is a coarse moduli space of $ \PP(a_0, \dots, a_n) $, which shows $ \PP(a_0, \dots, a_n) $ is also a projective stack.
	When $ n = 0, $ we have $ \PP(k) = B\mu_k $ and $ P(k) = \pt.$ When $ n = 1, $ we have
	\begin{equation*}
	P(a_0, a_1) \cong \PP^1.
	\end{equation*}
	Cyclic quotient singularities occur in coarse moduli spaces when $ n > 1. $ For example, $ P(1,1,2) $ is a singular surface isomorphic to a quadratic surface in $ \PP^3 $ with homogeneous coordinates $ [y_0, y_1, y_2, y_3] $ whose defining equation is
	\begin{equation*}
	y_1y_3 - y_2^2 = 0,
	\end{equation*} 
	as shown in \cite[Example 10.27]{harris1992algebraic}.
\end{ex}
Now we prove a fundamental fact on quotients by finite groups which should be known to experts but hard to find in the literature.  
\begin{thm}
	Let $ G $ be a finite group acting on a quasi-projective scheme $ X. $ Then there is a geometric quotient $ X/G $ of $ X $ by $ G $, which is also the coarse moduli space of the quotient stack $ [X/G]. $ If $ X $ is proper (or projective), so are $ [X/G] $ and $ X/G. $ 
\end{thm}
\begin{proof}
	Take an orbit $ Gx $ of $ G $ in $ X. $ Since $ Gx $ is finite, it is contained in an affine open subset of $ X $ by \cite[Proposition 3.3.36]{liu2006algebraic}. The existence of a geometric quotient $ X/G $ follows immediately from Proposition \ref{prop_geom_quot}. By Corollary \ref{cor_geom_quot_is_coarse_moduli_space}, $ X/G $ is also the coarse moduli space of $ [X/G] $. If $ X $ is proper, then the quotient stack $ [X/G] $ is proper since it's covered by $ X $, which forces its coarse moduli space $ X/G $ to be proper. If $ X $ is projective, then $ [X/G] $ is projective by definition, and $ X/G $ is projective by an explicit construction given on page 127 in \cite{harris1992algebraic}.
\end{proof}
\begin{ex}
	Let $ G $ be a finite subgroup of $ \SL(2,\CC) $. Then the quotient stack $ [\CC^2/G] $ is a separated DM stack but not proper. Its coarse moduli space is the same as the geometric quotient $ \CC^2/G $ which is a surface with ADE singularities.
\end{ex}
\begin{ex}\label{ex_K3_Fermat}
	Let $ X $ be a $ K3 $ surface in $ \PP^3 $ defined by the Fermat equation
	\begin{equation*}
	x_0^4 + x_1^4 + x_2^4 + x_3^4 = 0.
	\end{equation*}
	Let $ \mu_2 = \{1, \sigma\}$ act on $ X $ by
	\begin{equation*}
	\sigma \cdot [x_0,x_1,x_2,x_3] = [-x_0,-x_1,x_2,x_3].
	\end{equation*}
	This is a symplectic involution known as a \tb{Nikulin involution}. The homogeneous coordinate ring of $ X $ is 
	$$ S(X) = \CC[x_0,x_1,x_2,x_3]/(x_0^4 + x_1^4 + x_2^4 + x_3^4), $$ 
	whose invariant subring under the action of $ \mu_2 $ is
	\begin{equation*}
	S(X)^{\mu_2} = \CC[x_0^2,x_0x_1,x_1^2,x_2,x_3]/(x_0^4 + x_1^4 + x_2^4 + x_3^4).
	\end{equation*}
	The graded ring $ S(X)^{\mu_2} $ is not generated in degree one. Dropping its degree one component, we obtain a ring
	\begin{equation*}
	S(Y) = \CC[x_0^2,x_0x_1,x_1^2,x_2^2,x_2x_3,x_3^2]/(x_0^4 + x_1^4 + x_2^4 + x_3^4),
	\end{equation*}
	which is the homogeneous coordinate ring of a projective scheme $ Y $. By \cite[Exercise 10.25]{harris1992algebraic}, $ Y $ is the geometric quotient $ X/\mu_2. $ Note that
	\begin{equation*}
	Y = \Proj \frac{\CC[y_0, y_1, y_2, y_3, y_4, y_5]}{(y_0y_2 - y_1^2, y_3y_5 - y_4^2, y_0^2+y_2^2+y_3^2+y_5^2)}
	\end{equation*}
	is a surface in $ \PP^5 $ with eight singular points
	\begin{equation*}
	[1,\zeta_i,\zeta_i^2,0,0,0], \ [0,0,0,1,\zeta_i,\zeta_i^2], \quad \text{for} \ i = 1,2,3,4,
	\end{equation*}
	where $ \zeta_1,\zeta_2,\zeta_3,\zeta_4 $ are fourth roots of $ -1. $ The quotient stack $ [X/\mu_2] $ is a projective stack since its coarse moduli space $ X/\mu_2 $ is projective.
\end{ex}

\subsection{Quasi-coherent sheaves on quotient stacks}
In this section we study quasi-coherent sheaves on quotient stacks. Let $ G $ be a linear algebraic group. We will define sheaves on a stack $ \cX $ and $ G $-equivariant sheaves on a $ G $-scheme $ X $, and discuss their relation when $ \cX = [X/G]. $ We will also define $ G $-invariant sheaves on $ G $-schemes and discuss when they can lift to $ G $-equivariant sheaves. At the end of the section we will discuss various ways to construct sheaves on a quotient stack.

The following is an alternative form of Definition 2.1 in \cite[Chapter XIII]{arbarello2011geometry}.
\begin{defn}
	Let $ \cX $ be a stack over $ \CC $ or more generally, a category fibered in groupoids over $ (\Sch/\CC). $ A quasi-coherent sheaf $ \cE $ on $ \cX $ consists of the following data:
	\begin{enumerate}[font=\normalfont,leftmargin=*]
		\item For each object $ \eta $ in $ \cX $ over a scheme $ Y $: a quasi-coherent sheaf $ \cE(\eta) $ on $ Y. $
		\item For each arrow $ A_1: \eta_1 \to \eta_2 $ in $ \cX $ over a morphism $ a_1: Y_1 \to Y_2 $ in $ (\Sch/\CC) $: an isomorphism
		\begin{equation*}
		\cE(A_1): \cE(\eta_1) \xrightarrow{\sim} a_1^*\cE(\eta_2)
		\end{equation*}
		of quasi-coherent sheaves on $ Y_1, $ which satisfies the \tb{cocycle condition}, i.e., for a composition of arrows
		\begin{equation*}
		\eta_1 \xrightarrow{A_1} \eta_2 \xrightarrow{A_2} \eta_3 \quad \text{over} \quad Y_1 \xrightarrow{a_1} Y_2 \xrightarrow{a_2} Y_3,
		\end{equation*}
		we have
		\begin{equation*}
		\cE(A_2 \circ A_1) = a_1^* \cE(A_2) \circ \cE(A_1)
		\end{equation*}
		(here $ ``=" $ means equal up to a canonical isomorphism) in the following diagram of isomorphisms of quasi-coherent sheaves on $ Y_1: $
		\vspace{5pt}
		\begin{equation*}
		\begin{tikzcd}[column sep=1.5em,row sep=4em]
		\cE(\eta_1) 
		\ar[rr, "\cE(A_2 \circ A_1)", "\simeq"{below}] 
		\ar[dr, "\cE(A_1)"{below left}, "\simeq"{above, sloped}] 
		& & (a_2\circ a_1)^* \cE(\eta_3) \\
		& a_1^*\cE(\eta_2) 
		\arrow[ur, "a_1^*\cE(A_2)"{below right}, "\simeq"{above,sloped}] &,
		\end{tikzcd}\vspace{8pt}
		\end{equation*}
		where $ (a_2\circ a_1)^* \cE(\eta_3) \cong a_1^*a_2^*\cE(\eta_3) $ via the canonical pullback isomorphism.
	\end{enumerate}
	
	\vspace{8pt}
	If $ \cX $ is a locally noetherian algebraic stack over $ \CC, $ then coherent sheaves (resp. vector bundles) on $ \cX $ are defined by replacing ``quasi-coherent" with ``coherent" (resp. ``locally free") everywhere above and restricting the category $ (\Sch/\CC) $ to the category of locally noetherian schemes over $ \CC $. Let sheaves be either quasi-coherent sheaves, coherent sheaves, or vector bundles. A morphism $ \Phi: \cE \to \cF $ of sheaves on $ \cX $ is a collection of morphisms $ \Phi(\eta): \cE(\eta) \to \cF(\eta) $ of sheaves on $ Y $ for each object $ \eta $ in $ \cX $ over a scheme $ Y $ such that for each arrow
	\vspace{5pt}
	\begin{equation*}\label{eq_sheaf_morphism_on_stack}
	\begin{tikzcd}[column sep=3em,row sep=4em]
	\eta_1 \arrow[d,"A_1"{right, outer sep = 3pt}] \\
	\eta_2
	\end{tikzcd} 
	\quad \text{over} \quad
	\begin{tikzcd}[column sep=3em,row sep=4em]
	Y_1 \arrow[d,"a_1"{right, outer sep = 3pt}] \\
	Y_2
	\end{tikzcd}
	\quad, \quad \quad 
	\begin{tikzcd}[column sep=3em,row sep=4em]
	\cE(\eta_1) 
	\arrow[d, "\cE(A_1)"{left},  "\simeq"{above, sloped}] 
	\arrow[r, "\Phi(\eta_1)"{}] 
	&[-3pt] \cF(\eta_1) \arrow[d, "\cF(A_1)"{right},  "\simeq"{below, sloped}] \\
	a_1^*\cE(\eta_2) \arrow[r, "a_1^*\Phi(\eta_2)"{}] 
	& a_1^*\cF(\eta_2)
	\end{tikzcd}\vspace{8pt}
	\end{equation*}
	commutes as a diagram of morphisms of sheaves on $ Y_1 $.
\end{defn}
\begin{rmk}[Sheaves on stacks are a generalization of sheaves on schemes.]
	A sheaf $ \cE $ on a stack $ \cX $ contains huge amount of data which looks like abstract nonsense in category theory. But it's really a generalization of a sheaf $ E $ on a scheme $ X $ which gives a pullback $ f^*E $ for every morphism $ f: Y \to X. $ Suppose $ \cX = (\Sch/X) $ is represented by a scheme $ X. $ Then an object of $ (\Sch/X) $ over a scheme $ Y $ is precisely a morphism $ f: Y \to X, $ and an arrow $ A_1: f_1 \to f_2 $ in $ (\Sch/X) $ over a morphism $ a_1: Y_1 \to Y_2 $ in $ (\Sch/\CC) $ is exactly a commutative triangle:
	\vspace{5pt}
	\begin{equation}\label{eq_arrow_A1}
	\begin{tikzcd}[column sep=1.5em,row sep=3em]
	Y_1 \arrow{rr}{a_1} \arrow{dr}{f_1} & & Y_2 \arrow{dl}[swap]{f_2} \\
	& X &
	\end{tikzcd}\vspace{8pt}
	\end{equation}
	Fix a sheaf $ E $ on $ X. $ We define a sheaf $ \cE $ on $ (\Sch/X) $ as follows.  For a morphism $ f: Y \to X, $ put
	\begin{equation*}
	\cE(f) = f^*E.
	\end{equation*}
	For an arrow $ A_1 $ given by triangle (\ref{eq_arrow_A1}), let
	\begin{equation*}
	\cE(A_1): f_1^*E = (f_2 \circ a_1)^*E \xrightarrow{\sim} a_1^*f_2^*E
	\end{equation*}
	encode the canonical isomorphism of sheaves on $ Y_1. $ It follows that for a compositions of arrows
	\vspace{5pt}
	\begin{equation*}\label{eq_composition_of_arrows}
	\begin{tikzcd}[column sep=3em,row sep=3em]
	Y_1 \arrow{r}{a_1} \arrow{rd}{f_1} & Y_2 \arrow{r}{a_2} \arrow{d}[swap]{f_2} & Y_3 \arrow{ld}[swap]{f_3} \\
	& X &,
	\end{tikzcd}\vspace{8pt}
	\end{equation*}
	we have
	\vspace{5pt}
	\begin{equation*}
	\begin{tikzcd}[column sep=1.5em,row sep=4em]
	f_1^*E 
	\ar[rr, "\cE(A_2 \circ A_1)", "\simeq"{below}] 
	\ar[dr, "\cE(A_1)"{below left}, "\simeq"{above, sloped}] 
	& & (a_2\circ a_1)^* f_3^*E \\
	& a_1^*f_2^*E 
	\arrow[ur, "a_1^*\cE(A_2)"{below right}, "\simeq"{above,sloped}] &,
	\end{tikzcd}\vspace{8pt}
	\end{equation*}
	where $ (a_2\circ a_1)^* f_3^*E \cong a_1^*a_2^*f_3^*E $ canonically.  We see that a sheaf $ \cE $ on $ \cX = (\Sch/X) $ corresponds to a sheaf $ E $ on $ X. $ Take a morphism $ \vphi: E \to F $ of sheaves on $ X. $ Let $ \cE $ and $ \cF $ be the corresponding sheaves on $ (\Sch/X). $ Pulling back $ \vphi $ along triangle (\ref{eq_arrow_A1}) yields a commutative diagram:
	\begin{equation*}
	\begin{tikzcd}[column sep=4em,row sep=4em]
	f_1^*E \arrow[r, "\cE(A_1)", "\simeq"{below}] 
	\arrow[d, "f_1^*\vphi"{left, outer sep = 3pt}]
	& a_1^*f_2^*E \arrow[d, "a_1^*f_2^*\vphi"{outer sep = 3pt}] \\
	f_1^*F \arrow[r, "\cF(A_1)", "\simeq"{below}] & a_1^*f_2^*F
	\end{tikzcd}\vspace{8pt}
	\end{equation*}
	This defines a morphism $ \Phi: \cE \to \cF $ of sheaves on $ (\Sch/X) $ with $ \Phi(f) = f^*E $ for a morphism $ f: Y \to X. $ Therefore, the definition above generalizes the notion of a sheaf to stacks.
	
	Let $ \cX $ be a stack and $ Y $ be a scheme. Consider an object $ \eta $ in $ \cX(Y). $ It determines a morphism $ f_\eta: Y \to \cX, $ or a functor $ f_\eta: (\Sch/Y) \to \cX, $ to be precise. For a sheaf $ \cE $ on $ \cX, $ the sheaf $ \cE(\eta) $ on $ Y $ can be viewed as the pullback of $ \cE $ to $ Y, $ while for a morphism $ \Phi: \cE \to \cF $ of sheaves on $ \cX, $ the morphism $ \Phi(\eta): \cE(\eta) \to \cF(\eta) $ of sheaves on $ Y $ can be seen as the pullback of $ \Phi $ to $ Y, $ and we denote them by 
	\begin{equation*}
	f_\eta^* \cE = \cE(\eta) \quad \text{and} \quad f_\eta^* \Phi = \Phi(\eta)
	\end{equation*}
	as usual. Therefore, we can work with sheaves on stacks just like those on schemes.
\end{rmk}
\begin{ex}
	There is a structure sheaf $ \cO_\cX $ on a stack $ \cX, $ which assigns every object $ \eta $ of $ \cX $ over a scheme $ Y $ the structure sheaf $ \cO_Y, $ and assigns every arrow $ A_1: \eta_1 \to \eta_2 $ over a morphism $ a_1: Y_1 \to Y_2 $ the identity $ \id: \cO_{Y_1} \to a_1^*\cO_{Y_2} = \cO_{Y_1}. $ 
\end{ex}
Let $ X $ be a $ G $-scheme. Now we define $ G $-equivariant sheaves on $ X $ and $ G $-equivariant morphisms between them.
\begin{defn}\label{defn_equiv_sheaves}
	Let $ p_2: G \times X \to X $ be the projection to $ X, $ and let $ \sigma: G \times X \to X $ be the action of $ G $ on $ X. $ A \tb{$ G $-equivariant structure} of a quasi-coherent sheaf $ E $ on $ X $ is an isomorphism
	\begin{equation*}
	\phi: p_2^*E \xrightarrow{\sim} \sigma^*E
	\end{equation*}
	of quasi-coherent sheaves on $ G\times X $ which satisfies the \tb{cocycle condition}:
	\begin{equation*}
	(\mu \times \id_X)^* \phi = (\id_G \times \sigma)^* \phi \circ p_{23}^* \phi,
	\end{equation*}
	which is an identity of isomorphisms of quasi-coherent sheaves on $ G \times G \times X $ in the diagram
	\vspace{5pt}
	\begin{equation}\label{eq_equiv}
	\begin{tikzcd}[column sep=1.5em,row sep=3em]
	p_3^*E \arrow[rr, "(\mu \times \id_X)^* \phi", "\simeq"{below}] 
	\arrow[dr, "p_{23}^* \phi"{below left}, "\simeq"{above, sloped}] & & s^*E \\
	& r^*E \arrow[ur, "(\id_G \times \sigma)^* \phi"{below right}, "\simeq"{above, sloped}] &,
	\end{tikzcd}\vspace{8pt}
	\end{equation}
	where $ p_3, r, s $ are three morphisms from $ G \times G \times X $ to $ X $ defined by
	\begin{align*}
	p_3:\quad (g,h,x) & \quad \longmapsto \quad x, \\
	r: 	\quad (g,h,x) & \quad \longmapsto \quad hx, \\
	s: 	\quad (g,h,x) & \quad \longmapsto \quad ghx,
	\end{align*}
	and $ \mu \times \id_X,p_{23},\id_G \times \sigma $ are three morphisms from $ G \times G \times X $ to $ G \times X $ defined by
	\begin{align*}
	\mu \times \id_X:		\quad (g,h,x) & \quad \longmapsto \quad (gh,x), \\
	p_{23}: 				\quad (g,h,x) & \quad \longmapsto \quad (h,x), \\
	\id_G \times \sigma: 	\quad (g,h,x) & \quad \longmapsto \quad (g,hx).
	\end{align*}
	A quasi-coherent sheaf $ E $ on  $ X $ with a $ G $-equivariant structure $ \phi $ is called a \tb{$ G $-equivariant quasi-coherent sheaf} on $ X $ and is denoted by $ (E,\phi) $.
	
	\vspace{5pt}
	\tb{$ G $-equivariant coherent sheaves} on $ X $ are defined similarly. A \tb{$ G $-equivariant morphism} $ f: (E,\phi) \to (F,\psi) $
	is a morphism $ f: E \to F $ which is $ G $-equivariant, i.e., $ \psi \circ p_2^*f = \sigma^*f \circ \phi $ in the following diagram:
	\vspace{5pt}
	\begin{equation}\label{eq_equiv_morphism}
	\begin{tikzcd}[column sep=3em,row sep=3em]
	p_2^*E \arrow{r}{p_2^*f} \arrow[d, "\phi"{left, outer sep = 2pt}, "\simeq"{above, sloped}] 
	& p_2^*F \arrow[d, "\psi"{right, outer sep = 2pt}, "\simeq"{below, sloped}] \\
	\sigma^*E \arrow{r}{\sigma^*f} & \sigma^*F
	\end{tikzcd}\vspace{8pt}
	\end{equation}
\end{defn}
\begin{rmk}[Interpretation of $ G $-equivariant sheaves along fibers.]
	In Definition \ref{defn_equiv_sheaves}, the existence of an isomorphism $ \phi: p_2^*E \xrightarrow{\sim} \sigma^*E $ for a sheaf $ E $ on $ X $ implies that there exist isomorphisms 
	\begin{equation*}
	\phi^g: E \xrightarrow{\sim} g^*E
	\end{equation*}
	of sheaves on $ X $ for all $ g \in G $. Here $ \phi^g = i_g^*\phi $ (up to canonical isomorphism) is the pullback of $ \phi $ along the inclusion 
	\begin{equation*}
	i_g: X \to G \times X, \quad x \mapsto (g,x).
	\end{equation*}
	Moreover, the cocycle condition of $ \phi $ implies that the isomorphisms $ \phi^g $ for all $ g \in G $ satisfy the conditions
	\begin{equation*}\label{eq_2_cocycle}
	\phi^{gh} = h^*\phi^g \circ \phi^h
	\end{equation*}
	for all $ g,h \in G $, which is called the \tb{2-cocycle condition}. Indeed, pulling back diagram (\ref{eq_equiv}) along the inclusion
	\begin{equation*}
	i_{g,h}: X \to G \times G \times X, \quad x \mapsto (g,h,x) 
	\end{equation*}
	gives a commutative diagram
	\vspace{5pt}
	\begin{equation}\label{eq_equiv_local}
	\begin{tikzcd}[column sep=1.2em,row sep=3em]
	E \arrow[rr, "\phi^{gh}", "\simeq"{below, sloped}]
	\arrow[dr, "\phi^h"{below left}, "\simeq"{above, sloped}] & & (gh)^*E \\
	& h^*E \arrow[ur, "h^*\phi^g"{below right}, "\simeq"{above, sloped}] &,
	\end{tikzcd}\vspace{8pt}
	\end{equation}
	where $ (gh)^*E \cong h^*g^*E $ canonically.
	For a $ G $-equivariant morphism $  f: (E,\phi) \to (F,\psi) $, pulling back diagram (\ref{eq_equiv_morphism}) along the inclusion $ i_g: X \to G \times X $ gives a commutative diagram
	\vspace{5pt}
	\begin{equation*}\label{eq_equiv_morphism_local}
	\begin{tikzcd}[column sep=3em,row sep=3em]
	E \arrow{r}{f} \arrow[d, "\phi^g"{left}, "\simeq"{above, sloped}]
	& F \arrow[d, "\psi^g"{right}, "\simeq"{below, sloped}] \\
	g^*E \arrow{r}{g^*f} & g^*F.
	\end{tikzcd}\vspace{8pt}
	\end{equation*}
	For each point $ x \in X $, pulling back diagram (\ref{eq_equiv_local}) along the inclusion
	\begin{equation*}
	i_x: x \into X
	\end{equation*}
	gives a commutative diagram of vector spaces
	\vspace{5pt}
	\begin{equation*}\label{eq_equivariant_fiber_isom}
	\begin{tikzcd}[column sep=1.2em,row sep=3em]
	E_x \arrow[rr, "i_x^*\phi^{gh}", "\simeq"{below}] 
	\arrow[dr, "i_x^*\phi^h"{below left}, "\simeq"{above, sloped}] & & E_{ghx} \\
	& E_{hx} \arrow[ur, "i_{hx}^*\phi^g"{below right}, "\simeq"{above, sloped}] &,
	\end{tikzcd}\vspace{8pt}
	\end{equation*}
	where $ E_x = i_x^*E $ is the fiber of $ E $ at $ x $, etc. For every point $ x \in X $, there is a linear representation
	\begin{equation*}
	\phi_x: G_x \to \GL(E_x)
	\end{equation*}
	of the stabilizer $ G_x $ given by $ \phi_x(g) = i_x^*\phi^g = i_{g,x}^*\phi $ where $ i_{g,x}: x \into G \times X $ is the inclusion $ x \mapsto (g,x) $.
\end{rmk}
\begin{rmk}[Interpretation of $ G $-equivariant structures on vector bundles]
	Consider a $ G $-equivariant vector bundle $ (V, \phi) $ on $ X $. Let the same letter $ V $ denote its total space $ \Spec \Sym V^\vee $. Then the $ G $-equivariant structure $ \phi $ is the same as an action 
	\begin{equation*}
	\phi: G \to \Aut(V)
	\end{equation*}
	such that the projection $ \pi: V \to X $ is $ G $-equivariant, i.e., $ \pi \circ \phi(g) = g \circ \pi.$
	A $ G $-equivariant morphism $ f: (V, \phi) \to (W, \psi) $ between two $ G $-equivariant vector bundles is the same as a $ G $-equivariant morphism $ f: V \to W $
	of total spaces of vector bundles, i.e., a commutative diagram
	\vspace{5pt}
	\begin{equation}\label{eq_equiv_morphism_Vbundles}
	\begin{tikzcd}[column sep=3.3em,row sep=3.3em]
	V \arrow{r}{f} \arrow[d, "\phi(g)"{left}, "\simeq"{above, sloped}] & W \arrow[d, "\psi(g)"{right}, "\simeq"{below, sloped}] \\
	V \arrow{r}{f} & W
	\end{tikzcd}\vspace{8pt}
	\end{equation}
	for all $ g $ in $ G. $ Note that we don't require the morphism $ f $ to preserve fibers in diagram (\ref{eq_equiv_morphism_Vbundles}).
\end{rmk}
\begin{ex}
	The structure sheaf $ \cO_X $ carries a canonical $ G $-equivariant structure $ \phi_0 $, i.e., a $ G $-action
	\begin{equation*}
	\phi_0: G \to \Aut(X \times \CC), \quad g \mapsto (\phi_0(g): (x,v) \mapsto (gx,v))
	\end{equation*}
	on the total space $ X \times \CC $ of $ \cO_X. $
\end{ex}
\begin{rmk}
	We emphasize that being $ G $-equivariant is not a property of a sheaf $ E $ on $ X $, but rather an additional structure we put on $ E $ which is not always possible depending on both the action of $ G $ on $ X $ and the properties of $ E. $
\end{rmk}

The $ G $-action on the scheme $ X $ induces a right action
on the set of isomorphism classes in $ \Coh(X) $ by pullback, since 
$$ (gh)^*E \cong h^*g^*E $$
for all $ E $ in $ \Coh(X), $ and all $ g, h $ in $ G. $
\begin{defn}
	A sheaf $ E $ on $ X $ is called \tb{$ G $-invariant} if 
	$$ E \cong g^*E \quad \text{for all $ g $ in $ G $.} $$ 
	A $ G $-invariant sheaf $ E $ on $ X $ is said to be \tb{$ G $-linearizable} if $ E $ admits a $ G $-equivariant structure. 
\end{defn}
\begin{defn}
	A sheaf $ E $ on a scheme $ X $ is said to be \tb{simple} if its endomorphism ring $ \End(E) \cong \CC. $
\end{defn}
A $ G $-linearizable sheaf is $ G $-invariant, but the converse is not always true. We summarize \cite[Lemma 1]{ploog2007equivariant} and \cite[Propositions 2.4, 4.3 and 4.4]{kool2011fixed} in the following proposition.
\begin{prop}
	Let $ X $ be a projective scheme under an action of a linear algebraic group $ G. $ Let $ E $ be a $ G $-invariant simple sheaf on $ X. $ Then $ E $ is $ G $-linearizable if it satisfies either of the following conditions:
	\begin{enumerate}[font=\normalfont,leftmargin=2em]
		\item $ G $ is finite and $ H^2(G, \CC^*) = 0. $
		\item $ G $ is connected and $ \Pic(G) = 0. $
	\end{enumerate}
	Furthermore, if $ E $ is $ G $-linearizable, then the character group $ \widehat{G} $ acts freely and transitively on the set of $ G $-equivariant structures on $ E. $
\end{prop}
\begin{proof}
	Consider the second group cohomology $ H^2(G, \CC^*) $ of $ G $ with coefficients in $ \CC^* $ on which $ G $ acts trivially. We first recall an argument of Ploog in the proof of \cite[Lemma 1]{ploog2007equivariant}, which proves that the $ G $-invariant simple sheaf $ E $ determines a class $ [E] $ in $ H^2(G, \CC^*) $ such that $ E $ satisfies the $ 2 $-cocycle condition if and only if the class $ [E] $ is trivial. Note that Ploog's argument is valid for all linear algebraic groups, not only finite groups. Since $ E $ is $ G $-invariant, there exists a family 
	$$ \phi^g: E \xrightarrow{\sim} g^*E $$
	of isomorphisms of sheaves on $ X $ for all $ g $ in $ G. $
	Since $ E $ is simple, for all $ g $ and $ h $ in $ G, $ there is a unit
	\begin{equation*}
	u_{g,h} = (\phi^{gh})^{-1} \circ h^*\phi^g \circ \phi^h \in \CC^*.
	\end{equation*}
	The assignment $ (g,h) \mapsto u_{g,h} $ gives a $ 2 $-cocycle 
	$$ u: G^2 \to \CC^*, $$
	which is trivial if the family $ \{\phi^g\}_{g \in G} $ satisfies the $ 2 $-cocycle condition. Moreover, another family $ \psi^g: E \xrightarrow{\sim} g^*E $ of isomorphisms of sheaves on $ X $ for all $ g \in G $ determines a $ 2 $-cocycle $ v: G^2 \to \CC^* $ in the same class of $ u $ in $ H^2(G, \CC^*) $. The argument above proves that if $ G $ is finite and $ H^2(G, \CC^*) = 0, $ then $ E $ is $ G $-linearizable because when $ G $ is finite, the $ 2 $-cocyle condition of the family $ \{\phi^g\}_{g \in G} $ gives an isomorphism $ \phi: p_2^*E \xrightarrow{\sim} \sigma^*E $ which satisfies the cocycle condition.
	
	Now suppose $ G $ is connected and $ \Pic(G) = 0. $ Since $ E $ is simple and $ \Pic(G) = 0 $, the proof of \cite[Proposition 4.3]{kool2011fixed} shows that there exists an isomorphism
	\begin{equation*}
	\phi: p_2^*E \xrightarrow{\sim} \sigma^*E.
	\end{equation*}
	By \cite[Proposition 2.4]{kool2011fixed}, it's also enough to prove the family 
	$$ \phi^g = \phi|_{\{g\}\times X} $$
	for all $ g $ in $ G $ satisfies the $ 2 $-cocycle condition, which is indeed the case as shown in the proof of \cite[Proposition 4.4]{kool2011fixed}.
	
	The same argument as that in the proof of the second statement in \cite[Lemma 1]{ploog2007equivariant} proves the last statement.
\end{proof}
\begin{rmk}
	The second group cohomology $ H^2(G, \CC^*) $ of a group $ G $ is called the \tb{Schur multiplier}. It vanishes for cyclic groups $ \ZZ/n\ZZ $ and dihedral groups $ D_{2n+1} $ for $ n > 1. $ In these cases, a $ G $-invariant simple sheaf is always $ G $-linearizable independent of the group action and the sheaf. On the other hand, there are $ G $-invariant line bundles which are not $ G $-linearizable when $ G $ is as simple as $ \ZZ/2\ZZ \times \ZZ/2\ZZ. $ 
\end{rmk}
\begin{ex}
	Let $ E $ be an elliptic curve under the action of $ G = \ZZ/2\ZZ \times \ZZ/2\ZZ $ generated by two $ 2 $-torsion points on $ E. $ A line bundle of degree $ 2 $ on $ E $ is $ G $-invariant but never $ G $-linearizable: the geometric quotient $ E \to E/G $ is a morphism of degree $ 4 $ between two elliptic curves, and hence the degree of a $ G $-linearizable line bundle on $ E $ must be a multiple of $ 4 $.
\end{ex}
\begin{ex}
	Let $ X \to \PP^1 $ be an elliptic $ K3 $ surface under the action of $ G = \ZZ/2\ZZ \times \ZZ/2\ZZ $ generated by two $ 2 $-torsion sections. Let $ C_0 $ denote the zero section of $ X \to \PP^1 $. Then the line bundle $ \cO_X(2C_0) $ is $ G $-invariant but not $ G $-linearizable: if it were, then it would give rise to a $ G $-linearizable line bundle on an elliptic fiber of $ X, $ which is impossible by the previous example.
\end{ex}
\begin{notation}
	Let $ \Qcoh^G(X) $ and $ \Coh^G(X) $ denote the categories of $ G $-equivariant quasi-coherent sheaves and of $ G $-equivariant coherent sheaves on $ X $ respectively. Then $ \Qcoh^G(X) $ and $ \Coh^G(X) $ are abelian categories and $ \Qcoh^G(X) $ has enough injectives (see \cite[Proposition 5.1.1 and 5.1.2]{grothendieck1957quelques}). Let $ \Vect^G(X) $ denote the category of $ G $-equivariant vector bundles on $ X. $
\end{notation}

Now let's consider quotient stacks. Let $ \Sh $ denote either $ \Qcoh, \Coh $ or $ \Vect. $ Then we have the following:
\begin{prop}\label{prop_equiv_cat_sheaves}
	Let $ [X/G] $ be a quotient stack with the canonical morphism $ p: X \to [X/G] $. Then pulling back sheaves along $ p $ induces a canonical equivalence
	\begin{equation*}
	\begin{tikzcd}[column sep = 2em]
	T: \Sh([X/G]) \ar[r, "\simeq"] & \Sh^G(X)
	\end{tikzcd}	
	\end{equation*}
	between the category of sheaves on $ [X/G] $ and that of $ G $-equivariant sheaves on $ X $.
\end{prop}
This reduces to \cite[Example 7.21]{vistoli1989intersection} for a DM stack $ \cX = [X/G] $ where sheaves on $ \cX $ are equivalently defined as sheaves on the small \'etale site of $ \cX $. Since a complete proof was not given there or in any other literature (as far as we know), we will give a detailed proof in Appendix \ref{proof_equiv_cat_sheaves}.

When the quotient stack $ [X/G] $ is represented by a scheme $ Y, $ then we have a principal $ G $-bundle $ p: X \to Y, $ and Proposition \ref{prop_equiv_cat_sheaves} reduces to the following:
\begin{cor}\label{cor_equiv_cat_sheaves}
	Let $ G $ be a linear algebraic group and let $ p: X \to Y $ be a principal $ G $-bundle. Then pulling back sheaves along $ p $ induces a canonical equivalence
	\begin{equation*}
	\begin{tikzcd}[column sep = 2em]
	T: \Sh(Y) \ar[r, "\simeq"] & \Sh^G(X)
	\end{tikzcd}	
	\end{equation*}
	between the category of sheaves on $ Y $ and that of $ G $-equivariant sheaves on $ X $.
\end{cor}

Corollary \ref{cor_equiv_cat_sheaves} can be seen as a special case of descent theory for schemes. Consider a morphism $ f: X \to Y $ of schemes. There are projections $ \pi_i: X \times_Y X \to X $ for $ i = 1, 2, $ and $ \pi_{ij}: X \times_Y X \times_Y X \to X \times_Y X $ for $ 1 \leq i < j \leq 3. $ Let $ E $ be a sheaf on $ X. $ We say $ E $ has a descent datum $ \phi $ with respect to the morphism $ f: X \to Y $ if 
\begin{equation*}
\phi: \pi_1^*E \xrightarrow{\sim} \pi_2^*E
\end{equation*}
is an isomorphism of sheaves on $ X \times_Y X $ satisfying the cocycle condition, i.e.,
\begin{equation*}
\pi_{13}^* \phi = \pi_{23}^* \phi \circ \pi_{12}^* \phi,
\end{equation*}
which is an identity of isomorphisms of sheaves on $ X \times_Y X \times_Y X. $ The pairs $ (E, \phi) $ form a category, which we denote by $ \Sh(X, f). $ The descent theory for schemes says
\begin{prop}\label{prop_descent}
	Let $ f: X \to Y $ be a faithfully flat quasi-compact morphism of schemes. Pulling back sheaves along $ f $ gives an equivalence
	\begin{equation*}
	\begin{tikzcd}[column sep = 2em]
	\Sh(Y) \ar[r, "\simeq"] & \Sh(X, f)
	\end{tikzcd}	
	\end{equation*}
	between the category of sheaves on $ Y $ and that of sheaves on $ X $ with descent data. 
\end{prop}
Now, let $ p: X \to Y $ be a principal $ G $-bundle. Then it is a faithfully flat quasi-compact morphism of schemes. There are isomorphisms
\begin{align*}
u: G \times X \xrightarrow{\sim} X \times_Y X, \quad (g, x) \mapsto (x, gx)
\end{align*}
and
\begin{align*}
v: G \times G \times X \xrightarrow{\sim} X \times_Y X \times_Y X, \quad (g, h, x) \mapsto (x, hx, ghx)
\end{align*}
of schemes, so
\begin{equation*}
\pi_1 = p_2 \circ u^{-1}, \quad \pi_2 = \sigma \circ u^{-1},
\end{equation*}
and
\begin{equation*}
\pi_{13} = u \circ \mu \times \id_X \circ v^{-1}, \quad \pi_{23} = u \circ \id_G \times \sigma \circ v^{-1}, \quad \pi_{12} = u \circ p_{23} \circ v^{-1}.
\end{equation*}
For a sheaf $ E $ on $ X $, a descent datum $ \phi $ of $ E $ translates into a $ G $-equivariant structure $ u^*\phi $ on $ E $ by applying $ u^* $ to the isomorphism $ \phi: \pi_1^*E \xrightarrow{\sim} \pi_2^*E $ and applying $ v^* $ to the cocycle condition $ \pi_{13}^* \phi = \pi_{23}^* \phi \circ \pi_{12}^* \phi. $ Therefore, we have an isomorphism of categories
\begin{equation*}
\Sh(X,p) \longrightarrow \Sh^G(X), \quad (E, \phi) \longmapsto (E, u^*\phi).
\end{equation*}
Corollary \ref{cor_equiv_cat_sheaves} is then a special case of Proposition \ref{prop_descent}. 
\begin{rmk}
	Suppose $ G $ is a finite group. The categorical equivalence in Corollary \ref{cor_equiv_cat_sheaves} was described explicitly in \cite[Chapter II, Proposition 2]{mumford2008abelian}. Consider a pullback diagram:
	\vspace{5pt}
	\begin{equation*}\label{eq_can_morphism_to_scheme}
	\begin{tikzcd}[column sep=3em,row sep=3em]
	G \times X \arrow[r, "\sigma"{above, outer sep = 2pt}] \arrow[d, "p_2"{right, outer sep = 2pt}] & X \arrow[d, "p_2"{right, outer sep = 2pt}] \\
	X \arrow[r, "p"{above, outer sep = 2pt}] & Y
	\end{tikzcd}\vspace{8pt}
	\end{equation*}
	Pulling back a sheaf $ \cE $ on $ Y $ equips $ p^*\cE $ a canonical $ G $-equivariant structure
	\begin{equation*}
	\phi_\cE: p_2^* (p^*\cE) \xrightarrow{\sim} \sigma^* (p^*\cE).
	\end{equation*}
	For a $ G $-equivariant sheaf $ (E, \phi) $ on $ X, $ we can take the $ G $-invariant subsheaf $ (p_*E)^G $ which depends on $ \phi. $ Then the assignment $ (E, \phi) \mapsto (p_*E)^G $ is a categorical inverse of the assignment $ \cE \mapsto (p^*\cE, \phi_\cE). $ In other words, there are natural isomorphisms
	\begin{equation*}
	\cE \xrightarrow{\sim} \left(p_*p^*\cE\right)^G \quad \text{and} \quad (E,\phi) \xrightarrow{\sim} \left(p^*(p_*E)^G, \phi_{(p_*E)^G}\right)
	\end{equation*}
	of sheaves on $ Y $ and of $ G $-equivariant sheaves on $ X $ respectively.
\end{rmk}

Let's fix a quotient stack $ \cX = [X/G] $ with the canonical morphism $ p: X \to \cX. $ Now we can identify a sheaf $ \cE $ on $ \cX $ as a $ G $-equivariant sheaf $ (E, \phi) $ on $ X $ up to isomorphism. By abuse of notation, we will also write $ \cE = (E, \phi) $ when there is an isomorphism $ T(\cE) = (p^*\cE, \phi_\cE) \cong (E, \phi) $ of $ G $-equivariant sheaves on $ X. $
\begin{ex}
	The $ G $-equivariant sheaf $ (\cO_X, \phi_0) $ on $ X $ corresponds to the structure sheaf $ \cO_\cX $ on $ \cX, $ i.e., $ \cO_\cX = (\cO_X, \phi_0) $ where $ \cO_X \cong p^*\cO_\cX. $
\end{ex}
\begin{ex}
	Consider a special case: $ X = \pt $ and $ \cX = BG. $ A $ G $-equivariant sheaf on $ X $ is a pair $ (V, \rho) $ where $ V $ is a $ d $-dimensional vector space and $ \rho: G \to \GL(V) $ is a representation of $ G, $ and $ \Tot(V) = \Spec \CC[x_1, \dots, x_d] $ endows $ V $ a scheme structure. A $ G $-equivariant morphism $ f: (V_1, \rho_1) \to (V_2, \rho_2) $ of $ G $-equivariant sheaves on $ BG $ is an \tb{intertwining} map of representations of $ G $. Therefore, we have an equivalence of abelian categories 
	$$ \Coh(BG) \cong \Coh^G(\pt) \cong \Rep(G), $$
	where $ \Rep(G) $ is the category of representations of $ G. $
\end{ex}
There are plenty of sheaves on quotient stacks as shown by the following construction.
\begin{constrn}[Sheaves on a quotient stack]\label{constrn_equiv_sheaves}
	For the quotient stack $ \cX = [X/G] $, there is a 2-pullback diagram:
	\vspace{5pt}
	\begin{equation}\label{eq_flat_base_change}
	\begin{tikzcd}[column sep=3em,row sep=3em]
	X \arrow[r, "u"{above, outer sep = 2pt}] \arrow[d, "p"{right, outer sep = 2pt}] & \pt \arrow[d, "q"{right, outer sep = 2pt}] \\
	\cX \arrow[r, "v"{above, outer sep = 2pt}] & BG
	\end{tikzcd}\vspace{8pt}
	\end{equation}
	We will use it to construct sheaves on $ \cX $. 
	
	\vspace{10pt}
	\noindent{\textit{Pullback representations of $ G $ to sheaves on $ \cX. $}} 
	Consider a sheaf $ (V, \rho) $ on $ BG, $ i.e., $ V $ is a $ d $-dimensional vector space equipped with a representation $ \rho: G \to \GL(V) $ of $ G. $ Denote $ (V, \rho) $ by $ \rho $ for simplification. We have $ q^* \rho = V $ which forgets the $ G $-action on $ V. $ Pulling back $ V $ along $ u $ gives a vector bundle
	\begin{equation*}
	\cO_X \otimes V \cong \cO_X^{\oplus d}
	\end{equation*}
	on $ X, $ and pulling back $ \rho $ along $ v $ gives a vector bundle
	\begin{equation*}
	\cO_\cX \otimes \rho \cong (\cO_X \otimes V, \phi_\rho)
	\end{equation*}
	on $ \cX, $ where the $ G $-equivariant structure $ \phi_\rho: G \to \Aut(X \times V) $ is defined by
	\begin{equation*}
	\phi_\rho(g)(x,v) = (gx,\rho(g)(v)).
	\end{equation*}
	Here the symbol $ \otimes $ does not mean the tensor product of sheaves, but stands for the pullback of sheaves from $ \pt $ to $ X $ or from $ BG $ to $ \cX. $ 
	
	\vspace{10pt}
	\noindent{\textit{Twist sheaves on $ \cX $ by representations of $ G. $}} 
	More generally, for a sheaf $ \cE = (E, \phi) $ on $ \cX, $ we can twist it by $ \rho $ to form another sheaf
	\begin{equation*}
	\cE \otimes \rho = \cE \otimes_{\cO_\cX} \left(\cO_\cX \otimes \rho \right) = (E \otimes V, \phi \otimes \phi_\rho)
	\end{equation*}
	on $ \cX, $ where $E \otimes V \cong E \otimes_{\cO_X} {\cO_X^{\oplus d}} \cong E^{\oplus d} $
	is a sheaf on $ X $ with a $ G $-equivariant structure
	\begin{equation*}
	\phi \otimes \phi_\rho: p_2^*(E \otimes V) \xrightarrow{\sim} \sigma^*(E \otimes V).
	\end{equation*}
	When $ \rho = \rho_0: G \to \CC^* $ is the one-dimensional trivial representation of $ G, $ we have
	\begin{equation*}
	\cE \otimes \rho_0 \cong \cE.
	\end{equation*}
	
	\vspace{10pt}
	\noindent{\textit{Push-pull vector spaces to sheaves on $ \cX. $}} 
	Consider a finite group $ G. $ Take a vector space $ V. $ Since the morphism $ v: \cX \to BG $ in diagram (\ref{eq_flat_base_change}) is a flat morphism (of DM stacks), there is a natural isomorphism
	\begin{equation*}
	v^* q_* V \cong p_* u^* V,
	\end{equation*}
	where the pushforward 
	\begin{equation*}
	q_* V \cong V \otimes \rho_\reg 
	\end{equation*}
	is the tensor product of $ V $ (as a trivial representation of $ G $) and the regular representation $ \rho_\reg $ of $ G $. Letting $ V = \CC, $ we then have
	$$ p_* \cO_X \cong \cO_\cX \otimes \rho_\reg, \quad \text{and} \quad \left(p_* \cO_X\right)^G \cong \cO_\cX.$$
\end{constrn}
\begin{rmk}
	The assignment $ (\cE, \rho) \mapsto \cE \otimes \rho $ induces an action of the ring $ R(G) $ of representations of $ G $ on the Grothendieck group $ K(\cX) $ of coherent sheaves on $ \cX, $ which makes $ K(\cX) $ a module over $ R(G). $ We will talk about this again after we define $ K(\cX). $
\end{rmk}

\section{Orbifold Hirzebruch-Riemann-Roch}\label{sec_HRR_formula_stacks}

\subsection{$ K $-theory of quotient stacks}
In this section we review some basic facts about the algebraic $ K $-theory of quotient stacks.

Let $ \cX $ be a locally noetherian algebraic stack.
\begin{defn}
	The \tb{Grothendieck group} $ K_0(\cX) $ is the abelian group generated by isomorphism classes of coherent sheaves on $ \cX $ modulo the relations $ \cE - \cE' - \cE'' $
	for every short exact sequence 
	\begin{equation}\label{eq_ses_sheaves}
	0 \to \cE' \to \cE \to \cE'' \to 0
	\end{equation}
	in $ \Coh(\cX) $. Denote the class of a coherent sheaf $ E $ in $ K_0(\cX) $ by $ [\cE] $. Define similarly another Grothendieck group $ K^0(\cX) $ using vector bundles on $ \cX $: it is generated by isomorphism classes of vector bundles on $ \cX $ modulo the relations $ \cE - \cE' - \cE'' $	for every short exact sequence (\ref{eq_ses_sheaves}) in $ \Vect(\cX) $. Denote the class of a vector bundle $ \cV $ in $ K^0(\cX) $ by $ [\cV]^0 $.
\end{defn}
\begin{rmk}
	Both $ K_0(\cX) $ and $ K^0(\cX) $ are also called the \tb{$ K $-group} of $ \cX. $ There exist higher $ K $-groups $ K_i(\cX) $ and $ K^i(\cX) $ for $ i \geq 0, $ such that $ K_\bullet(\cX) $ and $ K^\bullet(\cX) $ behave like homology and cohomology respectively. But we will not use these higher $ K $-groups. We will use the simpler notation $ K(\cX) $ for $ K_0(\cX). $
\end{rmk}
\begin{rmk}
	Let $ \cX = [X/G] $ be a quotient stack. By definition, the scheme $ X $ of finite type and is in particular noetherian, so the algebraic stack $ \cX $ is also noetherian. Hence we have two abelian groups $ K(\cX) $ and $ K^0(\cX) $. We can also define two abelian group $ K_G(X) $ and $ K_G^0(\cX) $ using $ G $-equivariant coherent sheaves and $ G $-equivariant vector bundles on $ X $ respectively. By Proposition \ref{prop_equiv_cat_sheaves}, there are natural isomorphisms
	\begin{equation*}
	K_G(X) \cong K(\cX) \quad \text{and} \quad K_G^0(X) \cong K^0(\cX).
	\end{equation*}
	There is a natural map 
	$$ K^0(\cX) \to K(\cX), \quad [\cV]^0 \mapsto [\cV] $$ 
	which may not be an isomorphism unless $ \cX $ has the resolution property.
\end{rmk}

\begin{defn}
	We say $ \cX $ has the \tb{resolution property} if every coherent sheaf on $ \cX $ is a quotient of a vector bundle, i.e., if every coherent sheaf $ \cE $ on $ \cX $ admits a surjection $ \cV \to \cE $ in $ \Coh(\cX) $ where $ \cV $ is in $ \Vect(\cX) $.
\end{defn}
It's known that all factorial separated noetherian schemes have the resolution property. We will see that all smooth quotient stacks have the resolution property. Note that the smoothness of a quotient stack $ [X/G] $ is equivalent to that of the scheme $ X $ over any base scheme.
\begin{prop}\label{prop_smooth_quotient_stack}
	Let $ S $ be a scheme. A quotient stack $ [X/G] $ is smooth over $ S $ if and only if $ X $ is smooth over $ S $.
\end{prop}
\begin{proof}
	Recall that an algebraic stack $ \cX $ is smooth over $ S $ if there is a smooth morphism $ Y \to \cX $ from a smooth scheme $ Y $ over $ S. $ For the quotient stack $ [X/G] $, the canonical morphism $ X \to [X/G] $ is smooth since $ G $ is assumed to be smooth over $ S $ in Definition \ref{defn_quotient_stack}. It follows that $ [X/G] $ is smooth over $ S $ if and only if $ X $ is so.
\end{proof}
\begin{prop}
	Let $ k $ be a field. All smooth quotient stacks over $ k $ have the resolution property.
\end{prop}
\begin{proof}
	Let $ [X/G] $ be a quotient stack over $ k $. Recall that $ X $ is assumed to be separated and of finite type over $ k $ in Definition \ref{defn_quotient_stack}, so it is noetherian. By Proposition \ref{prop_smooth_quotient_stack}, since $ [X/G] $ is smooth over $ k, $ so is $ X. $ Hence $ X $ is a factorial separated noetherian scheme, which ensures that $ [X/G] $ has the resolution property by \cite{totaro2004resolution}.
\end{proof}
Let $ \cX = [X/G] $ be a smooth quotient stack. Then every coherent sheaf $ \cE $ on $ \cX $ has a locally free resolution $ \cE_{\boldsymbol{\cdot}} \to \cE \to 0 $. Since $ X $ is a regular scheme, by Hilbert's syzygy theorem, every coherent sheaf $ \cE $ has a locally free resolution of finite length. This proves the following:
\begin{lem}\label{prop_finite_resolution}
	Let $ \cX $ be a smooth quotient stack. Then every coherent sheaf $ \cE $ on $ \cX $ has a finite locally free resolution
	\begin{equation*}\label{eq_loc_free_resolution}
	0 \to \cE_n \to \cE_{n-1} \to \cdots \to \cE_1 \to \cE_0 \to \cE \to 0.
	\end{equation*}
\end{lem}
\begin{rmk}
	Let $ \cX = [X/G] $ be a quotient stack. The group $ K^0(\cX) $ always has a ring structure with the multiplication given by the tensor product of vector bundles and the unity $ 1 = [\cO_\cX]^0 $. We call $ K^0(\cX) $ the \tb{Grothendieck ring} of $ \cX. $ Suppose we have another quotient stack $ \cY = [Y/H] $, a morphism $ f_1: Y \to X $
	of schemes, and a group homomorphism $ f_2: H \to G $
	such that the morphism $ f_1 $ is compatible with the group homomorphism $ f_2, $ i.e.,
	\begin{equation*}
	f_1(h \cdot y) = f_2(h) \cdot f_1(y)
	\end{equation*}
	for all $ h $ in $ H $ and all $ y $ in $ Y. $ Then the pair $ (f_1,f_2) $ descends to a morphism
	\begin{equation*}
	f: \cY \to \cX
	\end{equation*}
	of quotient stacks. The morphism $ f $ induces a ring homomorphism
	\begin{equation*}
	f^*: K^0(\cX) \to K^0(\cY)
	\end{equation*}
	which maps the class $ [\cV]^0 $ of a vector bundle $ \cV $ on $ \cX $ to the class $ [f^*\cV]^0 $ of the pullback vector bundle $ f^*\cV $ on $ \cY. $
\end{rmk}
\begin{lem}
	Let $ \cX = [X/G] $ be a smooth quotient stack. Then the natural map $ K^0(\cX) \to K(\cX) $ is an isomorphism of abelian groups.
\end{lem}
\begin{proof}
	Take a coherent sheaf $ \cE = (E, \vphi) $ on $ \cX $. By Proposition \ref{prop_finite_resolution}, $ \cE $ has a finite locally free resolution $ \cE_{\boldsymbol{\cdot}} \to \cE \to 0 $, which is the same as a finite $ G $-equivariant locally free resolution $ E_{\boldsymbol{\cdot}} \to E \to 0 $ where each $ E_i $ has a $ G $-equivariant structure $ \vphi_i $. The class $ \sum (-1)^i [\cE_i]^0 $ in $ K^0(\cX) $ is independent of the resolution, and the map $ \cE \mapsto  \sum (-1)^i [\cE_i]^0 $ is additive on short exact sequences in $ \Coh(\cX) $ by modifying the argument in \cite[Appendix B.8.3]{fulton1998intersection} to make all morphisms $ G $-equivariant. Then the map $ [\cE] \mapsto  \sum (-1)^i [\cE_i]^0 $ is inverse to the natural map $ [\cV]^0 \mapsto [\cV] $.
\end{proof}
With the natural isomorphism $ K^0(\cX) \cong K(\cX) $, we can equip the group $ K(\cX) $ with a ring structure. Alternatively, we can define a multiplication directly on $ K(\cX) $ as follows:
\begin{defn}\label{defn_multiplication_sheaves}
	Let $ \cX $ be a smooth quotient stack. Define a multiplication on $ K(\cX) $ by
	\begin{equation}\label{eq_multiplication_sheaves}
	[\cE] [\cF] := \sum (-1)^i [\Tor_i^{\cO_\cX}(\cE,\cF)]
	\end{equation}
	for coherent sheaves $ \cE $ and $ \cF $ on $ \cX $ and extend linearly. 
\end{defn}
This is well-defined because each coherent sheaf on $ \cX $ has a finite resolution of locally free sheaves under the assumptions in Proposition \ref{prop_finite_resolution}.
\begin{rmk}
	Note that the product of $ [\cE] $ and $ [\cF] $ above corresponds to the class of the derived tensor product $ \cE \otimes^L \cF $ in the Grothendieck group of the derived category of coherent sheaves on $ \cX. $
\end{rmk}
The following proposition is an immediate consequence of Definition \ref{defn_multiplication_sheaves}.
\begin{prop}\label{lem_multiplication_sheaves}
	Let $ \cX  $ be a smooth quotient stack. Let $ \cE $ and $ \cF $ be coherent sheaves on $ \cX $. Then their product in $ K(\cX) $ is given by
	\begin{equation*}
	[\cE] [\cF] = \sum_i (-1)^i [\cE_i \otimes \cF]
	\end{equation*}
	for any finite locally free resolution $ \cE_{\boldsymbol{\cdot}} \to \cE \to 0 $.
	If either $ \cE $ or $ \cF $ is locally free, then $ [\cE] [\cF] = [\cE \otimes \cF] $ in $ K(\cX) $.
\end{prop}
The multiplication on $ K(\cX) $ defined by (\ref{eq_multiplication_sheaves}) is compatible with the multiplication on $ K^0(\cX) $ defined by the tensor product of vector bundles, and hence is
commutative and associative. Therefore, we have the following:
\begin{prop}\label{prop_isom_K_theory}
	Let $ \cX $ be a smooth quotient stack. The group $ K(\cX) $ has a ring structure with the multiplication given by  (\ref{eq_multiplication_sheaves}). The natural map $ K^0(\cX) \to K(\cX) $ is a ring isomorphism with the inverse given by $ [\cE] \mapsto  \sum_i (-1)^i [\cE_i]^0 $ for any finite locally free resolution $ \cE_{\boldsymbol{\cdot}} \to \cE \to 0 $.
\end{prop}
For a smooth quotient stack $ \cX, $ we also call $ K(\cX) $ the Grothendieck ring of $ \cX. $
\begin{notation}
	From now on, we will omit the superscript $ ``0" $ in the notation $ [\cE]^0 $ for the class of a vector bundle $ \cE $ on $ \cX $. 
	We denote the inverse map in Proposition \ref{prop_isom_K_theory} by
	\begin{equation*}
	\beta: K(\cX) \to K^0(\cX), \quad [\cE] \mapsto \sum_{i} (-1)^i[\cE_i].
	\end{equation*}
	A morphism $ f: \cY \to \cX $ of smooth quotient stacks induces a ring homomorphism
	\begin{equation*}
	f^{K}: K(\cX) \xrightarrow{\beta} K^0(\cX) \xrightarrow{f^*} K^0(\cY) \xrightarrow{\beta^{-1}} K(\cY).
	\end{equation*}
\end{notation}
The Grothendieck ring $ K^0(\cX) $ has an additional algebraic structure called a $ \lambda $-ring. We follow \cite[Section I.1]{fulton1985riemann} for the definition of a $ \lambda $-ring, which is called a pre-$ \lambda $-ring in some literatures, for example in \cite{yau2010lambda}.
\begin{defn}
	A \tb{$ \lambda $-ring} is a ring $ R $ with a group homomorphism
	\begin{align*}
	\lambda_t: (R, +) &\to R[[t]]^\times \\
	x & \mapsto \sum_{i \geq 0} \lambda^i(x) t^i
	\end{align*}
	such that $ \lambda^0(x) = 1 $ and $ \lambda^1(x) = x $ for all $ x \in R. $
	A \tb{positive structure} on a $ \lambda $-ring $ (R, \lambda_t) $ is a pair $ (\varepsilon, R_+) $ where $ \varepsilon: R \to \ZZ $ is a surjective ring homomorphism and $ R_+ $ is a subset of $ R $, which satisfies the following:
	\begin{enumerate}[font=\normalfont,leftmargin=*]
		\item $ R_+ $ is closed under addition and multiplication.
		\item $ R_+ $ contains the set $ \ZZ_+ $ of positive integers.
		\item $ R = R_+ - R_+, $ i.e., every element of $ R $ is a difference of two elements of $ R_+ $.
		\item If $ x \in R_+ $, then $ \varepsilon(x) = r > 0,\ \lambda^i(x) = 0 $ for $ i > r $, and $ \lambda^r(x) $ is a unit in $ R. $
	\end{enumerate}
	Elements of $ R_+ $ are said to be \tb{positive}. An element of $ R $ is \tb{non-negative} if it is either positive or zero. The set of non-negative elements is denoted by
	\begin{equation*}
	R_{\geq 0} = R_+ \cup \{0\}.
	\end{equation*}
\end{defn}
\begin{rmk}
	Consider a $ \lambda $-ring $ (R, \lambda_t). $ For all $ x, y \in R, $ we have
	\begin{equation*}
	\lambda_t(x+y) = \lambda_t(x)\lambda_t(y),
	\end{equation*}
	which is equivalent to
	\begin{equation*}
	\lambda^n(x+y) = \sum_{i=0}^{n} \lambda^i(x)\lambda^{n-i}(y)
	\end{equation*}
	for all integer $ n \geq 0.$ Note that we always have $ \lambda_t(0) = 1. $ Suppose $ (R, \lambda_t) $ has a positive structure $ (\varepsilon, R_+). $ Then we have $ \varepsilon(1) = 1 $ and $ \lambda_t(1) = 1 + t. $ On the other hand, we have $ \varepsilon(-1) = -1 $ and 
	\begin{equation*}
	\lambda_t(-1) = \lambda_t(1)^{-1} = 1 - t + t^2 - t^3 + \cdots.
	\end{equation*}
	From now on, we always assume there is a positive structure associated with a $ \lambda $-ring.
\end{rmk}

\begin{ex}
	The simplest example of a $ \lambda $-ring is $ \ZZ $ with $ \lambda_t(m) = (1+t)^m $ and the obvious positive structure $ (\varepsilon = \id: \ZZ \to \ZZ, \ZZ_+). $ In this case, we have 
	$$ \lambda^i(m) = \binom{m}{i} $$ 
	for $ m \in \ZZ $ and $ i \geq 0. $
\end{ex}

\begin{ex}[{\cite[Section V.1]{fulton1985riemann}}]
	Let $ X $ be a scheme. The Grothendieck ring $ K^0(X) $ of vector bundles on $ X $ is a $ \lambda $-ring: positive elements $ K^0(X)_+ $ are classes of vector bundles $ V $, taking exterior powers induces a map $ \lambda_t: K^0(X)_{\geq 0} \to K^0(X)[[t]]^\times $ given by
	\begin{equation*}
	\lambda_t(V) = \sum_{i \geq 0} \lambda^i (V) t^i = \sum_{i \geq 0} [\wedge^i V] t^i,
	\end{equation*}
	which is a semigroup homomorphism: for two vector bundles $ V $ and $ W,$
	\begin{equation*}
	\lambda_t(V + W) = \lambda_t(V)\lambda_t(W)
	\end{equation*}
	because 
	\begin{equation*}
	\wedge^n(V \oplus W) \cong \bigoplus_{i=0}^{n} \left(\wedge^i V \otimes \wedge^{n-i} W\right) \ \text{for all integers} \ n \geq 0.
	\end{equation*}
	Since each element in $ K^0(X) $ can be written as a difference $ [V] - [W] $ of two classes of vector bundles, we can extend the domain of $ \lambda_t $ to all of $ K^0(X) $, i.e.,
	\begin{equation*}
	\lambda_t(V-W) = \lambda_t(V)/\lambda_t(W).
	\end{equation*}
	The surjective ring homomorphism $ \varepsilon = \rk:  K^0(X) \to \ZZ $ is the \tb{rank} map: for an element $ x = [V] - [W], $
	$$ \rk(x) = \rk(V) - \rk(W). $$ 
	When $ X = \pt, $ we get back the simplest example of $ \ZZ. $ 
\end{ex}
\begin{prop}
	Let $ \cX = [X/G] $ be a quotient stack. The Grothendieck ring $ K^0(\cX) $ of vector bundles on $ \cX $ is a $ \lambda $-ring.
\end{prop}
\begin{proof}
	Consider a vector bundle $ \cV = (V, \phi) $ on $ \cX $. Exterior powers $ \wedge^i V $ of $ V $ carry canonical $ G $-equivariant structures $ \wedge^i \phi, $ so $ \wedge^i \cV $ is well-defined. Every element in $ K^0(\cX) $ is a difference $ [\cV] - [\cW] $ of two classes of equivariant vector bundles, so we can define the same $ \lambda_t $ and $ \varepsilon = \rk $ as those in the case of schemes. 
\end{proof}
\begin{ex}
	When $ X = \pt $, $ K^0(\cX) = K^0(BG) $ becomes the representation ring $ R(G), $ which is another example of $ \lambda $-rings (\cite{atiyah1969group}).
\end{ex}
\begin{defn}
	Let $ \cX $ be a quotient stack. Taking dual vector bundles respects exact sequences and hence defines an involution 
	$$ (\ \cdot \ )^\vee: K^0(\cX) \to K^0(\cX). $$
	When $ \cX $ is smooth, this gives an involution
	\begin{equation*}
	(\ \cdot \ )^\vee: K(\cX) \to K(\cX)
	\end{equation*}
	via the natural isomorphism $ K^0(\cX) \cong K(\cX). $
\end{defn}
\begin{rmk}
	The involution $ (\ \cdot \ )^\vee: K^0(\cX) \to K^0(\cX) $ is a ring automorphism, i.e., for all $ x $ and $ y $ in $ K^0(\cX), $ we have
	\begin{equation*}
	1^\vee = 1, \quad (x+y)^\vee = x^\vee + y^\vee, \quad \text{and} \quad (xy)^\vee = x^\vee y^\vee.
	\end{equation*}
	When $ \cX $ is smooth, $ K(\cX) $ is also a ring whose multiplication is given by (\ref{eq_multiplication_sheaves}) and hence the involution $ (\ \cdot \ )^\vee: K(\cX) \to K(\cX) $ is also a ring automorphism.
\end{rmk}
\begin{defn}
	Let $ \cX = [X/G] $ be a connected quotient stack. The \tb{$ K $-theoretic Euler class} of a non-negative element $ x = \cV $ in $ K^0(\cX) $ of rank $ r $ is defined by
	\begin{equation*}
	e^K(x) = \lambda_{t}(x^\vee) |_{t=-1} = \sum_{i=0}^{r} (-1)^i \lambda^i(x^\vee)  = \sum_{i=0}^{r} (-1)^i [\wedge^i \cV].
	\end{equation*}
\end{defn}
\begin{rmk}[Properties of the $ K $-theoretic Euler class]
	Let $ \cX $ be a connected quotient stack. The $ K $-theoretic Euler class gives a semigroup homomorphism 
	$$ e^K: (K^0(\cX)_{\geq 0}, +) \to (K^0(\cX), \times), $$ 
	i.e., for all non-negative elements $ x $ and $ y $ in $ K^0(\cX), $ we have
	\begin{equation*}
	e^K(x+y) = e^K(x)e^K(y).
	\end{equation*}
	In particular, $ e^K(0) = 1. $ By definition, we also have $ e^K(1) = 1 - 1 = 0. $ 
	Note that $ e^K $ cannot be defined on all of $ K^0(\cX) $: for example, $ e^K(-1) $ is undefined since $ \lambda_t(-1) $ is an infinite series. If $ x $ and $ y $ are positive elements of $ K^0(\cX) $ such that $ e^K(y) $ is a unit, then we define
	\begin{equation*}
	e^K(x-y) = e^K(x)/e^K(y).
	\end{equation*}
	For a morphism $ f: \cX \to \cY $ of connected quotient stacks,
	the $ K $-theoretic Euler class commutes with the induced ring homomorphism $ f^*: K^0(\cY) \to K^0(\cX), $
	i.e.,
	\begin{equation*}
	e^K(f^* y) = f^* e^K(y)
	\end{equation*}
	for all $ y \in K^0(\cY) $ when $ e^K(y) $ is defined.
\end{rmk}
\begin{rmk}[Geometric meaning of the $ K $-theoretic Euler class]
	Consider a vector bundle $ \cV = (V, \phi) $ of rank $ r $ on a quotient stack $ \cX = [X/G] $. Suppose the vector bundle $ \cV $ has a section $ s: \cO_\cX \to \cV $ ( i.e., a $ G $-equivariant section $ s: \cO_X \to V $) such that it cuts out a substack $ \cY = [Y/G] $, where $ Y $ is $ G $-invariant subscheme of $ X $ of codimension $ r. $ Let $ i: Y \to X $ denote the inclusion, and let $ s^\vee: V^\vee \to \cO_X $ denote the dual of $ s $. The Koszul complex
	\begin{equation*}
	0 \to \wedge^r V^\vee \to \wedge^{r-1} V^\vee \to \cdots \to V^\vee \xrightarrow{s^\vee} \cO_X \to i_*\cO_Y \to 0
	\end{equation*}
	on $ X $ is a $ G $-equivariant locally free resolution of the pushforward $ i_*\cO_Y, $ so we obtain a Koszul complex
	\begin{equation*}
	0 \to \wedge^r \cV^\vee \to \wedge^{r-1} \cV^\vee \to \cdots \to \cV^\vee \to \cO_\cX \to i_*\cO_\cY \to 0.
	\end{equation*}
	on $ \cX. $
	Therefore,
	\begin{equation}\label{eq_eK_V}
	[i_*\cO_\cY] = e^K(\cV) \in K_0(\cX),
	\end{equation}
	which explains the name ``$ K $-theoretic Euler class" since the usual Euler class $ e(\cV) $ is the class 
	$$ [\cY] = c_r(\cV) $$
	in the $ r $-th Chow group $ A^r(\cX) $ of $ \cX. $
	Suppose $ \cX $ and $ \cY $ are smooth. Taking the equivariant Chern class in Example 15.3.1 in \cite{fulton1998intersection}, we obtain a relation between the $ K $-theoretic Euler class and the usual Euler class:
	\begin{equation*}
	c_r(e^K(\cV)) = (-1)^{r-1} (r-1)! e(\cV) \in A^r(\cX).
	\end{equation*}
	By identifying $ K(\cX) \cong K^0(\cX) $ and $ K(\cY) \cong K^0(\cY), $ we have a $ K $-theoretic pullback 
	$$ i^K: K(\cX) \to K(\cY). $$ 
	Since the inclusion $ i: \cY \to \cX $ is proper, we also have a $ K $-theoretic pushward 
	$$ i_K: K(\cY) \to K(\cX). $$ 
	Note that $ i_! $ is a group homomorphism but not a ring homomorphism in general. Let $ 1_{K(\cY)} = [\cO_\cY] $ denote the unity in $ K(\cY). $ Then
	$$ i_K\left(1_{K(\cY)}\right) = [i_K\cO_\cY] = e^K(\cV) $$ 
	in $ K(\cX) $ since the morphism $ i $ is affine. 
	By (\ref{eq_eK_V}), it follows that\footnote{In general, $ i^K i_K (y) = e^K(\cN) \cdot y $ for every $ y \in K(\cY). $ See \cite{vezzosi2003higher} for a proof.}
	\begin{equation}\label{eq_eK_N}
	i^K i_K \left(1_{K(\cY)}\right) = e^K(\cN)
	\end{equation}
	in $ K(\cY), $ where $ \cN \cong  i^*\cV $ is the $ G $-equivariant normal bundle of $ Y $ in $ X. $ 
	Applying the $ K $-theoretic projection formula to (\ref{eq_eK_N}), we obtain the \tb{$ K $-theoretic self-intersection formula}
	\begin{equation*}
	i_K \left(e^K(\cN)\right) = i_K \left(1_{K(\cY)}\right) \cdot i_K \left(1_{K(\cY)}\right) 
	\end{equation*}
	in $ K(\cX), $ where $ \cdot $ is the multiplication on $ K(\cX) $ defined in (\ref{eq_multiplication_sheaves}).
\end{rmk}
We can use the $ K $-theoretic Euler class to compute the Grothendieck ring of a weighted projective stack $ \PP(a_0, \dots, a_n) $.
\begin{ex}[The Grothendieck ring of of a weighted projective stack]
	Let $ X = \CC^{n+1} $ with a closed subscheme $ p = \{0\} $, and let $ \cX = [X/\CC^*] $ where the $ \CC^* $-action on $ X $ is given by
	\begin{equation*}
	t (x_0, \dots, x_n) = (t^{a_0} x_0, \dots, t^{a_n} x_n)
	\end{equation*}
	for some positive integers $ a_0, \dots, a_n. $ Let $ U = X - p. $ We denote the weighted projective stack $ \PP(a_0, \dots, a_n) = [U/\CC^*] $ by $ \PP. $ Let $ i: p \into X $ and $ j: U \into X $ be the inclusions of $ p $ and $ U $ in $ X, $ which are $ G $-equivariant morphisms of schemes. Hence they induce inclusions 
	\begin{equation*}
	i: B\CC^* \to \cX \quad \text{and} \quad j: \PP \to \cX
	\end{equation*}
	of substacks in $ \cX. $
	Therefore, we have a right exact sequence
	\begin{equation}\label{eq_K_theory_right_exact}
	K(B\CC^*) \to K(\cX) \to K(\PP) \to 0
	\end{equation}
	of abelian groups where the first map is extension by zero and the second one is restriction. Since $ \cX $ and $ \PP $ are smooth quotient stacks with the resolution property, we can identify the two maps in (\ref{eq_K_theory_right_exact}) as the $ K $-theoretic pushforward and pullback
	$$ i_K: K(B\CC^*) \to K(\cX) \quad \text{and} \quad j^K: K(\cX) \to K(\PP) $$ 
	where $ K(\cX) $ is identified with $ K^0(\cX). $
	We know that 
	\begin{equation*}
	K(B\CC^*) \cong \ZZ[u, u^{-1}] \quad \text{and} \quad K(\cX) \cong \ZZ[x, x^{-1}],
	\end{equation*}
	where $ u $ is the class of the identity character 
	$$ \rho = \id: \CC^* \to \CC^*, $$ 
	and $ x $ is the class of the $ \CC^* $-equivariant line bundle $ \cO_\cX \otimes \rho. $
	The $ \CC^* $-equivariant tangent bundle $ TX \cong \cO_\cX \otimes (\rho^{a_0} \oplus \cdots \oplus \rho^{a_n}) $ on $ X $ has a canonical section 
	$$ s: \cO_X \to TX, \quad 1 \mapsto (x_0, \dots, x_n), $$ 
	which cuts out the origin $ p $ of $ X. $ Note that $ s $ corresponds to the morphism $ X \times \CC \to X \times \CC^{n+1} $ given by
	\begin{equation*}
	(x_0, \dots, x_n, z) \mapsto (x_0, \dots, x_n, x_0 z, \dots, x_n z),
	\end{equation*}
	which respects the $ \CC^* $-equivariant structures on $ \cO_X $ and $ TX $.
	Since
	\begin{equation*}
	[TX] = x^{a_0} + \cdots + x^{a_n} \in K(\cX),
	\end{equation*}
	the image of $ 1 = [\cO_{B\CC^*}] $ under the map $ i_K $ is
	\begin{equation*}
	i_K (1) = e^K(TX) = e^K \left(\sum_{i=0}^n x^{a_i}\right) = \prod_{i=0}^n e^K(x^{a_i}) = \prod_{i=0}^n \lambda_{-1}(x^{-a_i}) = \prod_{i=0}^n (1-x^{-a_i}) \in K(\cX).
	\end{equation*}
	In general, we have
	\begin{equation*}
	i_K (u^k) = x^k \prod_{i=0}^n (1-x^{-a_i}) \in K(\cX)
	\end{equation*}
	for all $ k \in \ZZ. $ Therefore,
	\begin{equation*}
	K(\PP) \cong \frac{\ZZ[x,x^{-1}]}{\inprod{(1-x^{-a_0})\cdots(1-x^{-a_n})}} \cong \frac{\ZZ[x]}{\inprod{(x^{a_0}-1)\cdots(x^{a_n}-1)}}
	\end{equation*}
	as a ring, where $ x $ is the class of the twisting sheaf
	$$ \cO_{\PP}(1) := j^*(\cO_\cX \otimes \rho) $$ 
	on $ \PP, $ and is invertible in $ K(\PP). $ Here each $ 1 - x^{-a_i} $ is the class of $ \cO_{\PP_i} $ in the short exact sequence
	\begin{equation*}
	0 \to \cO_\PP(-a_i) \to \cO_\PP \to \cO_{\PP_i} \to 0
	\end{equation*}
	of sheaves on $ \PP $, where $ \PP_i = \PP(a_0, \dots, \hat{a_i},\dots,a_n) $ is the substack in $ \PP $ cuts out by the coordinate $ x_i. $ Therefore, the relation
	\begin{equation}\label{eq_wps_relation}
	(1-x^{-a_0}) \cdots (1-x^{-a_n}) = 0
	\end{equation}
	in $ K^0(\PP) $ says
	\begin{equation*}
	[\cO_{\PP_0}] \cdots [\cO_{\PP_n}] = 0
	\end{equation*}
	in $ K^0(\PP) $, which reflects the fact that the intersection of $ \PP_0, \dots, \PP_n $ is empty.
\end{ex}
\begin{ex}[The $ K $-theoretic Euler class of $ T\PP $]
	A weighted projective stack $ \PP = \PP(a_0, \dots, a_n) $ has a tangent bundle $ T\PP $ determined by the Euler sequence
	\begin{equation*}
	0 \to \cO_\PP \to \oplus_{i=0}^n \cO_\PP(a_i) \to T\PP \to 0,
	\end{equation*}
	so we have
	\begin{equation*}
	[T\PP] = x^{a_0} + \cdots + x^{a_n} - 1
	\end{equation*}
	in $ K(\PP). $ Let's compute $ e^K(T\PP) = \lambda_t([T\PP]^\vee) |_{t=-1}. $ Here $ [T\PP]^\vee $ is the class of the cotangent bundle of $ \PP $ and is given by
	\begin{equation*}
	[T\PP]^\vee = x^{-a_0} + \cdots + x^{-a_n} - 1
	\end{equation*}
	in $ K(\PP). $ Since $ \lambda_t $ is multiplicative, we have
	\begin{equation*}
	\lambda_t([T\PP]^\vee) = \frac{\prod_i \lambda_t(x^{-a_i} )}{\lambda_t(1)} = \frac{\prod_i (1+x^{-a_i}t)}{1+t}.
	\end{equation*}
	The evaluation of $ \lambda_t([T\PP]^\vee) $ at $ t = -1 $ can be done by the L'H\^{o}pital's rule, i.e., by viewing $ \lambda_t([T\PP]^\vee) $ as a function in $ t $ and take the limit as $ t \to 0: $
	\begin{equation*}
	e^K(T\PP) = \lambda_t([T\PP]^\vee) |_{t=-1} = \sum_{i=0}^{n} \left(x^{-a_i} \prod_{j \neq i} (1-x^{-a_j})\right) = \sum_{i=0}^{n} \prod_{j \neq i} (1-x^{-a_j})
	\end{equation*}
	in $ K(\PP), $ where the last equality follows from the relation (\ref{eq_wps_relation}).
	Therefore, 
	$$ e^K(T\PP) = [\cO_{p_0}] + \cdots + [\cO_{p_n}] $$
	in $ K(\PP), $ where $ p_i = \PP(a_i) $ is a point in $ \PP $ cut out by the coordinates $ x_0, \cdots, \hat{x_i}, \cdots, x_n. $ On the other hand, from the Chern polynomial
	\begin{equation*}
	c_t(T\PP) = c_t\left(\oplus_{i=0}^n \cO_\PP(a_i)\right) = \prod_{i=0}^{n} (1+a_i h t)
	\end{equation*}
	with $ h = c_1(\cO_\PP(1)) $ in $ A^1(\PP), $ the usual Euler class is given by
	\begin{equation*}
	e(T\PP) = c_n(T\PP) = \left(\sum_{i=0}^{n} \prod_{j \neq i} a_i\right) h^n \in A^n(\PP).
	\end{equation*}
	The relation between the two Euler classes is
	\begin{equation*}
	c_n(e^K(T\PP)) = (-1)^{n-1} (n-1)! e(T\PP) \in A^n(\PP).
	\end{equation*}
	For example, if $ a_0 = \cdots = a_n = 1, $ then $ \PP = \PP^n $ is the usual projective space,
	\begin{equation*}
	K(\PP^n) \cong \frac{\ZZ[x,x^{-1}]}{\inprod{(1-x^{-1})^{n+1}}} \cong \frac{\ZZ[x]}{\inprod{(x-1)^{n+1}}},
	\end{equation*}
	and each $ [\cO_{p_i}] $ is the same as the class of any point $ p $ in $ \PP^n, $ so
	\begin{equation*}
	e^K(T\PP^n) = (n+1)[\cO_p] = (n+1)(1-x^{-1})^n \in K(\PP^n),
	\end{equation*}
	and
	\begin{equation*}
	e(T\PP) = (n+1)h^n \in A^n(\PP^n).
	\end{equation*}
	The relation between $ e^K(T\PP^n) $ and $ e(T\PP^n) $ becomes
	\begin{equation*}
	c_n((n+1)(1-x^{-1})^n) = (-1)^{n-1} (n-1)! (n+1)h^n \in A^n(\PP^n),
	\end{equation*}
	which can be verified directly by the identity
	\begin{equation*}
	c_n(v) = (-1)^{n-1} (n-1)! \ch_n(v) \in A^n(\PP^n)
	\end{equation*}
	for any $ v \in K(\PP^n) $ with vanishing $ c_i(v) $ for $ 0 < i < n, $ and computing
	\begin{equation*}
	\ch_n((n+1)(1-x^{-1})^n) = \left((n+1)(1-e^{-h})^n\right)_n = (n+1)h^n \in A^n(\PP^n),
	\end{equation*}
	where $ (\ \cdot\ )_n $ denotes taking the homogeneous part of degree $ n $ of an element in $ A(\PP_n) $.
\end{ex}
\begin{rmk}
	When $ \PP = \PP^n $ is the usual projective space, the Chern character (after tensored with $ \QQ $) gives a $ \QQ $-algebra isomorphism
	\begin{equation*}
	\ch: \frac{\QQ[x]}{\inprod{(x-1)^{n+1}}} \xrightarrow{\sim} \frac{\QQ[h]}{\inprod{h^{n+1}}}, \quad x \mapsto e^h
	\end{equation*}
	with the inverse given by sending $ h $ to $ \ln(x) = \sum_{i} (-1)^{i+1}(x-1)^i/i.$ The same map on the weighted projective stack $ \PP(a_0,\dots, a_n) $ doesn't give a $ \QQ $-algebra isomorphism between $ K(\PP)_\QQ $ and $ A(\PP)_\QQ \cong \QQ[h]/h^{n+1} $ since
	\begin{equation*}
	\frac{\QQ[x]}{\inprod{(x^{a_0}-1)\cdots(x^{a_n}-1)}} \cong \frac{\QQ[x]}{\inprod{(x-1)^{n+1}}} \times \frac{\QQ[x]}{\inprod{\prod_{i=0}^{n} \sum_{j=1}^{a_i}x^{a_i-j} }}
	\end{equation*}
	by the Chinese remainder theorem. In fact, there is an orbifold Chern character map
	\begin{equation*}
	\orbch: K(\PP) \to A(I\PP)_\CC
	\end{equation*}
	where $ I\PP $ is the inertia stack of $ \PP $, which becomes a $ \CC $-algebra isomorphism 
	\begin{equation*}
	\orbch: K(\PP)_\CC \to A(I\PP)_\CC
	\end{equation*}
	after tensored with $ \CC $.
\end{rmk}

\subsection{Inertia stacks of quotient stacks}
In this section we study the inertia stack associated to a quotient stack, which is the key ingredient for the orbifold HRR theorem. 
\begin{defn}\label{defn_inertia_stack}
	Let $ \cX = [X/G] $ be a quotient stack. The \tb{inertia scheme} of $ X $ under the action $ G $ is defined by the following pullback diagram:
	\vspace{5pt}
	\begin{equation}\label{eq_inertia_scheme}
	\begin{tikzcd}[column sep=4em,row sep=4em]
	I_G X \arrow{r} \arrow{d} & X \arrow{d}{\Delta} \\
	G \times X \arrow{r}{\alpha} & X \times X
	\end{tikzcd}\vspace{8pt}
	\end{equation} 
	Since $ \Delta $ and $ \alpha $ are $ G $-equivariant, they descend to morphisms $ \Delta_G: \cX \to [X \times X /G] $ and $ \alpha_G: [G \times X /G] \to [X \times X /G]. $ The \tb{inertia stack} of $ \cX $ is defined by the following 2-pullback diagram:
	\vspace{5pt}
	\begin{equation}\label{eq_inertia_stack}
	\begin{tikzcd}[column sep=4em,row sep=4em]
	I \cX \arrow{r} \arrow{d} & \cX \arrow{d}{\Delta} \\
	\left[G/G\right] \times \cX \arrow{r}{\alpha} & \cX \times \cX
	\end{tikzcd}\vspace{8pt}
	\end{equation}
\end{defn}
\begin{rmk}
	The inertia scheme $ I_G X $ can be identified as a subscheme of $ G \times X $ whose $ \CC $-points are
	$$ I_G X (\CC) = \{(g,x) \in G \times X \ | \ gx = x\}. $$
	The inertia stack $ I \cX $ can then be identified as a quotient stack $ [I_G X/G] $ where $ G $ acts on $ I_G X $ by
	$$ g(h,x) = (ghg^{-1}, gx) $$ 
	for $ g, h \in G $ and $ x \in X. $ For an algebraic stack $ \cX $, the inertia stack of $ \cX $ is defined by $ I\cX = \cX \times_\Delta \cX $ where $ \Delta: \cX \to \cX \times \cX $ is the diagonal morphism. For a quotient stack $ \cX $, this is equivalent to Definition \ref{defn_inertia_stack}. The inertia stack $ I\cX $ of a quotient stack $ \cX $ is independent of its presentation $ \cX = [X/G] $.
	The $ G $-equivariant projection $ I_G X \to X, (g,x) \mapsto x $ in (\ref{eq_inertia_scheme}) descends to the morphism $ I\cX \to \cX $ in (\ref{eq_inertia_stack}). In other words, 
	we have a 2-commutative cube
	\vspace{12pt}
	\begin{equation*}
	\begin{tikzcd}[row sep=2em, column sep=0.4em]
	& I_G X \arrow[dl] \arrow[rr] \arrow[dd] & & X \arrow[dl] \arrow[dd] \\
	I \cX \arrow[rr, crossing over] \arrow[dd] & & \cX \\
	& G \times X \arrow[dl] \arrow[rr] & & X \times X \arrow[dl] \\
	\left[G/G\right] \times \cX \arrow[rr] & & \cX \times \cX \arrow[from=uu, crossing over]\\
	\end{tikzcd}.
	\end{equation*}
\end{rmk}
Let $ \cX = [X/G] $ be a separated quotient DM stack. Recall from Construction \ref{constrn_equiv_sheaves} that we can twist sheaves on $ \cX $ by representations of $ G. $ This defines a (right) action of the representation ring $ R(G) = K^0(BG) $ on $ K(\cX), $ which makes $ K(\cX) $ a (right) module\footnote{In general, a morphism $ \cX \to \cY $ of algebraic stacks make $ K(\cX) $ a module over $ K^0(\cY) $ by tensoring a sheaf on $ \cX $ with the pullback of a vector bundle on $ \cY. $} over $ R(G), $ i.e.,
\begin{equation*}
[\cE] \cdot [\rho] = [\cE \otimes \rho]
\end{equation*}
for a sheaf $ \cE $ on $ \cX $ and a representation $ \rho $ of $ G. $
Tensoring by $ \CC, $ we then have a module $ K(\cX)_\CC $ over $ R(G)_\CC. $

Let $ \Supp K(\cX)_\CC $ denote the support of the $ R(G)_\CC $-module $ K(\cX)_\CC $. We will follow the treatment in \cite[Section 4.3]{edidin2013riemann} to decompose the inertia stack $ I\cX $ via $ \Supp K(\cX)_\CC. $ We first state a basic fact which combines two results in \cite[Proposition 2.5]{edidin2005nonabelian} and \cite[Remark 5.1]{edidin2000riemann}.
\begin{lem}
	Let $ G $ be a linear algebraic group. There is a bijection between the set of conjugacy classes of diagonalizable elements in $ G $ and the set of closed points in $ \Spec R(G)_\CC. $ If $ G $ acts on a scheme $ X $ with finite stabilizers, then the support of $ K(\cX)_\CC $ as a module over $ R(G)_\CC $ is finite.
\end{lem}
\begin{rmk}
	Recall that the linear algebraic group $ G $ is a subgroup of $ \GL(n,\CC). $ Each element $ g $ of $ G $ has a (unique) Jordan-Chevalley decomposition $ g = g_s + g_n $ where $ g_s \in G $ is semisimple, i.e., diagonalizable, $ g_n \in G $ is nilpotent, i.e., some power of it vanishes, and $ g_s g_n = g_n g_s. $ In particular, if $ G $ is finite, then each elements $ g $ of $ G $ is diagonalizable, i.e., $ g_n = 0 ,$ since $ g^k = 1 $ for some finite $ k $ and this implies the minimal polynomial of $ g $ divides $ x^k - 1 $ in $ \CC[x] $ and hence has distinct roots.
\end{rmk}
By the lemma above, $ K(\cX)_\CC $ is supported at a finite number of closed points $ p_0, \dots, p_l $ in $ \Spec R(G)_\CC $ which correspond to a finite number of conjugacy classes $ [g_0] = [1], \dots, [g_l] $ of diagonalizable elements in $ G. $ If $ [g] $ is one of these classes, then the fixed point locus $ X^g $ is non-empty, and hence it must be of finite order. On the contrary, if $ X^g $ is non-empty for an element $ g $, then the corresponding closed point $ \mathfrak{m}_g $ in $ \Spec R(G)_\CC $ is in the support of $ K(\cX)_\CC $. Therefore, we have the following:
\begin{prop}
	Let $ \cX = [X/G] $ be a separated quotient DM stack. Then the support of the $ R(G)_\CC $-module $ K(\cX)_\CC $ can be identified with a finite set of conjugacy classes of elements of finite order in $ G. $ Furthermore, $ X^g $ is non-empty if and only if $ [g] $ is one of those conjugacy classes.
\end{prop}
Now we assume $ \cX $ is connected. The we have a decomposition
\begin{equation*}
	I\cX \simeq \coprod_{i = 0}^{l} \{g_i\} \times \left[X^{g_i} / Z_{g_i}\right],
\end{equation*}
where $ \simeq $ means an equivalence of stacks as categories, $ g_0 = 1 $ is the identity in $ G, $ $ X^{g_i} $ is the fixed point locus in $ X $ under the action of $ g_i, $ and $ Z_{g_i} $ is the centralizer of $ g_i $ in $ G $. Note that each quotient stack $ \{g_i\} \times \left[X^{g_i} / Z_{g_i}\right] $ is not connected in general, but can be further decomposed. We can decompose $ I\cX $ into connected components as follows:
\begin{equation*}
	I\cX = \cX_0 \coprod I_t \cX, \quad I_t\cX = \coprod_{i=1}^{l}\coprod_{j=1}^{m_i} \{g_i\} \times [X_{ij}/G_{ij}], 
\end{equation*}
where $ \cX_0 = \{1\} \times \cX$, $ X_{ij} \subseteq X^{g_i} $ (i.e., $ g_i $ acts on $ X_{ij} $ trivially), and $ G_{ij} \subseteq Z_{g_i}. $ Note that the index $ g_i $ is important, as we will use it to define the orbifold Chern character and the orbifold Todd class later. Re-ordering these components by a single index set $ I $ of size
\begin{equation*}
	|I| = 1 + m_1 + \cdots + m_l,
\end{equation*}
we obtain the following:
\begin{prop}\label{prop_decomp_inertia}
	Let $ \cX = [X/G] $ be a connected separated quotient stack. Then its inertia stack has a decomposition
	\begin{equation*}
	I\cX \simeq \coprod_{i \in I} \cX_i
	\end{equation*}
	of connected components, where $ i $ runs over a finite index set $ I $ containing $ 0 $, and each $ \cX_i $ is a connected separated quotient stack of the form 
	\begin{equation*}
	\cX_i = \{g_i\} \times [X_i /G_i],
	\end{equation*}
	where $ X_i $ is a closed subscheme of $ X $ invariant under the action of a subgroup $ G_i \subset Z_{g_i} $ and fixed by an element $ g_i \in G_i $ with $ X_0 = X, $ $ G_0 = G $, and $ g_0 = 1. $
\end{prop}
\begin{rmk}
	The finite index set $ I $ in Proposition \ref{prop_decomp_inertia} is called the \tb{inertia index set} of $ \cX. $ The component $ \cX_0 \cong \cX $ and is called the \tb{distinguished component} and the other components $ \cX_i $ in $ I_t \cX $ are called \tb{twisted sectors}. If $ G $ acts on $ X $ freely, then there are no twisted sectors in $ I\cX $. 
\end{rmk}

Now we compute some examples of inertia stacks.
\begin{ex}
	Let $ \cX = BG = [\pt /G] $ be the classifying stack of a finite group $ G $. Then $ I_G X = G $, and the inertia stack of the classifying stack is $ I\cX = [G/G] $ where $ G $ acts on itself by conjugation. Let $ I = \{0, 1, \dots, n\} $ indexes the conjugacy classes $ [g_0], [g_1], \dots, [g_n] $ of $ G $ where $ g_0 = 1 $ is the identity in $ G $. Then we have a decomposition
	\begin{equation*}
	IBG = \coprod_{i \in I} \{g_i\} \times BZ_{g_i}.
	\end{equation*}
	Here $ BZ_{g_0} = BG $ is the distinguished component. If $ G $ is abelian, then every element $ g $ is a conjugacy class and each $ Z_{g} = G, $ so we have $$ IBG = G \times BG. $$
\end{ex}
\begin{ex}
	Let $ \cX = [\CC^2/\mu_n] $ where the $ \mu_n $-action on $ \CC^2 $ is given by $ \zeta \cdot (x_1, x_2) = (\zeta x_1, \zeta^{-1} x_2) $. Then $ I_G X = \left(\{1\} \times \CC^2\right) \coprod \left(\left(\mu_n-\{1\}\right) \times \{(0,0)\}\right), $ and the twisted sectors are
	\begin{equation*}
	I_t[\CC^2/\mu_n] = \left(\mu_n-\{1\}\right) \times B\mu_n.
	\end{equation*}
\end{ex}
\begin{ex}
	Let $ X = \CC^2 - \{(0,0)\} $, and let $ G = \CC^* $ act on $ X $ by $ \lambda \cdot (x_0, x_1) = (\lambda^2 x_0, \lambda^3 x_1) $. The quotient stack $ \cX = [X/G] $ is the weighted projective line $ \PP(2,3) $. We have $ K(BC^*)_\CC \cong \CC[x, x^{-1}], $ and
	\begin{equation*}
		K(\PP(2,3))_\CC \cong \frac{\CC[x, x^{-1}]}{\inprod{(1-x^{-2})(1-x^{-3})}}.
	\end{equation*}
	Therefore, we have
	\begin{equation*}
		\Supp K(\PP(2,3))_\CC = \{g_0, g_1, g_2, g_3\},
	\end{equation*}
	where $ g_0 = 1, g_1 = -1, g_2 = (-1+\sqrt{3}i)/2, $ and $ g_3 = (-1-\sqrt{3}i)/2 $. Now we compute the fixed point loci: $ X^{g_0} = X,\ X^{g_1} = \CC^* \times \{0\}, $ and $\ X^{g_2} = X^{g_3} = \{0\} \times \CC^* $. The centralizer $ Z_{g_i} = \CC^* $ for $ i = 0, 1, 2, 3. $ Therefore, we have a decomposition
	\begin{align*}
		I\PP(2,3) = \coprod_{i=0}^3 \{g_i\} \times [X^{g_i}/Z_{g_i}] \cong \{1\} \times \PP(2,3) \coprod \{-1\} \times B\mu_2 \coprod \left\{\frac{-1 \pm \sqrt{3}i}{2}\right\} \times B\mu_3.
	\end{align*}
\end{ex}

\vspace{1pt}

\subsection{Orbifold Chern character and orbifold Todd class}
In this section we define orbifold Chern characters and orbifold Todd classes of coherent sheaves on a connected separated smooth quotient DM stack.

We first review the Chow ring of a separated smooth quotient DM stack following \cite{edidin1998equivariant}.
\begin{defn-prop}
	Let $ \cX = [X/G] $ be a quotient stack. For each integer $ i \geq 0, $ the $ i $-th \tb{Chow group} of $ \cX $ is defined by
	$$ A^i(\cX) = A^i((X \times U)/G) $$ 
	where $ U $ is an open subscheme of codimension $ > i $ in an affine space $ \CC^n $ with a linear $ G $-action such that it restricts to a free $ G $-action on $ U $. The definition of $ A^i(\cX) $
	is independent of the choice of $ n, $ the $ G $-action on $ \CC^n $ and the open subscheme $ U. $ If $ \cX $ is a smooth quotient stack, then there is an intersection product 
	\begin{equation*}
	A^i(\cX) \times A^j(\cX) \to A^{i+j}(\cX)
	\end{equation*}
	for all non-negative integers $ i $ and $ j $ such that $ A^*(\cX) $ becomes a ring graded in $ \ZZ_{\geq 0}. $ We denote by $ A^*(\cX) $ or simply $ A(\cX) $ the \tb{Chow ring} of $ \cX. $ If $ \cX $ is a separated smooth quotient DM stack, then the proper pushforward 
	\begin{equation*}
	A(\cX)_\QQ \to A(X/G)_\QQ
	\end{equation*}
	(here the geometric quotient $ X/G $ is the same as the coarse moduli space of $ \cX $) 
	is an isomorphism of groups in each degree, which implies integral elements of each $ A^i(\cX)_\QQ $ are classes of substacks $ \cY = [Y/G] $ of $ \cX $ of codimension $ i $ represented by $ G $-invariant cycles $ Y $ on $ X $, and $ A^i(\cX) $ vanishes for $ i > \dim \cX = \dim X - \dim G. $ If $ \cX $ is a proper smooth quotient stack of dimension $ n $, then there is a group homomorphism
	\begin{equation*}
	\deg: A^n(\cX) \to \QQ
	\end{equation*}
	defined by
	\begin{equation*}
	\deg \, [Gx/G] = \frac{1}{|G_x|}
	\end{equation*}
	for a point $ [Gx/G] $ on $ \cX $ corresponding to an orbit $ Gx $ on $ X $ where $ G_x = \Stab_G(x). $ 
\end{defn-prop}
\begin{ex}
	The Chow ring of $ B\CC^* $ is given by
	\begin{equation*}
	A(B\CC^*) = \ZZ[[h]]
	\end{equation*}
	where $ h = c_1(\cO_{\PP^n}(1)) $ in $ A^1(B\CC^*) = A^1(\PP^n) $ for any $ n \geq 1. $
\end{ex}
\begin{ex}
	The Chow ring of the weighted projective stack $ \PP = \PP(a_0, \dots, a_n) $ is given by
	\begin{equation*}
	A(\PP) = \frac{\ZZ[h]}{(a_0\cdots a_nh^{n+1})}
	\end{equation*}
	where $ h = c_1(\cO_\PP(1)) $ in $ A^1(\PP). $
\end{ex}
\begin{rmk}
	If $ \cX $ is a separated smooth quotient DM stack, then so is its inertia $ I\cX $.
\end{rmk}
Let $ \cX $ be a connected separated smooth quotient DM stack. There is an orbifold Chern character map
$$ \orbch : K^0(\cX) \to A(I\cX)_{\CC}, $$
which is defined in \cite{toen1999theoremes}.
There is also an orbifold Todd class map
$$ \orbtd : K^0(\cX) \to A(I\cX)_{\CC}^\times, $$
which is given in \cite[Appendix A]{tseng2010orbifold} and \cite[Section 2.4]{iritani2009integral}. Since $ K(\cX) \cong K^0(\cX) $, we have two maps:
$$ \orbch : K(\cX) \to A(I\cX)_{\CC} \quad \text{and} \quad \orbtd : K(\cX) \to A(I\cX)_{\CC}^\times. $$ 
Now we give detailed constructions for these two maps for quotient DM stacks.

\begin{notation}\label{notation_quotient_stack}
	Let $ \cX = [X/G] $ be a connected separated smooth quotient DM stack. Recall the decomposition 
	$$ I\cX = \coprod_{i \in I} \cX_i = \coprod_{i \in I} \{g_i\} \times [X_i/G_i] $$ 
	in Proposition \ref{prop_decomp_inertia}, we have two rings
	\begin{equation*}
	K^0(I\cX) = \bigoplus_{i \in I} K^0(\cX_i) \quad \text{and} \quad A(I\cX) = \bigoplus_{i \in I} A(\cX_i).
	\end{equation*}
	Let $ \pi: I\cX \to \cX $ denote the natural morphism. For each $ i \in I, $ let $ p_i: X_i \to X $ and $ q_i: \cX_i \to \cX $ denote the inclusions of the subscheme $ X_i $ in $ X $ and the substack $ \cX_i $ in $ \cX $ respectively; also let $ r_i: \cX_i \to I\cX $ denote the inclusion of $ \cX_i $ in $ I\cX. $ Note that $ \pi \circ r_i = q_i $ for each $ i \in I. $
\end{notation}

We first prove a fundamental lemma that is crucial for the constructions of both the orbifold Chern character and the orbifold Todd class.
\begin{lem}\label{lem_eigenbundle}
	Let $ [X/G] $ be a quotient stack. Let $ g $ be an element in the center of $ G $ such that it has a finite order $ n $ and fixes $ X $. Then every $ G $-equivariant vector bundle $ \cV = (V, \phi) $ on $ X $ admits a $ G $-equivariant decomposition
	\begin{equation}\label{eq_equiv_decomp}
	\cV = \bigoplus_{\lambda \in \eig(\phi(g))} \cV_\lambda
	\end{equation}
	where $ \eig(\phi(g)) $ consists of eigenvalues of the linear automorphism $ \phi(g) $ on $ V $ which are $ n $th roots of unity, and each $ \cV_\lambda = (V_\lambda, \phi_\lambda) $ is a $ G $-equivariant vector bundle on $ X. $
\end{lem}
\begin{proof}
	Consider a $ G $-equivariant vector bundle $ (V, \phi) $ on $ X. $ This corresponds to an action
	\begin{equation*}
	\phi: G \to \Aut(V)
	\end{equation*}
	of $ G $ on $ V. $ By assumption, $ g^n = 1 $, so the linear automorphism $ \phi(g) $ on $ V $ satisfies $ (\phi(g))^n = \id_V, $ and hence each eigenvalue of $ \phi(g) $ is a $ n $th root. Since $ g $ acts on $ X $ trivially, there is an eigenbundle decomposition
	$$ V = \bigoplus_{\lambda \in \eig(\phi(g))} V_\lambda $$
	such that for each $ v \in V_\lambda, $
	\begin{equation*}
	\phi(g) (v) = \lambda v.
	\end{equation*}
	Take any $ h \in G. $ Since $ g $ is in the center of $ G, $ we have $ gh = hg. $ For any $ v \in V_\lambda, $ we have
	\begin{equation*}
	\phi(g) \left(\phi(h)v\right) = \phi(h) \left(\phi(g)v\right) = \phi(h) \left(\lambda v\right) = \lambda \left(\phi(h)v\right),
	\end{equation*}
	which implies each $ V_\lambda $ is an invariant subbundle of $ V $ under the action of $ G $, i.e., $ \phi $ restricts to an action
	\begin{equation*}
	\phi_\lambda: G \to \Aut(V_\lambda)
	\end{equation*}
	of $ G $ on each $ V_\lambda $. Therefore, we have the $ G $-equivariant decomposition of $ (V, \phi) $ claimed in the lemma.
\end{proof}
Applying Lemma \ref{lem_eigenbundle} to $ [X/G] = BG $, we obtain a result in representation theory.
\begin{cor}\label{cor_rep_decomp}
	Let $ G $ be a linear algebraic group and let $ g $ be an element in the center of $ G $. Then every representation $ \phi: G \to \GL(V) $ of $ G $ admits a decomposition
	\begin{equation}\label{eq_rep_decomp}
	\phi = \bigoplus_{\lambda \in \eig(\phi(g))} \phi_\lambda
	\end{equation}
	into subrepresentations $ \phi_\lambda: G \to \GL(V_\lambda). $
\end{cor}
\begin{rmk}
	In Lemma \ref{lem_eigenbundle} and Corollary \ref{cor_rep_decomp}, each $ \cV_\lambda $ or $ \phi_\lambda $ is not irreducible in general. If $ g = 1 $ is the identity of $ G $, then any decomposition (\ref{eq_equiv_decomp}) or (\ref{eq_rep_decomp}) is trivial.
\end{rmk}
\begin{defn}
	Let $ \cX = [X/G] $ be a quotient stack. Let $ g $ be an element in the center of $ G $ such that it has a finite order and fixes $ X $.
	The $ \boldsymbol{g} $\tb{-twisting morphism}
	\begin{equation*}
	\rho_g: K^0(\cX) \to K^0(\cX)_\CC
	\end{equation*}
	is defined by
	\begin{equation*}
	\rho_g ([\cV]) = \sum_{\lambda \in \eig(\phi(g))} \lambda [\cV_\lambda] 
	\end{equation*}
	for the decomposition (\ref{eq_equiv_decomp}) of a vector bundle $ \cV $ on $ \cX $ and extended linearly. 
\end{defn}
\begin{rmk}
	The twisting morphism $ \rho_g $ depends on the element $ g. $ If $ g = 1 $ is the identity of $ G, $ then $ \rho_g $ is the identity map.
\end{rmk}
\begin{lem}\label{lem_twisting_morphism}
	The twisting morphism $ \rho_g: K^0(\cX) \to K^0(\cX)_\CC $ is a ring homomorphism and induces a $ \CC $-algebra automorphism on $ K^0(\cX)_\CC. $
\end{lem}
\begin{proof}
	Consider a pair of vector bundles $ \cV = (V, \phi) $ and $ \cW = (W, \psi) $ on $ \cX. $ We have
	\begin{equation*}
	\cV \oplus \cW = (V \oplus W, \phi \oplus \psi) \quad \text{and} \quad \cV \otimes \cW = (V \otimes W, \phi \otimes \psi).
	\end{equation*}
	The linear automorphism $ (\phi \oplus \psi)(g) $ on $ \cV \oplus \cW $ has a set of eigenvalues
	\begin{equation*}
	\{\alpha + \beta: \alpha \in \eig(\phi(g)), \beta \in \eig(\psi(g))\},
	\end{equation*}
	the linear automorphism $ (\phi \otimes \psi)(g) $ on $ \cV \otimes \cW $ has a set of eigenvalues
	\begin{equation*}
	\{\alpha \beta: \alpha \in \eig(\phi(g)), \beta \in \eig(\psi(g))\},
	\end{equation*}
	and we have decompositions
	\begin{equation*}
	\cV \oplus \cW = \bigoplus_{\lambda \in \eig((\phi\oplus\psi)(g))} \left(\cV \oplus \cW\right)_\lambda = \bigoplus_{\substack{\alpha \in \eig(\phi(g)) \\ \beta \in \eig(\psi(g))} } \left(\cV_\alpha \oplus \cW_\beta\right)
	\end{equation*}
	and 
	\begin{equation*}
	\cV \otimes \cW = \bigoplus_{\lambda \in \eig((\phi\otimes\psi)(g))} \left(\cV \otimes \cW\right)_\lambda = \bigoplus_{\substack{\alpha \in \eig(\phi(g)) \\ \beta \in \eig(\psi(g))} } \left(\cV_\alpha \otimes \cW_\beta\right).
	\end{equation*}
	Therefore, we have
	\begin{equation*}
	\rho_g([\cV] + [\cW]) = \rho_g([\cV]) + \rho_g([\cW]) \quad \text{and} \quad \rho_g([\cV] [\cW]) = \rho_g([\cV]) \rho_g([\cW]).
	\end{equation*}
	By the linearity of $ \rho_g, $ we then have
	\begin{equation*}
	\rho_g(x+y) = \rho_g(x) + \rho_g(y) \quad \text{and} \quad \rho_g(xy) = \rho_g(x) \rho_g(y)
	\end{equation*}
	for all $ x $ and $ y $ in $ K^0(\cX). $
\end{proof}
\begin{defn}\label{defn_twisted_Chern}
	Let $ \cX = [X/G] $ be a separated smooth quotient DM stack. Let $ g $ be an element in the center of $ G $ such that it has a finite order and fixes $ X $. We define the $ \boldsymbol{g} $ \tb{-twisted Chern character} map by
	\begin{equation*}
	\ch^{\rho_g}: K^0(\cX) \xrightarrow{\rho_g} K^0(\cX)_\CC \xrightarrow{\ch} A(\cX)_\CC,
	\end{equation*}
	where $ \ch: K^0(\cX)_\CC \to A(\cX)_\CC $ is the Chern character map on the stack $ \cX $ tensored by $ \CC $.
\end{defn}
\begin{rmk}
	If $ g = 1 $ is the identity of $ G, $ then $ \rho_g $ is trivial, and hence the twisted Chern character map reduces to the ordinary Chern character map
	\begin{equation*}
	\ch: K^0(\cX) \to A(\cX)_\QQ
	\end{equation*}
	for the quotient stack $ \cX = [X/G], $ which is the $ G $-equivariant Chern character map
	\begin{equation*}
	\ch: K_G^0(X) \to A_G(X)_\QQ
	\end{equation*}
	for the $ G $-scheme $ X $.
\end{rmk}
Lemma \ref{lem_twisting_morphism} immediately implies the following:
\begin{lem}\label{lem_twisted_Chern}
	The twisted Chern character map $ \ch^{\rho_g}: K^0(\cX) \to A(\cX)_\CC $ is a ring homomorphism.
\end{lem}
\begin{ex}\label{ex_twisted_Chern_on_BG}
	Consider the classifying stack $ BG $ for a finite group $ G. $ Then $ K^0(BG) = K(BG) $ is the representation ring $ R(G). $ Take an element $ g $ in the center of $ G. $ Let's compute the $ g $-twisted Chern character
	\begin{equation*}
	\ch^{\rho_{g}}: R(G) \to A(BG)_\CC \cong \CC.
	\end{equation*}
	A representation $ \phi: G \to \GL(V) $ of $ G $ on a vector space $ V $ decomposes as
	\begin{equation*}
	\phi = \bigoplus_{\lambda \in \eig(\phi(g))} \phi_\lambda
	\end{equation*}
	with subrepresentations $ \phi_\lambda: G \to \GL(V_\lambda), $ so we have a twisting morphism
	\begin{equation*}
	\rho_{g}: R(G) \to R(G)_\CC, \quad \phi \mapsto \sum_{\lambda \in \eig(\phi(g))} \lambda \phi_\lambda.
	\end{equation*}
	The Chern character map $ \ch: R(G)_\CC \to \CC $ is a linear extension of the rank map, i.e., for a virtual representation $ \sum_i a_i \vphi_i $ of $ G $ where $ a_i \in \CC $ and representation $ \vphi_i: G \to \GL(V_i) $, we have
	\begin{equation*}
	\ch\left(\sum_i a_i \vphi_i\right) = \sum_i a_i \dim V_i
	\end{equation*}
	in $ \CC. $ Therefore, the twisted Chern character map is given by 
	\begin{equation*}
	\ch^{\rho_{g}}(\phi) = \ch\left(\rho_g(\phi)\right) = \sum_{\lambda \in \eig(\phi(g))} \lambda \dim V_\lambda = \tr(\phi(g)) = \chi_\phi(g)
	\end{equation*}
	in $ \CC, $ where $ \chi_\phi $ denotes the character of the representation $ \phi $ of $ G. $ The fact that the twisted Chern character map is a ring map implies that for any pair of representations $ \phi $ and $ \psi $ of $ G, $ we have
	\begin{equation*}
	\chi_{\phi \oplus \psi}(g) = \chi_{\phi}(g) + \chi_{\psi}(g), \quad \text{and} \quad \chi_{\phi \otimes \psi}(g) = \chi_{\phi}(g) \chi_{\psi}(g),
	\end{equation*}
	which is a basic fact in character theory.
\end{ex}

\begin{defn-prop}
	Let $ \cX = [X/G] $ be a connected separated smooth quotient DM stack. We have a finite inertia index set $ I $ and a decomposition 
	$$ I\cX \simeq \coprod_{i \in I} \cX_i = \coprod_{i \in I} \{g_i\} \times [X_i/G_i]. $$ 
	Each $ g_i $ is in the center of $ G_i $ (because $ G_i $ is contained in the centralizer $ Z_{g_i} $ of $ g_i $), it fixes $ X_i $ and has a finite order, so it induces a $ g_i $-twisting morphism
	\begin{equation*}
		\rho_{g_i}: K^0(\cX_i) \to K^0(\cX_i)_\CC,
	\end{equation*}
	and hence a $ g_i $-twisted Chern character map
	\begin{equation*}
		\ch^{\rho_{g_i}}: K^0(\cX_i) \to A(\cX_i)_\CC.
	\end{equation*}
	Let $ \rho: K^0(I\cX) \to K^0(I\cX)_\CC $ denote the \tb{twisting morphism} on $ K^0(I\cX) $, i.e.,
	\begin{equation*}
		\rho = \bigoplus_{i \in I} \rho_{g_i}.
	\end{equation*}
\end{defn-prop}

Now we have all the ingredients to define the orbifold Chern character map.
\begin{defn}\label{defn_orbifold_Chern}
	Let $ \cX = [X/G] $ be a connected separated smooth quotient DM stack. 
	The \tb{inertia Chern character} map is defined by the composition
	$$ \ch^\rho: K^0(I\cX) \xrightarrow{\rho} K^0(I\cX)_\CC \xrightarrow{\ch} A(I\cX)_\CC, $$
	i.e.,
	\begin{equation*}
	\ch^\rho(x) = \left(\ch(x_0), \ \bigoplus_{0 \neq i \in I} \ch^{\rho_{g_i}}(x_i)\right) \vspace{8pt}
	\end{equation*}
	for any $ x = \bigoplus_{i \in I} x_i $ in $ K^0(I\cX). $
	The \tb{orbifold Chern character} map on $ K^0(\cX) $ is defined by the composition
	\begin{equation*}
	\orbch : K(\cX) \xrightarrow{\beta} K^0(\cX) \xrightarrow{\pi^*} K^0(I\cX) \xrightarrow{\ch^\rho} A(I\cX)_\CC.
	\end{equation*}
\end{defn}
\begin{prop}\label{prop_orbch_ring_map}
	The orbifold Chern character map $ \orbch : K^0(\cX) \to A(I\cX)_{\CC} $ is a ring homomorphism, i.e.,
	\begin{equation*}
	\quad \orbch(x+y) = \orbch(x) + \orbch(y) \quad \text{and} \quad \orbch(xy) = \orbch(x) \orbch(y)
	\end{equation*}
	for all $ x $ and $ y $ in $ K(\cX). $
\end{prop}
\begin{proof}
	Lemma \ref{lem_twisted_Chern} implies the inertia Chern character map $ \ch^\rho: K^0(I\cX) \to  A(I\cX)_\CC $ is a ring homomorphism, where both $ K^0(I\cX) $ and $ A(I\cX)_\CC $ carry ring structures from the direct sum of the rings $ K^0(\cX_i) $ and $ A(\cX_i)_\CC $ over $ I. $ The orbifold Chern character $ \orbch : K(\cX) \to A(I\cX)_{\CC} $ is a ring homomorphism since it's a composition of three ring homomorphisms $ \beta, \pi^* $ and $ \ch^\rho $.
\end{proof}
\begin{rmk}[Explicit formulas for the orbifold Chern character]
	The orbifold Chern character of a vector bundle $ \cV = (V, \phi) $ on $ \cX $ is given by
	\begin{equation}\label{eq_orbch_concrete}
	\orbch (\cV) = \left(\ch(\cV),\ \bigoplus_{0 \neq i \in I} \ch^{\rho_{g_i}}(\cV_i)\right) = \left(\ch(\cV), \  \bigoplus_{0 \neq i \in I} \sum_{\lambda \in \eig(\phi_i(g_i))} \lambda \ch (\cV_{i,\lambda})\right), \vspace{5pt}
	\end{equation}
	where each $ \cV_i 
	= r_i^*\pi^*\cV = q_i^*\cV $ is the restriction of $ \cV $ from $ \cX $ to the substack $ \cX_i, $ with an eigenbundle decomposition $ \cV_i = \oplus_{\lambda \in \eig(\phi_i(g_i))} \cV_{i,\lambda} $ for $ 0 \neq i \in I $. For a sheaf $ \cE $ on $ \cX, $ choose any finite locally free resolution $ \cE_{\boldsymbol{\cdot}} \to \cE \to 0 $, then we have
	\begin{equation}\label{eq_orbch_sheaf}
	\orbch(\cE) = \sum_k (-1)^k \orbch(\cE_k).
	\end{equation}
\end{rmk}

\begin{rmk}[Care is required when working with sheaves.]
	Consider two sheaves $ \cE, \cF $ on $ \cX. $ Recall from Lemma \ref{lem_multiplication_sheaves} that their product is given by
	\begin{equation*}
	[\cE] [\cF] = \sum_k (-1)^k [\cE_k \otimes \cF]
	\end{equation*}
	for any finite locally free resolution $ \cE_{\boldsymbol{\cdot}} \to \cE \to 0. $
	The identity 
	$$ \orbch([\cE] [\cF]) = \orbch([\cE]) \orbch([\cF]) $$ 
	implies
	\begin{equation*}
	\orbch(\cE) \orbch(\cF) = \sum_k (-1)^k \orbch(\cE_k \otimes \cF),
	\end{equation*}
	but
	\begin{equation*}
	\orbch(\cE) \orbch(\cF) \neq \orbch(\cE \otimes \cF)
	\end{equation*}
	in general, since $ \cE_{\boldsymbol{\cdot}} \otimes \cF \to \cE \otimes \cF \to 0 $ may fail to be exact. However, if either $ \cE $ or $ \cF $ is locally free, then we have
	\begin{equation*}
	\orbch(\cE) \orbch(\cF) = \orbch(\cE \otimes \cF).
	\end{equation*}
	
	We also have a $ K $-theoretic pullback
	\begin{equation*}
	\pi^{K}: K(\cX) \xrightarrow{\beta} K^0(\cX) \xrightarrow{\pi^*} K^0(I\cX) \xrightarrow{\beta^{-1}} K(I\cX)
	\end{equation*}
	via the natural isomorphisms $ K(\cX) \xrightarrow{\sim} K^0(\cX) $ and $ K(I\cX) \xrightarrow{\sim} K^0(I\cX) $. For a sheaf $ \cE $ on $ \cX, $ its $ K $-theoretic pullback is then given by
	\begin{equation*}
	\pi^K [\cE] = \sum_k (-1)^k \pi^K[\cE_k] = \sum_k (-1)^k [\pi^*\cE_k]
	\end{equation*}
	for any finite locally free resolution $ \cE_{\boldsymbol{\cdot}} \to \cE \to 0 $ of $ \cE, $
	but
	$$ \pi^K[\cE] \neq [\pi^*\cE] $$ in general since the pullback functor $ \pi^*: \Coh(\cX) \to  \Coh(I\cX) $ is not exact (although $ \pi^*: \Vect(\cX) \to \Vect(I\cX) $ is exact) in general. 
\end{rmk}
Recall that for a separated smooth scheme $ X $ of finite type, the Chern character map $ \ch: K(X)_\QQ \xrightarrow{\sim} A(X)_\QQ $
is a $ \QQ $-algebra isomorphism after tensored with $ \QQ. $ A similar result holds for quotients stacks.
\begin{prop}\label{prop_orbch_isom}
	Let $ \cX $ be a connected separated smooth quotient DM stack. The orbifold Chern character map
	$$ \orbch : K(\cX)_\CC \xrightarrow{\sim} A(I\cX)_{\CC} $$
	is a $ \CC $-algebra isomorphism after tensored with $ \CC. $
\end{prop}
\begin{proof}
	Note that $ \cX $ is a separated DM stack of finite type over $ \CC $. By the discussion on page 9 in \cite{toen2000motives}, after tensored with $ \QQ, $ the orbifold Chern character $ \orbch: K(\cX)_\QQ \to A_\chi(\cX) $ is a $ \QQ $-algebra isomorphism onto $ A_\chi(\cX) $, the rational Chow ring with coefficients in the characters of $ \cX $. Since $ A_\chi(\cX)_\CC \cong A(I\cX)_{\CC}, $ the map $ \orbch $ becomes a $ \CC $-algebra isomorphism after tensored with $ \CC. $
\end{proof}
\begin{rmk}
	The ordinary Chern character $ \ch = \orbch_0: K^0(\cX)_\CC \to A(\cX)_\CC  $ is a surjective $ \CC $-algebra homomorphism, but not an isomorphism in general, as the following example shows.
\end{rmk}
\begin{ex}[The orbifold Chern character for $ BG $ is the character of representations of $ G $.]\label{ex_orbch_BG}
	Let $ \cX = BG $ for a finite group $ G $ with $ n $ conjugacy classes 
	$$ [g_0], [g_1], \dots, [g_{n-1}]. $$
	Recall that 
	$$ IBG = \coprod_{i=0}^{n-1} \{g_i\} \times BZ_{g_i}. $$ 
	Take a vector bundle $ \cV = (V, \phi) $ on $ BG, $ i.e., a linear representation $\phi: G \to \GL(V) $
	of $ G $. 
	On each component $ \{g_i\} \times BZ_{g_i} $, the representation $ \phi $ of $ G $ restricts to a representation
	\begin{equation*}
	\phi_i: Z_{g_i} \to \GL(V)
	\end{equation*}
	of $ Z_{g_i} $.
	Since each $ g_i $ is in the center of $ Z_{g_i}, $ by the computations in Example \ref{ex_twisted_Chern_on_BG}, we know that each $ g_i $-twisted Chern character map is
	\begin{equation*}
	\ch^{\rho_{g_i}}: R(Z_{g_i}) \to \CC, \quad \vphi \mapsto \chi_\vphi(g_i).
	\end{equation*}
	Therefore, the orbifold Chern character becomes
	\begin{alignat*}{3}
	\orbch: & R(G) & \ \longrightarrow & \ \bigoplus\limits_{i=0}^{n-1} R(Z_{g_i}) & \ \longrightarrow & \ \CC^n \\
	& \quad \phi & \ \longmapsto & \ (\phi_0, \dots, \phi_{n-1}) & \ \longmapsto & \ \left( \chi_\phi(g_0), \cdots, \chi_\phi(g_{n-1}) \right),
	\end{alignat*}
	where each $ \chi_\phi(g_i) = \chi_{\phi_i}(g_i) $ since $ \phi_i $ is a restriction of $ \phi $.
	Observe that the zeroth component 
	$$ \chi_\phi(g_0) = \chi_\phi(1) = \deg \phi = \dim V $$ 
	is the ordinary Chern character of $ \phi. $
	The space $ \CC^n $ can be identified with the space $ C(G) $ of class functions on $ G. $ 
	Therefore, the orbifold Chern character coincides with the character map
	\begin{align*}
	\orbch: R(G) & \rightarrow C(G) (\cong \CC^n) \\
	\phi & \mapsto \chi_\phi.
	\end{align*}
	This is an injective ring map, where irreducible representations map to irreducible characters which form a basis of $ C(G) $. Tensoring $ R(G) $ with $ \CC $, we obtain a $ \CC $-algebra isomorphism
	$$ \orbch: R(G)_\CC \xrightarrow{\sim} C(G). $$
	Compare with the ordinary Chern character map (i.e., the rank map)
	\begin{align*}
	\ch: R(G) \rightarrow \ZZ, \quad (V, \phi) \mapsto \dim V,
	\end{align*}
	which is surjective but not injective in general.
\end{ex}
\begin{ex}[The orbifold Chern character map on $ B\mu_n $ is the inverse discrete Fourier transform]\label{ex_discrete_FT}
	Take $ G = \mu_n $ in Example \ref{ex_orbch_BG}. Then its representation ring $ R(\mu_n) = \ZZ[x]/(x^n-1) $
	where $ x $ is the character
	$$ \mu_n \to \CC^*, \quad \omega = e^{2\pi i/n} \mapsto \omega $$
	of $ \mu_n. $ The orbifold Chern character for $ \mu_n $ is the ring map
	\begin{align*}
	\orbch: \frac{\ZZ[x]}{(x^n-1)} \to \CC^n, \quad x \mapsto (1, \omega, \dots, \omega^{n-1}).
	\end{align*}
	Tensoring $ R(\mu) $ with $ \CC, $ we obtain a $ \CC $-algebra isomorphism 
	$$ \orbch: \frac{\CC[x]}{(x^n-1)} \xrightarrow{\sim} \CC^n. $$ 
	Here the inverse map $ \orbch^{-1}: \CC^n \to \CC[x]/(x^n-1) $ can be identified with the discrete Fourier transform as follows. For each $ k = 0, \dots, n-1, $ put
	\begin{equation*}
	e_k = \orbch(x^k) = (1, \omega^k, \dots, \omega^{(n-1)k}).
	\end{equation*}
	Consider the weighted inner product on $ \CC^n $ defined by
	\begin{equation*}
	\inprod{a,b}_W = \frac{1}{n} \sum_{i=0}^{n-1} \cj{a}_i b_i
	\end{equation*}
	in $ \CC $ for all $ a = (a_0, \dots, a_{n-1}) $ and $ b = (b_0, \dots, b_{n-1}) $ in $ \CC^n. $ Then the elements $ e_0, \dots, e_{n-1} $ form an orthonormal basis, i.e., for all $ i,j = 0, \dots, n-1, $ we have
	\begin{equation*}
	\inprod{e_i,e_j}_W = \delta_{ij},
	\end{equation*}
	where $ \delta_{ij} $ is the Kronecker delta. Take any function $ f: \ZZ/n\ZZ \to \CC. $ Identifying $ f $ with an element $ (f(0), \dots, f(n-1)) $ in $ \CC^n, $ we then have the orthogonal decomposition
	\begin{equation*}
	f = \sum_{k=0}^{n-1} \inprod{e_k, f}_W e_k, 
	\end{equation*}
	where
	\begin{equation*}
	\inprod{e_k, f}_W = \frac{1}{n} \sum_{j=0}^{n-1} \cj{\omega}^{jk} f(j) = \frac{1}{n} \sum_{j=0}^{n-1} e^{-2\pi ijk/n} f(j)
	\end{equation*}
	for $ k = 0, \dots, n-1. $ This defines a map
	\begin{equation*}
	\psi: \CC^n \to \frac{\CC[x]}{(x^n-1)}, \quad f \mapsto \hat{f} = \sum_{k=0}^{n-1} \inprod{e_k, f}_W x^k,
	\end{equation*}
	which can be identified with the \tb{discrete Fourier transform} (up to a scaling depending on  the convention). The map $ \psi $ is indeed the inverse of $ \orbch $ since $ \psi(e_k) = x^k $ for all $ k = 0, \dots, n-1. $ Therefore, the orbifold Chern character for $ \mu_n $ is nothing but the inverse discrete Fourier transform, which maps an element
	\begin{equation*}
	f = f(0) + f(1) x + f(2) x^2 + \dots + f(n-1) x^{n-1}
	\end{equation*}
	in $ \CC[x]/(x^n-1) $ to an element $ \check{f} = (\check{f}(0), \dots, \check{f}(n-1)) $ in $ \CC^{n} $ where
	\begin{equation*}
	\check{f}(k) = \sum_{j=0}^{n-1} \omega^{jk} f(j) = \sum_{j=0}^{n-1} e^{2\pi ijk/n} f(j)
	\end{equation*}
	for $ k = 0, \dots, n-1. $
\end{ex}
\begin{ex}
	Let's compute the orbifold Chern character map for the weighted projective line $ \PP = \PP(2,3). $ Recall that its inertia stack is 
	\begin{equation*}
	I\PP(2,3) = \cX_0 \coprod \cX_1 \coprod \cX_2 \coprod \cX_3,
	\end{equation*}
	where each $ \cX_i = \{g_i\} \times [X_i/\CC^*], $
	\begin{equation*}
	g_0 = 1, \quad g_1 = -1, \quad g_2 = (-1+\sqrt{3}i)/2, \quad g_3 = (-1-\sqrt{3}i)/2
	\end{equation*}
	and
	\begin{equation*}
	X_0 = X = \CC^2 - \{0\}, \quad x_1 = \CC^* \times \{0\}, \quad X_2 = X_3 = \{0\} \times \CC^*.
	\end{equation*}
	Note that
	\begin{equation*}
	\cX_0 \simeq \PP(2,3), \quad \cX_1 \simeq B\mu_2, \quad \cX_2 \simeq \cX_3 \simeq B\mu_3.
	\end{equation*}
	We also computed
	\begin{equation*}
	K(\PP(2,3)) = \frac{\ZZ[x]}{\inprod{(x^2-1)(x^3-1)}}
	\end{equation*}
	where $ x = [\cO_\PP(1)]. $ The sheaf
	$$ \cO_\PP(1) = \cO_\PP \otimes \rho = (\cO_X, \phi) $$
	is obtained from twisting the structure sheaf $ \cO_\PP $ by the identity character $ \rho: \CC^* \to \CC^* $ which corresponds to a $ \CC^* $-equivariant line bundle $ (\cO_X, \phi) $ on $ X = \CC^2 - \{0\}. $ Here $ \phi $ is the $ \CC^* $ action
	\begin{equation*}
	\phi: \CC^* \to \Aut(X \times \CC), \quad t \mapsto \left(\phi(t): (x_0,x_1,v) \mapsto (t^2x_0,t^3x_1, tv)\right).
	\end{equation*}
	on the total space of the structure sheaf $ \cO_X. $ Now the $ K $-theoretic pullback of $ x $ from $ \PP(2,3) $ to $ I\PP(2,3) $ is given by
	\begin{equation*}
	\pi^K(x) = \left(x, [\cO_{X_1}, \phi_1], [\cO_{X_2}, \phi_2], [\cO_{X_3}, \phi_3]\right)
	\end{equation*}
	in $ K(I\PP(2,3)), $ where each $ (\cO_{X_i}, \phi_i) $ is a $ \CC^* $-equivariant line bundle on $ X_i $ with the $ \CC^* $ actions restricted from $ \phi $ and given by
	\begin{equation*}
	\phi_1: \CC^* \to \Aut(X_1 \times \CC), \quad t \mapsto \left(\phi_1(t): (x_0,0,v) \mapsto (t^2x_0,0, tv)\right),
	\end{equation*}
	\begin{equation*}
	\phi_2: \CC^* \to \Aut(X_2 \times \CC), \quad t \mapsto \left(\phi_1(t): (0,x_1,v) \mapsto (0,t^3x_1, tv)\right),
	\end{equation*}
	and $ \phi_3 = \phi_2. $ These $ \CC^* $ line bundles correspond to three representations
	\begin{equation*}
	\rho_1 = \cO_{\cX_1} \otimes \rho = \rho|_{\mu_2}, \quad \rho_2 = \cO_{\cX_2} \otimes \rho = \rho|_{\mu_3}, \quad \rho_3 = \cO_{\cX_3} \otimes \rho = \rho|_{\mu_3}.
	\end{equation*}
	Therefore, we have
	\begin{align*}
	\orbch(x) & = \ch^\rho(x,\rho_1,\rho_2,\rho_3) \\
	& = \left(\ch(x), \chi_{\rho_1}(g_1), \chi_{\rho_2}(g_2), \chi_{\rho_3}(g_3)\right) \\
	& = \left((1+h), g_1, g_2, g_3\right)
	\end{align*}
	in the Chow ring
	$$ A(I\PP(2,3))_\CC \cong \frac{\CC[h]}{(h^2)} \oplus \CC^3. $$
	Tensoring $ K(\PP(2,3)) $ with $ \CC, $ we see that the orbifold Chern character map
	\begin{align*}
	\orbch: \frac{\CC[x]}{\inprod{(x^2-1)(x^3-1)}} & \xrightarrow{\simeq} \frac{\CC[h]}{(h^2)} \oplus \CC^3 \\
	x & \mapsto \left(1+h, -1, \frac{-1+\sqrt{3}i}{2}, \frac{-1-\sqrt{3}i}{2}\right)
	\end{align*}
	is a $ \CC $-algebra isomorphism by the Chinese remainder theorem.
\end{ex}
Next we will define the orbifold Todd class map.
\begin{defn}
	Let $ \cX = [X/G] $ be a quotient stack. Let $ g $ be an element in the center of $ G $ such that it has a finite order and fixes $ X $. For the decomposition (\ref{eq_equiv_decomp}) of a vector bundle $ \cV $ on $ \cX, $ we write
	\begin{equation*}
	\cV^\fix = \cV_1 \quad \text{and} \quad \cV^\mov = \bigoplus_{1 \neq \lambda \in \eig(\phi(g))} \cV_\lambda.
	\end{equation*}
	We define two projections
	\begin{equation*}
	P_f: K^0(\cX) \to K^0(\cX) \quad \text{and} \quad P_m: K^0(\cX) \to K^0(\cX)
	\end{equation*}
	by setting
	\begin{equation*}
	P_f(x) = \cV^\fix - \cW^\fix \quad \text{and} \quad P_m(x) = \cV^\mov - \cW^\mov
	\end{equation*}
	for any element $ x = [\cV] - [\cW] $ in $ K^0(\cX) $ where $ \cV $ and $ \cW $ are vector bundles on $ \cX. $
\end{defn}
\begin{rmk}
	It's easy to check both $ P_f $ and $ P_m $ are indeed projections, i.e., they are linear and idempotent. But they are not multiplicative in general. For example, take $ x = [\cV] = [\cV_{-1}] $ and $ y = [\cW] = [\cW_{-1}], $ i.e., $ \cV $ and $ \cW $ only have the eigenvalue $ -1 $ under the $ g $ action. Then $ P_f(x) = 0 = P_f(y), $ $ P_m(x) = x, $ $ P_m(y) = y, $ and we have
	\begin{equation*}
	P_f(xy) = xy \neq P_f(x) P_f(y) \quad \text{and} \quad P_m(xy) = 0 \neq P_m(x) P_m(y).	
	\end{equation*}
\end{rmk}
\begin{defn}
	Let $ \cX = [X/G] $ be a connected separated smooth quotient DM stack. Let $ g $ be an element in the center of $ G $ such that it has a finite order and fixes $ X $. The $ \boldsymbol{g} $ \tb{-twisted Euler class} map is defined by
	\begin{equation*}
	e^{\rho_g}: P_m(K^0(\cX)) \xrightarrow{e^K} K^0(\cX) \xrightarrow{\ch^{\rho_g}} A(\cX)_\CC
	\end{equation*}
	where $ P_m(K^0(\cX)) $ is the image of the projection $ P_m. $
\end{defn}
\begin{lem}
	The $ g $-twisted Euler class $ e^{\rho_g} $ is multiplicative and maps into the units of $ A(\cX)_\CC, $ i.e., we have a group homomorphism
	\begin{equation*}
	e^{\rho_g}: \left(P_m(K^0(\cX)), +\right) \to (A(\cX)_\CC^\times, \times).
	\end{equation*}
\end{lem}
\begin{proof}
	It suffices to consider an element $ x = P_m(\cV) $ in $ K^0(\cX) $ for a vector bundle $ \cV = (V, \phi) $ on $ \cX. $
	By the splitting principle, we can assume that each eigenbundle $ \cV_{\lambda} = (V_\lambda, \phi_\lambda) $ decomposes into line bundles 
	\begin{equation*}
	\cV_{\lambda} = \bigoplus_{j=1}^{\rk(\cV_{\lambda})} \cL_{\lambda,j}.
	\end{equation*}
	Let $ x_{\lambda,j} = [\cL_{\lambda,j}] $ in $ K^0(\cX) $. Then we have
	\begin{equation*}
	x = \sum_{1 \neq \lambda \in \eig(\phi(g))} \sum_{j=1}^{\rk(\cV_{\lambda})} x_{\lambda,j}
	\end{equation*}
	in $ K^0(\cX) $. Since $ e^K $ is multiplicative, we have that
	\begin{align*}
	e^{\rho_g}(x) & = \ch^{\rho_g} \left(\prod_{\lambda \neq 1,j} e^K(x_{\lambda,j})\right) = \left(\ch \circ \rho_g\right) \prod_{\lambda \neq 1,j} \left(1-x_{\lambda,j}^\vee\right) \\
	& = \ch \left(\prod_{\lambda \neq 1,j} \left( 1 - \lambda^{-1} x_{\lambda,j}^\vee \right)\right) = \prod_{\lambda \neq 1,j} \left( 1 - \lambda^{-1} e^{-h_{\lambda,j}} \right)
	\end{align*}
	in $ A(\cX)_\CC, $ where $ h_{\lambda,j} = c_1(x_{\lambda,j}). $ Since the degree zero component in $ e^{\rho_g}(x) $ is nonzero, it is invertible in $ A(\cX)_\CC. $
\end{proof}
\begin{defn}\label{defn_twisted_Todd}
	Let $ \cX = [X/G] $ be a connected separated smooth quotient DM stack. Let $ g $ be an element in the center of $ G $ such that it has a finite order and fixes $ X $. We define the $ \boldsymbol{g} $\tb{-twisted Todd class} map 
	$$ \td^{\rho_g}: K^0(\cX) \to A(\cX)_\CC^\times $$
	by setting
	\begin{equation*}
	\td^{\rho_g}(x) = \frac{\td(P_f(x))}{e^{\rho_g}(P_m(x))}
	\end{equation*}
	for any $ x $ in $ K^0(\cX). $
\end{defn}
\begin{rmk}
	Applying Definition \ref{defn_twisted_Todd} to each substack $ \cX_i = \{g_i\} \times [X_i/G_i] $ of $ \cX $, we then have a $ g_i $-twisted Todd class map
	\begin{equation*}
		\td^{\rho_{g_i}}: K^0(\cX_i) \to A(\cX_i)_\CC^\times.
	\end{equation*}
\end{rmk}
Now we can define the orbifold Todd class map.
\begin{defn}
	Let $ \cX = [X/G] $ be a connected separated smooth quotient DM stack in Notation \ref{notation_quotient_stack}. Define projections $ IP_f $ and $ IP_m $ on $ K^0(I\cX) $ by the direct sums
	\begin{equation*}
	IP_f = \bigoplus_{i \in I} P_f \quad \text{and} \quad IP_m = \bigoplus_{i \in I} P_m.
	\end{equation*}
	Define the \tb{inertia Euler class} map by the composition
	\begin{equation*}
	e^\rho: IP_m(K^0(I\cX)) \xrightarrow{e^K} K^0(I\cX) \xrightarrow{\ch^\rho} A(I\cX)_\CC^\times.
	\end{equation*}
	The \tb{inertia Todd class} map
	$$ \td^\rho: K^0(I\cX) \to A(I\cX)_{\CC}^\times $$
	is defined by
	\begin{equation*}
	\td^\rho(x) = \frac{\td(IP_f(x))}{e^\rho(IP_m(x))} = \left(\td(x_0),\  \bigoplus_{0 \neq i \in I} \frac{\td(P_f(x_i))}{e^{\rho_{g_i}}(P_m(x_i))}\right) \vspace{5pt}
	\end{equation*}
	for any $ x = \oplus_i x_i $ in $ K^0(I\cX) $. The \textbf{orbifold Todd class} map on $ K^0(\cX) $ is defined by the composition
	\begin{equation*}
	\orbtd : K(\cX) \xrightarrow{\beta} K^0(\cX) \xrightarrow{\pi^*} K^0(I\cX) \xrightarrow{\td^\rho} A(I\cX)_\CC^\times.
	\end{equation*}
\end{defn}
The following result is immediate.
\begin{prop}
	The orbifold Todd class map $ \orbtd : K(\cX) \to A(I\cX)_\CC^\times $ is multiplicative, i.e.,
	\begin{equation*}
	\orbtd(x+y) = \orbtd(x) \orbtd(y)
	\end{equation*}
	for all $ x $ and $ y $ in $ K(\cX). $
\end{prop}
\begin{rmk}[Explicit formulas for the orbifold Todd class]
	The orbifold Todd class of a vector bundle $ \cV = (V, \phi) $ on $ \cX $ is given by
	\begin{align*}
	\orbtd (\cV) = \left(\td(\cV),\ \bigoplus_{0 \neq i \in I} \frac{\td(\cV_i^\fix)}{e^{\rho_{g_i}}(\cV_i^\mov)}\right) = \left(\td(\cV), \ \bigoplus_{0 \neq i \in I}  \frac{\td(\cV_i^\fix)}{\prod_{\lambda \neq 1, j } (1 - \lambda^{-1} e^{-h_{i,\lambda,j}})} \right), \vspace{5pt}
	\end{align*}
	where $ \lambda \in \eig(\phi_i(g_i)) $, and $ h_{i,\lambda,j} $'s are the Chern roots of each eigenbundle $ \cV_{i,\lambda} $ on $ \cX_i $ for $ 0 \neq i \in I $ and $ 1 \leq j \leq \rk(\cV_{i,\lambda}) $. For a sheaf $ \cE $ on $ \cX, $ choose any finite locally free resolution $ \cE_{\boldsymbol{\cdot}} \to \cE \to 0 $, then we have
	\begin{equation*}
	\orbtd(\cE) = \prod_k \orbtd(\cE_k)^{(-1)^k}.
	\end{equation*}
\end{rmk}
We are interested in the orbifold Todd class of the $ K $-group class corresponding to the tangent complex of $ \cX. $
\begin{defn}[The tangent complex of a quotient stack]
	Let $ \cX = [X/G] $ be a smooth quotient stack. The tangent bundle $ TX $ has a canonical $ G $-equivariant structure $ \tau $ defined as follows. Take any $ g \in G. $ The map $ x \mapsto gx $ gives a smooth automorphism on $ X $
	\begin{equation*}
	g: X \to X
	\end{equation*}
	which induces a linear isomorphism on the tangent spaces
	\begin{equation*}
	dg_x: T_x X \to T_{gx} X
	\end{equation*}
	at all points $ x \in X. $ Therefore, we have a $ G $-action
	\begin{equation*}
	\tau: G \to \Aut(TX), \quad g \mapsto (\tau(g): (x,v) \mapsto (gx,dg_x(v)).		
	\end{equation*}
	The $ G $-equivariant vector bundle $ (TX, \tau) $ corresponds to a vector bundle on $ \cX, $ which we may simply denote by $ TX. $ Take any $ x \in X. $ The $ G $-action on $ X $ also gives a smooth morphism
	\begin{equation*}
	\sigma_x: G \to X, \quad g \mapsto gx
	\end{equation*}
	which induces a linear map 
	\begin{equation}\label{eq_lieG_Tx}
	L_x : \mathfrak{g} \to T_x X, \quad A \mapsto (d\sigma_x)_1(A)
	\end{equation}
	where $ \mathfrak{g} $ is the Lie algebra of $ G. $ There is a \tb{fundamental vector field} $ s_A $ associated to every $ A \in \mathfrak{g}, $ i.e., a section $ s_A: X \to TX $ with values 
	\begin{equation*}
	s_A(x) = L_x (A)
	\end{equation*}
	in $ T_xX $ for every $ x \in X. $ The \tb{tangent complex} of $ \cX $ is defined by a 2-term complex
	\begin{equation*}
	T\cX = [(\cO_X \otimes \mathfrak{g},\phi_\Ad) \xrightarrow{L} (TX, \tau)]
	\end{equation*}
	concentrated in degree $ -1 $ and $ 0, $ where $ \Ad $ is the adjoint representation 
	\begin{equation*}
	\Ad: G \to \GL(\fg), \quad g \mapsto (\Ad(g): A \mapsto gAg^{-1})
	\end{equation*}
	of $ G $ and the morphism $ L $ is induced by the linear maps $ L_x $ in (\ref{eq_lieG_Tx}), i.e., $ L $ corresponds to the morphism
	\begin{equation*}
	X \times \fg \to \Tot(TX), \quad (x,A) \mapsto (x,s_A(x))
	\end{equation*}
	between total spaces of $ \cO_X \otimes \mathfrak{g} $ and $ TX $ which is $ G $-equivariant, as can be easily checked. To simplify notations, we may denote $ T\cX $ by
	\begin{equation*}
	T\cX = [\cO_X \otimes \mathfrak{g} \xrightarrow{L} TX].
	\end{equation*}
	The \tb{tangent space} of $ \cX $ at a point $ x \in \cX $ is given by
	\begin{equation*}
	T_x\cX = T_xX / L_x(\fg).
	\end{equation*}
	Note that the image $ L_x(\fg) \subset T_xX $ can be identified with the tangent space of the orbit $ Gx $ at $ x. $ When $ G $ is a finite group, then $ T\cX = TX $ is simply the $ G $-equivariant tangent bundle on $ X, $ and the tangent space $ T_x \cX = T_x X $ at all $ x \in X. $ 
\end{defn}
\begin{ex}
	Let $ X = \CC^{n+1} $ and $ U = X - \{0\}. $ Let's re-visit the example of a weighted projective stack 
	$$ \PP = \PP(a_0, \dots, a_n) = [U/\CC^*]. $$ 
	Let's compute the tangent complex of $ \PP. $ Recall that pulling back the identity character
	\begin{equation*}
	\rho: \CC^* \to \CC^*, \quad t \to t
	\end{equation*}
	of $ \CC^* $ from $ B\CC^* $ to $ \cX = [X/\CC^*] $ gives a vector bundle $ \cO_\cX \otimes \rho $ on $ \cX, $ which restricts to the twisting sheaf $ \cO_\PP(1) = j^* (\cO_\cX \otimes \rho) $ where $ j: \PP \to \cX $ is the inclusion map. The tangent bundle $ TU $ can be identified as
	\begin{equation*}
	TU \cong \bigoplus_{i=0}^n \cO_\PP(a_i).
	\end{equation*}
	The Lie algebra of $ \CC^* $ is $ \CC $ with the trivial adjoint representation. At a point $ x = (x_0, \dots, x_n) \in X, $ the linear map in (\ref{eq_lieG_Tx}) is an injection
	\begin{equation*}
	L_x: \CC \to T_x U, \quad 1 \mapsto (a_0 x_0, \dots, a_n x_n) \neq 0.
	\end{equation*}
	All these linear maps glue to the tangent complex
	\begin{equation*}
	T\PP = [\cO_\PP \xrightarrow{L} TU ] = [\cO_\PP \xrightarrow{L} \bigoplus_{i=0}^n \cO_\PP(a_i) ].
	\end{equation*}
	Since the morphism $ L $ is injective, the tangent complex $ T\PP $ is indeed a tangent bundle, which is the cokernel of $ L $ in the Euler sequence
	\begin{equation*}
	0 \to \cO_\PP \xrightarrow{L} \bigoplus_{i=0}^n \cO_\PP(a_i) \to T\PP \to 0.
	\end{equation*}
	The tangent space of $ \PP $ at $ x = (x_0, \dots, x_n) \in U $ is
	\begin{equation*}
	T_x \PP = \CC^{n+1}/\CC \cdot (a_0 x_0, \dots, a_n x_n) \cong \CC^n.
	\end{equation*}
\end{ex}

\vspace{10pt}
Now let's compute the orbifold Todd class of the $ K $-group class 
\begin{equation*}
[T\cX] = [TX] - [\cO_X \otimes \mathfrak{g}] = [TX, \tau] - [\cO_X \otimes \mathfrak{g}, \phi_\Ad]
\end{equation*}
of the tangent complex $ T\cX $ in $ K(\cX). $ 
By the definition of $ \orbtd $, we have
\begin{equation*}
\orbtd\left(T\cX\right) = \frac{\td(IP_f[\pi^* T\cX])}{e^\rho(IP_m[\pi^* T\cX])},
\end{equation*}
where $ IP_f[\pi^* T\cX] $ and $ IP_m[\pi^* T\cX] $ are computed in the following lemma.
\begin{lem}
	For a connected separated smooth quotient DM stack $ \cX, $ there is a splitting short exact sequence
	\begin{equation*}
	0 \to TI\cX \to \pi^*T\cX \to NI\cX \to 0
	\end{equation*}
	of 2-term complexes on the inertia stack $ I\cX $ such that
	\begin{equation*}
	[TI\cX] = IP_f[\pi^* T\cX] \quad \text{and} \quad [NI\cX] = IP_m[\pi^* T\cX]
	\end{equation*}
	in $ K^0(I\cX). $
\end{lem}
\begin{proof}
	Let $ I $ be the inertia index set of $ \cX $ in Notation \ref{notation_quotient_stack}. Take $ i \in I. $ Pulling back $ T\cX $ along $ q_i: \cX_i \to \cX $ gives a complex
	\begin{equation*}
	(T\cX)_i = q_i^*T\cX = [(\cO_{X_i} \otimes \mathfrak{g},(\phi_\Ad)_i) \to ((TX)_i, \tau_i)]
	\end{equation*}
	where $ (\phi_\Ad)_i $ and $ \tau_i $ are $ G_i $-equivariant structures on $  \cO_{X_i} \otimes \mathfrak{g} $ and $ (TX)_i = p_i^*TX $ induced from $ G$-equivariant structures $ \phi_\Ad $ and $ \tau $ respectively.
	We want to find the fixed and moved sub-complexes 
	\begin{equation*}
	(T\cX)_i^\fix \quad \text{and} \quad (T\cX)_i^\mov
	\end{equation*}
	of the tangent complex $ T\cX $ under the action of $ g_i $ (via linear operators $ \phi_\Ad(g_i) $ and $ \tau(g_i) $ on the two terms of $ T\cX $). To simplify notations, we write
	\begin{equation*}
	(T\cX)_i = [\cO_{X_i} \otimes \mathfrak{g} \to (TX)_i]
	\end{equation*}
	where $ G_i $-equivariant structures are understood. 
	The inclusion $ G_i \into G $ induces a splitting short exact sequence
	\begin{equation*}
	0 \to \fg_i \to \fg \to \fg/\fg_i \to 0
	\end{equation*}
	of linear representations of $ G_i $, where $ \fg_i $ carries the adjoint representation $$ \Ad_i: G_i \to \GL(\fg_i) $$
	of $ G_i $, and $ \fg $ carries the restriction $ \Ad|_{G_i} $ of the adjoint representation of $ G. $
	Since $ G_i $ is a subgroup of the centralizer $ Z_{g_i}, $ we have
	\begin{equation*}
	\Ad_i(g_i) = \id_{\fg_i}.
	\end{equation*}
	Therefore, we have a splitting short exact sequence
	\begin{equation*}
	0 \to \cO_{X_i} \otimes \fg_i \to \cO_{X_i} \otimes \fg \to \cO_{X_i} \otimes \fg/\fg_i \to 0
	\end{equation*}
	of $ G_i $-equivariant vector bundles on $ X_i, $ where
	\begin{equation*}
	\cO_{X_i} \otimes \fg_i = \left(\cO_{X_i} \otimes \fg\right)^\fix \quad \text{and} \quad \cO_{X_i} \otimes \fg/\fg_i = \left(\cO_{X_i} \otimes \fg\right)^\mov.
	\end{equation*}
	The inclusion $ p_i: X_i \into X $ gives a splitting short exact sequence
	\begin{equation*}
	0 \to TX_i \to p_i^* TX \to NX_i \to 0
	\end{equation*}
	of $ G_i $-equivariant vector bundles on $ X_i, $ where $ TX_i $ and $ NX_i $ are the tangent bundle and normal bundle of $ X_i $ with their canonical $ G_i $-equivariant structures. Since $ g_i $ acts trivially on $ X_i, $ we can identify
	\begin{equation*}
	TX_i = (TX)_i^\fix \quad \text{and} \quad NX_i = (TX)_i^\mov.
	\end{equation*}
	Therefore, we obtain a splitting short exact sequence 
	\begin{equation*}
	0 \to T\cX_i \to (T\cX)_i \to N\cX_i \to 0
	\end{equation*}
	of 2-term complexes on $ \cX_i $, where 
	\begin{equation*}
	T\cX_i = [\cO_{X_i} \otimes \fg_i \to TX_i] \quad \text{and} \quad N\cX_i = [\cO_{X_i} \otimes \fg/\fg_i \to NX_i]
	\end{equation*}
	are identified with $ (T\cX)_i^\fix $ and $ (T\cX)_i^\mov $ respectively.
	Define two complexes
	\begin{equation*}
	TI\cX = \coprod_{i \in I} T\cX_i \quad \text{and} \quad NI\cX = \coprod_{i \in I} N\cX_i
	\end{equation*}
	on $ I\cX. $ Since $ \pi^*T\cX = \coprod_{i \in I} \left(T\cX\right)_i, $ we have
	\begin{equation*}
	[TI\cX] = IP_f[\pi^* T\cX] \quad \text{and} \quad [NI\cX] = IP_m[\pi^* T\cX]
	\end{equation*}
	in $ K^0(I\cX) $ as claimed in the lemma.
\end{proof}
\begin{notation}[The orbifold Todd class of $ T\cX $]\label{notation_todd}
	Let $ \cX = [X/G] $ be a connected separated smooth quotient DM stack in Notation \ref{notation_quotient_stack}.
	We have
	\begin{equation*}
	\orbtd(T\cX) = \frac{\td(TI\cX)}{e^\rho(NI\cX)} = \left(\td(T\cX),\ \bigoplus_{0 \neq i \in I} \frac{\td(T\cX_i)}{e^{\rho_{g_i}}(N\cX_i)}\right)
	\end{equation*}
	in $ A(I\cX)_\CC^\times $, where
	\begin{equation*}
	\td(T\cX) = \frac{\td(TX)}{\td(\cO_X \otimes \fg)}
	\end{equation*}
	is the Todd class of the tangent complex of $ \cX $ in $ A(\cX)_\QQ^\times $. To simplify notations, we write
	\begin{equation*}
	\orbtd_\cX = \frac{\td_{I\cX}}{e^\rho_{I\cX}} = \left( \td_\cX, \bigoplus_{0 \neq i \in I} \frac{\td_{\cX_i}}{e^\rho_{\cX_i}} \right), \vspace{5pt}
	\end{equation*}
	where $ \td_{I\cX} = \td(TI\cX), \ e^\rho_{I\cX} = e^\rho(NI\cX), \ \td_{\cX_i} = \td(T\cX_i) $, and $ e^\rho_{\cX_i} = e^{\rho_{g_i}}(N\cX_i) $ for each $ i \neq 0. $
	If $ G $ is a finite group, then $ \fg $ vanishes and hence we have
	\begin{equation}\label{eq_orbtd_finite_G}
	\orbtd_\cX = \left(\td_X, \bigoplus_{0 \neq i \in I} \frac{\td_{X_i}}{e^{\rho_{g_i}}(NX_i)}\right) = \left(\td_X, \bigoplus_{0 \neq i \in I} \frac{\td_{X_i}}{\prod_{\lambda, j} \left(1 - \overline{\lambda} e^{-h_{i,\lambda,j}}\right)}\right),
	\end{equation}
	where the second identity is from the $ G_i $-equivariant decomposition
	\begin{equation*}
	NX_i = \bigoplus_{\lambda \in \eig(\tau_i(g_i))} N_{i,\lambda}
	\end{equation*}
	of the normal bundle $ NX_i $ for each $ 0 \neq i \in I $, where $ \tau_i $ is the $ G_i $-equivariant structure on $ NX_i $ induced from the $ G $-equivariant tangent bundle $ (TX, \tau) $ on $ X, $ $ \eig(\tau_i(g_i)) $ contains nontrivial $ \ord(g_i) $-th roots of unity, and $ h_{i,\lambda,1}, \dots, h_{i,\lambda,\rk(N_{i,\lambda})} $ are the Chern roots of each eigenbundle $ N_{i,\lambda} $ for each $ 0 \neq i \in I $.
\end{notation}

\subsection{Orbifold Mukai pairing and Orbifold HRR formula}\label{sec_orbifold_HRR}
In this section we will define an orbifold Mukai pairing and derive a new HRR formula for a pair of coherent sheaves on a connected proper smooth quotient DM stack.

\vspace{4pt}

We begin with the definition of the orbifold Euler characteristic.
\begin{defn}\label{defn_Euler_char}
	Let $ \cX = [X/G] $ be a proper quotient stack. The \tb{orbifold Euler characteristic} $ \chi(\cX, \ {\cdot}\ ): K(\cX) \to \ZZ $ is defined by
	\begin{equation*}
	\chi(\cX, \cE) = \sum_i (-1)^i \dim H^i(\cX,\cE)
	\end{equation*}
	for a coherent sheaf $ \cE $ on $ \cX $ and extended linearly. 
\end{defn}
\begin{notation}
	Let $ \cX = [X/G] $ be a proper quotient stack of dimension $ n $ with its structure morphism $ f: \cX \to \pt $. We denote the pushforward 
	$$ f_*: A(\cX)_\CC \to A(\pt)_\CC = \CC $$ 
	by the integral symbol $\int_\cX $, i.e., for an element $ cv \in A(\cX)_\CC $ with $ c \in \CC $ and $ v \in A(\cX) $, we have
	\begin{equation*}
	\int_\cX cv = f_*(cv) =  c\deg (v_n),
	\end{equation*}
	where $ v_n $ is the degree $ n $ component of $ v. $
	We also have a linear map
	\begin{equation*}
	\int_{I\cX} : A(I\cX)_\CC \to \CC	
	\end{equation*}
	defined by $ \int_{I\cX} = \sum_{i \in I} \int_{\cX_i} $. 
\end{notation}
\begin{rmk}
	If $ G $ is a finite group, then we can identify $ A(\cX) $ with $ A(X)^G, $ and for an element $ v \in A(X)^G, $ we have 
	\begin{equation*}
	\int_{\cX} v = \frac{1}{|G|}\int_X v
	\end{equation*}
	in $ \QQ. $
\end{rmk}
Now we can state the orbifold HRR theorem.
\begin{thm}[Orbifold HRR]\label{thm_HRR}
	Let $ \cX $ be a connected proper smooth quotient DM stack. Then for all $ x \in K(\cX), $ we have
	\begin{equation}\label{eq_HRR}
	\chi(\cX, x) = \int_{I\cX} \orbch(x) \orbtd_\cX = \int_{I\cX} \orbch(x) \frac{\td_{I\cX}}{e^\rho_{I\cX}}.
	\end{equation}
\end{thm}
\begin{proof}
	Take an element $ x $ in $ K(\cX). $ By the natural isomorphism  $ K(\cX) \cong K^0(\cX) $, we can write $ x = [\cV] - [\cW] $ for two vector bundles $ \cV $ and $ \cW $ on $ \cX $. Formula (\ref{eq_HRR}) follows from Theorem 4.19 in \cite{edidin2013riemann} for vector bundles on $ \cX $ and the linearities of the orbifold Euler characteristic and the orbifold Chern character.
\end{proof}
\begin{rmk}
	\begin{enumerate}[font=\normalfont,leftmargin=*]
		\item If $ \cX $ is a quasi-projective stack, i.e., the coarse moduli space of $ \cX $ is a quasi-projective scheme, then Theorem \ref{thm_HRR} is the content of Part 1 of Theorem 4.10 in \cite{toen1999theoremes} for the morphism $ \cX \to \Spec \CC. $
		\item An analytic version of Theorem \ref{thm_HRR} was first proved in \cite{kawasaki1979riemann} for holomorphic vector bundles on a compact complex orbifold. This explains the word ``orbifold" in the orbifold Chern character, the orbifold Todd class, and the orbifold HRR theorem.
	\end{enumerate}
\end{rmk}
\begin{rmk}[The explicit orbifold HRR theorem]
	Let $ I $ be the finite inertia index set of $ \cX $ in Notations \ref{notation_quotient_stack} and \ref{notation_todd}. For a vector bundle $ V $ on $ \cX $, the orbifold HRR theorem becomes
	\begin{equation*}
	\chi(\cX, \cV) = \int_{\cX} \ch(V) \td_X + \sum_{0 \neq i \in I} \int_{\cX_i} \ch(\cV_i) \frac{\td_{\cX_i}}{e^{\rho}_{\cX_i}}
	\end{equation*}
	where $ \cV_i $ is the restriction of the vector bundle $ \cV $ on the substack $ \cX_i $ of $ \cX $ for each index $ i \neq 0. $ For a sheaf $ \cE $ on $ \cX, $ we can compute 
	$$ \chi(\cX, \cE) = \sum_k (-1)^k \chi(\cX, \cE_k)$$ 
	from a finite locally free resolution $ \cE_{\boldsymbol{\cdot}} \to \cE \to 0. $
\end{rmk}
\begin{defn}\label{defn_Euler_pairing}
	Let $ \cX = [X/G] $ be a proper quotient stack. The \tb{orbifold Euler pairing} 
	$$ \chi: K(\cX) \times K(\cX) \to \ZZ $$ 
	is defined by
	\begin{equation*}
	\chi(\cE,\cF) = \sum_i (-1)^i \dim \Ext^i(\cE,\cF)
	\end{equation*}
	for two coherent sheaves $ \cE = (E, \phi) $ and $ \cF = (F, \psi) $ on $ \cX $ and extended bilinearly.
\end{defn}
\begin{lem}\label{lem_Euler_pairing}
	Let $ \cX = [X/G] $ be a proper smooth quotient DM stack. Then for all $ x $ and $ y $ in $ K(\cX), $ we have
	\begin{equation*}
	\chi(x, y) = \chi(\cX, x^\vee y).
	\end{equation*}
	In particular, $ \chi(1, x) = \chi(\cX, x) $ for all $ x $ in $ K(\cX) $ where $ 1 = [\cO_\cX]. $
\end{lem}
\begin{proof}
	By the bilinearity of $ \chi, $ it suffices to consider $ x = [\cV] $ and $ y =[\cW] $ for vector bundles $ \cV = (V, \phi) $ and $ \cW = (W, \psi) $ on $ \cX. $ Then we have
	\begin{align*}
	\chi(x, y) & = \chi(\cV, \cW) = \sum_i (-1)^i \dim \Ext^i(V,W)^G \\
	& = \sum_i (-1)^i \dim H^i(X,V^\vee \otimes W)^G & \text{because $ V $ is locally free} \\
	& = \chi(\cX, \cV^\vee \otimes \cW) \\
	& = \chi(\cX, x^\vee y).
	\end{align*}
\end{proof}
For a separated smooth quotient DM stack $ \cX = [X/G] $, we can define its canonical line bundle via the cotangent complex $ \Omega_X \to \cO_X \otimes \fg^\vee $.
\begin{defn}
	Let $ \cX $ be a separated smooth quotient DM stack. The canonical line bundle of $ \cX $ is defined as the $ G $-equivariant line bundle
	\begin{equation*}
	\omega_\cX = \det(\Omega_X) \otimes \det(\cO_X \otimes \fg^\vee)^\vee
	\end{equation*}
	on $ X $ where the vector bundles $ \Omega_X $ and $ \cO_X \otimes \fg^\vee $ are equipped with natural $ G $-equivariant structures.
\end{defn}
Serre duality for projective schemes can be generalized to projective stacks. The following was proved in \cite[Appendix B]{bruzzo2015framed}.
\begin{thm}[Serre Duality]
	Let $ \cX $ be a smooth projective stack of dimension $ n $ over $ \CC $. Then for all coherent sheaves $ \cE $ and $ \cF $ on $ \cX, $ there is a natural isomorphism
	\begin{equation*}
	\Ext^i(\cE, \cF) \cong \Ext^{n-i}(\cF, \cE \otimes \omega_\cX)
	\end{equation*}
	of vector spaces for all $ i \in \ZZ $.
\end{thm}
\begin{rmk}
	Let $ \cX $ be a smooth projective stack. Serre duality implies that
	\begin{equation*}
	\chi(\cE,\cF) = (-1)^{\dim \cX} \chi(\cF,\cE \otimes \omega_\cX)
	\end{equation*}
	for a pair of coherent sheaves $ \cE $ and $ \cF $ on $ \cX, $ and hence
	\begin{equation*}
	\chi(x,y) = (-1)^{\dim \cX} \chi(y,x \cdot [\omega_\cX])
	\end{equation*}
	for all $ x $ and $ y $ in $ K(\cX). $
	Therefore, if $ \cX $ has a trivial canonical line bundle, the orbifold Euler pairing for $ \cX $ is symmetric when $ \dim \cX $ is even, and is anti-symmetric when $ \dim \cX $ is odd. 
\end{rmk}

Next we will define the orbifold Mukai vector and the orbifold Mukai pairing for quotient stacks. The Mukai vector and the Mukai pairing for schemes are reviewed in Appendix \ref{review_mukai}.

Let $ \cX $ be a connected separated smooth quotient DM stack. We first define the square root of the Todd class of the inertia stack $ I\cX. $
\begin{defn}
	Let $ \cX $ be a connected separated smooth quotient DM stack with an inertia index set $ I $ in Notations \ref{notation_quotient_stack} and \ref{notation_todd}.
	Take an index $ i \in I. $ We have $ \td_{\cX_i} = 1+v_i \in A(\cX_i)_\CC^\times $ for some $ v_i \in A^{\geq 1}(\cX_i). $ If $ d_i = \dim \cX_i $, then the square root of $ \td_{\cX_i} $
	$$ \sqrt{\td_{\cX_i}} = \sum_{k=0}^{d_i} \binom{1/2}{k} v_i^k = 1 + \frac{v_i}{2} - \frac{v_i^2}{8} + \cdots + \binom{1/2}{d_i} v_i^{d_i} $$
	is well defined in $ A(\cX_i)_\CC^\times. $ We define the square root of $ \td_{I\cX} $ by
	\begin{equation*}
		\sqrt{\td_{I\cX}} = \bigoplus_{i \in I} \sqrt{\td_{\cX_i}}.
	\end{equation*}
\end{defn}
\begin{defn}\label{defn_orbv}
	The \tb{orbifold Mukai vector} map 
	$$ \orbv: K(\cX) \to A(I\cX)_\CC $$
	is defined by
	\begin{equation*}
	\orbv(x) = \orbch(x) \sqrt{\td_{I\cX}} 
	\end{equation*}
	for all $ x $ in $ K(\cX). $ 
\end{defn}
\begin{rmk}
	For an element $ x $ in $ K(\cX), $ if $ \orbv(x) = (v,(v_i)) \in A(\cX)_\QQ \oplus A(I_t\cX)_\CC $, then $ v $ is called the \tb{Mukai vector} of $ x $ and each $ v_i $ for $ 0 \neq i \in I $ is called the $ \boldsymbol{i} $\tb{-th twisted Mukai vector} of $ x $.
\end{rmk}
\begin{rmk}
	Our definition of the orbifold Mukai vector map is different from the one in \cite{popa2017derived} where
	\begin{equation*}
	\orbv(x) = \orbch(x) \sqrt{\orbtd_{\cX}} = \orbch(x) \sqrt{\td_{I\cX}/e^\rho_{I\cX}}
	\end{equation*}
	in the cohomology $ H^*(I\cX)_\CC. $ We made a different (and arguably better) choice because taking the square root of $ e^\rho_{I\cX} $ is an unnecessary computation and it would also make the orbifold Mukai pairing more complicated. We also used the Chow ring $ A(I\cX)_\CC $ rather than $ H^*(I\cX)_\CC $ for the target of the orbifold Mukai vector map.
\end{rmk}
The next lemma follows immediately from Definition \ref{defn_orbv}.
\begin{lem}
	The orbifold Mukai vector map $ \orbv: K(\cX) \to A(I\cX)_\CC $ satisfies
	\begin{equation*}
	\orbv(x+y) = \orbv(x) + \orbv(y) \quad \text{and} \quad \orbv(xy) = \orbv(x) \orbch(y)
	\end{equation*}
	for all $ x $ and $ y $ in $ K(\cX). $ 
\end{lem}
There is an involution on the complex Chow ring $ A(I\cX)_\CC $.
\begin{defn}\label{defn_involution}
	Let $ \cX $ be a connected separated smooth quotient DM stack. We define an involution 
	$$ (\ \cdot\ )^\vee: A(I\cX)_\CC \to A(I\cX)_\CC $$ 
	by defining one on each component of $ A(I\cX)_\CC $. Take $ i \in I. $ An element $ u \in A(\cX_i)_\CC $ is a finite sum:
	\begin{equation*}
		u = \sum_{j,k} a_{jk} u_{jk},
	\end{equation*}
	where $ a_{jk} \in \CC $, and $ u_{jk} \in A^j(\cX_i) $ is the class of an irreducible closed substack of $ \cX_i $ of codimension $ j $. The involution on $ A(\cX_i)_\CC $ is defined by
	\begin{equation*}
		u^\vee = \sum_{j, k} (-1)^j \overline{a}_{jk} u_{jk},	
	\end{equation*}
	where $ \overline{a}_{jk} $ denotes the complex conjugate of $ a_{jk} $. For every element $ \orbv = \bigoplus_{i \in I} v_i $ in $ A(I\cX)_\CC, $ define
	\begin{align*}
		\orbv^\vee = \bigoplus_{i \in I} v_i^\vee.
	\end{align*}
\end{defn}
\begin{rmk}
	It is straightforward to check that $ (\ \cdot \ )^\vee: A(I\cX)_\CC \to A(I\cX)_\CC $ is a ring automorphism which also commutes with $ \sqrt{(\ \cdot \ )} $ when this is defined.
\end{rmk}
\begin{lem}\label{lem_orbch_orbv}
	Let $ \cX $ be a connected separated smooth quotient DM stack. Then for all $ x $ and $ y$ in $ K(\cX), $ we have
	\begin{equation*}
	\orbch(x^\vee) = \orbch(x)^\vee \quad \text{and} \quad \orbv(x^\vee) = \orbv(x)^\vee \sqrt{\td_{I\cX}/\td_{I\cX}^\vee}
	\end{equation*}
	in $ A(I\cX)_\CC $.
\end{lem}
\begin{proof}
	We first prove the first identity. By the linearity of $ \orbch $ and the involutions on $ K(\cX) $ and $ A(I\cX)_\CC, $ it suffices to check the case $ x = [\cV] $ is the class of a vector bundle $ \cV = (V, \phi) $ on $ \cX. $ We can write
	$$ \orbch\left(\cV^\vee\right) = \left( \ch(\cV^\vee),\ \bigoplus_{i \neq 0} \ch^{\rho_{g_i}}\left(\cV_i^\vee\right) \right). $$ 
	We have $ \ch(\cV^\vee) = \ch(\cV)^\vee $ by the property of the ordinary Chern character. For each $ 0 \neq i \in I, $ we have
	\begin{align*}
	\ch^{\rho_{g_i}} \left(\cV_i^\vee\right) & = \ch \left(\sum_{\lambda \in \eig(\phi_i(g_i))} \lambda^{-1} \, [\cV_{i,\lambda}]^\vee\right) \\
	& = \sum_{\lambda \in \eig(\phi_i(g_i))} \sum_{j=0}^{\dim \cX_i} (-1)^j \lambda^{-1} \ch_j \left(\cV_{i,\lambda}\right) = \left(\ch^{\rho_{g_i}} (\cV_i)\right)^\vee.
	\end{align*}
	Now the second identity follows from two identities:
	\begin{equation*}
	\orbv (x^\vee) = \orbch(x^\vee) \sqrt{\td_{I\cX}} = \orbch(x)^\vee \sqrt{\td_{I\cX}} \quad \text{and}\quad \orbv (x)^\vee = \orbch(x)^\vee \sqrt{\td_{I\cX}}^\vee
	\end{equation*}
\end{proof}

Assume $ \cX $ is proper. Now we define an orbifold Mukai pairing on $ A(I\cX)_\CC $.
\begin{defn}\label{defn_orbv_pairing}
	Let $ \cX $ be a connected proper smooth quotient DM stack. The \tb{orbifold Mukai pairing} 
	$$ \inprod{\cdot \ {,} \ \cdot }_{I\cX}: A(I\cX)_\CC \times A(I\cX)_\CC \to \CC $$ 
	is defined by
	\begin{equation}\label{eq_orbv_pairing}
	\inprod{\orbv,\orbw}_{I\cX} = \int_{I\cX} \frac{\orbv^\vee \orbw}{e^\rho_{I\cX}} \sqrt{\frac{\td_{I\cX}}{\td_{I\cX}^\vee}}
	\end{equation}
	for all $ \orbv $ and $ \orbw $ in $ A(I\cX)_\CC. $
\end{defn}
\begin{rmk}
	For a pair of vectors $ \orbv = (v,(v_i)) $ and $ \orbw = (w, (w_i)) $ in $ A(I\cX)_\CC, $ the orbifold Mukai pairing can be written as
	\begin{equation*}
	\inprod{\orbv,\orbw}_{I\cX} = \inprod{v,w}_\cX + \sum_{0 \neq i \in I} \inprod{v_i,w_i}_{\cX_i},
	\end{equation*}
	where the \tb{Mukai pairing}
	\begin{equation*}
	\inprod{v,w}_\cX = \int_\cX v^\vee w \sqrt{\td_\cX/\td_\cX^\vee}
	\end{equation*}
	and the $ \boldsymbol{i} $\tb{-th twisted Mukai pairing}
	\begin{equation*}
	\inprod{v_i,w_i}_{\cX_i} = \int_{\cX_i} \frac{v_i^\vee w_i}{e^\rho_{\cX_i}} \sqrt{\td_{\cX_i}/\td_{\cX_i}^\vee}
	\end{equation*}
	for each index $ 0 \neq i \in I. $
\end{rmk}
The following is immediate from the definition of the orbifold Mukai pairing.
\begin{lem}
	The orbifold Mukai pairing $ \inprod{\cdot \ {,} \ \cdot }_{I\cX}: A(I\cX)_\CC \times A(I\cX)_\CC \to \CC $ is sesquilinear, i.e.,
	\begin{equation*}
	\inprod{\orbv + \orbs, \orbw + \orbt}_{I\cX} = \inprod{\orbv, \orbw}_{I\cX} + \inprod{\orbv, \orbt}_{I\cX} + \inprod{\orbs, \orbw}_{I\cX} + \inprod{\orbs, \orbt}_{I\cX}
	\end{equation*}
	and
	\begin{equation*}
	\inprod{a\orbv, b\orbw}_{I\cX} = \cj{a}b \inprod{\orbv, \orbw}_{I\cX}
	\end{equation*}
	for all $ \orbv, \orbs, \orbw, \orbt $ in $ A(I\cX)_\CC $ and all $ a,b $ in $ \CC. $
\end{lem}
Now we can state the orbifold HRR formula for a pair of coherent sheaves in terms of the orbifold Euler pairing and the orbifold Mukai pairing.
\begin{thm}[Orbifold HRR Formula]
	Let $ \cX $ be a connected proper smooth quotient DM stack. Then for all $ x $ and $ y $ in $ K(\cX), $ we have
	\begin{equation}\label{eq_HRR2}
	\chi(x, y) = \inprod{\orbv(x),\orbv(y)}_{I\cX}.
	\end{equation}
\end{thm}
\begin{proof}
	Take any $ x $ and $ y $ in $ K(\cX). $ Formula (\ref{eq_HRR2}) is a result of the orbifold HRR theorem and the definition of the orbifold Mukai pairing:
	\begin{align*}
	\chi(x,y) & = \chi(\cX, x^\vee y) & \text{by Lemma \ref{lem_Euler_pairing}} \\
	& = \int_{I\cX} \orbch(x^\vee y) \orbtd_\cX & \text{by the orbifold HRR theorem \ref{thm_HRR}} \\
	& = \int_{I\cX} \orbch(x^\vee) \orbch(y) \orbtd_\cX & \text{since $ \orbch $ is a ring map by Proposition \ref{prop_orbch_ring_map}} \\
	& = \int_{I\cX} \frac{\orbv(x^\vee) \orbv(y)}{e^\rho_{I\cX}} & \text{by Definition \ref{defn_orbv}} \\
	& = \int_{I\cX} \frac{\orbv(x)^\vee \orbv(y)}{e^\rho_{I\cX}} \sqrt{\td_{I\cX}/\td_{I\cX}^\vee} & \text{by Lemma \ref{lem_orbch_orbv}} \\
	& = \inprod{\orbv(x),\orbv(y)}_{I\cX} & \text{by Definition \ref{defn_orbv_pairing}}
	\end{align*}
\end{proof}
The orbifold Euler pairing on $ K(\cX) $ can be extended to a sesquilinear form on $ K(\cX)_\CC $. Since $ \sqrt{\td_{I\cX}} $ is a unit in $ A(I\cX)_\CC, $ Proposition \ref{prop_orbch_isom} implies the following:
\begin{prop}
	Let $ \cX $ be a connected proper smooth quotient DM stack. The orbifold Mukai vector map $ \orbv: K(\cX)_\CC \to A(I\cX)_\CC $ is a linear isometry
	\begin{equation*}
	\orbv: \left(K(\cX)_\CC, \chi\right) \xrightarrow{\simeq} \left(A(I\cX)_\CC, \inprod{\cdot \ {,} \ \cdot }_{I\cX}\right).
	\end{equation*}
\end{prop}
Note that formula (\ref{eq_HRR2}) reduces to (\ref{eq_HRR}) when $ x = [\cO_\cX] = 1.$ Let's apply the orbifold HRR formula for $ \cX = [X/G] $ to two extreme cases: 1) $ G = \{1\} $, and 2) $ X = \pt $.
\begin{ex}
	Suppose $ G = \{1\} $ is the trivial group and $ \cX = X $ is a proper smooth scheme of dimension $ n $. Then $ \orbch, \orbtd, $ and $ \orbv $ become the usual Chern character $ \ch $, Todd class $ \td $, and Mukai vector $ v $ on schemes. It suffices to tensor $ A(X) $ with $ \QQ $ for the target of these maps.
	The orbifold Mukai pairing reduces to the Mukai pairing for schemes. 
	Formula (\ref{eq_HRR}) reduces to the HRR theorem for schemes:
	\begin{equation*}
	\chi(X, x) = \int_X \ch(x) \td_X
	\end{equation*}
	for all $ x $ in $ K(X). $
	Formula (\ref{eq_HRR2}) then becomes the HRR formula
	\begin{equation*}
	\chi(x, y) = \inprod{v(x),v(y)}
	\end{equation*}
	for all $ x $ and $ y $ in $ K(X) $. If $ c_1(X) $ vanishes (e.g., when $ X $ is a Calabi-Yau manifold), then the Mukai pairing can be defined integrally on $ A(X), $ i.e., we have a bilinear map
	\begin{equation*}
	\inprod{\cdot\ {,}\ \cdot}: A(X) \times A(X) \to \ZZ
	\end{equation*}	
	given by
	\begin{equation*}
	\inprod{v,w} = \int_X \sum_{i=0}^{n} (-1)^i v_i \sum_{j=0}^{n} w_j = \deg \left(\sum_{i=0}^{n} (-1)^i v_i w_{n-i}\right) \vspace{8pt}
	\end{equation*}
	in $ \ZZ $ for all $ v = \sum v_i $ and $ w = \sum w_i $ in $ A(X) = \oplus_{i=0}^n A^i(X) $.
\end{ex}
\begin{ex}
	Consider the other extreme case: $ \cX = BG $ for a finite group $ G $ with conjugacy classes $ [g_0], [g_1], \dots, [g_{n-1}] $. Hence 
	$$ K(BG) = K^0(BG) = R(G). $$ 
	The orbifold Euler characteristic of a sheaf $ \cV = (V, \phi) $ on $ BG $ is
	\begin{equation*}
	\chi(BG, \phi) = \dim V^G,
	\end{equation*}
	and the orbifold Euler pairing for sheaves $ \cV = (V, \phi) $ and $ \cW = (W, \psi) $ on $ BG $ is
	\begin{equation*}
	\chi(\phi, \psi) = \dim \Hom(V,W)^G.
	\end{equation*}
	The inertia stack of $ BG $ is
	\begin{equation*}
	IBG \simeq \coprod_{i = 0}^{n-1} \{g_i\} \times Z_{g_i}, \quad \text{and} \quad A(IBG)_\CC \cong \CC^n.
	\end{equation*}
	The orbifold Mukai vector and orbifold Chern character coincide to be
	\begin{equation*}
	\orbv = \orbch: R(G) \to \CC^n, \quad [V, \phi] \mapsto \chi_\phi,
	\end{equation*}
	i.e., the character map for representations of $ G. $ The involution on $ R(G) $ is taking dual representations of $ G, $ and the involution on $ A(IBG)_\CC = \CC^n $ is complex conjugation. The orbifold Mukai pairing is an inner product $ \inprod{\ {,}\ }_W $ on $ \CC^n $ weighted by the centralizers $ Z_{g_0}, \dots, Z_{g_{n-1}}, $ since formula (\ref{eq_orbv_pairing}) reduces to
	\begin{equation*}
	\inprod{a,b}_{IBG} = \int_{IBG} \overline{a} b = \sum_{i=0}^{n-1} \int_{BZ_{g_i}} \overline{a}_i b_i = \sum_{i=0}^{n-1} \frac{1}{|Z_{g_i}|} \overline{a}_i b_i = \inprod{a,b}_W,
	\end{equation*}
	for all $ a = (a_0, \dots, a_{n-1}) $ and $ b = (b_0, \dots, b_{n-1}) $ in $ \CC^n. $
	For a sheaf $ \cV = (V,\phi) $ on $ BG $, the orbifold HRR theorem reads
	\begin{align*}
	\chi(BG, \phi) = \dim V^G = \int_{IBG} \orbch(\cV) = \sum_{i=0}^{n-1} \int_{BZ_{g_i}} \chi_\phi(g_i) = \sum_{i=0}^{n-1} \frac{1}{|Z_{g_i}|} \chi_\phi(g_i) = \frac{1}{|G|} \sum_{g \in G} \chi_\phi(g),
	\end{align*}
	where the last equality holds because 
	$$ |[g]| \times |Z_{g}| = |G| $$ 
	for all $ g $ in $ G. $
	For a pair of sheaves $ \cV = (V, \phi) $ and $ \cW = (W, \psi) $ on $ BG $, the orbifold HRR formula reads
	\begin{align*}
	\chi(\phi, \psi) = & \dim \Hom(V,W)^G = \inprod{\orbv(\cV),\orbv(\cW)}_{IBG} = \inprod{\chi_\phi,\chi_\psi}_W \\
	& = \sum_{i=0}^{n-1} \frac{1}{|Z_{g_i}|} \overline{\chi_\phi(g_i)}\chi_\psi(g_i) = \frac{1}{|G|} \sum_{g \in G} \chi_\phi(g^{-1}) \chi_\psi(g),
	\end{align*}
	which is a standard result in representation theory. When $ \phi $ is an irreducible representation of $ G, $ then $ \chi(\phi, \psi) $ counts the copies of $ \phi $ in $ \psi. $
\end{ex}
\begin{ex}[The orbifold HRR formula for $ B\mu_n $ is Parseval's theorem.]\label{ex_descrete_Parseval}
	Consider $ \cX = B\mu_n. $ The orbifold Euler pairing on $ K(B\mu_n) = \ZZ[x]/(x^n-1) $ is given by
	\begin{equation*}
	\inprod{x^i, x^j} = \delta_{ij}
	\end{equation*}
	for $ i,j = 0, \dots, n-1 $ and extended linearly. The orbifold Mukai pairing on
	$$ A(IB\mu_n) \cong \CC^n $$ 
	is the weighted inner product
	\begin{equation*}
	\inprod{a,b}_W = \frac{1}{n} \sum_{i=0}^{n-1} \cj{a}_i b_i
	\end{equation*}
	in $ \CC $ for all $ a = (a_0, \dots, a_{n-1}) $ and $ b = (b_0, \dots, b_{n-1}) $ in $ \CC^n. $
	Recall that the orbifold Mukai vector is the inverse discrete Fourier transform
	\begin{equation*}
	\orbv: \CC[x]/(x^n-1) \cong \CC^n \to \CC^n, \quad f \mapsto \check{f},
	\end{equation*}
	and its inverse is the discrete Fourier transform $ f \mapsto \hat{f}. $
	The orbifold HRR formula (\ref{eq_HRR2}) now reads
	\begin{equation*}
	\inprod{f, g} = \inprod{\check{f},\check{g}}_W
	\end{equation*}
	for all $ f, g $ in $ \ZZ[x]/(x^n-1). $ Extending the orbifold Euler pairing to $ K(B\mu_n)_\CC = \CC[x]/(x^n-1) $ sesquilinearly, we have
	\begin{equation*}
	\inprod{a, b}_W = \inprod{\hat{a},\hat{b}}
	\end{equation*}
	in $ \CC $ for all $ a, b $ in $ \CC^n, $
	which is \tb{Parseval's theorem} for the discrete Fourier transform.
\end{ex}
\begin{rmk}[The numerical  Grothendieck ring and the numerical Chow ring]
	Let $ \cX $ be a proper quotient stack. The abelian group $ K(\cX) $ can be huge, e.g., its rank can be infinity, and the orbifold Euler pairing is in general degenerate, i.e., the map
	\begin{equation*}
	\Phi: K(\cX) \to \Hom(K(\cX), \ZZ), \quad x \mapsto \chi(x, \ \cdot\ )
	\end{equation*}
	has a non-trivial kernel. Therefore, it is useful to define the \tb{numerical Grothendieck group}
	\begin{equation*}
	N(\cX) = K(\cX)/\ker \Phi
	\end{equation*}
	of $ \cX. $ Similarly, one can define $ N^0(\cX) $ from $ K^0(\cX). $ If $ x $ is an element in $ K(\cX), $ let $ \gamma(x) $ denote the image of $ x $ in $ N(\cX). $ A morphism $ f: \cY \to \cX $ of proper quotient stacks induces a ring homomorphism
	\begin{equation*}
	f^N: N(\cX) \to N(\cY), \quad \gamma(x) \mapsto \gamma(f^K(x)).
	\end{equation*}
	Then we have a non-degenerate orbifold Euler pairing
	\begin{equation*}
	\chi: N(\cX) \times N(\cX) \to \ZZ
	\end{equation*}
	on $ N(\cX). $
	
	Suppose $ \cX $ is connected and smooth. Then $ N(\cX) \cong N^0(\cX) $ is a ring. In this case, we call both $ N(\cX) $ and $ N^0(\cX) $ the \tb{numerical Grothendieck ring} of $ \cX $. We define the \tb{numerical Chow ring} 
	\begin{equation*}
	R(\cX) = A(\cX)/\sim_\text{num}
	\end{equation*}
	of $ \cX $ where $ \sim_\text{num} $ denotes the numerical equivalence: for an element $ v $ in $ A^i(\cX), $ we say $ v \sim_\text{num} 0 $ if $ \deg (vw) = 0 $ for all elements $ w $ in $ A^{\dim \cX -i}(\cX). $ The inertia stack $ I\cX $ is a proper smooth quotient DM stack, so there is a numerical Chow ring $ R(I\cX) $ of $ I\cX. $ The orbifold Todd class $ \orbtd_\cX $ is now a unit in $ R(I\cX)_\CC. $ The orbifold Chern character map and the orbifold Mukai vector map descend to two injective homomorphisms
	\begin{equation*}
	\orbch: N(\cX) \to R(I\cX)_\CC \quad \text{and} \quad \orbv: N(\cX) \to R(I\cX)_\CC
	\end{equation*}
	of rings and abelian groups respectively, which become $ \CC $-algebra isomorphism and $ \CC $-linear isomorphism respectively after tensored with $ \CC. $ 
	The involution
	\begin{equation*}
	(\ {\cdot}\ )^\vee: R(I\cX)_\CC \to R(I\cX)_\CC
	\end{equation*}
	is defined in the same way as that in Definition \ref{defn_involution}. The orbifold Mukai pairing
	\begin{equation*}
	\inprod{\cdot \ {,} \ \cdot }_{I\cX}: R(I\cX)_\CC \times R(I\cX)_\CC \to \CC
	\end{equation*}
	is still defined by (\ref{eq_orbv_pairing}). Now the orbifold Mukai vector map is a linear isometry
	\begin{equation*}
	\orbv: \left(N(\cX)_\CC, \chi\right) \xrightarrow{\simeq} \left(R(I\cX)_\CC, \inprod{\cdot \ {,} \ \cdot }_{I\cX}\right),
	\end{equation*}
	where both sides carry non-degenerate bilinear forms. The orbifold HRR formula (\ref{eq_HRR2}) remains the same except that the elements $ x $ and $ y $ are now in $ N(\cX) $ with the corresponding orbifold Mukai vectors $ \orbv(x) $ and $ \orbv(y) $ in $ R(I\cX)_\CC. $ 
	
	A morphism $ f: \cY \to \cX $ of proper smooth quotient DM stacks gives a commutative diagram:
	\vspace{5pt}
	\begin{equation*}
	\begin{tikzcd}[column sep=3em,row sep=3.5em,every label/.append style={font=\normalsize}]
	N(\cX) \arrow[r, "f^N"{above, outer sep = 2pt}] \arrow[d, "\orbv"{left, outer sep = 2pt}] & N(\cY) \arrow[d, "\orbv"{right, outer sep = 2pt}] \\
	R(I\cX)_\CC \arrow[r, "f^R"{above, outer sep = 2pt}] & R(I\cY)_\CC
	\end{tikzcd}\vspace{15pt}
	\end{equation*}
\end{rmk}
\begin{ex}
	Let $ C $ be an elliptic curve. The Mukai vector coincides with the Chern character since $ \td_C = 1 $, and gives an isomorphism of abelian groups
	\begin{equation*}
	v: K(C) \xrightarrow{\simeq} A(C) \cong \ZZ \oplus \Pic(C), \quad x \mapsto (\rk(x), c_1(x)).
	\end{equation*}
	Hence the rank of $ K(C) $ is uncountable since $ \Pic(C) \cong C $. The HRR formula says
	\begin{equation*}
	\chi(x,y) = \inprod{v(x), v(y)} = \rk(x)c_1(y) - \rk(y)c_1(x)
	\end{equation*}
	in $ \ZZ. $ The Mukai vector map descends to an isomorphism
	\begin{equation*}
	v: N(C) \xrightarrow{\simeq} R(C) \cong \ZZ^2, \quad x \mapsto (\rk(x), \deg(x))
	\end{equation*}
	of abelian groups.
	Observe that the numerical Grothendieck ring $ N(C) $ is much smaller than $ K(X) $ and only carries the numerical information of a coherent sheaf, i.e., its rank and degree. The Euler pairing $ \chi $ on $ N(C) $ can be identified with a non-degenerate skew-symmetric form on $ \ZZ^2 $ represented by the matrix
	\begin{equation*}
	\begin{pmatrix}
	0 & 1 \\
	-1 & 0 \\
	\end{pmatrix}.
	\end{equation*}
\end{ex}

\vspace{2pt}

\section{Equivariant moduli theory}\label{sec_moduli_sheaves_stacks}
In this section we review the constructions of moduli spaces of sheaves on projective schemes and projective stacks, and define $ G $-equivariant moduli spaces of sheaves on a projective scheme under an action of a finite group $ G $.


\subsection{Moduli spaces of sheaves on projective schemes}\label{sec_moduli_sheaves_schemes}
In this section we review Simpson's construction of moduli spaces of stable sheaves on projective schemes over $ \CC $ in \cite{simpson1994moduli}. 

Let's first recall some notions which are used to define Hilbert polynomials of coherent sheaves on projective schemes.
\begin{defn}
	The dimension of a coherent sheaf $ E $ on a scheme $ X $ is the dimension of its support as a subscheme of $ X $, i.e.,
	\begin{equation*}
	\dim E = \dim \Supp E.
	\end{equation*}
\end{defn}
\begin{rmk}
	The dimension of a zero sheaf is $ -\infty $ since the dimension of the empty set is $ -\infty $.
\end{rmk}
\begin{defn}
	A polynomial $ P(z) $ in $ \QQ[z] $ is called \tb{numerical} if $ P(\ZZ) \subeq \ZZ. $
\end{defn}
\begin{defn}
	Let $ X $ be a projective scheme. A choice of an ample line bundle $ H $ on $ X $ is called a \tb{polarization} of $ X. $ A pair $ (X, H) $ where $ H $ is an ample line bundle on $ X $ is called a \tb{polarized} projective scheme.
\end{defn}
\begin{defn-prop}
	Let $ (X, H) $ be a polarized projective scheme of dimension $ d $. Let $ E $ be a coherent sheaf on $ X. $ There is a unique numerical
	polynomial 
	$$ P(E, z) = a_n(E) z^n + \dots + a_1(E)z + a_0(E) $$
	in $ \QQ[z] $ such that
	\begin{equation*}
	P(E, m) = \chi(X, E\otimes H^{\otimes m})
	\end{equation*}
	in $ \ZZ $ for all $ m $ in $ \ZZ. $ We call $ P(E) $ the \tb{Hilbert polynomial} of $ E $ with respect to $ H. $ If $ E = 0 $, then $ P(E) = 0 $; else, $ n = \dim E $. 
	If $ E \neq 0, $ then the \tb{reduced Hilbert polynomial} of $ E $ is defined by the monic polynomial\footnote{It can also be defined by $ P(E, z)/n!a_n(E) $ as in \cite[Definition 1.2.3]{huybrechts2010geometry}. This will give the same notion of stability for sheaves.}
	\begin{equation*}
	p(E, z) = \frac{P(E, z)}{a_n(E)}
	\end{equation*}
	in $ \QQ[z] $.
\end{defn-prop}
\begin{notation}\label{notation_rk_deg}
	Let $ (X, H) $ be a connected smooth polarized projective scheme of dimension $ d $. There are two rings associated with $ X $: the numerical Grothendieck ring $ N(X) $ and the numerical Chow ring $ R(X). $ There is a \tb{rank} map, i.e., a ring homomorphism
	\begin{equation*}
	\rk: N(X) \to \ZZ
	\end{equation*}
	defined by
	\begin{equation*}
	\rk(x) = \rk(V) - \rk(W)
	\end{equation*}
	for any element $ x = \gamma(V) - \gamma(W) $ in $ N(X).$
	There is also a \tb{degree} map, i.e., a group homomorphism 
	$$ \deg: N(X) \to \ZZ $$
	defined by
	\begin{equation*}
	\deg(x) = c_1(x) h^{d-1}
	\end{equation*}
	for all $ x $ in $ N(X), $ where $ h = c_1(H) $ in $ R^1(X) $. If $ E $ is torsion-free, then $ n = d $ and by the HRR theorem, we have 
	\begin{equation*}
	a_d(E) = \frac{\rk(E)h^d}{d!}.
	\end{equation*}
\end{notation}
\begin{rmk}
	Let $ (X, H) $ be a smooth polarized projective scheme. The Hilbert polynomial is additive on short exact sequences in $ \Coh(X) $ since the Euler characteristic $ \chi(X, \ \cdot\ ) $ is. Therefore, the Hilbert polynomial of a coherent sheaf
	only depends on its numerical class in $ N(X) $, and hence descends to an additive map
	\begin{equation*}
	P: N(X) \to \QQ[z], \quad \gamma(E) \mapsto P(E).
	\end{equation*}
	The Hilbert polynomial of a coherent sheaf $ E $ on $ X $ is determined by the Mukai vector $ v(E) $ in $ R(X)_\QQ $ since the map $ v: N(X) \to R(X)_\QQ$ is injective.
\end{rmk}
\begin{ex}
	Let $ X = \pt. $ A coherent sheaf on $ X $ is the same as a vector space $ V. $ There is only one ample line bundle $ H $ on $ X $, i.e., the one-dimensional vector space $ \CC. $ Therefore, the Hilbert polynomial of a vector space $ V $ is a constant
	\begin{equation*}
	P(V, z) = \dim V
	\end{equation*}
	and its reduced Hilbert polynomial is simply $ p(V, z) = 1. $
\end{ex}
\begin{ex}
	Let $ X = \PP^n $ with the canonical ample line bundle $ H = \cO_{\PP^n}(1). $ The Hilbert polynomial of a line bundle $ \cO_{\PP^n}(k) $ is given by
	\begin{equation*}
	P(\cO_{\PP^n}(k), z) = \binom{z+k+n}{n}.
	\end{equation*}
	This defines an additive map
	\begin{equation*}
	P: N(\PP^n) = \frac{\ZZ[x]}{\inprod{(x-1)^{n+1}}} \to \QQ[z], \quad \gamma(\cO_{\PP^n}(k)) = x^k \mapsto \binom{z+k+n}{n}.
	\end{equation*}
	In particular, the Hilbert polynomial of $ \cO_Y $ for a hypersurface $ Y $ of degree $ d $ in $ \PP^n $ is given by
	\begin{equation*}
	P(\cO_Y, z) = P(\cO_{\PP^n}, z) - P(\cO_{\PP^n}(-d), z) = \binom{z+n}{n} - \binom{z-d+n}{n}.
	\end{equation*}
	For example, a degree $ d $ smooth curve $ C $ in $ \PP^2 $ corresponds to the Hilbert polynomial
	\begin{equation*}
	P(\cO_C, z) = dz + 1 - g
	\end{equation*}
	where $ g $ is the genus of $ C. $
\end{ex}
\begin{ex}
	Let $ (X, H) $ be a polarized projective scheme. If $ Y $ is a subscheme of $ n $ points on $ X, $ then the Hilbert polynomial
	\begin{equation*}
	P(\cO_Y, z) = n
	\end{equation*}
	is a constant, and the reduced Hilbert polynomial $ p(\cO_Y, z) = 1. $
\end{ex}
\begin{ex}
	Let $ X $ be an elliptic curve in $ \PP^2 $ with the canonical ample line bundle $ H = \cO_X(1). $ Then $ \deg(H) = 3. $ Let $ L $ be a line bundle on $ X $ with degree $ d. $ By the HRR theorem for curves, we have
	\begin{equation*}
	\chi(X, L \otimes H^{\otimes m}) = \deg(L \otimes H^{\otimes m}) + 1 - 1 = 3m + d,
	\end{equation*}
	and hence
	\begin{equation*}
	P(L, z) = 3z + d.
	\end{equation*} 
\end{ex}
\begin{ex}
	Let $ (X, H) $ be a polarized $ K3 $ surface. By the HRR theorem, the Hilbert polynomial of a coherent sheaf $ E $ on $ X $ with Mukai vector $ v(E) = (r, c_1, s) $ is
	\begin{equation}\label{eq_Hilb_poly_K3}
	P(E, z) = \frac{rh^2}{2} z^2 + (c_1 h)z + r + s.
	\end{equation}
	Take an ample line bundle $ H = \cO_X(2) $ in Example \ref{ex_K3_Fermat}. Then we have $ h^2 = 16, $ and hence the Hilbert polynomial of an ideal sheaf $ I_Y $ of a subscheme $ Y $ of $ n $ points on $ X $ is
	\begin{equation*}
	P(I_Y, z) = 8 z^2 + 2-n.
	\end{equation*}
\end{ex}
\begin{rmk}
	For $ X = \PP^n, $ the Hilbert polynomial map $ P: N(\PP^n) \to \QQ[z] $ is injective, but this is not true in general. Let $ X \to \PP^1 $ be an elliptic $ K3 $ surface with a section $ S $ and a fiber $ F $ such that $ \Pic(X) = \ZZ[S] \oplus \ZZ[F]. $ Let $ s = c_1(S) $ and $ f = c_1(F) $. Then we have the intersection data
	\begin{equation*}
	s^2 = -2, \quad sf = 1, \quad f^2 = 0.
	\end{equation*}
	Choose a line bundle $ H $ with $ h = c_1(H) = s + 3f. $ Then we have 
	\begin{equation*}
	h^2 = 4, \quad hs = 1, \quad hf = 1.
	\end{equation*}
	Therefore, $ H $ is an ample line bundle by the Nakai-Moishezen criterion.
	Choose two numerical classes $ x $ and $ y $ in $ N(X) $ such that
	\begin{equation*}
	\ch(x) = (0,h,0) \quad \text{and} \quad \ch(y) = (0, 4f, 0)
	\end{equation*}
	in the numerical Chow ring $ R(X) \cong \ZZ \oplus \Pic(X) \oplus \ZZ. $ From equation (\ref{eq_Hilb_poly_K3}), we have
	\begin{equation*}
	P(x, z) = h^2 z = 4z \quad \text{and} \quad P(y, z) = (4fh)z = 4z.
	\end{equation*}
	Therefore, $ x $ and $ y $ map to the same Hilbert polynomial although they are different classes in $ N(X) $.
\end{rmk}
Now we recall the notion of pure sheaves.
\begin{defn}
	Let $ E $ be a nonzero coherent sheaf on a scheme $ X. $ We say the sheaf $ E $ is \tb{pure} if $ \dim F = \dim E $ for all nonzero subsheaves $ F \subset E $ on $ X. $ 
\end{defn}
\begin{rmk}
	The support of a pure sheaf is pure dimensional, but the converse is not true. A pure sheaf is torsion-free on its support. Indeed, pure sheaves are a generalization of torsion-free sheaves.
\end{rmk}
\begin{ex}
	Let $ X $ be a scheme of dimension $ d $.
	\begin{enumerate}[font=\normalfont,leftmargin=*]
		\item Every vector bundle on $ X $ is pure of dimension $ d $.
		\item Every torsion-free sheaf on $ X $ is pure of dimension $ d $.
		\item The structure sheaf $ \cO_Y $ of a subscheme $ Y $ of dimension $ n $ is pure of dimension $ n $.
		\item If $ d > 0, $ the sheaf $ \cO_X \oplus \cO_p $ where $ p $ is a point on $ X $ is not pure.
	\end{enumerate}
\end{ex}
We introduce a notation to compare polynomials in $ \QQ[z]. $
\begin{notation}
	Let $ p $ and $ q $ be two polynomials in $ \QQ[z]. $ The inequality $ p \leq q $ means $ p(z) \leq q(z) $ for $ z \gg 0. $
	Similarly, the strict inequality $ p < q $ means $ p(z) < q(z) $ for $ z \gg 0. $
\end{notation}
Now we can define Gieseker stability of pure sheaves.
\begin{defn}
	Let $ (X, H) $ be a polarized projective scheme. Let $ E $ be a coherent sheaf on $ X. $ We say $ E $ is \tb{$ H $-(semi)stable}\footnote{The $ H $-stability was first introduced by Gieseker in \cite{gieseker1977moduli} for vector bundles on surfaces, and it is commonly called Gieseker $ H $-stability or Gieseker stability.} if $ E $ is a nonzero pure sheaf and
	\begin{equation*}
	p(F) \ (\leq)\footnote{There are two statements here: the inequality $ \leq $ is used to define $ H $-semistable sheaves, and the strict inequality $ < $ is used to define $ H $-stable sheaves.} \ p(E)
	\end{equation*}
	for all proper nonzero subsheaves $ F \subset E $ on $ X. $ When the ample line bundle $ H $ is understood, $ H $-semistable (resp. $ H $-stable) sheaves are also called \tb{semistable} (resp. \tb{stable}) sheaves. We say $ E $ is \tb{strictly semistable} if it is semistable but not stable; in this case, there is a proper subsheaf $ F $ of $ E $ with $ p(F) = p(E), $ and we say $ F $ is a \tb{destablizing} subsheaf of $ E $ or destabilizes $ E. $
\end{defn}
\begin{rmk}
	Let $ E $ be a coherent sheaf on a polarized projective scheme $ (X, H). $ By \cite[Proposition 1.2.6]{huybrechts2010geometry}, to check the stability or semistability of $ E, $ it's sufficient to consider saturated subsheaves $ F \subset E $, i.e., subsheaves $ F \subset E $ with torsion-free quotient $ E/F $. Suppose $ X $ is connected and smooth. Then there is a well-defined rank $ \rk(F) $ for every coherent sheaf $ F $ on $ X. $
	If $ E $ is torsion-free, then all subsheaves $ F \subset E $ are torsion-free, and hence it's sufficient to consider subsheaves $ F \subset E $ with $ 0 < \rk(F) < \rk(E) $. In particular, a torsion-free sheaf of rank one is stable with respect to any polarization.
\end{rmk}
\begin{ex}
	Let $ X = \pt. $ Let $ V $ be a vector space of dimension $ d $. As a coherent sheaf on $ X, $ $ V $ is stable if $ d = 1 $ and strictly semistable if $ d > 1. $ 
\end{ex}
\begin{ex}
	Let $ X = \PP^1$ with $ H = \cO_{\PP^1}(1). $
	\begin{enumerate}[font=\normalfont,leftmargin=*]
		\item Every line bundle $ \cO_{\PP^1}(d) $ is stable.
		\item The sheaf $ \cO_{\PP^1}\oplus\cO_p $ for any point $ p $ is not semistable because it is not pure.
		\item Consider the vector bundle $ E = \cO_{\PP^1}(d_1) \oplus \cO_{\PP^1}(d_2) $ with Hilbert polynomial $ P_E(z) = 2z + d_1 + d_2 + 2. $ If $ d_1 \neq d_2, $ then $ E $ is not semistable because the subsheaf $ F = \cO_{\PP^1}(d) $ gives
		$$ p(F, z) = z + d + 1 > p(E, z) = z + (d_1 + d_2)/2 + 1 $$ 
		where $ d = \max(d_1,d_2). $ If $ d_1 = d_2 = d, $ then $ E $ is strictly semistable because the subsheaf $ F = \cO_{\PP^1}(d) $ destabilizes $ E, $ i.e., $ p_F = p_E. $ 
		\item In general, a vector bundle $ E = \cO_{\PP^1}(d_1) \oplus \cdots \oplus \cO_{\PP^1}(d_n) $ where $ d_1 \geq \cdots \geq d_n $ with $ n \geq 2 $ is never stable since $ \cO_{\PP^1}(d_1) $ destabilizes $ E, $ and is only semistable when all $ d_i  $'s are the same. 
	\end{enumerate} 
\end{ex}
\begin{ex}
	Let $ (X, H) $ be a polarized $ K3 $ surface. The structure sheaf $ \cO_Y $ of a subscheme $ Y $ consisting of $ n $ points is strictly semistable unless $ n = 1 $ (in which case it is stable). An ideal sheaf $ I_Y $ of such a subscheme $ Y $ is always stable since its only proper subsheaves of dimension $ 2 $ are ideal sheaves of $ Z $ of $ m > n $ points (containing those $ n $ points), and hence
	\begin{equation*}
	p(I_Z, z) = z^2 + \frac{4-2m}{h^2} < z^2 + \frac{4-2n}{h^2} = p(I_Y, z).
	\end{equation*}
\end{ex}
Now we can define moduli functors of semistable sheaves as in \cite{huybrechts2010geometry}, but we will define them as moduli stacks (i.e., categories fibered in groupoids over $ (\Sch/\CC) $) which have the advantage of keeping track of automorphisms of sheaves.
\begin{defn}\label{defn_moduli_stack_sheaves_on_schemes}
	Let $ (X, H) $ be a polarized projective scheme. Fix a numerical polynomial $ P $ of degree $ d \leq \dim X $ in $ \QQ[z]. $ The \tb{moduli stack} $ \cM_H(X,P) $ of $ H $-semistable sheaves on $ X $ with Hilbert polynomial $ P $ is a category over $ (\Sch/\CC) $ defined as follows.
	\begin{enumerate}[font=\normalfont,leftmargin=*]
		\item An object in $ \cM_H(X,P) $ over a scheme $ Y $ is (the isomorphism class of) a coherent sheaf $ E $ on $ X \times Y $ such that $ E $ is flat over $ Y $ and for all closed points $ y $ in $ Y, $ the restriction $ E_y $ of $ E $ on $ X \times \{y\} $ is semistable with Hilbert polynomial $ P $.
		\item An arrow $ A: E \to F $ in $ \cM_H(X,P) $ over a morphism $ a: Y \to Z $ is an isomorphism $ A: E \to (\id_X \times a)^*F $ of coherent sheaves on $ X \times Y. $
	\end{enumerate}
	The moduli stack $ \cM_H^s(X,P) $ of $ H $-stable sheaves on $ X $ with Hilbert polynomial $ P $ is defined similarly. 
\end{defn}
\begin{rmk}
	If we forget about the automorphisms of sheaves, then we can also view $ \cM_{H}^{(s)}(X,P) $ as a \tb{moduli functor}.
\end{rmk}
The moduli stack $ \cM_H^{(s)}(X,P) $ is an algebraic stack. We will often suppress the subscript $ H $ and simply write $ \cM_H^{(s)}(X,P) $ when the polarization $ H $ is understood.
\begin{rmk}
	The moduli stacks $ \cM_H(X,P) $ or $ \cM_H^s(X,P) $ may be empty. For example, take $ X = \PP^1 $ with the standard polarization $ H = \cO_{\PP^1}(1). $ 
	\begin{enumerate}[font=\normalfont,leftmargin=*]
		\item If $ P(z) = -1, $ then $ \cM(X,P) $ is empty since $ P_E(m) = h^0(\PP^1,E(m)) \geq 0 $ for $ m \gg 0 $ for every sheaf $ E $ on $ \PP^1 $.
		\item If $ P(z) = rz + a $ where $ a $ and $ r \geq 2 $ are integers. Then $ \cM^s(X,P) $ is empty since there are no stable vector bundles of rank $ \geq 2 $. Here we used two facts: all torsion-free sheaves on a curve are locally free, and all vector bundles on $ \PP^1 $ are direct sums of line bundles on $ \PP^1 $. Note that $ \cM(X,P) $ is non-empty if and only if $ a $ is a multiple of $ r, $ in which case $ \cM(X,P)(\pt) $ consists of only one semistable vector bundle of rank $ r $.
	\end{enumerate}
\end{rmk}
Recall the following:
\begin{defn}
	Let $ X $ be a projective scheme with a very ample line bundle $ H $. Let $ E $ be a coherent sheaf on $ X. $ Let $ m $ be an integer. We say $ E $ is \tb{$ m $-regular} (with respect to $ H $) if
	\begin{equation*}
	H^i(X, E \otimes H^{\otimes (m-i)}) = 0 \quad \text{for all} \ i > 0.
	\end{equation*}
\end{defn}
\begin{rmk}
	If a coherent sheaf $ E $ is $ m $-regular, then $ E \otimes H^{\otimes m} $ is globally generated with vanishing higher cohomologies and
	\begin{equation*}
	P(E, m) = h^0(X, E\otimes H^{\otimes m}).
	\end{equation*}
\end{rmk}

Fix a polarized projective scheme $ (X,H). $ Let's take a numerical polynomial $ P $ of degree $ d \leq \dim X $ in $ \QQ[z] $. Suppose there is some semistable sheaf $ E $ on $ X $ such that $ P(E, z) = P(z). $ By the boundedness of semistable sheaves, the triple $ (X,H,P) $ determines an integer $ m $ such that every semistable sheaf $ E $ with Hilbert polynomial $ P $ is $ m $-regular and hence 
$$ P(m) = P(E, m) = h^0(X, E \otimes H^{\otimes m}). $$ 
Let $ N = P(m) $. Define a vector bundle 
$$ F = \cO_X^{\oplus N} \otimes H^{\otimes -m} $$
on $ X $. Then there is a projective scheme $ \Quot_X(F,P) $ (called a \tb{Quot scheme}) whose closed points parametrize quotients of $ V, $ i.e., surjective morphisms of sheaves
\begin{equation}\label{eq_quot_scheme}
V \xrightarrow{q} E \to 0
\end{equation}
with $ P(E, z) = P(z) $ up to equivalence. Here we say $ (E, q) $ is equivalent to $ (E',q') $ if there is an isomorphism $ \phi: E \to E' $ such that $ \phi \circ q = q' $.
There is an open subscheme $ Q $ (resp. $ Q^s $) of $ \Quot_X(V,P) $ which parametrizes quotients (\ref{eq_quot_scheme}) where $ E $ is semistable (resp. stable).

Now we recall Harder-Narasimhan and Jordan-Hölder filtrations.
\begin{defn}
	Let $ E $ be a nonzero pure sheaf on $ X. $ A \tb{Harder-Narasimhan filtration} for $ E $ is a chain of subsheaves
	\begin{equation*}
		0 = E_0 \subset E_1 \subset \cdots \subset E_n = E 
	\end{equation*}
	such that each factor $ F_i = E_i/E_{i-1} $ is semistable and
	\begin{equation*}
		p(F_1) > p(F_2) > \cdots > p(F_n).
	\end{equation*}
\end{defn}
\begin{defn}
	Let $ E $ be a semistable sheaf on $ X. $ A \tb{Jordan-Hölder filtration} for $ E $ is a chain of subsheaves
	\begin{equation*}
		0 = E_0 \subset E_1 \subset \cdots \subset E_n = E 
	\end{equation*}
	such that the factors $ \gr_i(E) = E_i/E_{i-1} $ are stable with reduced Hilbert polynomial $ p(E). $
\end{defn}
\begin{rmk}
	\begin{enumerate}[font=\normalfont,leftmargin=*]
		\item Let $ E $ be a nonzero pure sheaf on $ X. $ There there exists a unique Harder-Narasimhan filtration $ 0 = E_0 \subset E_1 \cdots \subset E_n = E $ such that
		\begin{equation*}
			p(E_1) > p(E_2) > \cdots > p(E_n).
		\end{equation*} 
		If $ E $ is semistable, then $ E_1 = E. $ Suppose $ E $ is not semistable. Then for $ i = 1, \dots, n-1, $ each subsheaf $ E_i $ destablizes $ E, $ i.e., $ p(E_i) > p(E) $. We call the first subsheaf $ E_1 $ the \tb{maximal destablizing subsheaf} of $ E. $
		\item Let $ E $ be a semistable sheaf on $ X. $ Its Jordan-Hölder filtration exists but is not unique. However, the factors $ \gr_i(E) $ are unique up to ordering, and hence the sheaf $ \gr(E) = \bigoplus \gr_i(E) $ is unique. 
	\end{enumerate}
\end{rmk}
\begin{defn}
	Two semistable sheaves $ E $ and $ E' $ on $ X $ are said to be \tb{$ S $-equivalent} if $ \gr(E) \cong \gr(E'). $ 
\end{defn}

The following result was first obtained in \cite{simpson1994moduli} and later summarized in \cite[Chapter 4]{huybrechts2010geometry}. 
\begin{thm}\label{thm_coarse_moduli_space_stable_sheaves_on_schemes}
	Let $ (X, H) $ be a polarized projective scheme. Fix a numerical polynomial $ P $ in $ \QQ[z]. $ There is a projective scheme $ M_H(X,P) $ that universally corepresents the moduli functor $ \cM_H(X,P) $. Closed points in $ M_H(X,P) $ are in bijection with $ S $-equivalent classes of semistable sheaves with Hilbert polynomial $ P $. There is an open subscheme $ M_H^s(X,P) \subset M_H(X,P) $ that universally corepresents the moduli functor $ \cM_H^s(X,P). $ Moreover, the natural morphism
	\begin{equation*}
	\pi: Q^s \to M_H^s(X,P)
	\end{equation*}
	is a geometric quotient and is also a principal $ \PGL(N, \CC) $-bundle.
\end{thm}
\begin{defn}
	The scheme $ M_H(X,P) $ (resp. $ M_H^s(X,P) $) in Theorem \ref{thm_coarse_moduli_space_stable_sheaves_on_schemes} is called the \tb{Gieseker moduli space} of $ H $-semistable (resp. $ H $-stable) sheaves on $ X $ with Hilbert polynomial $ P $.
\end{defn}
\begin{rmk}
	Theorem \ref{thm_coarse_moduli_space_stable_sheaves_on_schemes} gives a compactification $ M_H(X,P) $ of the moduli space $ M_H^s(X,P) $ when $ M_H^s(X,P) $ is not empty.
\end{rmk}
Let $ (X, H) $ be a smooth polarized projective scheme. We can also fix a numerical class instead of a Hilbert polynomial. Let $ x $ be an element in $ N(X). $ Then the moduli stack (or functor) $ \cM_H^{(s)}(X,x) $ is defined similarly to $ \cM_H^{(s)}(X,P) $. Suppose there is some semistable sheaf $ E $ on $ X $ with numerical class $ x. $ Then $ x $ determines a numerical polynomial $ P(z) $ in $ \QQ[z]. $ There are finitely many numerical classes $ y $ in $ N(X) $ such that $ P(y, z) = P(z) $ since the scheme $ Q $ is of finite type. Let $ Q^s(x) $ denote the subscheme of $ Q^s $ determined by the element $ x. $
Then we have the following:
\begin{cor}
	Let $ (X, H) $ be a smooth polarized projective scheme. Choose an element $ x $ in the numerical Grothendieck group $ N(X). $ There is a projective scheme $ M_H(X,x) $ that universally corepresents the moduli functor $ \cM_H(X,x) $. Closed points in $ M_H(X,x) $ are in bijection with $ S $-equivalent classes of semistable sheaves with numerical class $ x $. There is an open subscheme $ M_H^s(X,x) \subset M_H(X,x) $ that universally corepresents the moduli functor $ \cM_H^s(X,x). $ Moreover, the natural morphism
	\begin{equation*}
	\pi: Q^s(x) \to M_H^s(X,x)
	\end{equation*}
	is a geometric quotient and is also a principal $ \PGL(N, \CC) $-bundle.
\end{cor}
Let's turn the discussion above into a definition.
\begin{defn}
	Let $ (X, H) $ be a polarized projective scheme. Let $ x $ be an element in the numerical Grothendieck group $ N(X) $. The projective scheme $ M_H(X, x) $ is called the Gieseker moduli space of $ H $-semistable sheaves on $ X $ with numerical class $ x $. The quasi-projective scheme $ M_H^s(X, x) $ is called the Gieseker moduli space of $ H $-stable sheaves on $ X $ with numerical class $ x $.
\end{defn}
\begin{rmk}
	As discussed before, the moduli space $ M_H(X, x) $ may be empty. If we choose $ x = -1 = -[\cO_X]$ in $ N(X), $ then the moduli space $ M_H(X, -1) $ is empty for an obvious reason.
\end{rmk}
The following is a trivial example.
\begin{ex}
	Let $ X = \pt. $ Then $ N(X) = K(X) = \ZZ. $ Since $ \dim \pt = 0, $ a numerical polynomial $ P $ of degree zero in $ \QQ[z] $ is a constant, i.e., $ P(z) = n $ for some integer $ n. $ The Hilbert polynomial of coherent sheaves on $ X $ (i.e., finite-dimensional vector spaces) descends to an inclusion
	\begin{equation*}
	P: \ZZ \into \QQ[z], \quad n \mapsto n.
	\end{equation*}
	In other words, fixing an element $ x $ in $ K(\pt) $ is the same as fixing a constant numerical polynomial $ P $ in $ \QQ[z]. $ Choose $ P(z) = n \geq 0 $. Note that every vector space is $ m $-regular for any integer $ m $ for a trivial reason. Then the scheme $ \Quot_\pt(\CC^n, n) $ parametrizes surjective linear maps
	\begin{equation*}
	\CC^n \xrightarrow{q} \CC^n \to 0
	\end{equation*}
	up to equivalence. Here the $ \CC^n $ in the target is viewed as a semistable sheaf on $ \pt. $ Note that every quotient map $ q $ is an element in $ \GL(n, \CC), $ and equivalent to the identity map. In other words, $ \Quot_\pt(\CC^n, n) = \pt. $ The moduli stacks $ \cM(\pt,n) $ for non-negative $ n $'s are summarized as follows.
	\begin{enumerate}[font=\normalfont,leftmargin=*]
		\item If $ n = 0, $ then the moduli stack $ \cM(\pt,0) = \cM^s(\pt,0) = \pt $.
		\item If $ n = 1, $ then the moduli stack $ \cM(\pt,1) = \cM^s(\pt,1) = [\pt/\CC^*] = B\CC^* $ has a coarse moduli space $ M(\pt,1) = M^s(\pt,1) = \pt. $
		\item If $ n > 1, $ then the moduli stack $ \cM(\pt,n) = [\pt/GL(n, \CC)] = B\GL(n, \CC) $ has a coarse moduli space $ M(\pt,n) = \pt, $ but $ \cM^s(\pt, n) $ is empty. Note that the moduli stack $ \cM(\pt,n) $ encodes the automorphism group $ \GL(n, \CC) $ of the vector space $ \CC^n $, which is a semistable sheaf on $ X = \pt. $
	\end{enumerate}
\end{ex}
Now we recall the slope and its corresponding $ \mu $-stability of a coherent sheaf on a smooth projective scheme.
\begin{defn}
	Let $ (X, H) $ be a connected smooth polarized projective scheme in Notation \ref{notation_rk_deg}. We define a \tb{slope} map\footnote{This is neither additive nor multiplicative.}
	\begin{equation*}
	\mu: N(X) \to \QQ \cup \{\infty\}
	\end{equation*}
	by setting
	\begin{equation*}
	\mu(x) = \left \{
	\begin{array}{ll}
	\frac{\deg(x)}{\rk(x)} & \mbox{if $\rk(x) \neq 0$},\\
	\infty & \mbox{otherwise}.
	\end{array}
	\right.
	\end{equation*}
	If $ h = c_1(H) $ in $ R^1(X), $ then we will use the notation $ \mu_h $ instead of $ \mu $ when we want to emphasize the dependence of the slope map on $ h. $
\end{defn}
\begin{rmk}
	Let $ E $ be a nonzero coherent sheaf on $ X $. The slope $ \mu(E) $ is a notion similar to but in general coarser than the reduced Hilbert polynomial $ p(E) $.
\end{rmk}
\begin{defn}
	Let $ (X, H) $ be a connected smooth polarized projective scheme. We say a coherent sheaf $ E $ on $ X $ is \tb{$ \mu $-(semi)stable} if $ E $ is a nonzero torsion-free sheaf and
	\begin{equation*}
	\mu(F) \ (\leq) \ \mu(E)
	\end{equation*}
	for all subsheaves $ F $ of $ E $ with $ 0 < \rk(F) < \rk(E). $
	We also refer to \tb{$ \mu $-stability} as \tb{slope stability}.
\end{defn}
\begin{rmk}
	\begin{enumerate}[font=\normalfont,leftmargin=*]
		\item Every torsion-free sheaf of rank one is slope stable.
		\item Tensoring a nonzero torsion-free sheaf $ E $ by a line bundle $ L $ shifts its slope:
		\begin{equation*}
		\mu(E \otimes L) = \frac{\deg(E \otimes L)}{\rk(E)} = \frac{\deg(E)+\rk(E)\deg(L)}{\rk(E)} = \mu(E) + \deg(L).
		\end{equation*}
		Therefore, tensoring by line bundles preserves slope (semi)stability.
		\item The notions of stability and $ \mu $-stability are equivalent for torsion-free sheaves on a smooth polarized projective curve $ (C, H) $. Suppose $ C $ has genus $ g. $ Let $ E $ be a torsion-free sheaf on $ C. $ By the HRR theorem, the Hilbert polynomial of $ E $ is
		\begin{equation*}
		P(E, z) = \rk(E) \deg(H) z + \deg(E) + \rk(E)(1-g),
		\end{equation*}
		and hence the reduced Hilbert polynomial of $ E $ is
		\begin{equation*}
		p(E, z) = z + \frac{\mu(E)+ 1-g}{\deg(H)}.
		\end{equation*}
		For a subsheaf $ F $ of $ E $ with $ 0 < \rk(F) < \rk(E), $ we see that
		\begin{equation*}
		\mu(F)\ (\leq)\ \mu(E)\ \Leftrightarrow\ p(F)\ (\leq)\ p(E).
		\end{equation*}
	\end{enumerate}
\end{rmk}
In general, we have the following:
\begin{lem}[{\cite[Lemma 1.2.13 and Lemma 1.2.14]{huybrechts2010geometry}}]\label{lem_slope_Gieseker}
	Let $ (X, H) $ be a connected smooth polarized projective scheme. Let $ E $ be a nonzero torsion-free sheaf $ E $ on $ X $. Then we have a chain of implications:
	\begin{equation*}
	E \ \text{is}\ \mu\text{-stable}\ \Rightarrow\ E\ \text{is stable}\ \Rightarrow\ E\ \text{is semistable}\ \Rightarrow\ E \ \text{is}\ \mu\text{-semistable}.
	\end{equation*}
	If $ \gcd(\rk(E),\deg(E)) = 1,$ then we also have the implication
	\begin{equation*}
	E \ \text{is}\ \mu\text{-semistable}\ \Rightarrow\ E \ \text{is}\ \mu\text{-stable}.
	\end{equation*}
\end{lem}
\begin{rmk}
	Lemma \ref{lem_slope_Gieseker} tells that if $ x $ is an element in $ N(X) $ with 
	$$ \rk(x) > 0 \quad \text{and} \quad \gcd(\rk(x),\deg(x)) = 1, $$
	then $ M(X,x) = M^s(X,x). $ In particular, if $ \rk(x) = 1, $ then $ M(X,x) = M^s(X,x). $
\end{rmk}
\begin{rmk}
	In general, we expect the moduli space $ M(X,x) $ or $ M^s(X,x) $ to consist of singular connected components of various dimensions even if $ X $ is a smooth projective variety. The situation is much better when $ X $ has nice geometric properties, for example, when $ X $ is a $ K3 $ surface.
\end{rmk}
Consider a $ K3 $ surface $ X $. Since the Mukai vector map $ v: N(X) \to R(X) $ is an isomorphism of abelian groups, fixing an element $ x $ in $ N(X) $ is the same as fixing a vector $ v $ in $ R(X) $ for moduli spaces of stable sheaves on $ X $. Recall a few notions.
\begin{defn}
	A vector $ v $ in $ R(X) $ is \tb{primitive} if it is not a multiple of another vector $ w $ in $ R(X) $. An element $ x $ in $ N(X) $ is said to be primitive if it is not a multiple of another element $ y $ in $ N(X). $
\end{defn}
\begin{defn}
	The \tb{ample cone} of $ X $ is a convex cone defined by
	$$ \Amp(X) = \{C \in R^1(X)_\RR: C = \sum a_i [D_i] \ \text{where } a_i > 0 \ \text{and} \ D_i \ \text{is an ample divisor} \}. $$
\end{defn}
It's shown in \cite{yoshioka1996chamber} that an element $ x $ in $ N(X) $ or equivalently, its Mukai vector in $ R(X) $ determines a countable locally finite set of hyperplanes in $ \Amp(X) $, the ample cone of $ X, $ which are called \emph{walls}. We will call them \emph{$ x $-walls} to emphasize their dependence on the element $ x $.
\begin{defn}
	Let $ x $ be an element in $ N(X) $. An ample line bundle $ H $ on $ X $ is said to be \emph{$ x $-generic}, if the divisor class $ h = c_1(H) \in \Amp(X) $ does not lie on any of the $ x $-walls.
\end{defn}

\begin{lem}\label{lem_slope_semistable_implies_stable}
	Let $ (X, H) $ be a polarized $ K3 $ surface. Let $ x $ be an element in $ N(X) $ with Mukai vector $ v = (r, c_1, s) $ in $ R(X) $ where $ r > 0 $. Suppose $ (r, c_1) $ is primitive and $ H $ is $ x $-generic. If a torsion-free sheaf $ E $ on $ X $ with numerical class $ x $ is $ \mu $-semistable, then it is $ \mu $-stable. 
\end{lem}
\begin{proof}
	Take a $ \mu $-semistable sheaf $ E $ on $ X $ with numerical class $ x. $ Suppose it is not $ \mu $-stable. Then there is a subsheaf $ F \subset E $ with $ 0 < \rk(F) < \rk(E) $ such that $ \mu(F) = \mu(E), $ i.e.,
	\begin{equation*}
	\frac{\deg(F)}{\rk(F)} = \frac{\deg(E)}{\rk(E)}.
	\end{equation*}
	Since $ H $ is $ x $-generic, we must have $ \rk(E) c_1(F) = \rk(F) c_1(E) $ in $ R^1(X). $ This implies $ (r, c_1) = r/\rk(F)(\rk(F),c_1(F)), $ which is impossible since $ (r, c_1) $ is primitive.
\end{proof}
More generally, fixing a primitive element $ x $ in $ N(X) $ such that the polarization is $ x $-generic guarantees that semistability is the same as stability. We summarize the results in \cite[Corollary 4.6.7 and Remark 4.6.8]{huybrechts2010geometry} (applied to $ K3 $ surfaces) together with \cite[Corollary 10.2.1 and Proposition 10.2.5]{huybrechts2016lectures} in the following proposition.
\begin{prop}\label{prop_K3_smoothness_and_dim}
	Let $ (X, H) $ be a polarized $ K3 $ surface. Let $ x $ be an element in $ N(X) $ with Mukai vector $ v  = (r,c_1,s) $ in $ R(X) $ and degree $ d = \deg(x) $. Then $ M_H^s(X,x) $ is either empty or a smooth quasi-projective scheme of dimension $ n = 2 - \inprod{v}^2. $ Suppose either of the following conditions is satisfied:
	\begin{enumerate}[label=\textnormal{(\arabic*)},font=\normalfont,leftmargin=2em]
		\item \label{prop_K3_smoothness_and_dim_condition_1} $ \gcd(r, d, s) = 1.$
		\item $ x $ is primitive and $ H $ is $ x $-generic.
	\end{enumerate}
	Then $ M_H(X, x) = M_H^s(X, x), $ and it is either empty or a smooth projective scheme of dimension $ n. $ Moreover, $ M_H(X, x) $ is a fine moduli space if \ref{prop_K3_smoothness_and_dim_condition_1} holds.
\end{prop}
Recall that the Hilbert scheme $ \Hilb^n(X) $ for a $ K3 $ surface $ X $ is an irreducible symplectic manifold of dimension $ 2n, $ and $ \Hilb^n(X) $ is deformation equivalent to $ \Hilb^n(Y) $ for any other $ K3 $ surface $ Y. $ Now we state the following classical result.
\begin{thm}[{\cite[Theorem 8.1]{yoshioka2001moduli}}]\label{thm_yoshioka}
	Let $ (X, H) $ be a polarized $ K3 $ surface. Let $ x $ be an element in $ N(X) $ with $ \rk(x) > 0 $. Assume $ x $ is primitive and $ H $ is $ x $-generic. Then $ M_H(X, x) = M_H^s(X, x), $ which is non-empty if and only if $ \inprod{v(x)^2} \leq 2. $
	If $ M_H(X, x) $ is non-empty, then it is an irreducible symplectic manifold of dimension $ n = 2 - \inprod{v^2} $ deformation equivalent to $ \Hilb^{n/2}(X) $.
\end{thm}
We will generalize the two results above to equivariant moduli spaces of stable sheaves on a $ K3 $ surface with symplectic automorphisms in Sections \ref{sec_HRR_K3/G} and \ref{sec_proof_main_thm} respectively.


\subsection{Moduli spaces of sheaves on projective stacks}\label{sec_moduli_sheaves_stacks2}
In this section we review Nironi's construction of moduli spaces of stable sheaves on projective stacks over an algebraically closed field in \cite{nironi2009moduli}. 


Let $ S $ be a scheme. We first recall the notion of tame stacks.
\begin{defn}[{\cite[Definition 3.1]{abramovich2008tame}}]
	Let $ \cX $ be an algebraic stack locally of finite presentation over $ S $ with a finite inertia, i.e., the canonical morphism $ I\cX \to \cX $ is finite. Let $ \pi: \cX \to Y $ denote the coarse moduli space of $ \cX $. The stack $ \cX $ is \tb{tame} if the pushforward functor $ \pi_*: \Qcoh(\cX) \to \Qcoh(Y) $ is exact.
\end{defn}
\begin{rmk}
	Every separated DM stack locally of finite type over $ \CC $ is tame. In particular, every separated quotient DM stack over $ \CC $ is a tame DM stack.
\end{rmk}

A closely related notion is good moduli spaces of algebraic stacks.
\begin{defn}[{\cite[Definition 4.1]{alper2013good}}]
	Let $ \cX $ be an algebraic stack over $ S. $ Let $ f: \cX \to Y $ be a morphism from $ \cX $ to an algebraic space $ Y. $ We say $ f: \cX \to Y $ is a \tb{good moduli space} if the morphism $ f $ is quasi-compact, the pushforward functor $ f_*: \Qcoh(\cX) \to \Qcoh(Y) $ is exact, and the natural morphism $ \cO_Y \to f_*\cO_\cX $ is an isomorphism.
\end{defn}
\begin{ex}
	The quotient stack $ [\CC/\CC^*] $ (with the obvious action) has no coarse moduli space, but has a good moduli space $ [\CC/\CC^*] \to \pt. $
\end{ex}
\begin{rmk}
	For a tame algebraic stack over $ S $, its coarse moduli space is always a good moduli space.
\end{rmk}

Let $ S = \Spec k $ for a field $ k. $ Recall the notion of generating sheaves on tame DM stacks over $ k $ introduced in \cite{olsson2003quot}. 
\begin{defn-prop}[{\cite[Theorem 5.2]{olsson2003quot}}]
	Let $ \cX $ be a tame DM stack over $ k $ with a coarse moduli space $ \pi: \cX \to Y $. A locally free sheaf $ \cV $ on the stack $ \cX $ is a \tb{generating sheaf} if it satisfies either of the following equivalent conditions:
	\begin{enumerate}[font=\normalfont,leftmargin=*]
		\item For every geometric point $ x $ in $ \cX, $ the local representation $ \phi_x: G_x \to \GL(V_x) $ of the automorphism group $ G_x $ on the fiber $ V_x $ of $ \cV $ at $ x $ contains every irreducible representation of $ G_x $.
		\item For every quasi-coherent sheaf $ \cF $ on $ \cX, $ the natural morphism
		\begin{equation*}
		\left(\pi^*\pi_*(\cF \otimes \cV^\vee)\right) \otimes \cV \to \cF
		\end{equation*}
		is surjective.
	\end{enumerate}
\end{defn-prop}
\begin{ex}
	For a scheme $ X $ over $ k, $ every vector bundle on $ X $ is a generating sheaf. 
\end{ex}
\begin{ex}
	Let $ \cX = BG $ for a finite group $ G $. Then the regular representation $ \rho_\reg $ of $ G $ is a generating sheaf on $ BG. $
\end{ex}
\begin{ex}
	For a quotient stack $ \cX = [X/G] $ where $ G $ is a finite group, the vector bundle $ \cO_\cX \otimes \rho_\reg $ is a generating sheaf on $ \cX $. We will consider this example again in the next section.
\end{ex}
\begin{rmk}
	In general we do not know whether there are generating sheaves on tame DM stacks. However, they do exist on projective stacks. Generating sheaves are a key ingredient in the construction of moduli spaces of sheaves on projective stacks.
\end{rmk}

\begin{defn}[{\cite[Definition 2.20]{nironi2009moduli}}]\label{defn_proj_stack}
	Let $ \cX $ be a separated quotient DM stack over $ k $. We say $ \cX $ is a \tb{projective} (resp. quasi-projective) stack if it is a tame DM stack and its coarse moduli space is a projective (resp. quasi-projective) scheme.
\end{defn}
\begin{rmk}
	If the base field $ k = \CC $ or has characteristic zero, then a separated quotient DM stack over $ k $ is always a tame DM stack. Therefore, Definition \ref{defn_proj_stack} reduces to Definition \ref{defn_projective_stack_over_C} when $ k = \CC $.
\end{rmk}
Now we assume the base field $ k $ is algebraically closed. Let $ \cX $ be a projective stack over $ k $ with a coarse moduli space $ \pi: \cX \to Y. $ Since $ \pi: \cX \to Y $ is a proper morphism and $ \cX $ is tame, we have an exact functor
\begin{equation*}
\pi_*: \Coh(\cX) \to \Coh(Y).
\end{equation*}
Therefore, for every coherent sheaf $ \cE $ on $ \cX $, its pushforward $ \pi_*\cE $ is a coherent sheaf on $ Y $ and we have
\begin{equation*}
H^i(\cX, \cE) = H^i(Y, \pi_*\cE)
\end{equation*}
for all $ i \geq 0 $. Furthermore, since $ \cX $ is a quasi-compact quotient stack, it has a generating sheaf by \cite[Theorem 5.5]{olsson2003quot}.  Since $ Y $ is a projective scheme, it has ample line bundles. To facilitate further discussions, we make the following definition.
\begin{defn}
	Let $ \cX $ be a projective stack over $ k $ with a coarse moduli space $ \pi: \cX \to Y. $ 
	Let $ \cV $ be a generating sheaf on $ \cX $, and let $ L $ be an ample line bundle on $ Y. $ Let $ \cH = \pi^*L. $ The pair $ (\cH, \cV) $ is called a \tb{polarization} of $ \cX. $ The triple $ (\cX, \cH, \cV) $ is called a \tb{polarized} projective stack.
\end{defn}
Let $ \cE $ be a nonzero coherent sheaf on $ \cX. $ The dimension of $ \cE $ is defined to be the dimension of its support, as in the case of schemes. By \cite[Proposition 3.6]{nironi2009moduli}, we have
\begin{equation*}
\dim \cE = \dim \pi_*(\cE \otimes \cV^\vee).
\end{equation*}
By the projection formula for DM stacks, we have
\begin{equation*}
\chi(\cX, \cE \otimes \cV^\vee \otimes \cH^{\otimes m}) = \chi(Y, \pi_*(\cE \otimes \cV^\vee) \otimes L^{\otimes m}) = P(\pi_*(\cE \otimes \cV^\vee), m)
\end{equation*}
in $ \ZZ $ for all $ m $ in $ \ZZ. $ Therefore, the assignment $ m \mapsto \chi(\cX, \cE \otimes \cV^\vee \otimes \cH^{\otimes m}) $ is a numerical polynomial of degree $ n = \dim \cE. $ Now we can define a modified Hilbert polynomial of the sheaf $ \cE $ on the stack $ \cX $ in a similar fashion to the usual Hilbert polynomial of a sheaf on schemes.
\begin{defn-prop}
	Let $ (\cX, \cH, \cV) $ be a polarized projective stack over $ k $. Let $ \cE $ be a coherent sheaf on $ \cX. $ The \tb{modified Hilbert polynomial} of $ \cE $ with respect to the polarization $ (\cH, \cV) $ is the unique numerical polynomial
	\begin{equation*}
	\orbP(\cE, z) = a_n(\cE) z^n + \cdots + a_1(\cE) + a_0(\cE)
	\end{equation*}
	in $ \QQ[z] $ such that
	\begin{equation*}
	\orbP(\cE, m) = \chi(\cX, \cE \otimes \cV^\vee \otimes \cH^{\otimes m})
	\end{equation*}
	in $ \ZZ $ for all $ m $ in $ \ZZ. $ If $ \cE = 0, $ then $ \orbP(\cE, z) = 0; $ else, $ n = \dim \cE. $
	If $ \cE \neq 0, $ then the reduced Hilbert polynomial of $ \cE $ is defined by the monic polynomial
	\begin{equation*}
	\orbp(\cE, z) = \frac{\orbP(\cE, z)}{a_n(\cE)}
	\end{equation*}
	in $ \QQ[z]. $
\end{defn-prop}
\begin{rmk}
	The modified Hilbert polynomial is additive on short exact sequences on $ \Coh(\cX) $ because the functor $ \Coh(\cX) \to \Coh(Y), \ \cE \to \pi_*(\cE \otimes \cV^\vee) $ is a composition of exact functors and hence is exact.
\end{rmk}
\begin{ex}\label{ex_BG_reg}
	Let $ \cX = BG $ for a finite group $ G. $ 
	Then $ (\cX, \rho_\reg, \cO_{\cX}) $ is a polarized projective stack where $ \rho_\reg $ is the regular representation of $ G, $ and $ \cO_{\cX} = \rho_0 $ is a trivial one-dimensional representation of $ G. $
	Let $ \cE $ be a coherent sheaf on $ BG, $ i.e., a representation $ \rho: G \to \GL(V) $ of $ G. $ Then the modified Hilbert polynomial of $ \rho $ is
	\begin{equation*}
	\orbP(\rho) = \chi(BG, \rho \otimes \rho_\reg^\vee) = \chi(\rho_\reg, \rho) = \dim V = \chi(\pt, V) = P(V).
	\end{equation*}
	Therefore, the modified Hilbert polynomial of the representation $ \rho $ of $ G $ is the same as the Hilbert polynomial of the underlying space $ V $ of $ \rho $. Later we will see this is true for quotient stacks by finite groups in general.
\end{ex}
We define pure sheaves on projective stacks as usual.
\begin{defn}
	A nonzero coherent sheaf $ \cE $ on $ \cX $ is pure if $ \dim \cF = \dim \cE $ for all nonzero subsheaves $ \cF \subset \cE $ on $ \cX. $ 
\end{defn}
\begin{defn}
	Let $ (\cX, \cH, \cV) $ be a polarized projective stack over $ k $. A coherent sheaf $ \cE $ on $ \cX $ is \tb{$ (\cH, \cV) $-(semi)stable} if it is a nonzero pure sheaf and
	$$ \orbp(\cF) \ \ (\leq)\ \ \orbp(\cE) $$ 
	for all proper nonzero subsheaves $ \cF \subset \cE $. When the choice of $ \cV $ is understood, $ (\cH, \cV) $-stability is also called \tb{$ \cH $-stability}; if both $ \cH $ and $ \cV $ are understood, $ (\cH, \cV) $-stability is also called \tb{stability}.
\end{defn}
Now we define moduli stacks of semistable sheaves on the projective stack $ \cX $ as in the case of projective schemes.
\begin{defn}\label{defn_moduli_stack_sheaves_on_stacks}
	Let $ (\cX, \cH, \cV) $ be a polarized projective stack over $ k $. Fix a numerical polynomial $ P $ of degree $ d \leq \dim \cX $ in $ \QQ[z]. $ The moduli stack $ \cM_{\cH,\cV}(\cX,P) $ of $ \cH $-semistable sheaves on $ \cX $ with modified Hilbert polynomial $ P $ is a category over $ (\Sch/k) $ defined as follows.
	\begin{enumerate}[font=\normalfont,leftmargin=*]
		\item An object in $ \cM_{\cH,\cV}(\cX,P) $ over a scheme $ Y $ is (the isomorphism class of) a coherent sheaf $ \cE $ on $ \cX \times Y $ such that $ \cE $ is flat over $ Y $ and for all closed points $ y $ in $ Y, $ the restriction $ \cE_y $ of $ E $ on $ X \times \{y\} $ is $ \cH $-semistable with modified Hilbert polynomial $ P $.
		\item An arrow $ A: \cE \to \cF $ in $ \cM_{\cH,\cV}(\cX,P) $ over a morphism $ a: Y \to Z $ is an isomorphism $ A: \cE \to (\id_\cX \times a)^*\cF $ of coherent sheaves on $ \cX \times Y. $
	\end{enumerate}
	The moduli stack $ \cM_{\cH,\cV}^s(\cX,P) $ of $ \cH $-stable sheaves on $ \cX $ with modified Hilbert polynomial $ P $ is defined similarly. For an element $ x $ in the numerical Grothendieck group $ N(\cX), $ we define two moduli stacks $ \cM_{\cH,\cV}(\cX,x) $ and $ \cM_{\cH,\cV}^s(\cX,x) $ similarly to $ \cM_{\cH,\cV}(\cX,P) $ and $ \cM_{\cH,\cV}^s(\cX,P) $ respectively.
\end{defn}
We will omit the subscript $ \cV $ and simply write $ \cM_{\cH}(X,x) $ and $ \cM_{\cH}^s(X,x) $ if the generating sheaf $ \cV $ is understood. If both $ \cH $ and $ \cV $ are understood, we may simply write $ \cM(X,x) $ and $ \cM^s(X,x) $ in place of $ \cM_{\cH,\cV}(X,x) $ and $ \cM_{\cH,\cV}^s(X,x) $. The following definition was used in \cite{nironi2009moduli} although not specifically stated.
\begin{defn}
	Let $ (\cX, \cH, \cV) $ be a polarized projective stack over $ k $ with a coarse moduli space $ \pi: \cX \to Y $. Let $ \cE $ be a coherent sheaf on $ \cX. $ Let $ r $ be an integer. We say $ \cE $ is $ r $-regular if the coherent sheaf $ \pi_*(\cE \otimes \cV^\vee) $ on $ Y $ is $ r $-regular.
\end{defn}
Fix a numerical polynomial $ P $ of degree $ d \leq \dim \cX $ in $ \QQ[z]. $ Suppose the moduli stack $ \cM_\cH(\cX,P) $ is non-empty. By \cite[Theorem 5.1]{nironi2009moduli}, there is an integer $ m $ such that all semistable sheaves with modified Hilbert polynomial $ P $ are $ m $-regular. Let $ N = P(m). $ Let $ \cF = \cV^{\oplus N} \otimes \cH^{\otimes -m}. $ By \cite[Proposition 4.20]{nironi2009moduli}, as in the case of schemes, there is a projective scheme $ \Quot_\cX(\cF,P) $ whose closed points parametrize pairs $ (\cE, q) $ in surjective morphisms of sheaves
\begin{equation*}\label{eq_quot_scheme_on_stacks}
\cF \xrightarrow{q} \cE \to 0
\end{equation*}
with $ \orbP(\cE) = P $ up to equivalence. As in the case of schemes, two pairs $ (\cE, q) $ and $ (\cE', q') $ are equivalent if there is an isomorphism $ \phi: \cE \to \cE' $ such that $ \phi \circ q = q'. $
Semistable sheaves on $ \cX $ also have JH filtrations and $ S $-equivalent classes. Now we combine the results in \cite[Theorems 5.1, 6.20 and 6.21]{nironi2009moduli} in the following theorem.
\begin{thm}
	Let $ (\cX, \cH, \cV) $ be a polarized projective stack over $ k $. Choose a numerical polynomial $ P $ of degree $ d \leq \dim \cX $ in $ \QQ[z]. $ 
	\begin{enumerate}[font=\normalfont,leftmargin=*]
		\item Suppose the moduli stack $ \cM_{\cH,\cV}(\cX,P) $ is not empty. Then there exists an integer $ m $ such that semistable sheaves with modified Hilbert polynomial $ P $ are $ m $-regular. Denote 
		$$ \cF = \cV^{\oplus {N}} \otimes \cH^{\otimes -m}, $$
		where $ N = P(m) $ is a positive integer. 
		Then $ \cM_{\cH,\cV}(\cX,P) $ is a quotient stack over $ k $ of the form
		\begin{equation*}
		\cM_{\cH,\cV}(\cX,P) = [Q/\GL(N, k)],
		\end{equation*}
		where $ Q $ is an open subscheme of the projective scheme $ \Quot_\cX(\cF,P) $.
		There is a good moduli space 
		$$ \cM_{\cH,\cV}(\cX,P) \to M_{\cH,\cV}(\cX,P), $$
		where $ M_{\cH,\cV}(\cX,P) $ is a projective scheme parametrizing $ S $-equivalent classes of $ (\cH,\cV) $-semistable sheaves on $ \cX $ with modified Hilbert polynomial $ P $, which is a categorical quotient of $ Q $.
		\item Suppose the substack $ \cM_{\cH,\cV}^s(\cX,P) \subset \cM_{\cH,\cV}(\cX,P) $ is not empty. Then there is an open subscheme $ Q^s \subset Q $ such that
		\begin{equation*}
		\cM_{\cH,\cV}^s(\cX,P) = [Q^s/\GL(N, k)]
		\end{equation*}
		is a quotient substack of $ \cM_{\cH,\cV}(\cX,P) $ with a coarse moduli space 
		$$ \cM_{\cH,\cV}^s(\cX,P) \to M_{\cH,\cV}^s(\cX,P), $$
		where $ M_{\cH,\cV}^s(X,P) $ is an open subscheme of $ M_{\cH,\cV}(\cX,P) $ parametrizing $ (\cH,\cV) $-stable sheaves on $ \cX $ with numerical class $ x $, which is a geometric quotient of $ Q^s $ such that the morphism 
		$$ Q^s \to M_{\cH,\cV}^s(\cX,P) = Q^s/\PGL(N,k)$$ 
		is a principal $ \PGL(N,k) $-bundle.
	\end{enumerate}
\end{thm}
Since the moduli stacks $ \cM_{\cH,\cV}(\cX,P) $ and $ \cM_{\cH,\cV}^s(\cX,P) $ are algebraic stacks of finite type over $ k $, there are finitely many elements $ x $ in the numerical Grothendieck group $ N(\cX) $ such that $ \orbP(x) = P. $ Therefore, we have the following:
\begin{cor}\label{cor_coarse_moduli_space_sheaves_on_stacks}
	Let $ (\cX, \cH, \cV) $ be a polarized projective stack over $ k $. Choose an element $ x $ in the numerical Grothendieck group $ N(\cX). $
	\begin{enumerate}[font=\normalfont,leftmargin=*]
		\item Suppose the moduli stack $ \cM_{\cH,\cV}(\cX,x) $ is not empty. Then there exists an integer $ m $ such that semistable sheaves with modified Hilbert polynomial $ \orbP(x) $ are $ m $-regular. Denote 
		$$ \cF = \cV^{\oplus {N}} \otimes \cH^{\otimes -m}, $$ 
		where $ N = \orbP(x, m) $ is a positive integer. 
		Then $ \cM_{\cH,\cV}(\cX,x) $ is a quotient stack  over $ k $ of the form
		\begin{equation*}
		\cM_{\cH,\cV}(\cX,x) = [Q(x)/\GL(N,k)],
		\end{equation*}
		where $ Q(x) $ is an open subscheme of the projective scheme $ \Quot_\cX(\cF,\orbP(x)) $.
		There is a good moduli space 
		$$ \cM_{\cH,\cV}(\cX,x) \to M_{\cH,\cV}(\cX,x), $$
		where $ M_{\cH,\cV}(\cX,x) $ is a projective scheme parametrizing $ S $-equivalent classes of $ (\cH,\cV) $-semistable sheaves on $ \cX $ with numerical class $ x $, which is a categorical quotient of $ Q(x) $.
		\item Suppose the substack $ \cM_{\cH,\cV}^s(\cX,x) \subset \cM_{\cH,\cV}(\cX,x) $ is not empty. Then there is an open subscheme $ Q^s(x) \subset Q(x) $ such that
		\begin{equation*}
		\cM_{\cH,\cV}^s(\cX,x) = [Q^s(x)/\GL(N,k)]
		\end{equation*}
		is a quotient substack of $ \cM_{\cH,\cV}(\cX,x) $ with a coarse moduli space 
		$$ \cM_{\cH,\cV}^s(\cX,x) \to M_{\cH,\cV}^s(\cX,x), $$
		where $ M_{\cH,\cV}^s(X,x) $ is an open subscheme of $ M_{\cH,\cV}(\cX,x) $ parametrizing $ (\cH,\cV) $-stable sheaves on $ \cX $ with numerical class $ x $, which is a geometric quotient of $ Q^s(x) $ such that the morphism 
		$$ Q^s(x) \to M_{\cH,\cV}^s(\cX,x) = Q^s(x)/\PGL(N,k)$$ 
		is a principal $ \PGL(N,k) $-bundle.
	\end{enumerate}
\end{cor}
\begin{defn}
	In the setting of Corollary \ref{cor_coarse_moduli_space_sheaves_on_stacks}, the projective scheme $ M_{\cH,\cV}(\cX,x) $ is called the Gieseker moduli space of $ (\cH,\cV) $-semistable sheaves on $ \cX $ with numerical class $ x $ in $ N(\cX) $; the quasi-projective scheme $ M_{\cH,\cV}^s(\cX,x) $ is called the Gieseker moduli space of $ (\cH,\cV) $-stable sheaves on $ \cX $ with numerical class $ x $ in $ N(\cX) $.
\end{defn}
\begin{ex}
	Let $ \cX = B\mu_n $ over $ \CC $ with a coarse moduli space $ \pi: B\mu_n \to \pt. $ Consider the polarized projective stack $ (\cX, \rho_\reg, \rho_0) $ in Example \ref{ex_BG_reg}. Recall that 
	$$ N(B\mu_n) = K(B\mu_n) \cong \frac{\ZZ[x]}{(x^n-1)}. $$
	Then $ [\rho_0] = 1 $ and $ [\rho_\reg] = 1 + x + \cdots + x^{n-1} $ in $ N(B\mu_n). $
	As computed in Example \ref{ex_BG_reg}, the Hilbert polynomial of a representations $ \rho $ of $ \mu_n $ is simply $ \deg(\rho), $ and hence it descends to an additive map
	\begin{equation*}
	\orbP: \frac{\ZZ[x]}{(x^n-1)} \to \QQ[z], \quad [\rho] = f(x) \mapsto \orbP(\rho) = \deg(\rho) = f(1).
	\end{equation*}
	Take a constant numerical polynomial $ P(z) = d $ for some integer $ d \geq 0. $ The scheme $ \Quot_{B\mu_n}(\rho_\reg^{\oplus d}, d) $ parametrizes surjective morphisms
	\begin{equation*}
	\rho_\reg^{\oplus d} \xrightarrow{q} \rho \to 0
	\end{equation*}
	of sheaves on $ B\mu_n $ up to equivalence. A morphism $ q: \rho_\reg^{\oplus d} \to \rho $ is the same as an intertwining map between the two representations $ \rho_\reg^{\oplus d} $ and $ \rho $ of $ \mu_n $, i.e., a $ \mu_n $-equivariant linear map between their underlying spaces. Suppose we fix a representation $ \rho $ of degree $ d. $ Let $ \Hom_{\mu_n}(\rho_\reg^{\oplus d}, \rho) $ denote the set of intertwining maps between $ \rho_\reg^{\oplus d} $ and $ \rho $. Since $ \rho_\reg^{\oplus d} = \Ind \rho_o^{\oplus d}, $ the induced representation from $ d $-copies of trivial representation $ \rho_0 $ of $ B\mu_n, $ we have a set bijection
	\begin{equation*}
	\Hom_{\mu_n}(\rho_\reg^{\oplus d}, \rho) \cong \Hom(\rho_0^{\oplus d}, \rho_0^{\oplus d}) = \Hom(\CC^d, \CC^d) 
	\end{equation*}
	by the Frobenius reciprocity, so there is only one surjective morphism $ q:\rho_\reg^{\oplus d} \to \rho $ up to equivalence. 
	Note that every element 
	\begin{equation*}
	a_0 + a_1x+ \cdots + a_{n-1}x^{n-1}
	\end{equation*}
	in $ N(B\mu_n) $ is mapped to $ d $ if and only if $ a_0 + a_1 + \cdots + a_{n-1} = d. $ Therefore, we have a decomposition
	\begin{equation*}
	\Quot_{B\mu_n}(\rho_\reg^{\oplus d}, d) = \coprod_{[\rho]: \deg(\rho) = d} Q(\rho)
	\end{equation*}
	where each $ Q(\rho) $ is a point, i.e., the isomorphism class $ [\rho] $ of a representation $ \rho $ of $ \mu_n. $ This implies that each moduli space $ M(B\mu_n,[\rho]) $ is a point and $ M(B\mu_n,d) $ counts the number of isomorphism classes $ [\rho] $ such that $ \deg(\rho) = d, $ which is 
	\begin{equation*}
	l_{n,d} = \binom{n+d-1}{n-1},
	\end{equation*}
	the number of collections of non-negative integers $ a_0, a_1, \dots, a_{n-1} $ with $ \sum_{i=0}^{n-1} a_i = d. $ The moduli stacks $ \cM(B\mu_n,[\rho]) $ for different numerical classes $ [\rho] $ in $ N(B\mu_n) $ are summarized as follows.
	\begin{enumerate}[font=\normalfont,leftmargin=*]
		\item If $ d = 0, $ then the moduli stack $ \cM(B\mu_n,0) = \cM^s(B\mu_n,0) = \pt $.
		\item Let $ d = 1. $ Then the moduli stack 
		$$ \cM(B\mu_n,x^i) = \cM^s(B\mu_n,x^i) = B\CC^* \times \{x^i\} $$
		has a coarse moduli space $ M(B\mu_n,x^i) = M^s(B\mu_n,x^i) = \{x^i\} $ for $ i = 0, 1, \dots, n-1. $ So the moduli stack
		$$ \cM(B\mu_n,1) = \cM^s(B\mu_n,1) = \coprod_{i=0}^{n-1} \cM(B\mu_n,x^i) = B\CC^* \times \mu_n $$ 
		has a coarse moduli space $ M(B\mu_n,1) = M^s(B\mu_n,1) = \mu_n. $ Note that $ M^s(B\mu_n,1) $ is a trivial principal $ \mu_n $-bundle over $ M^s(\pt,1) = \pt. $
		\item If $ d > 1, $ then the moduli stack
		$$ \cM(B\mu_n,d) = \coprod_{[\rho]: \deg(\rho) = d} B\GL(d) \times \{[\rho]\} $$ has a coarse moduli space $ M(B\mu_n,d) $ which consists of $ l_{n,d} $ points, but $ \cM^s(B\mu_n, d) $ is empty since every representation $ \rho $ of $ \mu_n $ has a destabilizing sub-representation.
	\end{enumerate}
\end{ex}

Now we generalize the smoothness criteria of moduli spaces of sheaves on projective schemes in \cite[Proposition 10.1.11]{huybrechts2016lectures} to projective stacks.
\begin{prop}\label{prop_smoothness_criterion}
	Let $ (\cX, \cH, \cV) $ be a smooth polarized projective stack over $ k $. Choose an element $ x $ in $ N(\cX) $ such that the moduli stack $ \cM = \cM_{\cH,\cV}(\cX,x) $ is non-empty. Let $ t $ be a closed point in the moduli space $ M = M_{\cH,\cV}(\cX, x) $ corresponding to a stable sheaf $ \cE $ in the groupoid $ \cM(k) $.
	\begin{enumerate}[label=(\arabic*),font=\normalfont,leftmargin=*]
		\item \label{smoothness1} There is a natural isomorphism $ T_{t} M \cong \Ext^1(\cE,\cE) $.
		\item \label{smoothness2} If $ \Ext^2(\cE,\cE) = 0, $ then $ M $ is smooth at $ t. $
		\item \label{smoothness3} If $ \Pic(\cX) $ is smooth and the trace map $ \Ext^2(\cE,\cE) \to H^2(\cX, \cO_\cX) $ is injective, i.e., $ \Ext^2(\cE,\cE)_0 = 0, $ then $ M $ is smooth at $ t. $
	\end{enumerate}
\end{prop}
\begin{proof}
	Let $ M^s = M_{\cH,\cV}^s(\cX, x). $ By the first result in Corollary \ref{cor_coarse_moduli_space_sheaves_on_stacks}, the element $ x $ in $ N(\cX) $ determines an integer $ m $ such that all semistable sheaves with modified Hilbert polynomial $ \orbP(x) $ are $ m $-regular. In particular, $ \cE $ is $ m $-regular, which means $ \pi_*(\cE \otimes \cV^\vee) $ is $ m $-regular. Let $ N = \orbP(x, m), $ and let $ \cF = \cV^{\oplus N} \otimes \cH^{\otimes -r}.$
	The second result of Corollary \ref{cor_coarse_moduli_space_sheaves_on_stacks} says that there is an open subscheme $ Q^s \subset \Quot_\cX(\cF,\orbP(x)) $ such that the geometric quotient
	\begin{equation*}
	\pi: Q^s \to M^s
	\end{equation*}
	is a principal $ \PGL(N,k) $-bundle. Consider a point $ q $ in the fiber $ \pi^{-1}(t), $ which corresponds to a short exact sequence
	\begin{equation}\label{eq_quot_ses}
	0 \to \cK \xrightarrow{s} \cF \xrightarrow{q} \cE \to 0
	\end{equation}
	in $ \Coh(\cX) $. Let $ \cO_q = \PGL(N,k) \cdot q $ denote the orbit of $ q $ in $ Q^s. $ Then we have
	\begin{equation*}
	T_t M \cong T_q R^s/T_q \cO_q,
	\end{equation*}
	where the tangent space $ T_q Q^s \cong \Hom(\cK, \cE) $. Now we apply $ \Hom(\ {\cdot} \ ,\cE) $ to (\ref{eq_quot_ses}) and get a long exact sequence
	\begin{equation*}
	0 \to \End(\cE) \to \Hom(\cF,\cE) \xrightarrow{s^\sharp} \Hom(\cK,\cE) \to \Ext^1(\cE,\cE) \to \Ext^1(\cF,\cE) \to \cdots
	\end{equation*}
	We claim that 
	$$ \Ext^i(\cF,\cE) = 0 \quad \text{for all}\ i \geq 1. $$
	This is a consequence of the $ m $-regularity of $ \cE. $ Take an integer $ i \geq 1. $ Since $ \cF $ is locally free, we have
	\begin{align*}
	\Ext^i(\cF,\cE) & = H^i(\cX, \cE \otimes \cF^\vee) \\
	& = H^i(\cX, (\cE \otimes \cV^\vee \otimes \cH^{\otimes m})^{\oplus N}) \\
	& = \bigoplus_{k=1}^N H^i(\cX, \cE \otimes \cV^\vee \otimes \cH^{\otimes m}) \\
	& = \bigoplus_{k=1}^N H^i(Y, \pi_*(\cE \otimes \cV^\vee) \otimes L^{\otimes m}) & \text{since $ \cX $ is tame}\\
	& = 0 & \text{since $ \cE $ is $ m $-regular}\\
	\end{align*}
	The claim is proved. Note that $ T_q \cO_q \cong \img(s^\sharp) \subset \Hom(\cK, \cE). $  Therefore, we have
	\begin{equation*}
	\Ext^1(\cE,\cE) \cong \Hom(\cK, \cE)/ \img(s^\sharp) \cong T_t M,
	\end{equation*}
	which proves \ref{smoothness1}.
	From the vanishing of $ \Ext^2(\cF, \cE), $ we also have
	\begin{equation*}
	\Ext^1(\cK,\cE) \cong \Ext^2(\cE, \cE),
	\end{equation*}
	where $ \Ext^1(\cK,\cE) $ is the space of obstructions to deform $ \cE $ in $ M. $ Hence \ref{smoothness2} is proved. Since $ \cX $ is smooth and has the resolution property, there is a determinant morphism
	\begin{equation*}
	\det: M \to \Pic(\cX), \quad \cE \to \det(\cE),
	\end{equation*}
	where $ \det(\cE) $ is well-defined via a finite locally resolution of $ \cE. $ Since $ \det(\cE) $ is a line bundle, we have an isomorphism of $ k $-vector spaces:
	\begin{equation*}
	H^2(\cX, \cO_\cX) \cong \Ext^2(\det(\cE),\det(\cE))
	\end{equation*}
	Therefore, we can identify the trace map as the following map:
	\begin{equation*}
	\tr: \Ext^2(\cE,\cE) \to \Ext^2(\det(\cE),\det(\cE))
	\end{equation*}
	Consider a local Artinian $ k $-algebra $ (A, \frak{m}) $ where $ \frak{m} $ is the maximal ideal in $ A $. Let $ I $ be an ideal in $ A $ such that $ I \cdot \frak{m} = 0.$ Then the trace map induces a map
	\begin{equation*}
	\tr \otimes_k I: \Ext^2(\cE,\cE) \otimes_k I \to \Ext^2(\det(\cE),\det(\cE)) \otimes_k I
	\end{equation*}
	which sends the obstruction $ \ob(\cE,A) $ to lift an $ A/I $-flat deformation of $ \cE $ in $ M $ to the obstruction $ \ob(\det(\cE), A) $ to lift an $ A/I $-flat deformation of $ \det(\cE) $ in $ \Pic(\cX) $. Suppose the trace map is injective. Then the $ k $-linear map $ \tr \otimes_k I $ is injective. Therefore, the smoothness of $ \Pic(\cX) $ implies $ \ob(\det(\cE), A) = 0 $, and hence $ \ob(\cE,A) = 0. $ This proves \ref{smoothness3}.
\end{proof}
\begin{ex}
	Let $ M^s_g(r,d) $ denote the (coarse) moduli space of stable vector bundles of rank $ r $ and degree $ d $ on a smooth complex projective curve $ C $ of genus $ g. $ Note that vector bundles are the same as torsion-free sheaves on a curve. Suppose it is non-empty (which is the case whenever $ g \geq 2 $). Then we know that $ M^s_g(r,d) $ is smooth and that
	\begin{equation*}
	\dim M^s_g(r,d) = r^2(g-1) + 1.
	\end{equation*}
	This is a degenerate case of Proposition \ref{prop_smoothness_criterion} where $ \cX = C $ and $ \cV = \cO_C $. Smoothness is automatic since $ \dim C = 1. $ To compute the dimension of $ M^s_g(r,d), $ choose any $ E $ in $ M^s_g(r,d), $ and compute
	\begin{align*}
	\dim \Ext^1(E, E) & = 1 - \chi(E, E) \\ 
	& = 1 - \chi(C, E^\vee \otimes E) & \text{by Lemma \ref{lem_Euler_pairing}} \\
	& = 1 - r^2\chi(C, \cO_C) - r\deg(E^\vee \otimes E) & \text{by the HRR theorem}\\
	& = 1 + r^2(g-1). 
	\end{align*}
	Note this is independent of the degree $ d. $
	It's in general a difficult question whether the moduli space of stable sheaves on a variety $ X $ is irreducible or not, since we only know it's a quasi-projective scheme by the GIT construction. Fortunately, this is indeed the case for curves, which was first shown in \cite[Theorem 8.1]{seshadri1967space}. When $ g = 1, $ i.e., $ C $ is an elliptic curve, then a moduli space $ M^s_1(r,d) $ is either empty (when $ \gcd(r,d) \neq 1 $), or isomorphic to $ C $ itself (when $ \gcd(r,d) = 1 $) which was proved in \cite[Theorem 7]{atiyah1957vector}.
\end{ex}

\subsection{Equivariant moduli spaces of sheaves on projective schemes}
In this section we consider the following:
\begin{situation}\label{situation_finite_group}
	Let $ X $ be a smooth projective scheme, and let $ G $ be a finite group acting on $ X $. Then we have a commutative triangle
	\medskip
	\begin{equation*}
	\begin{tikzcd}[column sep=1.6em,row sep=3.2em]
	& X \arrow{ld}[swap]{p} \arrow{rd}{q} & \\
	\cX \arrow{rr}{\pi} & & Y,
	\end{tikzcd}\medskip
	\end{equation*}
	where $ p $ and $ q $ are finite, $ \pi $ is proper, $ \cX = [X/G] $ is a projective stack, $ Y = X/G $ is a projective scheme which is both the geometric quotient of $ X $ by $ G $ and the coarse moduli space of $ \cX $. Let $ L $ be an ample line bundle on $ Y $.
	Put
	$$ \cH = \pi^*L \quad \text{and} \quad H = q^*L. $$
	Then $ H $ is an ample line bundle on $ X $ since it's a pullback of an ample line bundle along a finite morphism.
	The line bundle $ \cH $ on $ \cX $ corresponds to an $ G $-equivariant ample line bundle $ (H, \phi) $ on $ X. $ Since it's a pullback from a line bundle on the coarse moduli space $ M, $ all local representations 
	\begin{equation*}
	\phi_x: G_x \to \GL(H_x)
	\end{equation*}
	of the stabilizers $ G_x $ on the fibers of $ H $ at every point $ x \in X $ are trivial.
	Consider the regular representation $ \rho_\reg $ of $ G $. Pulling it back along $ \cX \to BG $ gives a canonical generating sheaf 
	$$ \cV_\reg = \cO_\cX \otimes \rho_\reg $$
	on the stack $ \cX. $ Then we have a polarized projective scheme $ (X, H) $ and a polarized projective stack $ (\cX, \cH, \cV_{\reg}) $.
\end{situation}
By the categorical equivalence $ \Coh(\cX) \cong \Coh^G(X), $ we define the stability of $ G $-equivariant sheaves on $ X $ as the stability of the corresponding sheaves on $ \cX. $
\begin{defn}\label{stability_equivariant_sheaves}
	Let $ \cX = [X/G] $ in Situation \ref{situation_finite_group}. A $ G $-equivariant sheaf $ (E,\phi) $ on $ X $ is \tb{$ \cH $-semistable }(resp. \tb{$ \cH $-stable}) if the corresponding sheaf $ \cE $ on $ \cX $ is $ (\cH, \cV_\reg) $-semistable (resp. $ (\cH, \cV_\reg) $-stable).
\end{defn}
Recall that given an element $ x $ in $ N(\cX), $ there is a projective scheme $ M_\cH(\cX,x) $ parametrizing $ S $-equivalent classes of Gieseker $ (\cH, \cV_\reg) $-semistable sheaves on $ \cX $, which compactifies the quasi-projective scheme $ M_\cH^s(\cX,x) $ parametrizing Gieseker $ (\cH, \cV_\reg) $-stable sheaves on $ \cX $. These will be $ G $-equivariant moduli spaces of sheaves on $ X. $
\begin{defn}
	Let $ \cX = [X/G] $ in Situation \ref{situation_finite_group}. Let $ x $ be an element in $ N(\cX). $ We call the projective (resp. quasi-projective) scheme $ M_\cH(\cX,x) $ (resp. $ M_\cH^s(\cX,x) $) the \tb{$ G $-equivariant moduli space} of $ \cH $-semistable (resp. $ \cH $-stable) sheaves on $ X $ with numerical class $ x $.
\end{defn}
\begin{ex}
	Consider Example \ref{ex_K3_Fermat}. We have $ G = \mu_2, $
	\begin{equation*}
	X = \Proj \CC[x_0,x_1,x_2,x_3]/(x_0^4+x_1^4+x_2^4+x_3^4), \quad \text{and}
	\end{equation*}
	\begin{equation*}
	X/G = \Proj \CC[y_0, y_1, y_2, y_3, y_4, y_5]/(y_0y_2 - y_1^2, y_3y_5 - y_4^2, y_0^2+y_2^2+y_3^2+y_5^2).
	\end{equation*}
	We can choose an ample line bundle $ L $ to be the restriction of $ \cO_{\PP^5}(1) $ on $ X/G. $ Then $ H $ is isomorphic to  $ \cO_{X}(2), $ the restriction of  $ \cO_{\PP^3}(2) $ on $ X. $ The line bundle $ \cH = (H, \phi) $ on $ \cX $ restricts to trivial $ \mu_2 $-representations on each of the eight fixed points of $ \mu_2 $ on $ X. $ The regular representation of $ \mu_2 $ is 
	$$ \rho_\reg \cong \rho_0 \oplus \rho_1, $$
	where $ \rho_0 $ is the trivial irreducible representation of $ \mu_2 $ and $ \rho_1 $ is the only non-trivial irreducible representation of $ \mu_2 $. Then the generating sheaf
	\begin{equation*}
	\cV = \cO_\cX \oplus (\cO_\cX \otimes \rho_1)
	\end{equation*}
	is the direct sum of the structure sheaf on $ \cX $ and its twist by $ \rho_1 $.
\end{ex}

The following lemma identifies the modified Hilbert polynomials for sheaves $ \cE $ on $ \cX $ with the Hilbert polynomials for the pullback $ E = p^*\cE $ on $ X $.
\begin{lem}
	For a coherent sheaf $ \cE = (E, \phi) $ on $ \cX $ and an integer $ m, $
	\begin{equation*}
		H^i(X, E \otimes H^{\otimes m}) \cong H^i(\cX, \cE \otimes \cV_\reg \otimes \cH^{\otimes m})
	\end{equation*}
	for all $ i \geq 0, $ and hence
	\begin{equation*}
		\chi(X, E \otimes H^{\otimes m}) = \chi(\cX, \cE \otimes \cV_\reg \otimes \cH^{\otimes m}).
	\end{equation*}
\end{lem}
\begin{proof}
	The argument here is due to Promit Kundu. Take a coherent sheaf $ \cE = (E, \phi) $ on $ \cX. $ Fix an integer $ m $ and an integer $ i \geq 0. $	Since the morphism $ p: X \to \cX $ is finite and affine, the pushforward functor $ p_*: \Coh(X) \to \Coh(\cX) $
	is exact. Note that $ p_*\cO_X \cong \cO_\cX \otimes \rho_\reg = \cV_\reg $.  Applying the projection formula to $ p: X \to \cX $, we have
	\begin{equation*}
		p_*p^*(\cE \otimes \cH^{\otimes m}) \cong \cE \otimes \cH^{\otimes m} \otimes p_*\cO_X \cong \cE \otimes \cH^{\otimes m} \otimes \cV_\reg,
	\end{equation*}
	and hence
	\begin{align*}
		H^i(X, E \otimes H^{\otimes m}) \cong H^i(X, p^*(\cE \otimes \cH^{\otimes m})) \cong H^i(\cX, \cE \otimes \cV_\reg \otimes \cH^{\otimes m}).
	\end{align*}
\end{proof}

\begin{cor}\label{cor_P_orb_equals_P}
	Let $ \cE = (E, \phi) $ be a coherent sheaf on $ \cX $. Then
	$$ \orbP(\cE) = P(E). $$
\end{cor}
\begin{proof}
	Note that $ \rho_\reg $ is self-dual and hence $ \cV_\reg \cong \cV_\reg^\vee. $ We then have
	\begin{align*}
		P(E, m) = \chi(X, E \otimes H^{\otimes m}) = \chi(\cX, \cE \otimes \cV_\reg \otimes \cH^{\otimes m}) = \orbP(\cE, m)
	\end{align*}
	in $ \ZZ $ for all $ m $ in $ \ZZ $.
\end{proof}
\begin{rmk}
	An alternative definition of stability of $ G $-equivariant sheaves on $ X $ was introduced in \cite{amrutiya2015moduli} as follows: a $ G $-equivariant sheaf $ (E, \phi) $ on $ X $ is (semi)stable if $ E $ is a nonzero pure sheaf and $ p_F(z) \ (\leq)\ p_E(z) $ for all nonzero proper $ G $-equivariant subsheaves $ (F, \phi) \subset (E, \phi) $. Corollary \ref{cor_P_orb_equals_P} implies that this definition is equivalent to Definition \ref{stability_equivariant_sheaves}.
\end{rmk}
Next we want to compare the (semi)stabilities of $ G $-equivariant coherent sheaves $ (E, \phi) $ on $ X $ and the (semi)stabilities of their underlying sheaves $ E $ on $ X. $
\begin{prop}\label{prop_Gieseker_relation}
	Let $ \cX = [X/G] $ in Situation \ref{situation_finite_group}. Let $ \cE = (E, \phi) $ be a nonzero pure sheaf on $ \cX $. Then $ \cE $ is $ \cH $-semistable if and only if $ E $ is $ H $-semistable. If $ E $ is $ H $-stable, then $ \cE $ is $ \cH $-stable.
\end{prop}
\begin{proof}
	Since $ \cE $ is a nonzero pure sheaf, by Remark 3.3 in \cite{nironi2009moduli}, $ E $ is also a nonzero pure sheaf. Suppose $ \cE $ is not (semi)stable. Then there is a proper nonzero subsheaf $ \cF \subset \cE $ such that
	\begin{equation*}
		\orbp(\cF)\ (\geq) \ \orbp(\cE)
	\end{equation*}
	in $ \QQ[z] $. Since the morphism $ p: X \to \cX $ is smooth, the pullback functor $ p^*: \Coh(\cX) \to \Coh(X) $ is exact, which implies $ F = p^*\cF \subset E $ is also a proper nonzero subsheaf. By Corollary \ref{cor_P_orb_equals_P}, we have
	$$ p(F) \ (\geq) \ p(E). $$
	Therefore $ E $ is not (semi)stable. Hence we have proved two implications:
	\begin{equation*}
		\text{$ E $ is (semi)stable} \quad \Rightarrow \quad \text{$ \cE $ is (semi)stable}.
	\end{equation*}
	Suppose $ E $ is not semistable. We use an argument similar to that in the proof of \cite[Lemma 2.7]{biswas1997parabolic} to show that $ \cE $ is not semistable. Let $ E_1 $ be the maximal destabilizing subsheaf of $ E $, so in particular we have
	$$ p(E_1) > p(E). $$ 
	Since $ H $ is $ G $-invariant, i.e., $ g^*H \cong H $ for all $ g $ in $ G, $ the Hilbert polynomial $ P(g^*F) $ is the same as $ P(F) $ for all $ g $ in $ G. $ Therefore, for each $ g $ in $ G, $ $ g^*E_1 $ is the maximal destabilizing subsheaf of $ g^*E \cong E $, and the uniqueness of Harder-Narasimhan filtrations implies that
	$$ g^*E_1 \cong E_1. $$
	So $ E_1 $ is a $ G $-invariant subsheaf of $ E. $ Then we have a $ G $-equivariant subsheaf $ \cE_1 = (E_1, \psi) $ of $ \cE = (E, \phi) $ where $ \psi $ is the restriction of the $ G $-equivariant structure of $ E $ on $ E_1. $ Since
	\begin{equation*}
		\orbp(\cE_1) = p(E_1) > p(E) = \orbp(\cE),
	\end{equation*}
	$ \cE_1 $ destabilizes $ \cE $. In other words, we have proved the implication
	\begin{equation*}
		\text{$ \cE $ is semistable} \quad \Rightarrow \quad \text{$ E $ is semistable}.
	\end{equation*}
\end{proof}

\begin{rmk}
	Let's summarize the results in Proposition \ref{prop_Gieseker_relation}:
	\vspace{5pt}
	\begin{equation*}\label{diag_Gieseker_relation}
		\begin{tikzcd}[column sep=2.4em,row sep=3.6em]
			\text{$ p^*\cE $ is stable} \arrow[r, Rightarrow, start anchor = {[xshift = 3.3ex]}, end anchor = {[xshift = -3.3ex]}] \arrow[d, Rightarrow, start anchor = {[yshift = -2.5ex]}, end anchor = {[yshift = 2.5ex]}] & \text{$ \cE $ is stable} \arrow[d, Rightarrow, start anchor = {[yshift = -2.5ex]}, end anchor = {[yshift = 2.5ex]}] \\
			\text{$ p^*\cE $ is semistable} \arrow[r, Leftrightarrow, start anchor = {[xshift = 1ex]}, end anchor = {[xshift = -1ex]}] & \text{$ \cE $ is semistable}
		\end{tikzcd}\vspace{10pt}
	\end{equation*}
	In general, the stability of a sheaf $ \cE = (E, \phi) $ on $ \cX $ is weaker than the stability of the sheaf $ E $ on $ X $. For example, consider $ G = \{1,\sigma\} \cong \ZZ/2\ZZ $ acting on $ X $ nontrivially. Let $ \cE = \cO_p $ be the structure sheaf of a non-orbifold point $ p \in \cX $, so $ E = p^*\cE = \cO_{Gx} $ where $ Gx = \{x, \sigma x\} $ is a free orbit on $ X $, i.e., $ E = \cO_x \oplus \cO_{\sigma x} $. Note that the $ G $-equivariant structure $ \phi $ on $ E $ is from the regular representation of $ G $. Therefore, $ \cE $ is stable, but $ E $ is strictly semistable.
\end{rmk}

\section{Equivariant moduli theory on $ K3 $ surfaces}\label{sec_moduli_sheaves_K3/G}
\vspace{5pt}
In this section we study equivariant moduli spaces of sheaves on $ K3 $ surfaces with symplectic automorphisms. We first introduce a notation.
\begin{notation}
	Let $ X $ be a $ K3 $ surface. Let $ G $ be a finite subgroup of the symplectic automorphisms of $ X $. Then the quotient stack $ \cX = [X/G] $ is a connected smooth projective stack over $ \CC. $ The quotient stack $ \cX $ is called a model of $ [K3/G] $ and is denoted by $ \cX = [K3/G] $.
\end{notation}
Take a model $ \cX = [K3/G] $ in Situation \ref{situation_finite_group}. Recall that we have a line bundle $ \cH = (H,\phi) $ on $ \cX $ which descends to an ample line bundle on the surface $ X/G. $
We first record a result for later use:
$$ \chi(\cX, \cO_\cX) = 2. $$ 
This follows from a direct computation:
\begin{enumerate}[font=\normalfont,leftmargin=*]
	\item $ h^0(\cX, \cO_\cX) = h^0(X, \cO_X) = 1 $ since $ H^0(X,\cO_X) = H^0(X,\cO_X)^G $.
	\item $ h^1(\cX, \cO_\cX) = 0 $ since $ H^1(X,\cO_X) = 0 $ and hence $ H^1(X,\cO_X)^G = 0 $.
	\item $ h^2(\cX, \cO_\cX) = h^0(\cX, \cO_\cX) = 1 $ by Serre duality for projective stacks.
\end{enumerate}
An immediate consequence is that for every line bundle $ \cL $ on $ \cX $,
\begin{equation*}
\chi(\cL,\cL) = \chi(\cX, \cL^\vee \otimes \cL) = \chi(\cX, \cO_\cX) = 2,
\end{equation*}
and hence $ \inprod{\orbv(\cL)^2}_{I\cX} = 2. $

Let's set up a few notations for later use.
\begin{notation}\label{notation_K3/G}
	Let $ \cX = [K3/G] $ in Situation \ref{situation_finite_group}. Label the conjugacy classes in $ G $ by 
	$$ [g_0], [g_1], \dots, [g_l] $$
	where $ [g_0] = [1] $ is trivial. We choose a representative element $ g_i $ in each class once and for all. For each $ g \in G, $ denote its order by $ n_g. $ 
	For every point $ x $ in the set $ X^g $ fixed by an element $ g \in G $, the linear automorphism 
	$$ dg_x: T_xX \to T_xX $$ 
	has two eigenvalues $ \lambda_{g,x} $ and $ \lambda_{g,x}^{-1} $ which are primitive $ n_g $-th roots of unity.
	Each element $ 1 \neq g \in G $ determines a finite non-empty fixed point set $ X^g $ of size
	$$ s_g = |X^g|. $$ 
	For each $ 1 \leq i \leq l, $ we have
	\begin{equation*}
	X^{g_i} = \{x_{i,1},x_{i,2},\dots,x_{i,s_{g_i}}\}.
	\end{equation*}
	Note that $ X^{g_1}, \cdots,  X^{g_l} $ are not distinct when $ l > 1 $. Indeed, for any triple $ i,j,k, $
	\begin{equation*}
	g_j = g_i^k \quad \Rightarrow \quad X^{g_i} \subset X^{g_j}. 
	\end{equation*}
	The union of these $ X^{g_i} $ consist of $ r $ non-free $ G $-orbits on $ X $, and they correspond to the orbifold points on $ \cX, $ and the singular points on $ Y $: 
	\begin{equation*}
		Gx_1, \dots, Gx_r \ \text{on} \ X \quad \leftrightarrow \quad \{p_1, \dots, p_r\} \ \text{on} \ \cX \quad \leftrightarrow \quad y_1, \dots, y_r \ \text{on} \ Y.
	\end{equation*}
	Here a point $ x_k $ is chosen in each $ G $-orbit on $ X $ once and for all. Take an element $ g_i \neq 1. $ The finite set $ X^{g_i} $ is invariant under the centralizer $ Z_{g_i} $ of $ g_i $. The number of orbits on $ X^{g_i} $ under the $ Z_{g_i} $-action is
	\begin{equation*}
	m_i = |X^{g_i}/Z_{g_i}|.
	\end{equation*} 
	Label these $ Z_{g_i} $-orbits by
	\begin{equation*}
	O_{i,1}, O_{i,2}, \dots, O_{i,m_i}.
	\end{equation*}
	If $ G $ is abelian, then each $ O_{ij} = Z_{g_i} \cdot x_{k_{ij}} $ for a unique index $ 1 \leq k_{ij} \leq r. $ Denote the total number of these $ Z_{g_i} $-orbits $ O_{ij} $ on $ X^{g_i} $ by
	$$ m = m_1 + \cdots + m_l. $$ 
	For each point $ x_{ij} \in X^{g_i} $, define
	\begin{equation*}
	G_{ij} = \Stab_{Z_{g_i}}(x_{ij}) = Z_{g_i} \cap \Stab_G(x_{ij}).
	\end{equation*}
	Note that each $ G_{ij} $ is a subgroup of $ Z_{g_i} $ with order
	$$ |G_{ij}| = |Z_{g_i}|/|O_{ij}|. $$
	The inertia stack of $ \cX $ is given by $ I\cX = \cX \coprod I_t \cX, $ where
	\begin{equation*}
	I_t \cX = \coprod_{i=1}^l \{g_i\} \times \left[ X^{g_i}/Z_{g_i} \right] = \coprod_{i=1}^l \coprod_{j=1}^{m_i} \{g_i,x_{ij}\} \times BG_{ij}.
	\end{equation*}
\end{notation}
\begin{ex}
	If $ G = \mu_2 = \{1, \sigma\}$ is generated by a Nikulin involution i.e., a symplectic involution $ \sigma, $ then there are eight fixed points of $ \sigma $,
	\begin{equation*}
	x_1, \cdots, x_8,
	\end{equation*}
	which correspond to eight $ A_1 $ singularities on $ Y. $ So $ l = 1, m_1 = 8, $ and each orbit $ O_{1,j} = \{x_j\} $ with stabilizer $ G_{1,j} = \mu_2. $ The twisted sectors are
	\begin{equation*}
	I_t \cX = \coprod_{j=1}^{8} \{\sigma,x_j\} \times B\mu_2.
	\end{equation*}
\end{ex}
\begin{ex}
	Let $ G = \mu_4 = \{1, g_1, g_2, g_3\} $ with $ g_i = g_1^i $ for $ i = 1,2,3. $ The fixed point set of the whole group $ G $ is 
	$$ X^G = \{x_1,x_2,x_3,x_4\}. $$ There are four more points 
	\begin{equation*}
	x_5, x_6, x_5', x_6'
	\end{equation*}
	fixed by $ g_2. $ The elements $ g_1 $ and $ g_3 $ swap $ x_j $ and $ x_j' $ for $ j = 5,6. $ 
	So the fixed point sets of $ g_1, g_2, $ and $ g_3 $ are:
	\begin{align*}
	X^{g_1} & = X^{g_3} = \{x_1, x_2, x_3, x_4\}. \\
	X^{g_2} & = \{x_1, x_2, x_3, x_4, x_5, x_6, x_5', x_6'\}.
	\end{align*}
	Under the action of $ G, $ the two sets $ X^{g_1} $ and $ X^{g_3} $ have $ m_1 = m_3 = 4 $ orbits:
	\begin{align*}
	\{x_1\}, \{x_2\}, \{x_3\}, \{x_4\},
	\end{align*}
	and $ X^{g_2} $ has $ m_2 = 2 $ orbits:
	\begin{align*}
	\{x_5, x_5'\}, \{x_6, x_6'\}.
	\end{align*} 
	These are six non-free orbits on $ X $ which descend to singularities
	$ 4A_3 + 2A_1$
	on $ Y $. There are fourteen twisted sectors:
	\begin{equation*}
	I_t \cX = \left(\{g_1,g_2,g_3\} \times X^G \times B\mu_4\right) \coprod \left(\{g_2\} \times \{x_5,x_6\} \times B\mu_2\right).
	\end{equation*}
\end{ex}

\vspace{1pt}

\subsection{Orbifold HRR formula}\label{sec_HRR_K3/G}
In this section we work out an explicit orbifold HRR formula for $ \cX = [K3/G] $ in Notation \ref{notation_K3/G}. 
Recall that the Chern character map on the $ K3 $ surface $ X $ gives a ring isomorphism 
\begin{equation*}
\ch: N(X) \to R(X) = R^0(X) \oplus R^1(X) \oplus R^2(X) \cong \ZZ \oplus \Pic(X) \oplus \ZZ,
\end{equation*}
where $ \Pic(X) $ is a free $ \ZZ $-module with rank $ \rho(X) $ which depends on the surface $ X. $ 
The complex numerical Chow ring is then identified as
$$ R(I\cX)_\CC = R(\cX) \oplus R(I_t\cX)_\CC \cong R(X)^G \oplus \CC^m \cong \ZZ \oplus \Pic(X)^G \oplus \ZZ \oplus \CC^m. $$
Here $ R(\cX) $ is not tensored by $ \CC $, $ R(X)^G $ denotes the $ G $-invariant part of $ R(X) $ under the induced action of $ G $ on $ R(X) $, and $ \Pic(X)^G $ is generated by isomorphism classes of $ G $-invariant line bundles.

Now let's compute the orbifold Chern character map
\begin{equation*}
\orbch: N(\cX) \to R(X)^G \oplus \CC^m.
\end{equation*}
Take a vector bundle $ \cV = (V, \phi) $ on $ \cX $. For each point $ x_{k_{ij}}, $ the equivariant structure $ \phi $ on $ V $ restricts to a representation
\begin{equation*}
\phi_{ij}: G_{ij} \to \GL(V_{ij})
\end{equation*}
of the group $ G_{ij} $ on the fiber $ V_{ij} $ of $ V $. By formula (\ref{eq_orbch_concrete}), we have
\begin{equation*}
\orbch(\cV) = \left(\ch(V) , \left(\chi_{\phi_{ij}}(g_i)\right) \right)
\end{equation*}
where $ 1 \leq i \leq l, 1 \leq j \leq m_i. $
For a sheaf $ \cE = (E, \phi) $ on $ \cX, $ we resolve it by any finite locally free resolution $ \cE_{\boldsymbol{\cdot}} \to \cE \to 0 $ with each $ \cE_k = (E_k, \phi_k) $, and hence we have
\begin{equation*}
\orbch(\cE) = \sum_{k} (-1)^k \orbch(\cE_k).
\end{equation*}

Since $ G $ is a finite group, the Todd class $ \td_\cX = \td_X. $ By formula (\ref{eq_orbtd_finite_G}), we have
\begin{equation*}
\orbtd_\cX = \frac{\td_{I\cX}}{e^\rho_{I\cX}} = \left(\td_X,\ \bigoplus_{i,j} \frac{1}{e^{\rho_{g_i}}(T_{ij} X)}\right) = \left((1,0,2),\ \bigoplus_{i,j} \frac{1}{2-2\Real(\lambda_{ij})}\right),
\end{equation*}
where $ T_{ij} X $ is the tangent space of $ X $ at $ x_{ij} $, $ \lambda_{ij} $ is either one of the two conjugate eigenvalues of $ g_i $ at the point $ x_{ij} $,
\begin{equation*}
\td_{I\cX} = \left((1,0,2), 1, \dots, 1 \right), \quad \text{and} \quad e^\rho_{I\cX} = \left((1,0,0),\  \bigoplus_{i,j} (2-2\Real(\lambda_{ij}))\right).
\end{equation*}
Note that $ \td_{I\cX} = \td_{I\cX}^\vee $ since $ \td_X = \td_X^\vee. $

The orbifold Mukai vector map
\begin{equation*}
\orbv: N(\cX) \to R(X)^G \oplus \CC^m
\end{equation*}
is given by 
\begin{equation*}
\orbv(\cE) = \orbch(\cE) \sqrt{\td_{I\cX}} = \orbch(\cE) \left((1,0,1),1, \dots, 1\right) = \left(v(E), \left(\chi_{\phi_{ij}}(g_i)\right) \right)
\end{equation*}
for a sheaf $ \cE = (E, \phi) $ on $ \cX $, where 
$$ v(E) = (\rk(E), \ch_1(E), \rk(E) + \ch_2(E)) $$ 
in $ R(X)^G $ is the Mukai vector of the sheaf $ E $ on $ X $.

The orbifold Mukai pairing 
\begin{equation*}
\inprod{\cdot \ {,}\ \cdot}_{I\cX}: R(I\cX)_\CC \times R(I\cX)_\CC \to \CC
\end{equation*}
is then given by
\begin{align}\label{eq_orbv_pairing_K3/G}
\inprod{\orbv,\orbw}_{I\cX} & = \frac{\inprod{v,w}_X}{|G|} + \frac{1}{2} \sum_{i=1}^{l} \sum_{j=1}^{m_i} \frac{\overline{v}_{ij} w_{ij}}{|G_{ij}|(1-\Real(\lambda_{ij}))}
\end{align}
for all $ \orbv = (v,(v_{ij})) $ and $ \orbw = (w,(w_{ij})) $ in $ R(I\cX)_\CC \cong R(X)^G \oplus \CC^m $, where 
$$ \inprod{v,w}_X = \int_X v^\vee w = v_0 w_2 - v_1 w_1 + v_2 w_0 $$
is the integral Mukai pairing of $ v $ and $ w $ in $ R(X). $

For every pair of sheaves $ \cE = (E, \phi) $ and $ \cF = (F, \psi) $ on $ \cX, $ the orbifold HRR formula (\ref{eq_HRR2}) reads
\begin{equation*}
\chi(\cE,\cF) = \inprod{\orbv(\cE),\orbv(\cF)}_{I\cX} = \frac{\inprod{v(E),v(F)}_X}{|G|} + \frac{1}{2} \sum_{i=1}^l \sum_{j=1}^{m_i} \frac{\chi_{\phi_{ij}}(g_i^{-1}) \chi_{\psi_{ij}}(g_i)}{|G_{ij}|(1-\Real(\lambda_{ij}))};
\end{equation*}
in particular, if $ \cE = \cF, $ then
\begin{equation}\label{eq_orbv_squared}
\chi(\cE,\cE) = \inprod{\orbv(\cE)^2}_{I\cX} = \frac{\inprod{v(E)^2}_X}{|G|} + \frac{1}{2} \sum_{i=1}^l \sum_{j=1}^{m_i} \frac{|\chi_{\phi_{ij}}(g_i)|^2}{|G_{ij}|(1-\Real(\lambda_{ij}))}.
\end{equation}
\begin{rmk}
	Consider any coherent sheaf $ \cE $ on $ \cX. $ Since $ \chi(\cE,\cE) $ is an integer, the double summation in (\ref{eq_orbv_squared}) must be a rational number although each term is irrational in general. Moreover, $ \chi(\cE,\cE) $ is always an even integer as we will see shortly.
\end{rmk}
We first generalize Proposition \ref{prop_K3_smoothness_and_dim} for $ K3 $ surfaces to $ [K3/G]. $
\begin{prop}\label{prop_K3/G_smoothness_dim}
	Let $ \cX = [K3/G] $ in Situation \ref{situation_finite_group}. Let $ x $ be an element in $ N(\cX) $ with orbifold Mukai vector $ \orbv $ in $ R(I\cX)_\CC $. Then $ M^s(\cX,x) $ is either empty or a smooth quasi-projective scheme with
	\begin{equation*}
	\dim M^s(\cX,x) = 2 - \inprod{\orbv^2}_{I\cX}.
	\end{equation*}
	Denote the Mukai vector of $ y = p^N x $ in $ N(X) $ by $ v = (r, c_1, s) $ in $ R(X)^G $ with $ d = \deg(y). $ Suppose either of the following conditions is satisfied:
	\begin{enumerate}[font=\normalfont,leftmargin=2em]
		\item $ \gcd(r, d, s) = 1 $.
		\item $ y $ is primitive and $ H $ is $ y $-generic.
	\end{enumerate}
	Then $ M(\cX,x) = M^s(\cX,x). $ If $ M(\cX,x) $ is not empty, then it is a smooth projective scheme of dimension $ 2 - \inprod{\orbv^2}_{I\cX}. $
\end{prop}
\begin{proof}
	Suppose $ M_\cH^s(\cX,x) $ is not empty. Then it is a quasi-projective scheme by construction. Take a stable sheaf $ \cE $ on $ \cX $ with numerical class $ x. $ Then we have $ \End(\cE) \cong \CC. $
	In particular, the inclusion map
	$$ H^0(\cX, \cO_\cX) \cong \CC \to \Ext^0(\cE, \cE) $$
	is an isomorphism of $ \CC $-vector spaces.
	Since the canonical sheaf $ \omega_\cX \cong \cO_\cX, $ Serre duality for smooth projective stacks gives an isomorphism
	\begin{equation*}
	\Ext^2(\cE, \cE) \xrightarrow{\sim} H^2(\cX, \cO_\cX).
	\end{equation*}
	The Picard scheme $ \Pic(\cX) = \Pic^G(X) $ is a group scheme, which is always smooth in characteristic zero. Therefore, by Proposition \ref{prop_smoothness_criterion}, $ M^s(\cX,x) $ is smooth at the closed point $ t $ corresponding to $ \cE. $ So we can compute its dimension at the tangent space $ T_tM^s(\cX,x) \cong \Ext^1(\cE,\cE) $. Recall the orbifold Euler pairing
	\begin{align*}
	\chi(\cE, \cE) = \sum_{i=0}^2 \dim \Ext^i(\cE, \cE) = 2 - \dim \Ext^1(\cE, \cE).
	\end{align*}
	By the orbifold HRR formula (\ref{eq_HRR2}), we have
	\begin{equation*}
	\dim \Ext^1(\cE, \cE) = 2 - \chi(\cE, \cE) = 2 - \inprod{\orbv(\cE)^2}_{I\cX}.
	\end{equation*}
	Now suppose either of the two conditions in Proposition \ref{prop_K3/G_smoothness_dim} is satisfied. If $ M(\cX,x) $ is empty, there is nothing to prove. Suppose it's not empty. Take any semistable sheaf $ \cF = (F, \psi) \in M(\cX,x) $. By Proposition \ref{prop_Gieseker_relation}, the coherent sheaf $ F $ on $ X $ is also semistable. By Proposition \ref{prop_K3_smoothness_and_dim}, $ F $ is also stable, and hence $ \cE $ is stable as well again by Proposition \ref{prop_Gieseker_relation}. The last statement follows since $ M(\cX,x) $ is a projective scheme by construction.
\end{proof}
\begin{rmk}
	The results in Proposition \ref{prop_K3/G_smoothness_dim} do not say that the moduli space $ M^s(\cX,x) $ is irreducible, which will be proved later in Theorem \ref{thm_K3/G} under a slightly stronger assumption, i.e., $ r >0, $ and either $ d > 0 $ or $ \gcd(r,d) = 1. $
\end{rmk}
Let's apply Proposition \ref{prop_K3/G_smoothness_dim} to a few examples. 
\begin{ex}
	Let $ \cL $ be a line bundle on $ \cX $ with numerical class $ \gamma(\cL) $ in $ N(\cX). $ Since $ \inprod{\orbv(\cL)^2} = 2, $ Proposition \ref{prop_K3/G_smoothness_dim} implies
	\begin{equation*}
	\dim  M(\cX, \gamma(\cL)) = 0,
	\end{equation*}
	since there are no strictly semistable sheaves of rank one. Indeed, each moduli space $ M(\cX, \gamma(\cL)) $ is a point in the smooth Picard scheme $ \Pic(\cX) $ of $ \cX $.
\end{ex}
\begin{ex}
	Let $ \cO_p $ be the structure sheaf of a generic point $ p $ on $ \cX $ corresponding to a free orbit $ Gx $ on $ X. $ Then $ \cO_p $ is a stable sheaf on $ \cX $ with an orbifold Mukai vector
	\begin{equation*}
	\orbv(\cO_p) = (0,0,|G|,0, \dots, 0)
	\end{equation*}
	in $ R(I\cX)_\CC. $ Therefore, $ \inprod{\orbv(\cO_p)^2} = 0, $ which implies
	\begin{equation*}
	\dim  M(\cX, \gamma(\cO_p)) = 2.
	\end{equation*}
	Indeed, we can identify
	$$ M(\cX, \gamma(\cO_p)) \cong \GHilb X, $$ the $ G $-Hilbert scheme of free orbits on $ X. $ Hence $ M = M(\cX, \gamma(\cO_p)) $ is a $ K3 $ surface which is the minimal resolution of the surface $ X/G. $ Moreover, $ M $ is a fine moduli space of stable sheaves on $ \cX $.
\end{ex}

\begin{ex}
	Consider the tangent bundle $ T\cX =(TX, \tau) $ on $ \cX. $ It's known that the tangent bundle $ TX $ is $ \mu $-stable (and hence stable) with respect to any polarization $ H $ on $ X $ (see for example, \cite[Proposition 4.5]{huybrechts2016lectures}). So the $ G $-equivariant tangent bundle $ T\cX $ is also stable. Since $ c_2(TX) = 24, $ the Mukai vector of $ TX $ is 
	$$ v(TX) = (2, 0, -22) $$ 
	in $ R(X) $ with $ \inprod{v(TX)^2}_X = -88, $ which implies
	\begin{equation*}
	\dim M^s(X,\gamma(TX)) = 90.
	\end{equation*}
	However, 
	$$ M(X,\gamma(TX)) \neq M^s(X,\gamma(TX)) $$ 
	because there exist strictly semistable sheaves of the form
	\begin{equation*}
	E = I_Z \oplus I_Z \quad \text{with} \quad v(E) = v(TX),
	\end{equation*}
	where $ Z $ is a zero-dimensional subscheme of $ X $ of length $ 12. $ 
	The twisted Mukai vectors of $ TX $ are 
	$$ \chi_{\tau_{ij}}(g_i) = \lambda_{ij} + \lambda_{ij}^{-1} = 2 \Real(\lambda_{ij}) $$ 
	in $ \RR $ for all $ 1 \leq i \leq l, 1 \leq j \leq m_i. $ Therefore, the orbifold Mukai vector of $ T\cX $ is
	\begin{equation*}
	\orbv(T\cX) = (2,0,-22,2 \Real(\lambda_{ij}))
	\end{equation*}
	in $ R(I\cX)_\RR, $ and
	\begin{equation*}
	\inprod{\orbv(T\cX)^2}_{I\cX} = -\frac{88}{|G|} + \sum_{i=1}^l \sum_{j=1}^{m_i} \frac{2r_{ij}^2}{|G_{ij}|(1-r_{ij})}
	\end{equation*}
	where each $ r_{ij} = \Real(\lambda_{ij}). $ Suppose $ G = \mu_2. $ Then we have
	\begin{equation*}
	\orbv(T\cX) = (2,0,-22,-2,\dots,-2)
	\end{equation*}
	in $ R(I\cX) = R(X)^G \oplus \ZZ^8 $ with 
	$$ \inprod{\orbv(T\cX)^2}_{I\cX} = -40. $$ 
	Therefore,
	\begin{equation*}
	\dim  M^s(\cX, \gamma(T\cX)) = 42.
	\end{equation*}
	We claim 
	$$ M(\cX, \gamma(T\cX)) \neq M^s(\cX, \gamma(T\cX)). $$ 
	This is expected since the Mukai vector $ v(TX) $ is not primitive in $ R(X) $. Indeed, there are strictly semistable sheaves on $ \cX $ of the form
	\begin{equation*}
	\cE = I_{\cZ_1} \oplus I_{\cZ_2} \quad \text{with} \quad \orbv(\cE) = \orbv(T\cX),
	\end{equation*}
	where $ \cZ_1 = [Z_1/\mu_2] $ is an object in the equivariant Hilbert scheme $ \Hilb_\cX(6, 0) $, and $ \cZ_2  = [Z_2/\mu_2] $ is another one in $ \Hilb_\cX(10, -e_1-\cdots-e_8) $ such that the lengths of $ Z_1 $ and $ Z_2 $ are both $ 12. $ Refer to Section \ref{sec_Hilb} for the notations for equivariant Hilbert schemes of points on $ \cX $.
\end{ex}
The following result is immediate from the orbifold HRR formula for $ [K3/G] $.
\begin{prop}
	Let $ \cX = [K3/G] $ in Situation \ref{situation_finite_group} with Notation \ref{notation_K3/G}. Then
	\begin{equation}\label{eq_K3/G_identity}
	\frac{1}{|G|} + \frac{1}{4} \sum_{i=1}^l \sum_{j=1}^{m_i} \frac{1}{|G_{ij}|(1-\Real(\lambda_{ij}))} = 1.
	\end{equation}
\end{prop}
\begin{proof}
	Consider $ \cE = \cO_\cX $ in (\ref{eq_orbv_squared}). The result follows from the fact $ \chi(\cX,\cO_\cX) = 2. $ 
\end{proof}

A result of Nikulin \cite{nikulin1979finite} says that if a $ K3 $ surface admits a symplectic automorphism of finite order $ n $, then $ n \leq 8. $ In \cite[Section 1]{mukai1988finite}, Mukai proved that the number of points $ f_n $ fixed by a symplectic automorphism only depends on its order $ n $, and computed all such numbers in Table \ref{table_number_fixed_points} on page \pageref{table_number_fixed_points}.

We now use equation (\ref{eq_K3/G_identity}) to reproduce these numbers $ f_n = |X^G|, $ where $ G $ is the cyclic group generated by a symplectic automorphism of order $ n. $

\begin{cor}\label{cor_table_fixed_points}
	Let $ \cX = [K3/G] $ in Situation \ref{situation_finite_group} where $ G \cong \ZZ/n\ZZ $. If $ n $ is a prime $ p $, then the number of fixed points of $ G $ is 
	$$ |X^G| = \frac{24}{p+1}. $$ 
	If $ n = 4, 6, $ and $ 8, $ then $ |X^G| = 4, 2 $ and $ 2 $ respectively. Thus we reproduce Table \ref{table_number_fixed_points} on page \pageref{table_number_fixed_points}. 
\end{cor}
\begin{proof}
	Let $ g $ be a generator of $ G, $ i.e., 
	$$ G = \{1, g, g^2, \dots, g^{n-1}\}. $$
	Then $ X^G = X^g $. 
	For each $ 1 \leq i \leq n-1, $ let $ s_i $ denote the number of fixed points of $ g^i. $ 
	Suppose $ n $ is a prime. Then every nontrivial $ g^i $ has the same fixed point set $ X^G $, and hence $ s_i = s_1 = |X^G| $ for all $ 1 \leq i \leq n-1. $ Since $ G $ is abelian, each centralizer $ Z_{g^i} = G $, so each stack 
	$$ [X^{g^i}/Z_{g^i}] = [X^G/G]. $$ 
	Therefore, in equation (\ref{eq_K3/G_identity}), the number of orbits $ m_i = |X^G| $ for all $ 1 \leq i \leq n-1 $, and each stabilizer $ G_{ij} $ has size $ |G_{ij}| = n. $
	Note that for each $ 1 \leq k \leq n-1, $
	\begin{equation*}
	\lambda_{k} = \exp(2k\pi i/n)
	\end{equation*}
	is an eigenvalue of $ g^k. $ 
	Equation (\ref{eq_K3/G_identity}) then becomes
	\begin{equation}\label{eq_K3/G_identity_prime_order}
	\frac{1}{n} + \frac{|X^G|}{4n} \sum_{k=1}^{n-1} \frac{1}{1-\cos(2k\pi /n)} = 1.
	\end{equation}
	Recall a trigonometric identity
	\begin{equation}\label{eq_trig_identity}
	\sum_{k=1}^{n-1} \frac{1}{1-\cos(2k\pi /n)} = \frac{n^2-1}{6}.
	\end{equation}
	Plugging (\ref{eq_trig_identity}) into (\ref{eq_K3/G_identity_prime_order}) yields
	\begin{equation*}
	|X^G| = \frac{24}{n+1}.
	\end{equation*}
	This gives $ f_2,f_3,f_5,f_7 $ in Table \ref{table_number_fixed_points}. The remaining three cases can be computed directly. For a non-prime $ n, $ we use the same argument, but now for each $ 1 \leq i \leq n-1, $ we have
	\begin{equation*}
	\sum_{j=1}^{m_i} \frac{1}{|G_{ij}|} = \sum_{j=1}^{m_i} \frac{|O_{ij}|}{|G|} = \frac{s_i}{n}.
	\end{equation*}
	Therefore, equation (\ref{eq_K3/G_identity}) becomes
	\begin{equation}\label{eq_K3/G_identity_non_prime_order}
	\frac{1}{n} + \frac{1}{4n} \sum_{k=1}^{n-1} \frac{s_k}{1-\cos(2k\pi /n)} = 1.
	\end{equation}
	Note that for each $ 1 \leq k \leq n-1, $ the element $ g^k $ has order 
	$$ n_k = \frac{n}{\gcd(n,k)}, $$ 
	so the number of fixed points $ s_k = f_{n_k}$ when $ n_k $ is a prime. Now applying (\ref{eq_K3/G_identity_non_prime_order}) to $ n = 4,$ we get
	\begin{align*}
	\frac{1}{4} + \frac{1}{16} \left(|X^G| + \frac{f_2}{2} + |X^G|\right) = 1,
	\end{align*}
	which gives $ |X^G| = 4. $ This implies all groups $ G \cong \ZZ/4\ZZ $ acting faithfully and symplectically on $ X $ have the same number of fixed points, and gives $ f_4 = 4 $ in Table \ref{table_number_fixed_points}.
	Carrying out the same procedure for $ n = 6 $ and $ 8, $ we have
	\begin{align*}
	& \frac{1}{6} + \frac{1}{24} \left(\frac{f_6}{1/2} + \frac{f_3}{3/2} + \frac{f_2}{2} + \frac{f_3}{3/2} + \frac{f_6}{1/2}\right) = 1, \ \text{and}\\
	& \frac{1}{8} + \frac{1}{32} \left(\frac{f_8}{1-1/\sqrt{2}} + f_4 + \frac{f_8}{1+1/\sqrt{2}} + \frac{f_2}{2} + \frac{f_8}{1+1/\sqrt{2}} + f_4 + \frac{f_8}{1-1/\sqrt{2}}\right) = 1,
	\end{align*}
	which gives $ f_6 = f_8 = 2, $ thus completing Table \ref{table_number_fixed_points} for $ 2 \leq n \leq 8 $. 
\end{proof}
\begin{rmk}
	An explicit formula for $ f_n $ was derived in \cite[Proposition 1.2]{mukai1988finite} as
	\begin{equation}\label{eq_number_fixed_points}
	f_n =  \frac{24}{n} \prod_{p | n} \left(1+\frac{1}{p}\right)^{-1}
	\end{equation}
	by the Lefschetz fixed point formula and the identity
	\begin{equation*}
	\sum_{k \in (\ZZ/n\ZZ)^\times} \frac{1}{(1-\lambda^k)(1-\lambda^{-k})} = \frac{n^2}{12} \prod_{p | n} \left(1-\frac{1}{p^2}\right),
	\end{equation*}
	where $ \lambda $ is a primitive $ n $th root of unity.
	It would be interesting to derive (\ref{eq_number_fixed_points}) directly using formula (\ref{eq_K3/G_identity_non_prime_order}). 
\end{rmk}


\subsection{Derived McKay correspondence}
Let $ p $ denote a generic point on $ \cX = [K3/G] $. Consider again the moduli space 
$$ M = M(\cX,\gamma(\cO_p)). $$
We have seen that $ M = \GHilb X, $ which is a fine moduli space. Therefore, the universal family on $ \cX \times M $ induces a Fourier-Mukai transform
\begin{equation*}
\Phi: \D(\cX) \xrightarrow{\sim} \D(M)
\end{equation*}
which is an equivalence between derived categories known as the derived McKay correspondence \cite{bridgeland2001mckay}. Passing to the numerical Grothendieck rings and the numerical Chow rings, we have a commutative diagram
\vspace{1pt}
\begin{equation*}\label{diag_Fourier_Mukai}
\begin{tikzcd}[column sep=2.5em,row sep=3em,every label/.append style={font=\normalsize}]
N(\cX) \arrow[r,"\Phi^N"{outer sep = 2pt}] \arrow[d, "\orbv"{left, outer sep = 2pt}]
& N(M) \arrow[d, "v"{right, outer sep = 2pt}] \\
R(I\cX)_\CC \arrow[r,"\Phi^R"{outer sep = 2pt}] & R(M)_\CC,
\end{tikzcd}\vspace{3pt}
\end{equation*}
where $ \Phi^N $ and $ \Phi^R $ are isomorphisms of rings and $ \CC $-algebra respectively, $ \orbv $ and $ v $ are injections of abelian groups, and each arrow is compatible with the canonical pairings on the sources and targets. For example, $ \Phi^N $ preserves orbifold Euler pairings:
\begin{equation*}
\chi(x, y) = \chi(\Phi^N(x),\Phi^N(y))
\end{equation*}
in $ \ZZ $ for all $ x $ and $ y $ in $ N(\cX), $ and $ \Phi^R $ preserves orbifold Mukai pairings:
\begin{equation}\label{eq_Fourier_Mukai}
\inprod{\orbv,\orbw}_{I\cX} = \inprod{\Phi^R(\orbv),\Phi^R(\orbw)}_M
\end{equation}
in $ \CC $ for all $ \orbv $ and $ \orbw $ in $ R(I\cX)_\CC. $ This implies the following:
\begin{prop}
	Let $ \cX = [K3/G] $ in Situation \ref{situation_finite_group}. Let $ \cE $ be a coherent sheaf on $ \cX $. Then the integer $ \inprod{\orbv(\cE)^2}_{I\cX} $ is even.
\end{prop}
\begin{proof}
	Consider the minimal resolution $ M \to X/G. $ We have 
	$$ \inprod{\orbv(\cE)^2}_{I\cX} = \inprod{\Phi^R(\orbv(\cE))^2}_M = \inprod{v(\Phi^N(\cE))^2}_M, $$ 
	which is even since the self-intersection on $ R^1(Y) $ is even for any $ K3 $ surface $ Y $.  
\end{proof}

\subsection{Bridgeland stability conditions}\label{sec_Bridgeland}
In this section we review Bridgeland stability conditions on triangulated categories, their constructions on surfaces, and induced stability conditions on a quotient stack $ \cX = [X/G] $ where $ X $ is a smooth projective variety and $ G $ is a finite group acting on $ X $.

We first review the Bridgeland stability conditions on triangulated categories introduced in \cite{bridgeland2007stability}. Fix a triangulated category $ \cD $. Let $ K(\cD) $ denote the Grothendieck group of $ \cD. $ Fix a finite rank lattice $ \Lambda $ and a surjective group homomorphism
\begin{equation*}
	v: K(\cD) \onto \Lambda.
\end{equation*}
For example, when $ \cD = \D(X) $ for a complex smooth projective variety $ X, $ we can choose $ \Lambda $ to be the numerical Grothendieck group $ N(X) $ of $ X. $ Let $ || \cdot || $ be a norm on $ \Lambda_\RR. $
\begin{defn}
	A \tb{stability condition} $ \sigma = (Z, \cP) $ on $ \cD $ consists of an additive map $ Z: K(\cD) \to \CC $ called the \tb{central charge}, and a slicing $ \cP $ of $ \cD $ satisfying the following conditions: 
	\begin{enumerate}[label=(\arabic*),font=\normalfont,leftmargin=2em]
		\item For each $ \phi $ in $ \RR, $ if a nonzero object $ E $ is in $ \cP(\phi), $ then
		\begin{equation*}
			Z(E) = m(E) e^{i\pi\phi}
		\end{equation*}
		for some $ m(E) > 0 $.
		\item (Support Property)
		\begin{equation*}
			\text{inf}\left\{\frac{|Z(E)|}{||v(E)||} : 0 \neq E \in \cP(\phi), \phi \in \RR \right\} > 0.
		\end{equation*}
	\end{enumerate}
\end{defn}
If there exists a stability condition on $ \cD $, then there are notions of stable objects and semistable objects in $ \cD. $
\begin{defn}
	Let $ \sigma = (Z, \cP) $ be a stability condition on $ \cD. $ For each $ \phi $ in $ \RR, $ the nonzero objects of $ \cP(\phi) $ are said to be \tb{$ \sigma $-semistable} of phase $ \phi $, and the simple objects of $ \cP(\phi) $ are said to be \tb{$ \sigma $-stable}.
\end{defn}

Alternatively, one can define stability conditions via a stability function on the heart $ \cA $ of a bounded $ t $-structure on $ \cD $, which is an abelian subcategory of $ \cD. $
\begin{defn}
	A \tb{stability function} on an abelian category $ \cA $ is an additive map $ Z: K(\cA) \to \CC $ such that for all nonzero objects $ E $ in $ \cA, $ we have
	\begin{equation*}
		Z(E) = m(E) e^{i\pi\phi(E)}
	\end{equation*}
	where $ m(E) > 0 $ and $ \phi(E) \in (0,1]. $
\end{defn}
\begin{defn}
	Let $ Z: K(\cA) \to \CC $ be a stability function on an abelian category $ \cA $. A nonzero objects $ E $ of $ \cA $ is said to be \tb{(semi)stable} with respect to $ Z $ if
	\begin{equation*}
		\phi(F) \ (\leq)\ \phi(E)
	\end{equation*} 
	for all nonzero proper subobjects $ F \subset E. $ A \tb{Harder-Narasimhan (HN) filtration} of a nonzero object $ E $ in $ \cA $ is a chain of subobjects
	\begin{equation*}
		0 = E_0 \subset E_1 \subset \cdots \subset E_n = E
	\end{equation*}
	such that each factor $ F_i = E_i/E_{i-1} $ is a semistable object in $ \cA $ and
	\begin{equation*}
		\phi(F_1) > \phi(F_2) > \cdots > \phi(F_n).
	\end{equation*}
	The stability function $ Z $ is said to have the \tb{HN property} if every nonzero object of $ \cA $ has an HN filtration.
\end{defn}
If $ \sigma = (Z, \cP) $ is a stability condition on $ \cD $, then the abelian subcategory $ \cA = \cP((0,1]) $ of $ \cD $ is the heart of the $ t $-structure $ \cP(>0). $ This yields two equivalent ways to impose stability conditions on $ \cD $.
\begin{prop}[{\cite[Proposition 5.3]{bridgeland2007stability}}]
	A stability condition $ \sigma = (Z, \cP) $ on $ \cD $ is the same as a pair $ (\cA, Z) $ where $ \cA $ is the heart of a bounded $ t $-structure on $ \cD $ and $ Z $ is a stability function on $ \cA $ with the HN property and the support property
	\begin{equation*}
		\mathrm{inf}\left\{\frac{|Z(E)|}{||v(E)||} : 0 \neq E \in \cA \ \mathrm{is \ semistable}\right\} > 0.
	\end{equation*}
\end{prop}
\begin{rmk}
	Let $ \sigma = (Z, \cP) $ be a stability condition on $ \cD. $ The category $ \cP(\phi) $ is abelian for all $ \phi \in \RR $. The support property ensures that each $ \cP(\phi) $ is of finite length, and hence every $ \sigma $-semistable object in $ \cP(\phi) $ has a finite Jordan-Hölder (JH) filtration whose factors are $ \sigma $-stable with the same phase $ \phi $. Two semistable objects in $ \cP(\phi) $ are said to be $ S $-equivalent if they have isomorphic JH factors.
\end{rmk}
Let $ X $ be a smooth projective variety. It's an open question whether there exist stability conditions on the derived category $ \D(X). $ The answer is affirmative when $ \dim X \leq 2. $
\begin{ex}
	Let $ X = \pt. $ Let $ \cA = \Vect_\CC $ be the abelian category of finite dimensional vector spaces over $ \CC. $ Then the additive map
	\begin{equation*}
	Z: K(\Vect_\CC) \to \ZZ, \quad V \mapsto i\dim(V)
	\end{equation*}
	is a stability function on $ \Vect_\CC $ and gives a stability condition on the derived category $ \D(\Vect_\CC). $ A nonzero vector space $ V $ is semistable of phase $ 1/2 $ and is stable if $ \dim(V) = 1. $
\end{ex}
\begin{ex}[{\cite[Example 5.4]{bridgeland2007stability}}]
	Let $ X $ be a smooth projective curve. Let $ \cA = \Coh(X). $ Then the additive map
	\begin{equation*}
	Z: K(X) \to \ZZ, \quad E \mapsto -\deg(E) + i\rk(E)
	\end{equation*}
	is a stability function on $ \Coh(X) $, and the pair $ (\Coh(X), Z) $ gives a stability condition on $ X $ which coincides with the slope stability.
\end{ex}

Now we recall the method of tilting to construct stability conditions on a surface.

\begin{constrn}[Bridgeland stability conditions on surfaces]\label{stability_conditions_on_surfaces}
	Let $ X $ be a smooth projective surface. Fix a divisor class $ \beta $ in $ R^1(X)_{\RR} $ and an ample divisor class $ \omega \in R^1(X)_{\RR}. $ Define two full additive subcategories of the abelian category $ \Coh(X) $ as follows:
	\begin{align*}
		\cT(\beta,\omega) & = \{E \in \Coh(X): T\ \text{has HN factors $ F $ with}\ \mu_\omega(F) > \beta \omega \}. \\
		\cF(\beta,\omega) & = \{E \in \Coh(X): F\ \text{has HN factors $ F $ with}\ \mu_\omega(F) \leq \beta \omega \}.
	\end{align*}
	Then $ (\cT(\beta,\omega), \cF(\beta,\omega)) $ is torsion pair. Note that $ \cT(\beta,\omega) $ contains torsion sheaves and each sheaf in $ \cF(\beta,\omega) $ is torsion-free. Define an abelian subcategory of $ \D(X) $ by
	\begin{equation*}
		\cA(\beta,\omega) = \{E \in \D(X): H^i(E) = 0 \ \text{for}\ i \neq -1, 0, H^{-1}(E) \in \cF(\beta,\omega), H^0(E) \in \cT(\beta,\omega) \}.
	\end{equation*}
	Then $ \cA(\beta,\omega) $ is the heart of a bounded $ t $-structure on $ \D(X). $
	Define an additive map
	\begin{equation*}
		Z_{\beta,\omega}: K(\cA(\beta,\omega)) \to \CC, \quad x \mapsto -\int_X e^{-\beta-i\omega} \ch(x).
	\end{equation*}
\end{constrn}
\begin{prop}[{\cite[Theorem 6.10]{macri2019lectures}}]\label{prop_stability_conditions_on_surfaces}
	Let $ X $ be a smooth projective surface. The pair $ \sigma_{\beta,\omega} = (Z_{\beta,\omega}, \cA(\beta,\omega)) $ in Construction \ref{stability_conditions_on_surfaces} gives a stability condition on $ \D(X). $
\end{prop}
\begin{rmk}
	Let $ X $ be a smooth projective variety. A stability condition $ \sigma = (Z, \cP) $ on $ \D(X) $ is called \emph{numerical} if the central charge $ Z: K(X) \to \CC $ factors through $ N(X). $ For example, the stability condition $ \sigma_{\beta,\omega} $ in Construction \ref{stability_conditions_on_surfaces} is numerical. Let $ \Stab(X) $ denote the space of all numerical stability conditions on $ X $. Suppose $ \Stab(X) $ is not empty. Then the map 
	$$ \Stab(X) \to \Hom(N(X), \CC), \quad (\cA, Z) \mapsto Z $$
	is a local homeomorphism which endows each connected component of $ \Stab(X) $ a structure of a complex manifold of dimension $ \rk(N(X)) $. When $ X $ is a $ K3 $ surface, there is a distinguished component which contains all stability conditions $ \sigma_{\beta,\omega} $ in Construction \ref{stability_conditions_on_surfaces}, and is denoted by $ \Stab^\dagger(X). $
\end{rmk}

Let $ X $ be a smooth projective variety. Let $ \sigma = (Z, \cP) $ be a numerical stability condition on $ \D(X) $. Fix an element $ x $ in $ N(X) $ and a phase $ \phi \in \RR $ such that $ Z(x) = m(x)e^{i\pi\phi} $. We denote by $ \cM_\sigma(X,x,\phi) $ (resp. $ \cM_\sigma^s(X,x,\phi) $) the moduli stack of $ \sigma $-semistable (resp. $ \sigma $-stable) objects in $ P(\phi) $ with numerical class $ x $. If we restrict the phase $ \phi $ to $ (0,2], $ then the element $ x $ uniquely determines $ \phi $, and in this case we will omit $ \phi $ from the notation  $ \cM_\sigma^{(s)}(X,x,\phi) $. 
The moduli stack $ \cM_\sigma(X,x) $ has a good moduli space $ M_\sigma(X,x) $ parametrizing $ S $-equivalent classes of $ \sigma $-semistable objects in $ \D(X) $ with numerical class $ x $, which becomes a coarse moduli space if $ \cM_\sigma(X,x) = \cM_\sigma^s(X,x) $.
\begin{thm}[{\cite[Theorem 21.24]{bayer2021stability}}]
	Let $ X $ be a smooth projective variety over $ \CC $. Let $ \sigma $ be a stability condition on $ \D(X) $, and let $ x $ be an element in $ N(X). $ Then $ \cM_\sigma(X,x) $ is an algebraic stack of finite type over $ \CC $ with a good moduli space $ M_\sigma(X,x) $ which is a proper algebraic space. If $ \cM_\sigma(X,x) = \cM_\sigma^s(X,x), $ then it is a $ \CC^* $-gerbe over its coarse moduli space $ M_\sigma(X,x) = M_\sigma^s(X,x). $
\end{thm}
Now we consider a $ K3 $ surface $ X. $ Under certain conditions, the moduli space $ M_\sigma(X,x) $ is an irreducible symplectic manifold equivalent to a Hilbert scheme of points on a $ K3 $ surface. Fix an element $ x $ in $ N(X). $ Then it determines a set of real codimension-one submanifolds with boundaries in the manifold $ \Stab(X) $ known as \tb{walls} as shown in \cite[Proposition 2.3]{bayer2014projectivity}.
\begin{defn}[{\cite[Definition 2.4]{bayer2014projectivity}}]
	Let $ X $ be a $ K3 $ surface. Let $ x $ be an element in $ N(X). $ A stability condition $ \sigma $ in $ \Stab(X) $ is said to be \tb{$ x $-generic} if it does not lie on any of the walls determined by $ x. $ 
\end{defn}
The following result is a generalization of Theorem \ref{thm_yoshioka} to Bridgeland moduli spaces of stable objects in $ \D(X) $.
\begin{thm}[{\cite[Theorem 1.1]{bottini2024stable}}]\label{thm_bottini}
	Let $ X $ be a $ K3 $ surface. Let $ x $ be an element in $ N(X) $ with a primitive Mukai vector $ v $ in $ R(X). $ Let $ \sigma $ be a stability condition in $ \Stab^\dagger(X) $ such that it is $ x $-generic. Then $ M_\sigma(X,x) =  M_\sigma^s(X,x), $ and 
	it is non-empty if and only if $ \inprod{v^2} \leq 2. $ In this case, $ M_\sigma(X,x) $ is an irreducible symplectic manifold of dimension $ 2 - \inprod{v^2} $ deformation equivalent to a Hilbert scheme of points on a $ K3 $ surface.
\end{thm}
Under some conditions, Bridgeland (semi)stable objects in $ \D(X) $ are the same as Gieseker (semi)stable sheaves in $ \Coh(X) $.
\begin{prop}\label{prop_Bridgeland_vs_Gieseker}
	Let $ (X, H) $ be a polarized $ K3 $ surface, and let $ h = c_1(H) $ in $ R^1(X) $. Let $ x $ be an element in $ N(X) $ with $ \rk(x) > 0 $. Then there exists a real number $ b $ such that for all real numbers $ t \gg 0 $, there is a stability condition $ \sigma_t = \sigma_{bh,th} $ in the distinguished component $ \Stab^\dagger(X) $ such that there is an isomorphism
	\begin{equation*}
		M_H^{(s)}(X,x) \cong M_{\sigma_t}^{(s)}(X, x)
	\end{equation*}
	between Gieseker and Bridgeland moduli spaces.
\end{prop}
\begin{proof}
	By Proposition \ref{prop_stability_conditions_on_surfaces}, there is a stability condition 
	$$ \sigma_t = \sigma_{bh,th} \in \Stab^\dagger(X) $$
	for all real numbers $ b $ and $ t > 0. $  
	Choose $ b < \mu_h(x)/h^2. $ Since $ \rk(x) > 0 $ and $ bh^2 < \mu_h(x), $ by \cite[Exercise 6.27]{macri2019lectures}, we can choose $ t_0 > 0 $ such that for all $ t \geq t_0, $ $ E $ is a  $ \sigma_t $-(semi)stable object in $ \D(X) $ with numerical class $ x $ is if and only if $ E $ is a $ H $-(semi)stable sheaf on $ X $ with numerical class $ x $. Therefore, we have
	\begin{equation*}
		M_{\sigma_t}^{(s)}(X, x) \cong M_H^{(s)}(X,x).
	\end{equation*}
\end{proof}
\begin{ex}[Hilbert schemes of points on $ K3 $ surfaces are Bridgeland moduli spaces.]
	Let $ (X,H) $ be a polarized $ K3 $ surface. Let $ I_Z $ be the ideal sheaf of a zero-dimensional subscheme $ Z \subset X $ of length $ n. $ Then the Mukai vector of $ I_Z $ is
	$$ v(I_Z) = (1,0,1-n) $$ 
	in $ R(X), $ the slope $ \mu_h(I_Z) = 0, $ and $ \gcd(r,d) = 1 $ since $ r = 1. $ Choose $ \beta = -h $ in $ R^1(X). $ Then we have $ \deg(\beta h) = -\deg(h^2) < 0 = \mu(I_Z).$ Choose a sufficiently large $ t > 0. $ Then we can identify the Hilbert scheme of $ n $ points on $ X $ as a Bridgeland moduli space:
	\begin{equation*}
	\Hilb^n(X) \cong M_H(X,\gamma(I_Z)) = M_H^s(X,\gamma(I_Z)) \cong M_{\sigma_t}^s(X, \gamma(I_Z))
	\end{equation*}
	where $ \sigma_t = \sigma_{-h,th} $ is a stability condition in the distinguished component $ \Stab^\dagger(X). $
\end{ex}
\begin{constrn}[Induced Bridgeland stability conditions on quotient stacks]\label{constrn_induced_stability}
	Consider a smooth projective variety $ X $ under an action of a finite group $ G $. Then we have a quotient stack $ \cX = [X/G] $ with its derived category 
	$$ \D(\cX) = \D(\Coh(\cX)) \cong \D(\Coh^G(X)). $$ 
	The canonical morphism $ p: X \to \cX $ is smooth, and hence induces an exact faithful functor 
	\begin{equation*}
		p^*: \D(\cX) \to \D(X)
	\end{equation*}
	and a ring homomorphism
	\begin{equation*}
		p^N: N(\cX) \to N(X), \quad \gamma(\cE) = \gamma(E,\phi) \mapsto \gamma(E)
	\end{equation*}
	for a $ G $-equivariant sheaf $ \cE = (E,\phi) $ on $ X. $
	The $ G $-action on $ X $ induces a $ G $-action on the manifold $ \Stab(X) $ via the auto-equivalence 
	$$ g^*: \D(X) \to \D(X) $$ 
	for each $ g $ in $ G. $
	Fix a $ G $-invariant stability condition $ \sigma = (Z, \cP) $
	on $ \D(X). $ By \cite[Section 2]{macri2009inducing}, it induces a stability condition $ \wtilde{\sigma} = (\wtilde{Z}, \wtilde{\cP}) $
	on $ \D(\cX) $ where the central charge $ \wtilde{Z} $ is a composition
	\begin{equation*}
		\wtilde{Z}: N(\cX) \xrightarrow{p^N} N(X) \xrightarrow{Z} \CC,
	\end{equation*}
	and the slicing $ \wtilde{\cP} $ is given by
	\begin{equation*}
		\wtilde{\cP}(\phi) = \{ \cE \in \D(\cX): p^*\cE \in \cP(\phi)\}
	\end{equation*}
	for all $ \phi \in \RR. $ Moreover, by \cite[Theorem 1.1]{macri2009inducing}, the assignment $ \sigma \mapsto \wtilde{\sigma} $ gives a closed embedding of complex manifolds 
	\begin{equation*}
		\Stab(X)^G \into \Stab(\cX).
	\end{equation*}
\end{constrn}
The following is an analogous result to Proposition \ref{prop_Gieseker_relation}.
\begin{prop}\label{prop_Bridgeland_relation}
	Let $ \cX = [X/G] $ in Construction \ref{constrn_induced_stability}. Let $ \sigma $ be a $ G $-invariant stability condition on $ \D(X). $ Let $ \cE $ be an object in $ \D(\cX) $ with $ E = p^*\cE $ in $ \D(X) $. Then $ \cE $ is $ \wtilde{\sigma} $-semistable if and only if $ E $ is $ \sigma $-semistable. If $ E $ is $ \sigma $-stable, then $ \cE $ is $ \wtilde{\sigma} $-stable.
\end{prop}
\begin{proof}
	This is a tautology by Construction \ref{constrn_induced_stability}.
\end{proof}
Suppose $ X $ is a $ K3 $ surface. If $ \sigma $ is $ G $-invariant stability condition in $ \Stab^\dagger(X), $ then the moduli stack $ \cM_{\wtilde{\sigma}}(\D(\cX), x) $ of $ \wtilde{\sigma} $-semistable complexes in $ \D(\cX) $ with numerical class $ x $ in $ N(\cX) $ is an algebraic stack and has a good moduli space $ M_{\wtilde{\sigma}}(\D(\cX), x) $.
\begin{thm}
	Let $ \cX = [X/G] $ where $ X $ is a $ K3 $ surface and $ G $ is a finite group acting on $ X $. Let $ \sigma $ be a $ G $-invariant stability condition in $ \Stab^\dagger(X). $ For every element $ x $ in $ N(\cX) $, the moduli stack $ \cM_{\wtilde{\sigma}}(\D(\cX), x) $ is an algebraic stack of finite type over $ \CC $ with a proper good moduli space.
\end{thm}
\begin{proof}
	This is {\cite[Theorem 3.22]{beckmann2022equivariant}} where $ X $ is a $ K3 $ surface.
\end{proof}
\begin{rmk}
	If there are no strictly $ \wtilde{\sigma} $-semistable complexes in $ \D(\cX) $ with numerical class $ x $, then the good moduli space $ M_{\wtilde{\sigma}}(\cX, x) $ is also a coarse moduli space.
\end{rmk}

\subsection{Proof of the main theorem}\label{sec_proof_main_thm}
In this section we will  prove the main theorem, i.e., Theorem \ref{thm_main} stated in the introduction.

Let $ \cX = [K3/G] $ in Situation \ref{situation_finite_group}. Recall that we have a line bundle $ \cH $ on $ \cX $ which descends to an ample line bundle on the projective surface $ X/G. $ The natural morphism $ p: X \to \cX $ pulls back $ \cH $ to a $ G $-invariant ample line bundle $ H $ on $ X. $ We have a commutative diagram 
\vspace{3pt}
\begin{equation*}
\begin{tikzcd}[column sep=2.5em,row sep=3em,every label/.append style={font=\normalsize}]
N(\cX) \arrow[r, "p^N"{outer sep = 2pt}] \arrow[d, "\orbv"{left, outer sep = 2pt}] 
& N(X) \arrow[d, "v"{right, outer sep = 2pt}, ] \\ 
R(I\cX)_\CC \arrow[r, "p^R"{outer sep = 2pt}] & R(X),
\end{tikzcd}\vspace{6pt}
\end{equation*}
where the pullbacks $ p^N $ and $ p^R $ are ring homomorphisms, the Mukai vector map $ v $ and the orbifold Mukai vector map $ \orbv $ are additive maps. Recall that
$$ R(I\cX)_\CC = R(\cX) \oplus R(I_t\cX)_\CC \cong R(X)^G \oplus \CC^m. $$
If $ G $ is non-trivial, then we observe the following:
\begin{enumerate}[font=\normalfont,leftmargin=*]
	\item $ \orbv $ is an injection and $ v $ is an isomorphism.
	\item $ p^N $ is not injective since the numerical classes of two $ G $-equivariant sheaves $ (E, \phi_1) $ and $ (E, \phi_2) $ on $ X $ are mapped to the same element $ \gamma(E) $ in $ N(X) $. $ p^N $ maps into the $ G $-invariant subspace $ N(X)^G $ of $ N(X) $, but may not map onto $ N(X)^G $ since there may be $ G $-invariant sheaves on $ X $ which are not $ G $-linearizable, i.e., which do not lift to $ G $-equivariant sheaves on $ X $.
	\item $ p^R $ is a projection $ (v,(v_{ij})) \mapsto v $ and hence is not injective. $ p^R $ maps onto the $ G $-invariant subspace $ R(X)^G $ of $ R(X) $.
\end{enumerate}
\vspace{8pt}

We first identify Gieseker moduli spaces of stable sheaves in $ \Coh(\cX) $ with Bridgeland moduli spaces of stable objects in $ \D(\cX) $.
\begin{lem}\label{lem_Bridgeland_vs_Gieseker_K3/G}
	Let $ \cX = [K3/G] $ in Situation \ref{situation_finite_group} with a $ G $-invariant ample line bundle $ H $ on $ X $. Choose an element $ x $ in $ N(\cX). $ Let $ y = p^N x $ in $ N(X)^G $ with Mukai vector $ v = (r, c_1, s) $ in $ R(X)^G $. Suppose $ r > 0 $ and either of the following conditions is satisfied:
	\begin{enumerate}[label=\textnormal{(\roman*)},font=\normalfont,leftmargin=2em]
		\item $ \gcd(r, d, s) = 1. $
		\item $ y $ is primitive and $ H $ is $ y $-generic.
	\end{enumerate}
	Then there exists an interval $ (a,b) $ in $ \RR $ such that for all $ t \in (a,b) $, there is a $ G $-invariant stability condition $ \sigma_t $ in $ \Stab^\dagger(X) $ which is $ y $-generic, and
	\begin{equation*}
		M_\cH^{s}(\cX, x) = M_\cH(\cX, x) \cong M_{\wtilde{\sigma}_t}(\cX,x) = M_{\wtilde{\sigma}_t}^{s}(\cX,x).
	\end{equation*}
\end{lem}
\begin{proof}
	Since $ r > 0 $, Proposition \ref{prop_Bridgeland_vs_Gieseker} applies, so there is a real number $ b $ and a real number $ t_0 > 0 $ such that for all $ t \geq t_0 $, there is a stability condition 
	$$ \sigma_t = \sigma_{bh,th} \in \Stab^\dagger(X) $$ 
	such that there is an isomorphism
	$$ M_H^{(s)}(X,x) \cong M_{\sigma_t}^{(s)}(X, x) $$ 
	between Gieseker and Bridgeland moduli spaces on $ X $.
	Moreover, these $ \sigma_t $ are $ G $-invariant since $ H $ is $ G $-invariant. Choose any $ t \geq t_0. $ Take an object $ \cE $ in $ \D(\cX) $ with numerical class $ x $ in $ N(\cX) $. Then we have a chain of equivalences:
	\begin{align*}
		& \cE\ \text{is } \wtilde{\sigma}_t \text{-semistable in } \D(\cX)\\
		\Leftrightarrow \quad & p^*\cE\ \text{is } \sigma_t\text{-semistable in } \D(X) & \text{by Proposition \ref{prop_Bridgeland_relation}}\\
		\Leftrightarrow \quad & p^*\cE \ \text{is } H\text{-semistable in } \Coh(X) & \text{by Proposition \ref{prop_Bridgeland_vs_Gieseker}}\\
		\Leftrightarrow \quad & \cE \ \text{is } \cH\text{-semistable in } \Coh(\cX) & \text{by Proposition \ref{prop_Gieseker_relation}}
	\end{align*}
	Hence there is an isomorphism 
	$$ M_\cH(\cX, x) \cong M_{\wtilde{\sigma_t}}(\cX,x) $$
	between Gieseker and Bridgeland moduli space on $ \cX $.
	
	Since the set of walls determined by $ y $ is locally finite, every interval in $ \RR $ intersects finitely many of them. Therefore, we can choose an interval $ (a,b) $ with $ a \geq t_0 $ such that no walls cross it. Then every $ t \in (a,b) $ determines a stability condition $ \sigma_t $ which is $ y $-generic. 
	
	By assumption, there are no strictly $ H $-semistable sheaves in $ \Coh(X) $ with numerical class $ y. $ By Proposition \ref{prop_Bridgeland_vs_Gieseker}, there are also no strictly $ \sigma_t $-semistable sheaves in $ \D(X) $ with numerical class $ y. $ Note that Proposition \ref{prop_Bridgeland_relation} gives an implication
	\begin{equation}\label{implication1}
		p^*\cE\ \text{is } \sigma_t\text{-stable in } \D(X) \ \Rightarrow \ \cE \ \text{is } \wtilde{\sigma}_t\text{-stable in } \D(\cX),
	\end{equation}
	and Proposition \ref{prop_Gieseker_relation} gives another implication
	\begin{equation}\label{implication2}
		p^*\cE\ \text{is } H\text{-stable in } \Coh(X) \ \Rightarrow \ \cE \ \text{is } \cH\text{-stable } \Coh(\cX).
	\end{equation}
	By the implications (\ref{implication1}) and (\ref{implication2}), we deduce there are neither strictly $ \wtilde{\sigma}_t $ semistable objects in $ \D(\cX) $ nor strictly $ \cH $-semistable sheaves in $ \Coh(\cX) $ with numerical class $ x $. This completes the proof.
\end{proof}

Now we consider the minimal resolution $ M \to X/G, $ which gives a derived McKay correspondence
\begin{equation*}
	\Phi: \D(\cX) \xrightarrow{\sim} \D(M).
\end{equation*}
Recall that we have a space $ \Stab(\cX) $ of stability conditions on the derived category $ \D(\cX), $ as well as a space $ \Stab(M) $ of stability conditions on $ \D(M). $ The two spaces $ \Stab(\cX) $ and $ \Stab(M) $ are complex manifolds which carry distinguished components $ \Stab^\dagger(\cX) $ and $ \Stab^\dagger(M) $ respectively. The derived equivalence $ \Phi $ induces an isomorphism of complex manifolds
\begin{equation*}
	\Phi^S: \Stab(\cX) \xrightarrow{\sim} \Stab(M), \quad (Z, \cP) \mapsto (Z \circ \Phi^{-1}, \Phi \circ \cP),
\end{equation*}
where the slicing $ \Phi \circ \cP $ is defined by $ (\Phi \circ \cP)(\phi) = \Phi(\cP(\phi)) $ for each $ \phi \in \RR. $ We don't know whether $ \Phi^S $ preserves distinguished components. But the following result suffices for our purpose.  
\begin{lem}[{\cite[Proposition 6.1]{beckmann2022equivariant}}]\label{lem_distinguished_FM}
	Let $ \cX = [K3/G] $ in Situation \ref{situation_finite_group} with the minimal resolution $ M \to X/G $ where $ M = \GHilb X $. Let $ \sigma $ be a $ G $-invariant stability condition in the distinguished component $ \Stab^\dagger(X) $. Let $ \wtilde{\sigma} $ denote the induced stability condition on $ \D(\cX) $. Then the derived equivalence $ \Phi: \D(\cX) \xrightarrow{\sim} \D(M) $ gives a stability condition $ \Phi^S(\wtilde{\sigma}) $ in the distinguished component $ \Stab^\dagger(M). $
\end{lem}

The derived equivalence $ \Phi $ preserves Bridgeland stabilities.
\begin{lem}\label{lem_isom_Bridgeland_FM}
	Let $ \tau $ be a stability condition on $ \D(\cX) $. Let $ x $ be an element in $ N(\cX). $ Then there is an isomorphism of Bridgeland moduli stacks
	\begin{equation*}
		\Phi: M_\tau^{(s)}(\cX,x) \xrightarrow{\simeq} M_{\Phi^S(\tau)}^{(s)}(M,\Phi^N(x)).
	\end{equation*}
\end{lem}
\begin{proof}
	This is a tautology since the derived equivalence $ \Phi $ maps $ \tau $-(semi)stable objects in $ \D(\cX) $ to $ \Phi^S(\tau) $-(semi)stable objects in $ \D(M) $.
\end{proof}
Now we are ready to prove our main theorem. Let's restate it here.
\begin{thm}\label{thm_K3/G}
	Let $ \cX = [K3/G] $ in Situation \ref{situation_finite_group}. Let $ x $ be an element in $ N(\cX) $ with $ y = p^N x $ in $ N(X)^G. $ Suppose $ y $ is primitive with $ \rk(y) > 0 $ and $ H $ is $ y $-generic.
	Then $ M(\cX,x) = M^s(\cX,x), $ which is non-empty if and only if $ \inprod{\orbv(x)^2}_{I\cX} \leq 2. $ If $ M(\cX,x) $ is non-empty, then it is an irreducible symplectic manifold of dimension $ n = 2 - \inprod{\orbv(x)^2}_{I\cX} $ deformation equivalent to $ \Hilb^{n/2}(X) $.
\end{thm}
\begin{proof}
	Since $ y $ is primitive and $ H $ is $ y $-generic, by Proposition \ref{prop_K3/G_smoothness_dim}, we have
	$$ M(\cX,x) = M^s(\cX,x), $$ 
	which is a smooth projective scheme of dimension $ n = 2 - \inprod{\orbv(x)^2}_{I\cX} $ if not empty. Lemma \ref{lem_Bridgeland_vs_Gieseker_K3/G} tells that there is an interval $ (a,b) $ such that for each $ t \in (a,b), $ there exists a $ G $-invariant stability condition $ \sigma_t $ in $ \Stab^\dagger(X) $ which induces a stability condition $ \wtilde{\sigma}_t $ on $ \D(\cX) $, and there is an isomorphism
	\begin{equation}\label{eq_first_isom}
		M(\cX,x) \cong M_{\wtilde{\sigma}_t}(\cX,x)
	\end{equation}
	between Gieseker and Bridgeland moduli spaces on $ \cX $.
	Take any $ t \in (a,b). $ By Lemma \ref{lem_isom_Bridgeland_FM}, there is an isomorphism
	\begin{equation}\label{eq_second_isom}
		M_{\wtilde{\sigma}_t}(\cX,x) \cong M_{\Phi^S(\wtilde{\sigma}_t)}(M,\Phi^N(x))
	\end{equation}
	between Bridgeland moduli spaces on $ \cX $ and $ M $.
	
	Since $ y $ is primitive in $ N(X), $ the element $ x $ is also primitive in $ N(\cX) $. Therefore, under the isomorphism 
	$$ \Phi^N: N(\cX) \xrightarrow{\sim} N(M), $$ 
	the transformed element $ \Phi^N(x) $ is also primitive in $ N(M) $.
	
	By Lemma \ref{lem_distinguished_FM}, the stability condition $ \Phi^S(\wtilde{\sigma}_t) $ lies in the distinguished component $ \Stab^\dagger(M) $ because $ \sigma_t $ is in $ \Stab^\dagger(X) $. Now, if we vary $ t $ in $ (a,b), $ the stability conditions $ \sigma_t $ form a real curve in the manifold $ \Stab(X)^G $, which maps to another one in $ \Stab(\cX) $ via the embedding 
	$$ \Stab(X)^G \into \Stab(\cX), \quad \sigma_t \mapsto \wtilde{\sigma}_t. $$ 
	Under the isomorphism 
	$$ \Phi^S: \Stab(\cX) \xrightarrow{\sim} \Stab(M), $$ 
	the stability conditions $ \Phi^S(\wtilde{\sigma}_t) $ trace a real curve in the manifold $ \Stab(M) $ as $ t $ varies in $ (a,b). $ Choose some $ c \in (a,b) $ such that $ \Phi^S(\wtilde{\sigma}_c) $ doesn't lie on any of the walls determined by $ \Phi^N(x). $ Now we have a primitive Mukai vector $ \Phi^N(x) $ in $ N(M) $ and a stability condition $ \Phi^S(\wtilde{\sigma}_c) $ in the distinguished component $ \Stab^\dagger(M) $ which is $ \Phi^N(x) $-generic. By Theorem \ref{thm_bottini}, the Bridgeland moduli space 
	$$ M_{\Phi^S(\wtilde{\sigma}_c)}(M,\Phi^N(x)) $$ 
	is non-empty if and only if $ \inprod{v(\Phi^N(x))^2}_M \leq 2, $ and is an irreducible symplectic manifold of dimension $ 2 - \inprod{v(\Phi^N(x))^2}_M  $ deformation equivalent to a Hilbert scheme of points on a $ K3 $ surface if not empty. By formula (\ref{eq_Fourier_Mukai}) and the two isomorphisms (\ref{eq_first_isom}) and (\ref{eq_second_isom}), the same results hold for the Gieseker moduli space $ M_\cH(\cX,x) $.
\end{proof}

\begin{rmk}
	Theorem \ref{thm_K3/G} reduces to Theorem \ref{thm_yoshioka} when $ G $ is trivial. On the other hand, a $ G $-equivariant Hilbert scheme of points on $ X $ (see \cite{bryan2022g}) is a special case of $ M(\cX, x) $ in Theorem \ref{thm_K3/G} where $ x = \gamma(I_\cZ) $ for a zero-dimensional substack $ \cZ $ of $ \cX $.
\end{rmk}

\subsection{Equivariant Hilbert schemes of points}\label{sec_Hilb}
In this section we define equivariant Hilbert schemes of points on a projective scheme under a finite group action, and apply Theorem \ref{thm_K3/G} in the rank one case to obtain a result on equivariant Hilbert schemes of points on a $ K3 $ surface.

Fix a projective stack $ \cX = [X/G] $ where $ X $ is a projective scheme and $ G $ is a finite group acting on $ X $. 

\begin{defn}
	Let $ N_0(\cX) $ denote the subgroup of $ N(\cX) $ generated by numerical classes $ \gamma(\cO_\cY) $ in $ N(\cX) $ for all zero-dimensional substack $ \cY $ of $ \cX $. For an element $ \delta $ in $ N_0(\cX), $ define a \tb{Hilbert functor}
	\begin{align*}
	\Hilb_\cX^\delta: (\Sch/\CC) & \to (\text{Set}) \\
	T & \mapsto 
	\left\{ 
	\begin{matrix}
	\text{closed substacks $ \cZ \subset \cX \times T $ such that} \\  
	\text{$ \cZ \to T $ is of finite presentation, flat, and proper} \\
	\text{with $ \gamma(\cO_{\cZ_t}) = \delta $ in $ N_0(\cX) $ for all $ t \in T $}
	\end{matrix}
	\right\}.
	\end{align*}
\end{defn}
\vspace{5pt}

\begin{rmk}
	Take an element $ \delta $ in $ N_0(\cX) $. The Hilbert functor $ \Hilb_\cX^\delta $ is a special case of a \tb{Quot functor} in \cite{olsson2003quot}. It is represented by a projective scheme $ \Hilb_\cX(\delta) $ which parametrizes zero-dimensional substacks of $ \cX $ with numerical class $ \delta. $
\end{rmk}

In the same spirit of how equivariant moduli spaces of sheaves are defined, we make the following definition.
\begin{defn}
	Let $ \cX = [X/G] $ where $ X $ is a projective scheme and $ G $ is a finite group acting on $ X. $ Let $ \delta $ be an element in $ N_0(\cX). $ The projective scheme $ \Hilb(\cX, \delta) $ is called the \tb{$ G $-equivariant Hilbert scheme} of points on $ X $ with numerical class $ \delta $.
\end{defn}

A zero-dimensional substack $ \cY = [Y/G] $ of $ \cX $ determines a short exact sequence 
\begin{equation*}
0 \to I_\cY \to \cO_\cX \to \cO_\cY \to 0
\end{equation*}
of coherent sheaves on $ \cX, $ where $ I_\cY $ is the ideal sheaf of $ \cY $, which corresponds to the $ G $-equivariant ideal sheaf $ I_Y $ of the $ G $-invariant subscheme $ Y $ of $ X. $ Therefore, we have the following:
\begin{lem}\label{lem_Hilb_isom}
	Let $ \cX = [X/G] $ where $ X $ is a projective scheme and $ G $ is a finite group acting on $ X $. For every element $ \delta $ in $ N_0(\cX), $ there is a natural isomorphism
	\begin{equation*}
	\Hilb_\cX(\delta) \cong M(\cX, 1 - \delta) 
	\end{equation*}
	of projective schemes, where $ 1 = \gamma(\cO_\cX) $ in $ N(\cX). $
\end{lem}

Now we prove a result which was implied in \cite[Section 5]{bryan2022g}.
\begin{lem}\label{lem_basis}
	Let $ \cX = [X/G] $ where $ X $ is an irreducible projective scheme and $ G $ is a finite group acting on $ X. $ Suppose there are finitely many non-free $ G $-orbits on $ X $. Let $ p $ be a generic point on $ \cX $ corresponding to a free $ G $-orbit $ Gx $, and let $ p_1, \dots, p_r $ be the orbifold points on $ \cX $ corresponding to the non-free $ G $-orbits $ Gx_1, \dots, Gx_r. $ For each $ 1 \leq i \leq r, $ label the irreducible representations of the finite group $ G_i = \Stab_G(x_i) $ by
	$$ \rho_{i,0},\ \rho_{i,1}\ \dots,\ \rho_{i,n_i}, $$ 
	where $ \rho_{i,0} $ is trivial. Then 
	$ N_0(\cX) $ is a free abelian group of rank
	\begin{equation*}
	r_0 = 1 + n_1 + \cdots + n_r
	\end{equation*}
	with a basis
	\begin{equation*}
	B = \{ b, b_{ij}: 1 \leq i \leq r, 1 \leq j \leq n_i \},
	\end{equation*}
	where $ b = \gamma(\cO_p) $ in $ N_0(\cX) $, and
	\begin{equation*}
	b_{ij} = \gamma(\cO_{p_i} \otimes \rho_{ij})
	\end{equation*}
	in $ N_0(\cX) $ for each $ 1 \leq i \leq r, 1 \leq j \leq n_i $.
\end{lem}
\begin{proof}
	Take an element $ \delta \in N_0(\cX). $ Then it can be written as
	\begin{equation*}
	\delta = n\gamma(\cO_p) + \sum_{i=1}^r \sum_{j=0}^{n_i} m_{ij} \gamma(\cO_{p_i} \otimes \rho_{ij}),
	\end{equation*}
	where $ n $ and $ m_{ij} $'s are integers. Fix an index $ 1 \leq i \leq r $. 
	We have two natural morphisms
	\begin{equation*}
	u_i: G_ix \to p = [G_ix/G_i] \xrightarrow{\sim} \pt \quad \text{and} \quad v_i: x_i \to p_i = [x_i/G_i] \xrightarrow{\sim} BG_i.
	\end{equation*}
	Therefore, we have two natural equivalences of abelian categories
	\begin{equation*}
	\Coh^{G_i}(G_ix) \cong \Coh(p) \quad \text{and} \quad \Coh^{G_i}(x_i) \cong \Coh(p_i),
	\end{equation*}
	where both $ \Coh(p) $ and $ \Coh(p_i) $ are subcategories of $ \Coh(\cX). $
	Under these equivalences, we can identify
	\begin{align*}
	\cO_p = (\cO_{G_i}, \phi_i), \quad \cO_{G_i} = u_i^*\cO_p, \quad \phi_i \cong \rho_{i,\reg}
	\end{align*}
	where $ \rho_{i,\reg} $ is the regular representation of $ G_i $, and
	\begin{equation*}
	\cO_{p_i}  = (\cO_{x_i}, \psi_i), \quad \cO_{x_i} = v_i^*\cO_{p_i}, \quad \psi_i \cong \rho_{i,0}.
	\end{equation*}
	Tensoring $ \cO_{p_i} $ by $ \rho_{i,\reg} $, i.e., pulling back the sheaf $ \rho_{i,\reg} $ from $ BG_i $ to $ p_i $, we then have another identification
	\begin{equation*}
	\cO_{p_i} \otimes \rho_{i,\reg} = (\cO_{x_i}^{\oplus |G_i|}, \vphi_i), \quad \cO_{x_i}^{\oplus |G_i|} = v_i^*\left(\cO_{p_i} \otimes \rho_{i,\reg}\right), \quad \vphi_i \cong \rho_{i,\reg}.
	\end{equation*}
	Therefore, $ \cO_p \cong \cO_{p_i} \otimes \rho_{i,\reg} $ as representations of $ G_i, $ which implies 
	\begin{equation*}
	b = \gamma(\cO_p) = \gamma(\cO_{p_i} \otimes \rho_{i,\reg})
	\end{equation*}
	in $ N_0(\cX), $ and hence there is a relation
	\begin{equation}\label{eq_relation}
	b = \gamma(\cO_{p_i}) + \sum_{j=1}^{n_i} \deg(\rho_{ij}) b_{ij}
	\end{equation}
	in $ N_0(\cX).$ The elements $ b $ and $ b_{ij} $'s are independent in $ N_0(\cX), $ and by the relation (\ref{eq_relation}), every element $ \delta $ in $ N_0(\cX) $ can now be uniquely written as
	\begin{equation*}
	\delta =  nb + \sum_{i=1}^r \sum_{j=1}^{n_i} m_{ij}b_{ij},
	\end{equation*}
	where $ n $ and $ m_{ij} $'s are integers.
\end{proof}
\begin{rmk}
	In the setup of Lemma \ref{lem_basis}, each $ (n,m) = (n,(m_{ij})) $ in $ \ZZ \oplus \ZZ^{n_0-1} $ corresponds to a numerical class
	\begin{equation*}
	\delta(n,m) = nb + \sum_{i=1}^r \sum_{j=1}^{n_i} m_{ij}b_{ij}
	\end{equation*}
	in $ N_0(\cX), $ which determines a $ G $-equivariant Hilbert scheme 
	\begin{equation*}
	\Hilb_\cX(n,m) := \Hilb_\cX(\delta(n,m))
	\end{equation*}
	of points on $ X $ with numerical class $ \delta(n,m) $.
	Take an object $ \cY = [Y/G] $ in $ \Hilb_\cX(n,m). $ Then the $ G $-invariant zero-dimensional subscheme $ Y $ of $ X $ has length 
	\begin{equation*}
	l(n,m) = h^0(Y, \cO_Y) = n|G| + \sum_{i=1}^r \sum_{j=1}^{n_i} |G/G_i| m_{ij} \deg(\rho_{ij})
	\end{equation*}
	in $ \ZZ, $ where $ |G/G_i| $ counts the number of points in each orbit $ Gx_i. $ Note that $ \Hilb_\cX(n,m) $ is empty if $ l(n,m) < 0 $.
\end{rmk}
By Lemma \ref{lem_Hilb_isom} and Theorem \ref{thm_K3/G} where $ \rk(x) = 1 $, we obtain the following:
\begin{cor}\label{cor_K3/G}
	Let $ \cX = [K3/G] $ in Situation \ref{situation_finite_group}.
	Let $ (n,m) $ be a vector in $ \ZZ \oplus \ZZ^{r_0-1} $ where $ r_0 = \rk(N_0(\cX))$. Let $ \orbv(n,m) = \orbv(1 - \delta(n,m)) $ in $ R(I\cX)_\CC. $ Then $ \Hilb_\cX(n,m) $ is non-empty if and only if $ \inprod{\orbv(n,m)^2}_{I\cX} \leq 2 $. If $ \Hilb_\cX(n,m) $ is non-empty, then it is an irreducible symplectic manifold of dimension $ d =  2 - \inprod{\orbv(n,m)^2}_{I\cX} $ deformation equivalent to $ \Hilb^{d/2}(X) $.
\end{cor}
\begin{rmk}
	Consider $ \cX = [K3/G] $. It's proved in \cite[Proposition 5.1]{bryan2022g} that $ \Hilb_\cX(n,m) $ is birational to $ \Hilb^{d/2}(M) $ and hence is also deformation equivalent to it, where $ d $ is computed via the Nakajima quiver varieties associated to the Dykin diagrams of the ADE root systems at the orbifold points on $ \cX $. Recall that at each orbifold point $ p_i = [x_i/G_i], $ the local repesentation $ G_{i} \to \GL(T_{x_i}X) $ corresponds to a geometric quotient $ \CC^2/G_i $ with a ADE singularity at the origin. Let $ \Lambda_i $ denote the free abelian group generated by the exceptional divisors $ E_{i,1}, \dots, E_{i,n_i} $ with an intersection form
	\begin{equation*}
	\mu_i(\ {,} \ ): \Lambda_i \times \Lambda_i \to \ZZ.
	\end{equation*}
	Then the dimension $ d $ of $ \Hilb_\cX(n,m) $ is given by
	\begin{equation*}
	d = 2n + \sum_{i=1}^{r} \mu_i(D_i^2),
	\end{equation*}
	where for each $ i = 1, \dots, r, $ $ D_i = \sum_{j=1}^{n_i} m_{ij} E_{ij} $ in $ \Lambda_i. $ In particular, if $ G \cong \ZZ/2\ZZ, $ then there are eight orbifold points and
	\begin{equation}\label{eq_BG}
	d = 2(n - m_1^2 - \dots m_8^2).
	\end{equation}
	We will recover this formula in the next example.
\end{rmk}
\begin{ex}
	Let $ G = \{1, \sigma\} \cong \ZZ/2\ZZ $, i.e., $ \sigma $ is a Nikulin involution. Then we have eight fixed points $ x_1, \dots, x_8 $ on $ X $ which correspond to eight orbifold points $ p_1, \dots, p_8 $ on $ \cX, $ where
	$$ p_i = [Gx_i/G] = [x_i/\mu_2] \cong B\mu_2 $$
	for $ 1 \leq i \leq 8. $
	The abelian group $ N_0(\cX) $ is of rank $ n_0 = 9 $ with a canonical basis
	\begin{equation*}
	b = \gamma(\cO_p), \ b_1 = \gamma(\cO_{p_1} \otimes \rho_1),\ \dots,\ b_8 = \gamma(\cO_{p_8} \otimes \rho_1),
	\end{equation*}
	where $ \rho_1 $ is the non-trivial irreducible representation of $ \mu_2. $
	Take a vector
	$$ (n,m) = (n,m_1, \dots, m_8) \ \in \ZZ^9. $$ 
	We have an element
	\begin{equation*}
	\delta(n,m) = nb + m_1 b_1 + \dots + m_8 b_8 \ \in N_0(\cX)
	\end{equation*}
	with length
	\begin{equation*}
	l = l(n,m) = 2n + m_1 + \dots + m_8 \ \in \ZZ.
	\end{equation*}
	Now we compute the dimension $ d = d(n,m) $ of $ \Hilb_\cX(n,m) $. 
	The orbifold Mukai vector map
	\begin{equation*}
	\orbv: N(\cX) \to R(X)^G \oplus \ZZ^8
	\end{equation*}
	lands into $ R(X)^G \oplus \ZZ^8 \subset R(X)^G \oplus \CC^8 $ since for every vector bundle $ \cE = (E, \phi) $ on $ \cX, $ on a twisted sector $ \{\sigma,x_i\} \times B\mu_2 $ for $ 1 \leq i \leq 8, $ the twisted Mukai vector of $ \cE $
	$$ v_i(\cE) = \chi_{\phi_i}(\sigma) = \tr(\phi_i(\sigma))\ \in \ZZ $$
	since $ \phi_i(\sigma) $ can only have eigenvalues $ \pm 1. $ 
	Therefore, we have an orbifold Mukai vector
	\begin{equation*}
	\orbv(n,m) = \orbv(1) - n\orbv(b) - m_1 \orbv(b_1) - \dots - m_8 \orbv(b_8) \ \in R(X)^G \oplus \ZZ^8,
	\end{equation*}
	where:
	\begin{align*}
	\orbv(1) & = (1,0,1,1,\dots,1) \\
	\orbv(b) & = (0,0,2,0,\dots,0)
	\end{align*}
	Next, we use $ K $-theoretic techniques to compute each $ \orbv(b_i). $ For a coherent sheaf $ \cE $ on $ \cX, $ denote the numerical class of $ e^K(\cE) $ by $ e^N(\cX). $ Recall formula (\ref{eq_eK_V}), which in our case says
	\begin{equation*}
	\gamma(\cO_{p_i}) = e^N(\cE)
	\end{equation*}
	in $ N_0(\cX) $ for a rank two vector bundle $ \cE = (E,\phi) $ on $ \cX $. 
	Note that the fiber $ E_{x_i} $ is the tangent space $ T_{x_i}X: $ we have a $ G $-equivariant Koszul complex 
	\begin{equation*}
	0 \to \wedge^2 E^\vee \xrightarrow{t} E^\vee \xrightarrow{s} I \to 0
	\end{equation*}
	on $ X $ where $ s: E \to \cO $ vanishes at $ x_i, $ which restricts to a right exact sequence
	\begin{equation*}
	\wedge^2 E_{x_i}^\vee \xrightarrow{t_{x_i}} E_{x_i}^\vee \xrightarrow{s_{x_i}} I_{x_i}/I_{x_i}^2 \to 0
	\end{equation*}
	of vector spaces on $ x_i $, where $ t_{x_i} = 0 $, so we have
	\begin{equation*}
	E_{x_i} = \left(I_{x_i}/I_{x_i}^2\right)^\vee = T_{x_i}X.
	\end{equation*}
	Now we use formula (\ref{eq_orbch_sheaf}) to compute the $ i $-th component of $ \orbch(\cO_{p_i}) $ as
	\begin{align*}
	\ch^{\rho_{-1}} \left((e^N(E))_{x_i}\right) & = \ch^{\rho_{-1}} (e^N(E_{x_i})) = \ch^{\rho_{-1}} (e^N(T_{x_i}X)) = \ch^{\rho_{-1}} (e^N(2\rho_1)) \\
	& = \ch^{\rho_{-1}} (e^N(\rho_1))^2 = \ch^{\rho_{-1}} (1-\rho_1^\vee))^2 = \ch^{\rho_{-1}} (1-\rho_1)^2 \\
	& = \chi_{(1-\rho_1)^2}(-1) \\
	& = 4,
	\end{align*}
	where $ e^N(E) = e^N(\cE) $ is the $ G $-equivariant Euler class of $ E $.
	Therefore,
	\begin{align*}
	\orbv(b_i) = \orbch(\cO_{p_i}) = (0,0,1, 4e_i),
	\end{align*}
	where $ e_i $ is the $ i $-th element in the standard basis of $ \ZZ^8 $ for $ 1 \leq i \leq 8. $
	Now we have
	\begin{equation*}
	\orbv(n,m) = \left(1,0,1-l,1+4m_1,\dots,1+4m_8\right).
	\end{equation*}
	Recall the orbifold Mukai pairing (\ref{eq_orbv_pairing_K3/G}), which gives
	\begin{align*}
	\inprod{\orbv^2}_{I\cX} & = \frac{\inprod{v^2}_X}{|G|} + \frac{1}{8} \sum_{i=1}^{8} | v_{i} |^2 \ \in \ZZ.
	\end{align*}
	Now we have
	\begin{equation*}
	d(n,m)= 2(n - m_1^2 - \cdots - m_8^2) \ \in 2\ZZ.
	\end{equation*}
	Note that this matches formula (\ref{eq_BG}).
	
	Suppose $ \Hilb_\cX(n, m) $ is non-empty. Then we must have:
	\begin{align*}
	l & = 2n + m_1 + \cdots + m_8 \geq 0. \\
	d & = 2n - 2m_1^2 - \cdots - 2m_8^2 \geq 0.
	\end{align*}
	An immediate consequence is
	\begin{equation*}
	0 \leq n \leq l.
	\end{equation*}
	Since $ m_i + m_i^2 \geq 0, $ the second inequality $ d \geq 0 $ dominates. For each length $ l \geq 0 $, and each $ 0 \leq n \leq l, $ let's compute the number of solutions $ s $ in the set 
	\begin{equation*}
	S(l,n) = \{ (x_1, \dots, x_8) \in \ZZ^8: \sum_{i=1}^8 x_i = l - 2n, \ \sum_{i=1}^8 x_i^2 \leq n \},
	\end{equation*}
	and determine the corresponding equivariant Hilbert schemes. 
	Let's consider a few cases.
	\begin{enumerate}[font=\normalfont,leftmargin=*]
		\item Set $ l = 0. $ Then $ n = 0. $ We then have
		\begin{align*}
		& S(0,0) = \{0\},\ s = 1,\ d = 0,
		\end{align*}
		which gives
		\begin{align*}
		\Hilb_\cX(0,0) & = \{\cO_\cX\} = \pt.
		\end{align*}
		\item Set $ l = 1. $ Then $ 0 \leq n \leq 1. $ We then have:
		\begin{align*}
		& S(1,0) = \emptyset. \\
		& S(1,1) = \{-e_i: 1 \leq i \leq 8\}, \ s = 8,\ d = 0.
		\end{align*}
		Indeed, for each $ 1 \leq i \leq 8, $ we have:
		\begin{align*}
		& \Hilb_\cX(0,e_i) = \emptyset. \\
		& \Hilb_\cX(1,-e_i) = \{p_i\} = \pt.
		\end{align*}
		\item Set $ l = 2. $ Then $ 0 \leq n \leq 2. $ We have:
		\begin{align*}
		& S(2,0) = \emptyset. \\
		& S(2,1) = \{0\},\ s = 1,\ d = 2. \\
		& S(2,2) = \{-e_i-e_j: 1 \leq i \neq j \leq 8\}, \ s = \binom{8}{2} = 28, \ d = 0.
		\end{align*}
		The corresponding Hilbert schemes are:
		\begin{align*}
		& \Hilb_\cX(0,e_i+e_j) = \emptyset, \ & & \text{for} \ 1 \leq i, j \leq 8. \\
		& \Hilb_\cX(1,0) = \GHilb X = K3. \\
		& \Hilb_\cX(2,-e_i-e_j) = \{p_i,p_j\} = \pt, \ & & \text{for} \ 1 \leq i \neq j \leq 8.
		\end{align*}
		\item Set $ l = 3. $ Then $ 0 \leq n \leq 3. $ We have:
		\begin{align*}
		& S(3,0) = \emptyset.\\
		& S(3,1) = \{e_i: 1 \leq i \leq 8\},\ s = 8,\ d = 0. \\
		& S(3,2) = \{-e_i: 1 \leq i \leq 8\}, \ s = 8, \ d = 2. \\
		& S(3,0) = \{-e_i-e_j-e_k: 1 \leq i \neq j \neq k \neq i \leq 8\}, \ s = \binom{8}{3} = 56, \ d = 0.
		\end{align*}
		The corresponding Hilbert schemes are:
		\begin{align*}
		& \Hilb_\cX(0,e_i+e_j+e_k) = \emptyset, \ & & \text{for} \ 1 \leq i, j, k \leq 8. & \\
		& \Hilb_\cX(1,e_i) = \pt, \ & & \text{for} \ 1 \leq i \leq 8. & \\
		& \Hilb_\cX(2,-e_i) = K3, \ & & \text{for} \ 1 \leq i \leq 8.  & \\
		& \Hilb_\cX(3,-e_i-e_j-e_k) = \{p_i,p_j,p_k\} = \pt, \ & & \text{for } 1 \leq i \neq j \neq k \neq i \leq 8. & 
		\end{align*}
		Let's describe $ \Hilb_\cX(1,e_i) $ and $ \Hilb_\cX(2,-e_i). $ Take a point in $ \Hilb_\cX(1,e_i). $ Then it is a zero-dimensional substack 
		$ \cY_i = [Y_i/G] $
		of $ \cX $ supported at $ p_i $ such that 
		\begin{equation*}
		\gamma(\cO_{\cY_i}) = \gamma(\cO_{p_i} \otimes \rho_{i,\reg}) + \gamma (\cO_{p_i} \otimes \rho_1) = \gamma (\cO_{p_i} \otimes \rho_0) + 2\gamma (\cO_{p_i} \otimes \rho_1).
		\end{equation*}
		In some local coordinate system $ (x,y) $ at $ x_i $, this is a
		$ \mu_2 $-representation:
		\begin{equation*}
		\CC[x,y]/(x^2,xy,y^2) = \CC \oplus \CC x \oplus \CC y 
		\end{equation*}
		since the $ \mu_2 $-action near $ x_i $ is generated by the involution
		\begin{equation*}
		(x,y) \mapsto (-x, -y).
		\end{equation*}
		We also have
		$$ \Hilb_\cX(2,-e_i) \cong \Hilb_\cX(1,0) = \GHilb X = K3, $$ 
		because each point in $ \Hilb_\cX(2,-e_i) $ is a pair 
		$ \{p,p_i\}, $
		where $ p $ is a generic point on $ \cX. $ 
	\end{enumerate}
\end{ex}

\begin{rmk}
	Note that the results here agree with those obtained in \cite[Examples 4.5 and 4.6]{kamenova2022symplectic} where the authors described the connected components in the fixed loci $ \Hilb^l(X)^{G} $ for $ l = 2 $ and $ l = 3 $, where $ G \cong \ZZ/2\ZZ $ acts on $ \Hilb^l(X) $ by pulling back structure sheaves $ \cO_Y $ of zero-dimensional subschemes $ Y \subset X. $ In other words, there are identifications between $ G $-equivariant Hilbert schemes of points on $ X $ and the connected components in the fixed loci of the induced $ G $-actions on Hilbert schemes of points on $ X $.
\end{rmk}

\appendix
\section{Proof of Proposition \ref{prop_equiv_cat_sheaves}}\label{proof_equiv_cat_sheaves}

Let's fix a quotient stack $ \cX = [X/G]. $ Consider the 2-pullback diagram
\vspace{10pt}
\begin{equation*}\label{eq_can_morphism_to_stack}
\begin{tikzcd}[column sep=3em,row sep=3em]
G \times X \arrow[r, "\sigma"{outer sep = 2pt}] \arrow[d, "p_2"{right, outer sep = 2pt}] & X \arrow[d, "p"{right, outer sep = 2pt}] \\
X \arrow[r, "p"{outer sep = 2pt}] & \cX,
\end{tikzcd}\vspace{10pt}
\end{equation*}
where $ p: X \to \cX = [X/G] $ is the canonical morphism in the groupoid $ \Hom(X, \cX) $ given by the object
\vspace{10pt}
\begin{equation}\label{eq_can_obj}
\eta_1 = \left(
\begin{tikzcd}[column sep=3em,row sep=3em]
G \times X \arrow[r, "\sigma"{outer sep = 2pt}] \arrow[d, "p_2"{right, outer sep = 2pt}] & X \\
X
\end{tikzcd}\right) \vspace{10pt}
\end{equation}
of $ \cX $ over $ X. $

We first define the functor $ T $ on the objects of $ \Sh(\cX). $ Let $ \cE $ be a sheaf on $ \cX. $ Its pullback $ p^*\cE = \cE(\eta_1) $ on $ X $ comes with an isomorphism
\begin{equation*}\label{eq_can_equiv_structure}
\phi_\cE: p_2^* (p^*\cE) \xrightarrow{\sim} \sigma^* (p^*\cE)
\end{equation*}
of sheaves on $ G \times X, $ which is determined by $ \cE $ as follows.
There is another object
\vspace{10pt}
\begin{equation*}\label{eq_obj_over_GxX}
\eta_2 = \left(
\begin{tikzcd}[column sep=3em,row sep=3em]
G^2 \times X \arrow[r, "s"{outer sep = 2pt}] \arrow[d, "p_{23}"{right, outer sep = 2pt}] & X \\
G \times X
\end{tikzcd}\right) \vspace{8pt}
\end{equation*}
of $ \cX $ over $ G \times X. $ There are two arrows
\vspace{5pt}
\begin{equation*}
\begin{tikzcd}[column sep=5em, every label/.append style={font=\normalsize}]
\eta_2 \arrow[r, bend left=40, "A_1"{above, outer sep=3pt}] \arrow[r, bend right=40, "A_2"{below, outer sep=3pt}] & \eta_1
\end{tikzcd}
\quad \text{over} \quad
\begin{tikzcd}[column sep=4em, every label/.append style={font=\normalsize}]
G \times X \arrow[r, bend left=40, "a_1 = p_2"{above, outer sep=3pt}] \arrow[r, bend right=40, "a_2 = \sigma"{below, outer sep=3pt}] & X,
\end{tikzcd} \vspace{5pt}
\end{equation*}
which are given by two commutative diagrams:
\vspace{5pt}
\begin{equation*}\label{eq_arrow_1}
A_1 :
\begin{tikzcd}[column sep=3em,row sep=4em]
G^2 \times X \arrow[rr, bend left=40, "s"{outer sep = 2pt}] \arrow[r, "\mu \times \id_X"{outer sep = 2pt}] \arrow[d, "p_{23}"{right, outer sep = 2pt}] & G \times X \arrow[d, "p_2"{right, outer sep = 2pt}] \arrow[r, "\sigma"{outer sep = 2pt}] & X \\
G \times X \arrow[r, "p_2"{outer sep = 2pt}]  & X &
\end{tikzcd}
\end{equation*}
\begin{equation*}\label{eq_arrow_2}
A_2 :
\begin{tikzcd}[column sep=3em,row sep=4em]
G^2 \times X \arrow[rr, bend left=40, "s"{outer sep = 2pt}]  \arrow[r, "\id_G \times \sigma"{outer sep = 2pt}] \arrow[d, "p_{23}"{right, outer sep = 2pt}] & G \times X \arrow[d, "p_2"{right, outer sep = 2pt}] \arrow[r, "\sigma"{outer sep = 2pt}] & X \\
G \times X \arrow[r, "\sigma"{outer sep = 2pt}]  & X &
\end{tikzcd} \vspace{15pt}
\end{equation*}
The data of $ \cE $ include two isomorphisms of sheaves on $ G \times X $:
\vspace{3pt}
\begin{equation*}
\cE(A_1): \cE(\eta_2) \xrightarrow{\sim} p_2^* \cE(\eta_1) \quad \text{and} \quad
\cE(A_2): \cE(\eta_2) \xrightarrow{\sim} \sigma^* \cE(\eta_1) \vspace{3pt}
\end{equation*}
Therefore, we have an isomorphism
\vspace{3pt}
\begin{equation*}
\phi_\cE = \cE(A_2) \circ \cE(A_1)^{-1}: p_2^* \cE(\eta_1) \xrightarrow{\sim} \sigma^* \cE(\eta_1)
\end{equation*}
of sheaves on $ G \times X. $ In the next three pages, we check that $ \phi_\cE $ satisfies the cocycle condition in a $ G $-equivariant structure. There is a third object
\vspace{3pt}
\begin{equation*}\label{eq_obj_over_GxGxX}
\eta_3 = \left(
\begin{tikzcd}[column sep=3em,row sep=3em]
G^3 \times X \arrow[r, "t"{outer sep = 2pt}] \arrow[d, "p_{234}"{right, outer sep = 2pt}] & X \\
G^2 \times X
\end{tikzcd}\right) \vspace{8pt}
\end{equation*}
of $ \cX $ over $ G \times X, $ where $ t $ maps $ (g_1, g_2, g_3, x) \mapsto g_1 g_2 g_3 x. $ There are three arrows
\vspace{5pt}
\begin{equation*}
\begin{tikzcd}[column sep=7em, every label/.append style={font=\small}]
\eta_3 
\arrow[r, bend left=40, shorten <= 10pt, shorten >= 9pt, "B_1"{above, outer sep = 2pt}] 
\arrow[r, shorten <= 7pt, shorten >= 7pt, "B_2"{above}] 
\arrow[r, bend right=40, shorten <= 7pt, shorten >= 7pt, "B_3"{above}] & \eta_2
\end{tikzcd}
\quad \text{over} \quad
\begin{tikzcd}[column sep=5em, every label/.append style={font=\small}]
G^2 \times X 
\arrow[r, bend left=40, shorten <= 7pt, shorten >= 7pt, "b_1"{outer sep = 2pt}] 
\arrow[r, shorten <= 3pt, shorten >= 5pt, "b_2"{xshift = -1.5 pt}] 
\arrow[r, bend right=40, shorten <= 10pt, shorten >= 8pt, "b_3"] 
& G \times X,
\end{tikzcd}
\ \
\begin{tikzcd}[column sep=5em, every label/.append style={font=\small}]
\arrow[draw = none, r, shift left = 23pt, "b_1 = p_{23}"{xshift = -12 pt, yshift = 0 pt}] 
\arrow[draw = none, r, "b_2 = \mu \times \id_X"{yshift = -0.5 pt}] 
\arrow[draw = none, r, shift right = 23 pt, "b_3 = \id_G \times \sigma"{yshift = 0 pt}]
& {},
\end{tikzcd}\vspace{8pt}
\end{equation*}
which are described by three commutative diagrams:
\vspace{5pt}
\begin{equation*}
B_1 :
\begin{tikzcd}[column sep=4em,row sep=4.2em]
G^3 \times X \arrow[rr, bend left=40, "t"{outer sep = 2pt}]
\arrow[r, "\mu \times \id_{G \times X}"{outer sep = 2pt}]  \arrow[d, "p_{234}"{right, outer sep = 2pt}] 
& G^2 \times X \arrow[d, "p_{23}"{right, outer sep = 2pt}] \arrow[r, "s"{outer sep = 2pt}] & X \\
G^2 \times X \arrow[r, "p_{23}"{outer sep = 2pt}] & G \times X &
\end{tikzcd}
\end{equation*}
\begin{equation*}
B_2 :
\begin{tikzcd}[column sep=4em,row sep=4.2em]
G^3 \times X \arrow[rr, bend left=40, "t"{outer sep = 2pt}] \arrow[r, "\id_G \times \mu \times \id_X"{outer sep = 2pt}] \arrow[d, "p_{234}"{right, outer sep = 2pt}] & G^2 \times X \arrow[d, "p_{23}"{right, outer sep = 2pt}] \arrow[r, "s"{outer sep = 2pt}] & X \\
G^2 \times X \arrow[r, "\mu \times \id_X"{outer sep = 2pt}] & G \times X &
\end{tikzcd}
\end{equation*}
\begin{equation*}
B_3 :
\begin{tikzcd}[column sep=4em,row sep=4.2em]
G^3 \times X \arrow[rr, bend left=40, "t"{outer sep = 2pt}] \arrow[r, "\id_{G \times G} \times \sigma"{outer sep = 2pt}] \arrow[d, "p_{234}"{right, outer sep = 2pt}] & G^2 \times X \arrow[d, "p_{23}"{right, outer sep = 2pt}] \arrow[r, "s"{outer sep = 2pt}] & X \\
G^2 \times X \arrow[r, "\id_G \times \sigma"{outer sep = 2pt}] & G \times X &
\end{tikzcd} \vspace{10pt}
\end{equation*}
The six compositions $ A_i \circ B_j $ for $ i = 1, 2 $ and $ j = 1, 2, 3 $ give three different morphisms from $ \eta_3 $ to $ \eta_1 $ as follows:
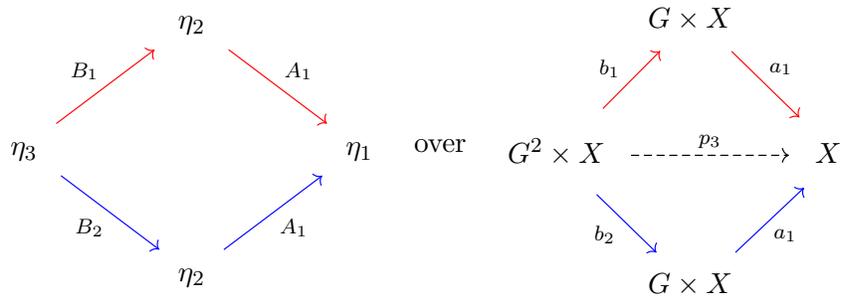
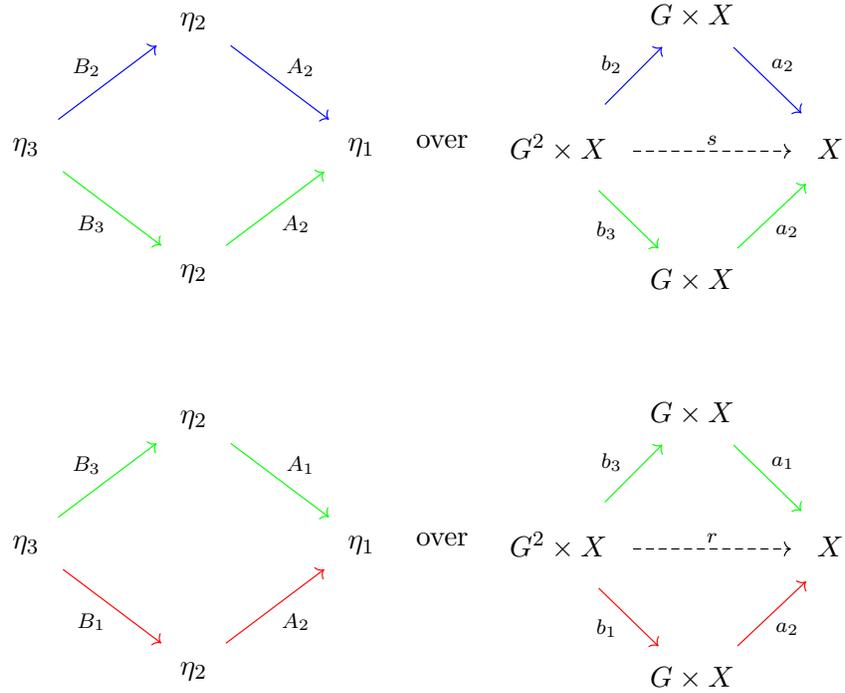
\begin{equation*}
\begin{tikzcd}[column sep=4em,row sep=3em]
& \eta_2 \arrow[dr, draw = red, shorten <= 5pt, shorten >= 5pt, "A_1"{}] & \\
\eta_3 \arrow[ur, draw = red, shorten <= 5pt, shorten >= 5pt, "B_1"{}] 
\arrow[dr, draw = blue, shorten <= 5pt, shorten >= 5pt, "B_2"{below left}] & & \eta_1 \\
& \eta_2 \arrow[ur, draw = blue, shorten <= 5pt, shorten >= 5pt, "A_1"{below right}]
\end{tikzcd}
\quad \text{over} \quad
\begin{tikzcd}[column sep=2em,row sep=3em]
&[-15pt] G \times X \arrow[dr, draw = red, shorten <= 8pt, shorten >= 5pt, shift left, "a_1"{}] & \\
G^2 \times X \arrow[ur, draw = red, shorten <= 10pt, shorten >= 5pt, "b_1"{}] 
\arrow[dr, draw = blue, shorten <= 10pt, shorten >= 5pt, "b_2"{below left}] 
\arrow[rr, shorten <= 5pt, shorten >= 5pt, dashed, "p_3"{}] & & X \\
& G \times X \arrow[ur, draw = blue, shorten <= 8pt, shorten >= 5pt, shift right, "a_1"{below right}] &
\end{tikzcd}
\end{equation*}
\begin{equation*}
\begin{tikzcd}[column sep=4em,row sep=3em]
& \eta_2 \arrow[dr, draw = blue, shorten <= 5pt, shorten >= 5pt, "A_2"{}] & \\
\eta_3 \arrow[ur, draw = blue, shorten <= 5pt, shorten >= 5pt, "B_2"{}] 
\arrow[dr, draw = green, shorten <= 5pt, shorten >= 5pt, "B_3"{below left}] & & \eta_1 \\
& \eta_2 \arrow[ur, draw = green, shorten <= 5pt, shorten >= 5pt, "A_2"{below right}]
\end{tikzcd}
\quad \text{over} \quad
\begin{tikzcd}[column sep=2em,row sep=3em]
&[-15pt] G \times X \arrow[dr, draw = blue, shorten <= 8pt, shorten >= 5pt, shift left, "a_2"{}] & \\
G^2 \times X \arrow[ur, draw = blue, shorten <= 10pt, shorten >= 5pt, "b_2"{}] 
\arrow[dr, draw = green, shorten <= 10pt, shorten >= 5pt, "b_3"{below left}] 
\arrow[rr, shorten <= 5pt, shorten >= 5pt, dashed, "s"{}] & & X \\
& G \times X \arrow[ur, draw = green, shorten <= 8pt, shorten >= 5pt, shift right, "a_2"{below right}] &
\end{tikzcd}
\end{equation*}
\begin{equation*}
\begin{tikzcd}[column sep=4em,row sep=3em]
& \eta_2 \arrow[dr, draw = green, shorten <= 5pt, shorten >= 5pt, "A_1"{}] & \\
\eta_3 \arrow[ur, draw = green, shorten <= 5pt, shorten >= 5pt, "B_3"{}] 
\arrow[dr, draw = red, shorten <= 5pt, shorten >= 5pt, "B_1"{below left}] & & \eta_1 \\
& \eta_2 \arrow[ur, draw = red, shorten <= 5pt, shorten >= 5pt, "A_2"{below right}]
\end{tikzcd}
\quad \text{over} \quad
\begin{tikzcd}[column sep=2em,row sep=3em]
&[-15pt] G \times X \arrow[dr, draw = green, shorten <= 8pt, shorten >= 5pt, shift left, "a_1"{}] & \\
G^2 \times X \arrow[ur, draw = green, shorten <= 10pt, shorten >= 5pt, "b_3"{}] 
\arrow[dr, shorten <= 10pt, draw = red, shorten >= 5pt, "b_1"{below left}] 
\arrow[rr, shorten <= 5pt, shorten >= 5pt, dashed, "r"{}] & & X \\
& G \times X \arrow[ur, draw = red, shorten <= 8pt, shorten >= 5pt, shift right, "a_2"{below right}] &
\end{tikzcd}\vspace{3pt}
\end{equation*}

Applying the cocycle condition of the sheaf $ \cE $ to compositions of arrows with the same color yields three commutative diagrams of isomorphisms of sheaves on $ G \times G \times X $:
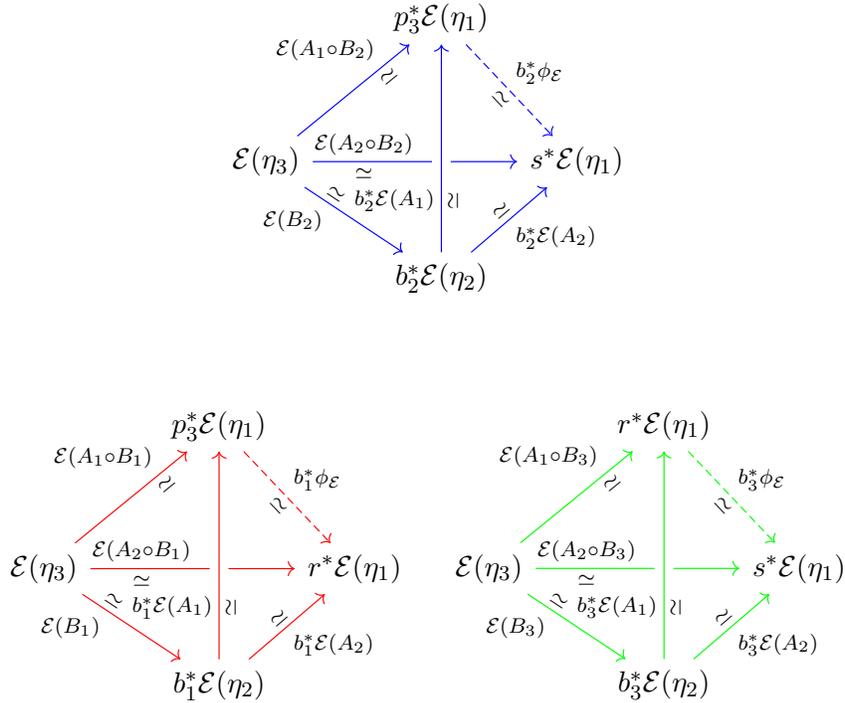
\vspace{10pt}
\begin{equation*}
\begin{tikzcd}[column sep=2.4em,row sep=3.2em, draw = blue]
& p_3^* \cE(\eta_1) 
\arrow[dr, dashed, "\simeq"{below, sloped}, "b_2^* \phi_{\cE}"{above right}] &[-20pt] \\
\cE(\eta_3) 
\arrow[ur, "\simeq"{below, sloped, near end}, "\cE(A_1 \circ B_2)"{above left, near end}] 
\arrow[rr, "\simeq"{below, sloped, near start}, "\cE(A_2 \circ B_2)"{above, near start}] 
\arrow[dr, "\simeq"{above, sloped, near start}, "\cE(B_2)"{below left, near start}] 
& & s^* \cE(\eta_1) \\[-10pt]
& b_2^* \cE(\eta_2) 
\arrow[uu, crossing over, "\simeq"{below, sloped, near start}, "b_2^* \cE(A_1)"{left, near start}] 
\arrow[ur, "\simeq"{above, sloped}, "b_2^* \cE(A_2)"{below right}] &
\end{tikzcd}
\end{equation*}
\begin{equation*}
\begin{tikzcd}[column sep=2.4em,row sep=3.2em, draw = red]
& p_3^* \cE(\eta_1) 
\arrow[dr, dashed, "\simeq"{below, sloped}, "b_1^* \phi_{\cE}"{above right}] &[-20pt] \\
\cE(\eta_3) 
\arrow[ur, "\simeq"{below, sloped, near end}, "\cE(A_1 \circ B_1)"{above left, near end}] 
\arrow[rr, "\simeq"{below, sloped, near start}, "\cE(A_2 \circ B_1)"{above, near start}] 
\arrow[dr, "\simeq"{above, sloped, near start}, "\cE(B_1)"{below left, near start}] 
& & r^* \cE(\eta_1) \\[-10pt]
& b_1^* \cE(\eta_2) 
\arrow[uu, crossing over, "\simeq"{below, sloped, near start}, "b_1^* \cE(A_1)"{left, near start}] 
\arrow[ur, "\simeq"{above, sloped}, "b_1^* \cE(A_2)"{below right}] &
\end{tikzcd}
\quad
\begin{tikzcd}[column sep=2.4em,row sep=3.2em, draw = green]
& r^* \cE(\eta_1) 
\arrow[dr, dashed, "\simeq"{below, sloped}, "b_3^* \phi_{\cE}"{above right}] &[-20pt] \\
\cE(\eta_3) 
\arrow[ur, "\simeq"{below, sloped, near end}, "\cE(A_1 \circ B_3)"{above left, near end}] 
\arrow[rr, "\simeq"{below, sloped, near start}, "\cE(A_2 \circ B_3)"{above, near start}] 
\arrow[dr, "\simeq"{above, sloped, near start}, "\cE(B_3)"{below left, near start}] 
& & s^* \cE(\eta_1) \\[-10pt]
& b_3^* \cE(\eta_2) 
\arrow[uu, crossing over, "\simeq"{below, sloped, near start}, "b_3^* \cE(A_1)"{left, near start}] 
\arrow[ur, "\simeq"{above, sloped}, "b_3^* \cE(A_2)"{below right}] &
\end{tikzcd} \vspace{15pt}
\end{equation*}
Therefore, we have
\begin{align*}
b_3^*\phi_\cE \circ b_1^*\phi_\cE & =  \left(\cE(A_2 \circ B_3) \circ \cE(A_1 \circ B_3)^{-1}\right) \circ \left(\cE(A_2 \circ B_1) \circ \cE(A_1 \circ B_1)^{-1}\right) \\
& = \cE(A_2 \circ B_3) \circ \left(\cE(A_1 \circ B_3)^{-1} \circ \cE(A_2 \circ B_1)\right) \circ \cE(A_1 \circ B_1)^{-1} \\
& = \cE(A_2 \circ B_2) \circ \id_{\cE(\eta_3)} \circ \cE(A_1 \circ B_2)^{-1} \\
& = b_2^*\phi_\cE ,
\end{align*}
and hence we have a commutative diagram
\vspace{5pt}
\begin{equation*}
\begin{tikzcd}[column sep=1.5em,row sep=4em]
p_3^*\cE(\eta_1) \arrow[rr, draw = blue, "b_2^* \phi_\cE"{outer sep = 2pt}, "\simeq"{below}] 
\arrow[dr, draw = red,  "b_1^* \phi_\cE"{below left}, "\simeq"{above, sloped}] & & s^*\cE(\eta_1) \\
& r^*\cE(\eta_1) \arrow[ur, draw = green, "b_3^* \phi_\cE"{below right}, "\simeq"{above, sloped}] &,
\end{tikzcd}\vspace{8pt}
\end{equation*}
where $ \cE(\eta_1) = p^*\cE. $ This is precisely the cocycle condition of isomorphism $ \phi_\cE $. Therefore, the assignment
\begin{equation*}
T(\cE) = (p^*\cE, \phi_\cE)
\end{equation*}
is a $ G $-equivariant sheaf on $ X, $ i.e., an object of $ \Sh^G(X). $

Next, we define the functor $ T $ on the morphisms of $ \Sh(\cX). $ Let $ \Phi: \cE \to \cF $ be a morphism of sheaves on $ \cX. $ Its pullback $ p^*\Phi = \Phi(\eta_1) $ is a morphism
$$ \Phi(\eta_1): \cE(\eta_1) \to  \cF(\eta_1) $$
of sheaves on $ X $. It's compatible with $ G $-equivariant structures $ \phi_{\cE} $ and $ \phi_{\cF} $ as can be checked by the compatibility of $ \cE $ and $ \cF $ along the two arrows $ A_1 $ (over $ a_1 = p_2 $) and $ A_2 $ (over $ a_2 = \sigma $):
\vspace{5pt}
\begin{equation}\label{eq_triangular_prism1}
\begin{tikzcd}[column sep=6em,row sep=1.2em]
& \cF(\eta_2) \ar[dd, "\cF(A_1)"{left}, "\simeq"{above, sloped}] 
\ar[dddr, "\cF(A_2)", "\simeq"{below, sloped}] & \\
\cE(\eta_2) \ar[ur, "\Phi(\eta_2)"{near end, outer sep = 2pt}] 
\ar[dd, "\cE(A_1)"{left}, "\simeq"{above, sloped}] & & \\
& p_2^*\cF(\eta_1) \ar[dr, dashed, "\phi_\cF"{below, near start}, "\simeq"{near start, sloped}] & \\
p_2^*\cE(\eta_1) \ar[ur, "p_2^*\Phi(\eta_1)"{very near end, outer sep = 2pt}] 
\ar[dr, dashed, "\phi_\cE"{below, near start}, "\simeq"{near start, sloped}] & & \sigma^*\cF(\eta_1) \\
& \sigma^*\cE(\eta_1) 
\ar[uuul, leftarrow, crossing over, "\cE(A_2)"{above right, near start}, "\simeq"{below, sloped, near start}] 
\ar[ur, "\sigma^*\Phi(\eta_1)"] &
\end{tikzcd}\vspace{8pt}
\end{equation}
Therefore, 
$$ p^*\Phi: (p^*\cE, \phi_\cE) \to  (p^*\cF, \phi_\cF) $$
is a $ G $-equivariant morphism between $ G $-equivariant sheaves on $ X $. We define 
$$ T(\Phi) = p^*\Phi. $$ 

Now we show that the assignment
\begin{align*}
T: \Hom_{\Sh(\cX)}(\cE, \cF) & \to \Hom_{\Sh^G(X)}(T(\cE), T(\cF)) \\
\Phi & \mapsto p^*\Phi
\end{align*}
is a bijection. Take a $ G $-equivariant morphism $ \vphi: T(\cE) \to T(\cF). $ It suffices to show that there is a unique morphism $ \Phi: \cE \to \cF $ such that $ p^*\Phi = \vphi. $
Take any object
\vspace{8pt}
\begin{equation*}\label{eq_arbitrary_object}
\eta_Y = \left(
\begin{tikzcd}[column sep=3em,row sep=3em]
Z \arrow[r, "c_2"{outer sep = 2pt}] \arrow[d, "c_1"{right, outer sep = 2pt}] & X \\
Y
\end{tikzcd}\right) \vspace{8pt}
\end{equation*}
of $ \cX $ over a scheme $ Y. $ We will construct a morphism
\begin{equation*}
\vphi_Y: \cE(\eta_Y) \to \cF(\eta_Y)
\end{equation*}
such that the collection of morphisms $ \vphi_Y $ satisfies condition (\ref{eq_sheaf_morphism_on_stack}).
Pulling back $ \eta_Y $ along $ \pi: Z \to Y $ gives another object
\vspace{5pt}
\begin{equation*}
\eta_Z = \left(
\begin{tikzcd}[column sep=3em,row sep=3em]
G \times Z \arrow[r, "u"{outer sep = 2pt}] \arrow[d, "p_2"{right, outer sep = 2pt}] & X \\
Z
\end{tikzcd}\right), \vspace{8pt}
\end{equation*}
where $ u $ sends $ (g,z) \mapsto f(gz)=gf(z). $
There are two arrows
\vspace{10pt}
\begin{equation*}
\begin{tikzcd}[column sep=3em,row sep=1.5em]
& \eta_Y \\
\eta_Z \arrow[ur, "C_1"] \arrow[dr,"C_2"{below left}] & \\
& \eta_1
\end{tikzcd}
\quad \text{over} \quad
\begin{tikzcd}[column sep=3em,row sep=1.5em]
& Y \\
Z \arrow[ur, "c_1"] \arrow[dr,"c_2"{below left}] & \\
& X
\end{tikzcd},\vspace{10pt}
\end{equation*}
which are described by two commutative diagrams
\vspace{5pt}
\begin{equation*}
C_1 :
\begin{tikzcd}[column sep=3em,row sep=4em]
G \times Z \arrow[rr, bend left=40, "u"{outer sep = 2pt}] \arrow{r}{\sigma_Z} \arrow[d, "p_2"{right, outer sep = 2pt}] & Z \arrow[d, "c_1"{right, outer sep = 2pt}] \arrow[r, "c_2"{outer sep = 2pt}] & X \\
Z \arrow[r, "c_1"{outer sep = 2pt}]  & Y &
\end{tikzcd}
\end{equation*}
and
\begin{equation*}
C_2 :
\begin{tikzcd}[column sep=3em,row sep=4em]
G \times Z \arrow[rr, bend left=40, "u"{outer sep = 2pt}] \arrow[r, "\id_G \times c_2"{outer sep = 2pt}] \arrow[d, "p_2"{right, outer sep = 2pt}] & G \times X \arrow[d, "p_2"{right, outer sep = 2pt}] \arrow[r, "\sigma"{outer sep = 2pt}] & X \\
Z \arrow[r, "c_2"{outer sep = 2pt}]  & X &,
\end{tikzcd} \vspace{15pt}
\end{equation*}
where $ \sigma_Z: G \times Z \to Z $ denotes the free $ G $-action on $ Z. $ There exist unique morphisms $ \vphi_Z: \cE(\eta_Z) \to \cF(\eta_Z) $ and $ \psi_Y: c_1^*\cE(\eta_Y) \to c_1^*\cF(\eta_Y) $ of sheaves on $ Z $ such that the following diagram commutes:
\vspace{10pt}
\begin{equation}\label{eq_triangular_prism2}
\begin{tikzcd}[column sep=6em,row sep=1.2em]
& \cF(\eta_Z) \ar[dd, "\cF(C_1)"{left}, "\simeq"{above, sloped}] 
\ar[dddr, "\cF(C_2)"{near end}, "\simeq"{below, near end, sloped}] & \\
\cE(\eta_Z) \ar[ur, dashed, "\vphi_Z"{near end, outer sep = 2pt}] 
\ar[dd, "\cE(C_1)"{left}, "\simeq"{above, sloped}] & & \\
& c_1^*\cF(\eta_Y) \ar[dr, ""{below, sloped}, "\simeq"{near start, sloped}] & \\
c_1^*\cE(\eta_Y) \ar[ur, dashed, "\psi_Y"{very near end, outer sep = 2pt}] 
\ar[dr, ""{below, sloped}, "\simeq"{near start, sloped}] & & c_2^*\cF(\eta_1) \\
& c_2^*\cE(\eta_1) 
\ar[uuul, leftarrow, crossing over, "\cE(C_2)"{above right, near start}, "\simeq"{below, sloped, near start}] 
\ar[ur, "c_2^*\vphi"] &
\end{tikzcd}\vspace{15pt}
\end{equation}
Since $ c_2 $ is a $ G $-equivariant morphism of schemes, the sheaves $ c_2^*\cE(\eta_X) $ and $ c_2^*\cF(\eta_X) $ on $ Z $ can be endowed with canonical $ G $-equivariant structures, which are respected by the morphism $ c_2^*\vphi. $ Therefore, every sheaf in the above diagram can be endowed with a $ G $-equivariant structure and each morphism is a $ G $-equivariant morphism. Now we apply the descent theory 
for schemes to the $ G $-torsor $ c_1: Z \to Y $ with the following pullback diagram:
\vspace{10pt}
\begin{equation*}
\begin{tikzcd}[column sep=3em,row sep=3em]
G \times Z \arrow[r, "\sigma_Z"{outer sep = 2pt}] \arrow[d, "p_2"{right, outer sep = 2pt}] & Z \arrow[d, "c_1"{right, outer sep = 2pt}] \\
Z \arrow[r, "c_1"{outer sep = 2pt}]  & Y
\end{tikzcd}\vspace{10pt}
\end{equation*}
The $ G $-equivariant morphism $ \psi_Y: c_1^*\cE(\eta_Y) \to c_1^*\cF(\eta_Y) $ of $ G $-equivariant sheaves on $ Z $ descends to a morphism $ \vphi_Y: \cE(\eta_Y) \to \cF(\eta_Y) $ of sheaves on $ Y, $ i.e., 
$$ \psi_Y = c_1^*\vphi_Y. $$
When $ Y = X $ and $ \eta_Y = \eta_1, $ diagram (\ref{eq_triangular_prism2}) coincides with diagram (\ref{eq_triangular_prism1}) and hence $ \vphi_X = \vphi. $
Take any arrow $ A: \eta_{Y_1} \to \eta_{Y_2} $ over a morphism $ a: Y_1 \to Y_2 $ of schemes. We have a commutative diagram of arrows:
\vspace{10pt}
\begin{equation}\label{eq_4_arrows}
\begin{tikzcd}[column sep=3em,row sep=3em]
\eta_{Z_1} \arrow[r, "C"{outer sep = 2pt}] \arrow[d, "D"{right, outer sep = 2pt}] & \eta_{Z_2} \arrow[d, "B"{right, outer sep = 2pt}] \\
\eta_{Y_1} \arrow[r, "A"{outer sep = 2pt}]  & \eta_{Y_2}
\end{tikzcd}
\quad \quad \text{over} \quad \quad
\begin{tikzcd}[column sep=3em,row sep=3em]
{Z_1} \arrow[r, "c"{outer sep = 2pt}] \arrow[d, "d"{right, outer sep = 2pt}] & {Z_2} \arrow[d, "b"{right, outer sep = 2pt}] \\
{Y_1} \arrow[r, "a"{outer sep = 2pt}]  & {Y_2}
\end{tikzcd}\vspace{10pt}
\end{equation}
We also have a composition of arrows:
\vspace{3pt}
\begin{equation}\label{eq_2_arrows}
\begin{tikzcd}[column sep=2em,row sep=3em]
\eta_{Z_1} \arrow[r,"C"] & \eta_{Z_2} \arrow[r,"E"] & \eta_1
\end{tikzcd}
\quad \quad \text{over} \quad \quad
\begin{tikzcd}[column sep=2em,row sep=3em]
{Z_1} \arrow[r,"c"] & {Z_2} \arrow[r,"e"] & X
\end{tikzcd}\vspace{5pt}
\end{equation}
By construction, the morphisms $ \vphi_{Z_1} $ and $ c^*\vphi_{Z_2} $ of sheaves on $ Z_1 $ fit into the following commutative diagram:
\vspace{5pt}
\begin{equation*}
\begin{tikzcd}[column sep=3em,row sep=3em]
\cE(\eta_{Z_1}) \ar[d,"\cE(E \circ C)"{left},"\simeq"{above,sloped}] \ar[r,"\vphi_{Z_1}"] &  \cF(\eta_{Z_1}) \ar[d,"\cF(E \circ C)"{right},"\simeq"{below,sloped}] \\
c^*e^*\cE(\eta_{1}) \ar[d,"c^*\cE(E)^{-1}"{left},"\simeq"{above,sloped}] \ar[r,"c^*e^*\vphi"] & c^*e^*\cF(\eta_{1}) \ar[d,"c^*\cF(E)^{-1}"{right},"\simeq"{below,sloped}] \\
c^*\cE(\eta_{Z_2}) \ar[r,"c^*\vphi_{Z_2}"] &  c^*\cF(\eta_{Z_2}),
\end{tikzcd}\vspace{8pt}
\end{equation*}
and hence 
$$ c^*\vphi_{Z_2} \circ \cE(C) = \cF(C) \circ \vphi_{Z_1}. $$
Now we have the following commutative cube:
\vspace{5pt}
\begin{equation*}
\begin{tikzcd}[row sep=2.6em, column sep=1.2em]
& \cE(\eta_{Z_1}) \ar[dl,"\cE(C)"{above left},"\simeq"{below,sloped}] \ar[rr,"\vphi_{Z_1}"] 
\ar[dd,"\cE(D)"{left,near end},"\simeq"{above, near end, sloped}] & & \cF(\eta_{Z_1}) \ar[dl,"\cF(C)","\simeq"{above, sloped}] \ar[dd,"\cF(D)"{right},"\simeq"{below, sloped}] \\
c^*\cE(\eta_{Z_2}) \ar[rr, crossing over, "c^*\vphi_{Z_2}"{near end}] \ar[dd,"c^*\cE(B)"{left},"\simeq"{above, sloped}] & & c^*\cF(\eta_{Z_2}) \\
& d^*\cE(\eta_{Y_1}) \ar[dl, dashed, "d^*\cE(A)"{above left},"\simeq"{below, sloped}] \ar[rr, "d^*\vphi_{Y_1}"{near start}] & & d^*\cF(\eta_{Y_1}) \ar[dl, dashed, "d^*\cF(A)","\simeq"{above, sloped}] \\
c^*b^*\cE(\eta_{Y_2}) \ar[rr, "c^*b^*\vphi_{Y_2}"] & & c^*b^*\cF(\eta_{Y_2}) \ar[from=uu, crossing over, "c^*\cF(B)"{right,near start},"\simeq"{below, near start, sloped}]
\end{tikzcd},\vspace{10pt}
\end{equation*}
where the bottom is the pullback $ d^* $ of the commutative diagram
\vspace{10pt}
\begin{equation*}
\begin{tikzcd}[column sep=3em,row sep=4em]
\cE(\eta_{Y_1}) \arrow[r,"\vphi_{Y_1}"] \arrow[d,"\cE(A)"{left},"\simeq"{above, sloped}] & \cF(\eta_{Y_1}) \arrow[d,"\cF(A)"{right},"\simeq"{below, sloped}] \\
a^*\cE(\eta_{Y_2}) \arrow[r,"a^*\vphi_{Y_2}"]  & a^*\cF(\eta_{Y_2}).
\end{tikzcd}\vspace{8pt}
\end{equation*}
Now we define $ \Phi: \cE \to \cF $ by $ \Phi(\eta_Y) = \vphi_Y. $ Then $ \Phi $ satisfies condition (\ref{eq_sheaf_morphism_on_stack}), $ p^*\Phi = \Phi(\eta_1) = \vphi_X = \vphi $ and is unique by diagram (\ref{eq_triangular_prism2}).

Last, take any object $ (F, \phi) $ in $ \Sh^G(X). $ We will show that there is a sheaf $ \cE $ on $ \cX $ such that $ T(\cF) = (p^*\cE, \phi_\cE) = (\cE(\eta_1), \phi_\cE) \cong (F, \phi). $ For any object $ \eta_Y $ in diagram (\ref{eq_arbitrary_object}), the pullback $ c_2^*F $ of the sheaf $ F $ along $ c_2 $ descends to a sheaf $ \cE(\eta_Y) $ along $ c_1, $ i.e.,
\begin{equation*}
c_1^*\cE(\eta_Y) = c_2^*F.
\end{equation*} 
Consider the arrows in diagrams (\ref{eq_4_arrows}) and (\ref{eq_2_arrows}). The isomorphism $ \cE(A): \cE(\eta_{Y_1}) \xrightarrow{\sim} a^*\cE(\eta_{Y_2}) $ is uniquely determined by its pullback along $ d, $
\begin{equation*}
d^*\cE(A): d^*\cE(\eta_{Y_1}) = (e \circ c)^*F \xrightarrow{\sim} c^*e^*F = d^*a^*\cE(\eta_{Y_2}).
\end{equation*}
Then $ \cE $ satisfies the cocycle condition to be a sheaf on $ \cX. $ For the canonical object $ \eta_1 $ in diagram (\ref{eq_can_obj}), we have
\begin{equation*}
p_2^*\cE(\eta_1) = \sigma^*F,
\end{equation*}
which restricts on $ \{1\} \times X $ to an isomorphism $ f: \cE(\eta_1) \xrightarrow{\sim} F $. The $ G $-equivariant structure $ \phi_{\cE}: p_2^*\cE(\eta_1) \xrightarrow{\sim} \sigma^*\cE(\eta_1) $ is given by
\begin{equation*}
\phi_{\cE} = \cE(A_2) \circ \cE(A_1)^{-1} = \sigma^* f^{-1},
\end{equation*}
which makes the following diagram commutes:
\vspace{2pt}
\begin{equation*}
\begin{tikzcd}[column sep=3em,row sep=3em]
p_2^*\cE(\eta_1) \arrow[r,"p_2^*f"{outer sep = 2pt}] \arrow[d,"\phi_{\cE}"{left, outer sep = 2pt},"\simeq"{above, sloped}] & p_2^*F \arrow[d,"\phi"{right, outer sep = 2pt},"\simeq"{below, sloped}] \\
\sigma^*\cE(\eta_1) \arrow[r,"\sigma^*f"{outer sep = 2pt}]  & \sigma^*F
\end{tikzcd}\vspace{8pt}
\end{equation*}
Therefore, we have an isomorphism
\begin{equation*}
f: T(\cF) = (\cE(\eta_1), \phi_\cE) \xrightarrow{\sim} (F, \phi)
\end{equation*}
of $ G $-equivariant sheaves on $ X. $
\vspace{1pt}

\section{Mukai pairing and HRR formula}\label{review_mukai}

In this appendix we review a few notions in the intersection theory of schemes, including the Mukai vector and the Mukai pairing, and prove the Hirzebruch-Riemann-Roch (HRR) formula. 

Fix an algebraically closed field $ k $. Let $ X $ be a separated scheme of finite type over $ k $. We first recall the Euler characteristic which appears in the left-hand side of the HRR theorem.
\begin{defn}
	If $ X $ is proper, then the \tb{Euler characteristic} 
	$$ \chi(X, \ \cdot\ ): K(X) \to \ZZ $$
	is defined by
	\begin{equation*}
	\chi(X, E) = \sum_i (-1)^i \dim H^i(X, E)
	\end{equation*}
	for a coherent sheaf $ E $ on $ X $ and extended linearly.
\end{defn}
Recall the HRR theorem for schemes.
\begin{thm}[HRR]
	Let $ X $ be a proper smooth scheme. For all $ x $ in $ K(X), $ we have
	\begin{equation*}
	\chi(X, x) = \int_X \ch(x) \td_X
	\end{equation*}
	in $ \ZZ. $
\end{thm}
\begin{defn}
	If $ X $ is proper, then the \tb{Euler pairing}
	\begin{equation*}
	\chi: K(X) \times K(X) \to \ZZ
	\end{equation*}
	is defined by
	\begin{equation*}
	\chi(E,F) = \sum_{i} (-1)^i \dim \Ext^I(E,F)
	\end{equation*}
	for two coherent sheaves $ E $ and $ F $ on $ X $ and extended bilinearly.
\end{defn}
If $ X $ is smooth, then $ X $ has the resolution property and $ K(X) \cong K^0(X), $ so we can define an involution on $ K(X) $ via the involution on $ K^0(X) $ which is induced by taking dual vector bundles. 

Let $ X $ be a separated smooth scheme of finite type.
\begin{defn}
	The involution $ (\ \cdot \ )^\vee: K(X) \to K(X) $ is defined by the composition
	\begin{equation*}
	(\ \cdot \ )^\vee: K(X) \xrightarrow{\beta} K^0(X) \xrightarrow{(\ \cdot \ )^\vee} K^0(X) \xrightarrow{\beta^{-1}} K(X)
	\end{equation*}
	where $ \beta $ is the natural isomorphism $ K(X) \xrightarrow{\sim} K^0(X). $
\end{defn}
The following result establishes the relation between the Euler characteristic and the Euler pairing.
\begin{lem}
	Let $ X $ be a proper smooth scheme. For all $ x $ and $ y $ in $ K(X), $ we have
	\begin{equation*}
	\chi(x,y) = \chi(X, x^\vee y)
	\end{equation*}
	in $ \ZZ. $ In particular, $ \chi(1,x) = \chi(X,x) $ for all $ x $ in $ K(X) $ where $ 1 = [\cO_X] $.
\end{lem}
\begin{proof}
	By the bilinearity of $ \chi, $ it suffices to consider $ x = [V] $ and $ y =[W] $ for vector bundles $ V $ and $ W $ on $ X. $ Then we have
	\begin{align*}
	\chi(x, y) & = \sum_i (-1)^i \dim \Ext^i(V,W) \\
	& = \sum_i (-1)^i \dim H^i(X,V^\vee \otimes W) & \text{because $ V $ is locally free} \\
	& = \chi(X, V^\vee \otimes W) \\
	& = \chi(X, x^\vee y)
	\end{align*}
\end{proof}
Since $ X $ is smooth, the Chow groups $ A(X) = A^*(X) $ form a graded ring with an intersection product. We define an involution on the rational Chow ring $ A(X)_\QQ $ as follows.
\begin{defn}
	The involution $ (\ \cdot \ )^\vee: A(X)_\QQ \to A(X)_\QQ $ is defined by 
	\begin{equation*}
	v^\vee = \sum (-1)^i v_i
	\end{equation*}
	for all $ v = \sum_i v_i $ in $ A(X)_\QQ $ where each $ v_i \in A^i(X)_\QQ. $
\end{defn}
For a unit $ 1 + x = 1 + \sum_{i\geq 1} x_i $ in a graded ring $ R $ with $ x_i \in R^i $ for $ i \geq 1 $ such that $ x^n = 0 $ for some $ n \geq 1, $ we can define its square root by a Taylor series:
\begin{equation*}
\sqrt{1+x} = \sum_i \binom{1/2}{i} x^i = 1 + \frac{x}{2} - \frac{x^2}{8} + \cdots
\end{equation*}
A quick computation shows the following properties of the involution $ (\ \cdot \ )^\vee $ on $ A(X)_\QQ $.
\begin{lem}
	For all $ v $ and $ w $ in $ A(X)_\QQ $, we have 
	\begin{equation*}
	(vw)^\vee = v^\vee w^\vee \quad \text{and} \quad \sqrt{v^\vee} = \sqrt{v}^\vee
	\end{equation*}
	in $ A(X)_\QQ $ whenever $ \sqrt{v} $ is defined.
\end{lem}
The Todd class of $ X $
\begin{equation*}
\td_X = 1 + \frac{c_1}{2} + \frac{c_1^2+c_2}{12} + \cdots
\end{equation*}
is a unit in $ A(X)_\QQ $ with $ c_i = c_i(X) = c_i(TX) $ for $ i > 0 $,
and it has a well-defined square root
\begin{equation*}
\sqrt{\td_X} = 1 + \frac{c_1}{4} + \frac{c_1^2+4c_2}{96} + \cdots
\end{equation*}
in $ A(X)_\QQ^\times. $
\begin{defn}\label{defn_v}
	The \tb{Mukai vector} map 
	$$ v: K(X) \to A(X)_\QQ $$
	is defined by
	\begin{equation*}
	v(x) = \ch(x) \sqrt{\td_X}
	\end{equation*}
	for all $ x \in K(X). $
\end{defn}
The following properties of the Chern Character, the Mukai vector, and the involutions on $ K(X) $ and $ A(X)_\QQ $ can be easily checked.
\begin{lem}\label{lem_ch_v_inv}
	For all $ x $ and $ y $ in $ K(X), $ we have:
	\begin{enumerate}[font=\normalfont,leftmargin=2em]
		\item $ v(x+y) = v(x) + v(y) $
		\item $ v(xy) = v(x)\ch(y) $
		\item $ \ch(x^\vee) = \ch(x)^\vee $
		\item $ v(x^\vee) = v(x)^\vee \sqrt{\td_X/\td_X^\vee} $
	\end{enumerate}
\end{lem}
\begin{rmk}
	We have the identity
	$$ \sqrt{\td_X/\td_X^\vee} = e^{c_1(X)/2} $$
	in $ A(X)_\QQ $ by applying the splitting principle to $ \td_X $
	as in the proof of Lemma 5.41 in \cite{huybrechts2006fourier}.
\end{rmk}
If $ X $ is proper, then we can define a pairing on $ A(X)_\QQ $.
\begin{defn}\label{defn_v_pairing}
	Let $ X $ be a proper smooth scheme. The \tb{Mukai pairing}
	$$ \inprod{\cdot\ {,}\ \cdot}: A(X)_\QQ \times A(X)_\QQ \to \QQ $$ 
	is defined by
	\begin{equation*}
	\inprod{v,w} = \int_{X} v^\vee w \sqrt{\td_X/\td_X^\vee}
	\end{equation*}
	for all $ v $ and $ w $ in $ A(X)_\QQ. $
\end{defn}
The HRR theorem for schemes implies the following HRR formula (\ref{eq_HRR_schemes}) in terms of the Euler pairing and the Mukai pairing. Another version of formula (\ref{eq_HRR_schemes}) is formula (5.5) in \cite{huybrechts2006fourier} where the inputs are complexes of sheaves on the scheme $ X. $
\begin{thm}[HRR Formula]
	Let $ X $ be a proper smooth scheme. For all $ x $ and $ y $ in $ K(X), $ we have
	\begin{equation}\label{eq_HRR_schemes}
	\chi(x, y) = \inprod{v(x),v(y)}
	\end{equation}
	in $ \ZZ. $
\end{thm}
\begin{proof}
	Take any $ x $ and $ y $ in $ K(X). $ We then have
	\begin{align*}
	\chi(x,y) & = \chi(X, x^\vee y) & \text{by Lemma \ref{lem_Euler_pairing} for schemes} \\
	& = \int_{X} \ch(x^\vee y) \td_X & \text{by the HRR theorem for schemes} \\
	& = \int_{X} \ch(x^\vee) \ch(y) \td_X & \text{because $ \ch $ is a ring map} \\
	& = \int_{X} v(x^\vee) v(y) & \text{by Definition \ref{defn_v}} \\
	& = \int_{X} v(x)^\vee v(y) \sqrt{\td_{X}/\td_{X}^\vee} & \text{by Lemma \ref{lem_ch_v_inv}} \\
	& = \inprod{v(x),v(y)} & \text{by Definition \ref{defn_v_pairing}}
	\end{align*}
\end{proof}
For a separated smooth scheme $ X $ of finite type, the Chern character map 
$$ \ch: K(X) \to A(X)_\QQ $$ 
is a ring homomorphism and becomes a $ \QQ $-algebra isomorphism after tensored with $ \QQ, $ so the Mukai vector map $ v: K(X) \to A(X)_\QQ $
becomes an isomorphism of $ \QQ $-vector spaces after tensored with $ \QQ $, since $ \sqrt{\td_X} $ is a unit in $ A(\cX)_\QQ. $ If $ X $ is proper, then we can extend the Euler pairing to $ K(X)_\QQ, $ and hence we have two vector spaces $ K(X)_\QQ $ and $ A(\cX)_\QQ $ with bilinear forms. Therefore, we have the following:
\begin{prop}
	Let $ X $ be a proper smooth scheme. The Mukai vector map 
	$$ v: K(X) \to A(X)_\QQ $$ 
	induces a linear isometry
	\begin{equation*}
	v: \left(K(X)_\QQ, \chi\right) \xrightarrow{\simeq} \left(A(X)_\QQ, \inprod{\cdot\ {,}\ \cdot}\right).
	\end{equation*}
\end{prop}

\newcommand{\etalchar}[1]{$^{#1}$}

\end{document}